\documentclass[times, doublespace]{nmeauth}

\usepackage{graphicx}

\usepackage{amssymb}

\usepackage{amsmath}
\usepackage[hang,raggedright]{subfigure}
\usepackage{mathtools,stmaryrd}
\DeclarePairedDelimiter{\dbr}{\llbracket}{\rrbracket}

\usepackage{lineno}

\usepackage{xspace}

\usepackage[square, numbers, comma, sort&compress]{natbib}

\usepackage[squaren]{SIunits}

\usepackage{binary}
\usepackage{cite}
\usepackage{verbatim}  
\usepackage{color}

\usepackage{listings}
\usepackage{url}
\usepackage{math tools,cancel}
\usepackage{units}

\usepackage[mathscr]{euscript}
\DeclareSymbolFont{rsfs}{U}{rsfs}{m}{n}
\DeclareSymbolFontAlphabet{\mathscrsfs}{rsfs}

\lstdefinelanguage{C++enhanced}[]{C++}{%
  morekeywords={volDiagTensorField, volSymmTensor4thOrderField, fvVectorMatrix},
}
\definecolor{gray}{rgb}{0.5,0.5,0.5}
\usepackage{courier}
\usepackage{algorithm}
\usepackage[noend]{algpseudocode}
\makeatletter
\def\BState{\State\hskip-\ALG@thistlm}
\makeatother

\usepackage{pgfplots}
\pgfplotsset{compat=1.8}

\newcommand{\etc} {\textit{etc.}\@\xspace}
\newcommand{\eg} {\textit{e.g.}\@\xspace}
\newcommand{\ie} {\textit{i.e.}\@\xspace}

\newcommand{\etal} {et~al.\xspace}
\renewcommand*{\refname}{}

\newcommand{\ra}[1]{\renewcommand{\arraystretch}{#1}}

\begin{document}


\runningheads{Thirty years of the finite volume method for solid mechanics}{Thirty years of the finite volume method for solid mechanics}

\title{Thirty years of the finite volume method for solid mechanics}

\author{
P. Cardiff\affil{1}\corrauth,
I.\ Demird\v{z}i\'{c}\affil{2}
}

\address{
\affilnum{1}University College Dublin, Bekaert University Technology Centre, School of Mechanical and Materials Engineering, Belfield, Ireland \break
\affilnum{2}Ma\v{s}inski fakultet Sarajevo, Vilsonovo \v{s}etali\v{s}te 9, 71000 Sarajevo, Bosnia-Herzegovina
}

\corraddr{P. Cardiff, University College Dublin, Bekaert University Technology Centre, School of Mechanical and Materials Engineering, Belfield, Ireland. E-mail: philip.cardiff@ucd.ie}
\corraddr{I. Demird\v{z}i\'{c}, Ma\v{s}inski fakultet Sarajevo, Vilsonovo \v{s}etali\v{s}te 9, 71000 Sarajevo, Bosnia Herzegovina. E-mail: i.demirdzic@yahoo.com}

\begin{abstract}
Since early publications in the late 1980s and early 1990s, the finite volume method has been shown suitable for solid mechanics analyses.
At present, there are several flavours of the method, which can be classified in a variety of ways, such as grid arrangement (cell-centred vs staggered vs vertex-centred), solution algorithm (implicit vs explicit), and stabilisation strategy (Rhie-Chow vs Jameson-Schmidt-Turkel vs Godunov upwinding).
This article gives an overview, historical perspective, comparison and critical analysis of the different approaches where a close comparison with the \emph{de facto} standard for computational solid mechanics, the finite element method, is given.
The article finishes with a look towards future research directions and steps required for finite volume solid mechanics to achieve more widespread acceptance.
\end{abstract}


\maketitle



\section{Introduction}
``\emph{Now however I recognise the FVM/FEM dichotomy as being comparable with those between Protestant and Catholic, or Sunni and Shia. That is to say that it promotes needless conflict; and expense; and loss of opportunity.}'' \citep{Runchal2017}
With these words, Brian Spalding highlighted the, at times, unproductive nature of the debate around the relative merits of the finite volume (FVM) and finite element (FEM) methods.
Although accepted in the Computational Fluid Dynamics (CFD) field, there remains a reluctance and general confusion around the application of the finite volume method to solid mechanics.
The aim of this article is to clarify this confusion, by: providing an overview of the significant developments within the field; linking variants of the finite volume method for solid mechanics analyses; comparing finite volume methods with standard finite element methods; and, finally, identifying future directions for the field.
Building on the foundations of the finite difference method, the finite volume method is a generalisation in terms of geometry and topology: simple finite volume schemes reduce to finite difference schemes.
Whereas the finite difference method is based on nodal relations for \emph{differential} equations, the finite volume method balances forces acting on control volumes, directly discretising the \emph{integral} form of the conservation laws.
A number of prior developments within CFD provided the foundation for the earliest paper on the cell-centred finite volume method for solid mechanics by \citet{Demirdzic1988}.
In the subsequent three decades the finite volume method for solid mechanics has developed in a number of directions, differing in terms of discretisation, solution methodology and overall philosophy.
The varied approaches may be classified in a number of ways, including, for example:
\begin{itemize}
	\item \emph{Grid arrangement}: cell-centred \citep{Demirdzic1988, Demirdzic1993} vs vertex-centred \citep{Bailey1991, Fryer1991, Zienkiewicz1991, Cross1992, Onate1992, Cross1993, Fryer1993, Idelsohn1994, Onate1994} vs staggered-grid \citep{Beale1990a, Beale1990b, Bukhari1990, Hattel1990, Bukhari1991, Hattel1992, Hattel1993a}, as well as more recently face-centred \citep{Sevilla2018a, Sevilla2018b} and meshless \citep{Ebrahimnejad2014, Fallah2018a};
	\item \emph{Solution algorithm}: implicit \citep{Demirdzic1993, Bailey1995, Cardiff2016a} vs explicit (matrix-free) \citep{Kluth2010, Lee2013, Haider2017};
	\item \emph{Stabilisation approach}: Rhie-Chow \citep{Demirdzic1995, Cardiff2016b} vs Jameson-Schmidt-Turkel \citep{Aguirre2014} vs Godunov \emph{two-sided} upwinding \citep{LeVeque2004, Lee2013}.
\end{itemize}
There are countless other ways to classify the approaches, for example, based on force discretisation at the control volume face, however, the current article will base its discussion primarily around the three classification types above.

The essential characteristic of a finite volume method is the integration of the governing conservation equation over finite control volumes.
In this sense, the \emph{cell-centred} vs \emph{vertex-centred} vs \emph{staggered-grid} approaches primarily differ in \emph{how} these control volumes are constructed.
Cell-centred approaches use the primary mesh control volumes, whereas vertex-centred approaches construct a \emph{dual mesh} with control volumes surrounding the vertices of the primary mesh.
On the other hand, staggered-grid approaches create multiple secondary meshes, one for each scalar component of the primary solution vector, constructed about the primary mesh faces.
Despite their close relationship, as explored further in Section \ref{sec:compareFiniteVolumeVariants}, approaches based on differing grid arrangements often differ greatly in terms of philosophy: as noted by \citet{Baliga2009}, the cell-centred method is often thought of as a control-volume finite difference method, combining ideas borrowed from finite volume and finite difference methods; whereas, the vertex-centred approach is viewed as a control-volume finite element method, formulated by amalgamating concepts native to finite volume and finite element methods.

Regardless of the chosen grid arrangement, spatio-temporal integration of the governing equations may adopt an \emph{implicit} or \emph{explicit} solution algorithm.
Implicit algorithms are characterised by the solution of a linear system of equations and are unconditionally stable with respect to the time increment size. In contrast, explicit or matrix-free algorithms avoid the need to construct such a system of linear equations but the time increment size is restricted by the classic Courant-Friedrichs-Lewy constraint \citep{Courant1928}.
The choice of solution algorithm often depends on the problems of interest, with implicit methods favoured for elliptic and parabolic cases (steady-state and quasi-steady-state) and explicit methods for hyperbolic (high rates, wave propagation).
Once the grid arrangement and solution algorithm are selected, care must be taken in the construction of a stable discretisation that does not suffer from unphysical instabilities in the solution field.
To this end, a variety of stabilisation approaches have been proposed.
Generally, each approach can be viewed upon as adding some form of diffusion to the surface force discretisation with the purpose of quelling high-frequency oscillations.
Rhie-Chow-style stabilisation is common in the implicit approaches originating from the work of \citet{Demirdzic1995}, while Godunov upwinding and Jameson-Schmidt-Turkel approaches are more commonly seen in the explicit approaches, rooted in the solution of compressible gas flow Euler equations.

From afar, each of these main variants of finite volume method can seem quite distinct; however, upon closer inspective, they share many similarities in terms of discretisation and solution algorithm.
Nevertheless, the fragmented nature of the finite volume solid mechanics community is evident from many recent publications in the area, where authors, reviewers and editors tend to be unaware of developments within the field, for example, see \citep{Demirdzic2015}.
In addition, and more generally, there is a lack of awareness in the computational mechanics community around the capabilities of the finite volume method for solid mechanics.
Accordingly, a comparative and critical review of the finite volume method for solid mechanics is timely, relevant and essential for the future progress of the field.
Within this domain, there are a number of open questions:
What are the strengths and weaknesses of the various approaches?
Which approaches show the greatest potential and widest applicability?
Are there possibilities for the various methods to be combined to produce superior methods?
How do the finite volume approaches relate to finite element methods?
Are there directions of development which are missing?
This article will attempt to provide answers to these questions, as well as providing a unifying framework for the discussion of finite volume methods for solid mechanics, and their relationship with finite element methods.
Given the scale of the field, it is outside the scope of the article to provide an exhaustive review of all formulations; instead, the article aims to provide detailed analysis on the main variants of approach common in the literature.
The primary novel contributions of the current article are twofold:
(i) The first detailed review of the finite volume solid mechanics field is presented;
and 
(ii) The first analysis of similarities and differences between all main variants of the finite volume method for solid mechanics is detailed.
(ii) A unique analysis of the finite volume method for solid mechanics is performed, cumulating in the description of a unified approach.

The remainder of the article is constructed as follows:
Section 2 provides a chronological overview of the prominent finite volume developments for solid mechanics.
Section 3 compares and contrasts the three main flavours of the finite volume method for solid mechanics.
In Section 4, the finite volume method for solid mechanics is compared to the ``standard'' continuous Bubnov-Galerkin finite element method, highlighting similarities and differences in terms of discretisation, solution methodology and overall philosophy.
Section 5 reviews the variety of structural applications to which the finite volume method has been applied in its thirty year history.
The penultimate section briefly reviews software that use, or have previously used, the finite volume method for solid mechanics simulations.
Finally, Section 7 summaries the main conclusions of the article and considers the current challenges facing the field.
\section{History of the finite volume method for computational solid mechanics}
\label{sec:history}



The development of the finite volume method for solid mechanics has occurred independently in a number of forms, with finite difference methods, computational fluid mechanics algorithms, and finite element methods providing much of the inspiration.
This section gives an overview of the historical development of these finite volume methods.
The treatise is primarily partitioned based on the grid arrangement, where comments regarding the solution algorithm and stabilisation scheme are given where appropriate.
In-depth dissections of the technical details are left to Section \ref{sec:compareFiniteVolumeVariants}.
While the field is small and fragmented, there are a number of notable reviews worth mentioning, including those of \citet{Maneeratana2000b},
 \citet{Vaz2009} 
 and Cavalcante \etal \citep{Cavalcante2012c}.

\begin{figure}[htb]
   \centering
   \includegraphics[width=\textwidth]{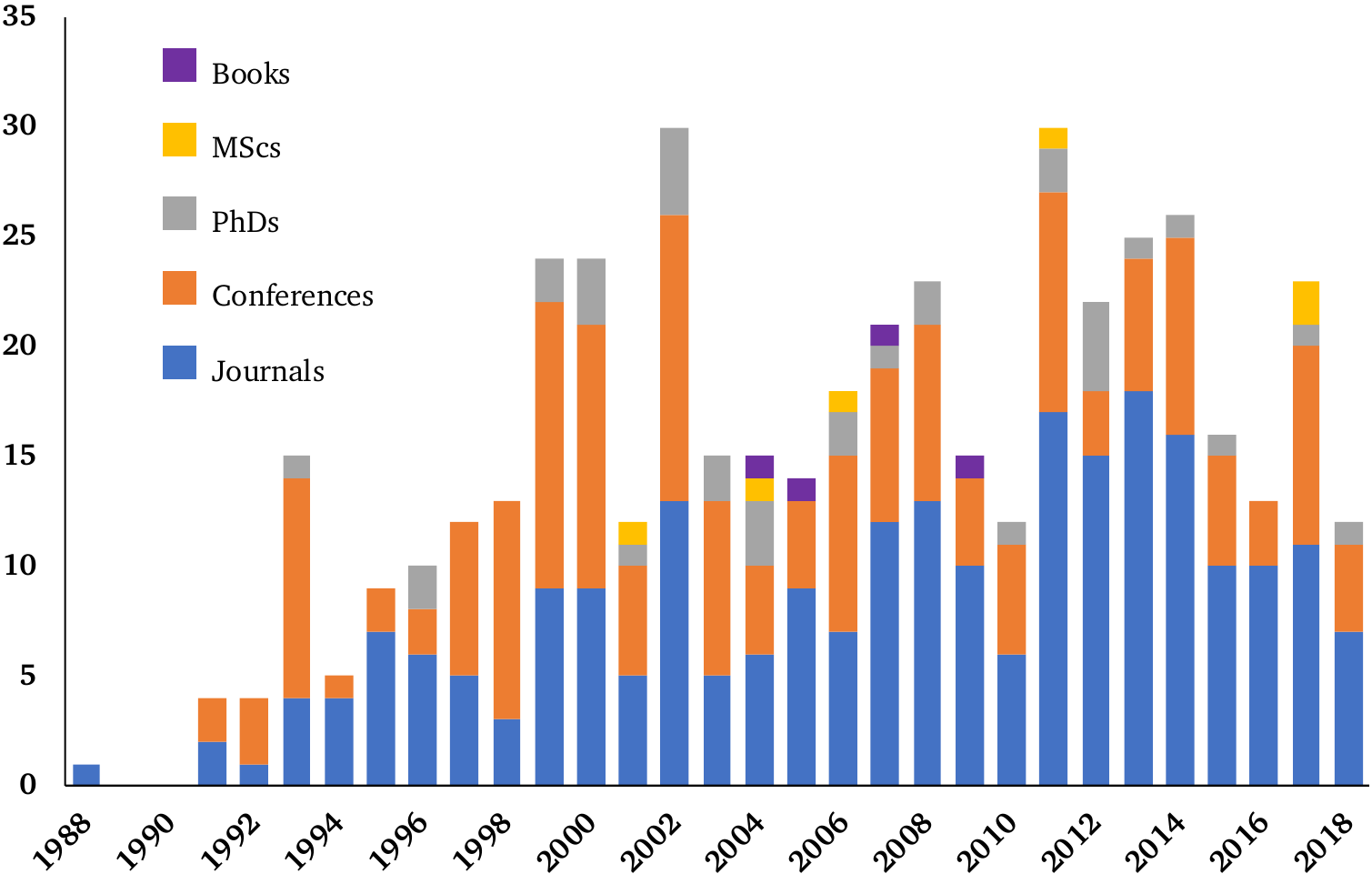} 
   \caption{Histogram of the finite volume method for solid mechanics publications to-date}
   \label{fig:publicationHistogram}
\end{figure}

Before delving into details of influential publications, it is insightful to first consider the literature landscape as a whole.
To this end, Figure \ref{fig:publicationHistogram} presents a histogram of the publications in the area to-date, separated into journal articles, conferences, Ph.D. theses, Masters theses and books.
Of course, exact records of each publication type are difficult to track and consequently the data should be taken as indicative.
To complement this, Figure \ref{fig:journalArticleCounts} lists the most popular international journals for publishing finite volume solid mechanics articles; only journals with greater than five articles have been included.
\begin{figure}
\centering
\begin{tikzpicture}
\begin{axis}[ 
xbar, xmin=0,
xlabel={Number of related articles},
symbolic y coords={%
    {Composite Structures},
    {International Journal of Fracture},
    {Composites Part B: Engineering},
    {Computer Modeling in Eng. and Sci.},
    {Engineering Fracture Mechanics},
    {International Journal of Plasticity},
    {Numerical Heat Transfer, Part B},
    {Computational Materials Science},
    {Computers \& Structures},
    {ASME Journal of Applied Mechanics},
    {Int. J. Solids Struct.},
    {Comput. Methods Appl. Mech. Eng.},
    {Applied Mathematical Modelling},
    {Journal of Computational Physics},
    {Int. J. Numer. Methods Eng.}
    },
ytick=data,
nodes near coords, 
nodes near coords align={horizontal},
ytick=data,
]
\addplot coordinates {
    (5,{Composite Structures})
    (5,{International Journal of Fracture})
    (6,{Composites Part B: Engineering})
    (6,{Computer Modeling in Eng. and Sci.})
    (7,{Engineering Fracture Mechanics})
    (8,{International Journal of Plasticity})
    (9,{Numerical Heat Transfer, Part B})
    (9,{Computational Materials Science})
    (9,{Computers \& Structures})
    (11,{ASME Journal of Applied Mechanics})
    (12,{Int. J. Solids Struct.}) 
    (15,{Comput. Methods Appl. Mech. Eng.})
    (17,{Applied Mathematical Modelling})
    (20,Journal of Computational Physics)
    (24,{Int. J. Numer. Methods Eng.}) 
    };
\end{axis} 
\end{tikzpicture}
\caption{Number of publications per journal related to computational solid mechanics using the finite volume method, based on the cited references, where only journals with greater than five articles have been counted}
\label{fig:journalArticleCounts}
\end{figure}
Furthermore, a table of the most cited articles related to the finite volume method for solid mechanics can be found in Appendix A.

\begin{figure}[htb]
   \centering
   \includegraphics[width=\textwidth]{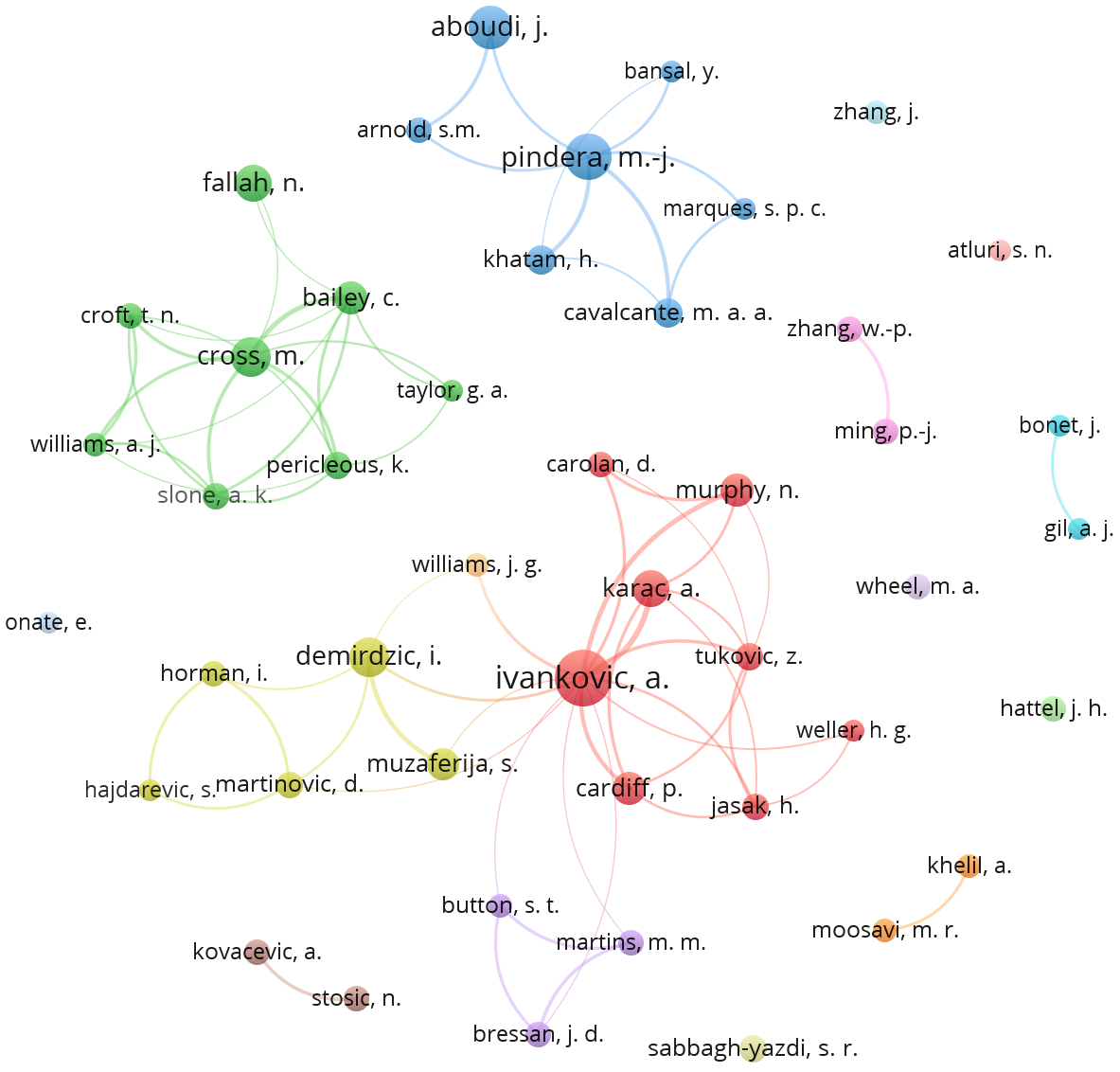} 
   \caption{Co-authorship network of the finite volume method for solid mechanics, generated using the VOSviewer software \citep{VanEck2007}. To maintain interpretability, only authors with five or more related peer-reviewed publications have been included in the network. The network provides insight into collaborations within the field but is not a direct measure of contribution.}
   \label{fig:co-authorshipMap}
\end{figure}

A perspective on the literature landscape is gained from the \emph{co-authorship network} presented in Figure \ref{fig:co-authorshipMap}, which has been generated using the VOSviewer software \citep{VanEck2007}.
A co-authorship network is a visual method to assess research collaborations within a field, as well as observe detached regions of research.
Referring to Figure \ref{fig:co-authorshipMap}, the three larger sub-networks broadly correspond to:
implicit cell-centred approaches indicated by the red sub-network centred on \emph{demird\v{z}i\'{c}} and \emph{ivankovi\'{c}};
implicit vertex-centred approaches indicated  by the green sub-network centred on \emph{cross}, \emph{bailey} and \emph{fallah};
and a specialised form of finite volume method for microstructural analysis indicated by the blue sub-network centred on \emph{aboudi}, \emph{pindera} and \emph{cavalcante}.
Although the co-authorship graph is not a direct measure of contribution to the field, it does provide insight into the collaboration between authors and subdivisions within the domain.
In the graph, individual authors are represented by filled circles, and a co-authored publication between two individual authors is symbolised by a line connecting the two filled circles;
the line thickness increases with increasing number of co-authored publications;
the size of an author's filled circle is directly proportional to the number of their related publications.
To maintain interpretability, only authors with five or more related publications have been included in the network, and only peer-reviewed publications have been incorporated.


Finally, before considering individual contributions, we will clarify the meaning of a number of terms in the current context to avoid any confusion.
For the domain spatial discretisation, this will be referred to interchangeably as the ``mesh'' or ``grid'';
each sub-domain in a finite volume mesh will be equivalently referred to as a ``control volume'', a ``cell'', or even an ``element'';
similarly, each sub-domain in the finite element mesh will be referred to as an ``element'' or ``cell'.
The term ``node'' indicates the mesh location where the solution variable is stored, which could be a vertex, face-centre or cell-centre depending on the method.

\subsection{Cell-centred approaches}
Thirty years ago, \citet{Demirdzic1988} proposed the first application of the cell-centred finite volume method in its modern form to solid mechanics.
Subsequently, cell-centred developments have primarily focussed on two relatively disconnected approaches:
\begin{itemize}
	\item \emph{Implicit} cell-centred approaches based on the original approach of \citet{Demirdzic1988}, and
	\item \emph{Explicit} Godunov-type cell-centred approaches stemming from the work of \citet{Trangenstein1991}.
\end{itemize}

The cell-centred approach takes its name from the dependent variable(s) residing at the cell centres (control volume centroids); equivalently, the approach has been termed the \emph{colocated},  \emph{co-located} or \emph{collocated} finite volume method, as the dependent variables share their location at the cell centres/centroids.

\begin{figure}[htb]
	\centering
	\subfigure[2-D structured quadrilateral mesh from \citet{Demirdzic1988}]
	{
		\includegraphics[width=0.45\textwidth]{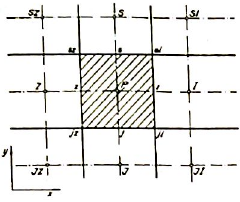}
	}
	\subfigure[3-D general convex polyhedral control volume from \citet{Demirdzic1994a}]
	{
		\includegraphics[width=0.45\textwidth]{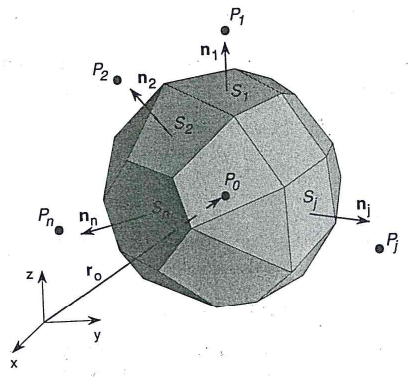}
	}
	\caption{Original 2-D structured quadrilateral mesh of \citet{Demirdzic1988} (left) and the subsequent generalisation to 3-D unstructured convex polyhedra \citep{Demirdzic1994a} (right)}
	\label{fig:controlVolumeFig1Demirdzic1994}
\end{figure}

\paragraph{Implicit cell-centred methods}
Considering first the implicit methods: 
in the original approach of \citet{Demirdzic1988}, a structured rectangular 2-D method was applied to the simulation of thermal deformations in welded workpieces (Fig. \ref{fig:controlVolumeFig1Demirdzic1994}(a)).
The displacement was assumed to vary linearly between computational nodes, and the small strain material behaviour was described by the \emph{Duhamel-Neumann} form of Hooke's law.
A distinguishing feature of the proposed solution algorithm was the partitioning of the surface force term into a compact-stencil \emph{implicit} term and a larger stencil \emph{explicit} term.
As a result, the linear momentum \emph{vector} equation was temporarily decoupled into three \emph{scalar} component equations that were independently solved, where outer fixed-point/Gauss-Seidel/Picard iterations provided the required coupling.
This form of solution methodology is termed a \emph{segregated} approach, as the governing conservation of linear momentum equation is segregated into three scalar equations during solution.
This style of implicit segregated solution algorithm was inspired by the methods adopted in similar CFD procedures,
where the restrictive computer memory sizes available at the time necessitated memory efficient procedures.

The original 2-D method was later generalised to 3-D convex polyhedral cells by \citet{Demirdzic1994a} (Fig. \ref{fig:controlVolumeFig1Demirdzic1994}(b)).
This form of the cell-centred approach has since been extended to deal with a wide variety of solid and multi-physics phenomena, including:
\begin{itemize}
	\item
	Elastoplasticity \citep{Demirdzic1993, Dioh1994, Rente2000, Maneeratana2000b, Teskeredzic2002, Teskeredzic2004, Basic2005, Martins2010, Bressan2010, Leonard2012, Bressan2015, Cardiff2016a, Martins2016, Martins2017, Cardiff2017a, Bressan2017}
	and viscoelasticity \citep{Zarrabi1999, Zarrabi2000, Dzaferovic2000, Demirdzic2005, Das2012, Safari2016};
	\item
	Thermo-elasticity \citep{Demirdzic1994b, Demirdzic1996a, Demirdzic1988, Demirdzic1994a, Osman2011, Junior2013, Bibin2013} and hygro-thermo-elasticity \citep{Martinovic1998, Demirdzic2000, Martinovic2000, Martinovic2001, Martinovic2002, Horman2008, Martinovic2008, Horman2009, Horman2010a, Horman2012, Fu2018, Holzmann2018};
	\item
	Poro-elasticity \citep{Tang2015, Bryant2015, Cardiff2015a, Cardiff2015b, Lee2015, Elsafti2016, Manchanda2017, Asadollahi2017};
	\item Anisotropy \citep{Fainberg1996, Martinovic1998, Horman1999, Demirdzic2000, Martinovic2001, Martinovic2008, Horman2009, Horman2010a, Horman2012, Cardiff2014a, Golubovic2017a} and
	heterogeneous material properties \citep{Tukovic2012, Carolan2013a, Cardiff2012d, Bryant2015, Manchanda2017}.
	\item Incompressibility and quasi-incompressibility \citep{Greenshields1999a, Greenshields1999b, Bijelonja2005a, Greenshields2005, Bijelonja2006, Giannopapa2006, Giannopapa2008, Bijelonja2017};
	\item Contact mechanics \citep{Jasak2000b, Cardiff2011c, Cardiff2011d, Cardiff2012a, Cardiff2014b, Cardiff2016a};
	\item Finite strains and rotations \citep{Maneeratana1999a, Maneeratana1999b, Maneeratana1999c, Maneeratana2000a, Maneeratana2000b, Basic2002, Bijelonja2002, Bijelonja2005a, Tukovic2007a, Tukovic2007b, Cardiff2012c, Cardiff2013, Cardiff2014a, Cardiff2014d, Cardiff2016a, Liu2018};
	\item Fracture mechanics \citep{Ivankovic1993, Ivankovic1994, Ivankovic1997a, Ivankovic1998, Ivankovic1999, Ivankovic2001a, Stylianou2002a, Stylianou2002b, Ivankovic2004a, Rager2005, Murphy2005, Murphy2006, Tropsa2006, Tukovic2010, Karac2011, Carolan2013a, Cardiff2015a, Lee2015, Bryant2015, Safari2016, Manchanda2017};
	\item Casting, melting, solidification and residual stresses \citep{Ivankovic1997b, Tropsa2000, Teskeredzic2002, Teskeredzic2004, Sato2006, Teskeredzic2015a, Teskeredzic2015b};
	\item Fluid-solid interaction \citep{Henry1993a, Henry1993b, Demirdzic1995, Greenshields1999a, Greenshields1999b, Greenshields2000, Ivankovic2001b, Schafer2001b, Schafer2002, Ivankovic2002b, Torlak2002a, Torlak2002b, Kovacevic2004b, Greenshields2005, Stosic2005, Shaw2005, Torlak2006, Giannopapa2006, Kovacevic2006, Kovacevic2007, Stosic2007, Giannopapa2008, Papadakis2008, Tukovic2007b, Jasak2007, Kanyanta2009a, Jagad2011, Wiedemair2012, Habchi2013, Tukovic2014, Smith2014, Sekutkovski2016, Jagad2016, Cardiff2017b, Jagad2017, Jagad2018};
	\item Beams, plates and shells \citep{Demirdzic1997b, Torlak2002a, Torlak2002b, Fallah2004, Fallah2006a, Fallah2006b, Fallah2006c, Fallah2006d, Hatami2006, Torlak2006, Isic2007a, Isic2007b, Isic2007c, Isic2008, Das2012, Fallah2013a, Fallah2014a, Fallah2014b, Jing2016, Fallah2017a, Fallah2017b, Golubovic2017a, Mohebi2017, Fallah2018b, Amraei2018, Tukovic2019};
	\item Solid-electrostatic interaction \citep{Das2011a, Das2011b} and wave propagation \citep{Dioh1994, Dioh1995a, Dioh1995b, Oosterkamp2000};
\end{itemize}


Of particular note are the developments of \citet{Weller1998} and \citet{Jasak2000a}, where the same implicit cell-centred form has been demonstrated to be well-suited to running on distributed-memory supercomputers.
The \emph{domain decomposition method} was used where the solution domain is decomposed into a number of sub-domains, each solved on a separate Central Processing Unit (CPU) core; necessary coupling between the sub-domains is performed using a message-passing protocol, initially Parallel Virtual Machine, but Message Passing Interface in later publications.
The \citet{Weller1998} and \citet{Jasak2000a} implementations form a component of the popular open-source C++ library OpenFOAM, formerly commercial software FOAM.
Many of the subsequent developments in the implicit cell-centred field have been based on the OpenFOAM platform, for example, \citep{Tukovic2007a, Tukovic2007b, Tukovic2010, Suvanjumrat2011, Cardiff2012a, Tukovic2012, Cardiff2014a, Tang2015, Cardiff2016a, Cardiff2016b, Elsafti2016, Holzmann2018}.
It should, however, be noted that the explicit Godunov-type approaches have also been implemented in OpenFOAM \citep{Haider2017, Haider2018}.

The essence of the implicit cell-centre finite volume method has been the use of a displacement approach combined with a segregated algorithm;
however, alternative solution methodologies have also been developed, including geometric multi-grid procedures \citep{Fainberg1996, Demirdzic1997c, Ivankovic1997a, Schafer2001b, Schafer2002}, block-coupled algorithms \citep{Das2011a, Das2011b, Cardiff2014e, Cardiff2016b, Gonzalez2018}, hybrid/mixed pressure-displacement formulations \citep{Henry1993a, Henry1993b, Basic2005, Fowler2003, Bijelonja2005a, Bijelonja2006, Bijelonja2011a, Bijelonja2017}, Aitken acceleration \citep{Tukovic2014, Tang2015, Gonzalez2018}, and curvilinear formulations \citep{Oliveira1999, Rente2000}.
In addition, fourth-order accuracy variants have been proposed \citep{Demirdzic2016} as well as novel gradient and tractions calculation methods \citep{Fallah2008a, Jagad2011, Tukovic2014, Nordbotten2014, Nordbotten2015, Cardiff2016a, Keilegavlen2017, Tukovic2018a, Gonzalez2018}.
Apart from the standard \emph{continuum} approaches, a number of authors have proposed implicit cell-centred finite volume methods for beams, plates and shells \citep{Demirdzic1997b, Torlak2002a, Torlak2002b, Hatami2006, Torlak2006, Isic2007a, Fallah2008b, Das2012, Jing2016, Golubovic2017a, Golubovic2017b, Fallah2018b, Amraei2018}.

\begin{figure}[htb]
	\centering
	\subfigure[2-D unstructured cell-centred polygonal mesh from \citet{Kluth2010}]
	{
		\includegraphics[height=0.35\textwidth]{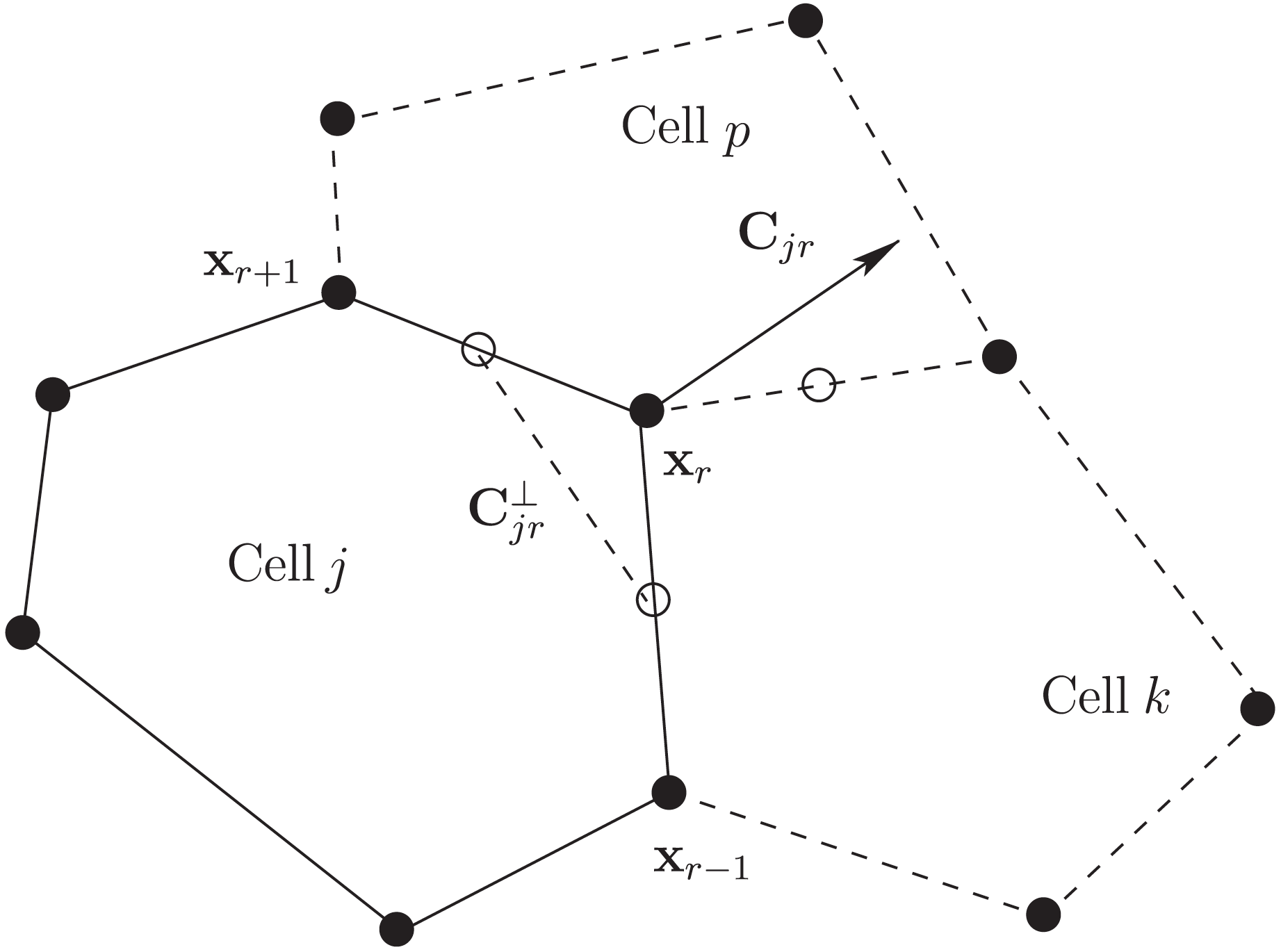}
	} \quad\quad
	\subfigure[2-D/3-D unstructured vertex-centred polygonal/polyhedral mesh from \citet{Aguirre2015}]
	{
		\includegraphics[height=0.35\textwidth]{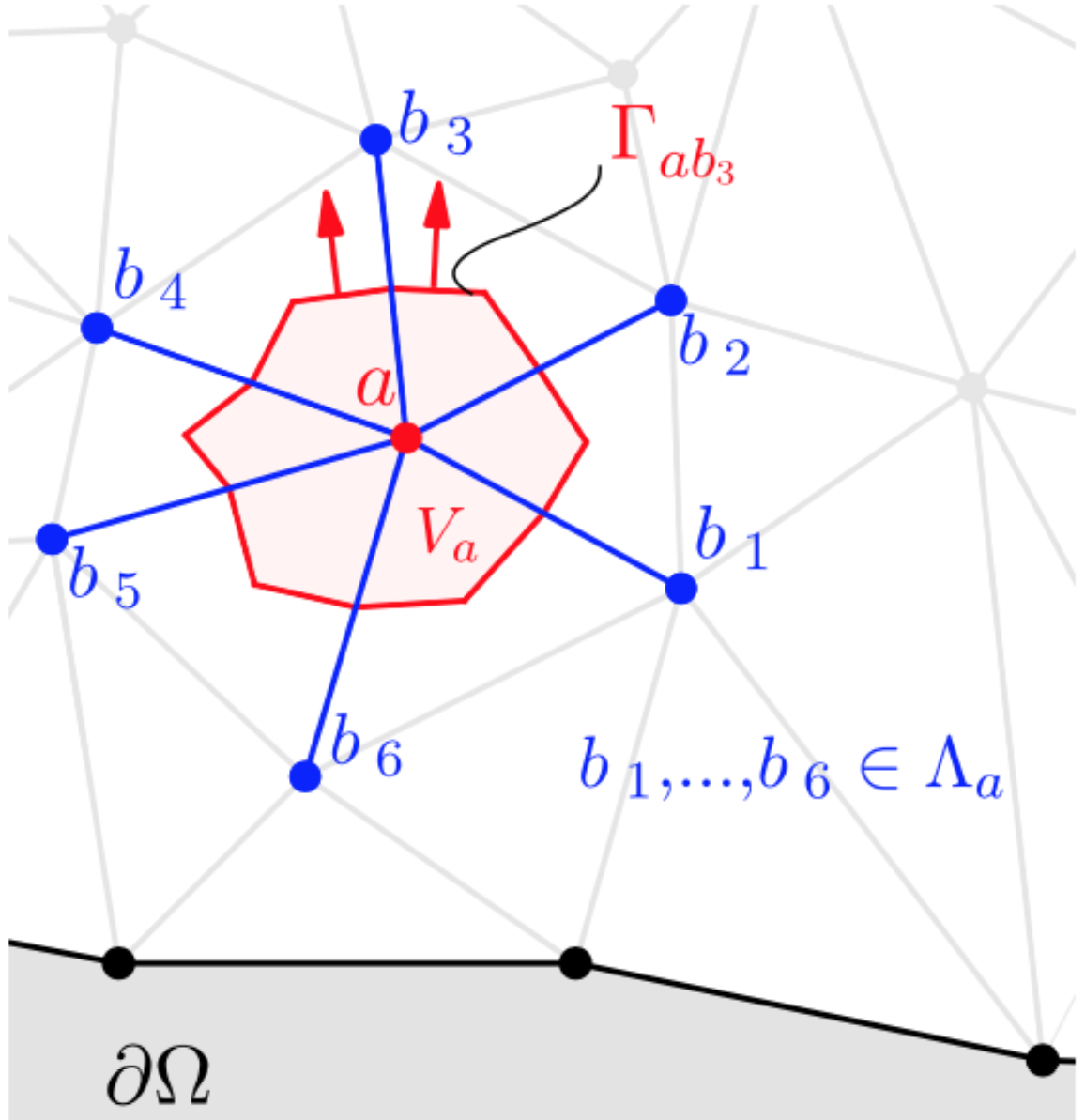}
	}
	\caption{Forms of mesh employed by \citet{Kluth2010} and \citet{Aguirre2015}}
	\label{fig:godunovMeshes}
\end{figure}

\paragraph{Explicit cell-centred methods} \label{sec:GodunovHistory}
Explicit Godunov-type finite volume methods were first proposed for the solution of hyperbolic problems characterised by waves and shocks \citep{Godunov1959, Godunov1962, LeVeque2004}, and have been popularised for the solution of Euler compressible gas flow equations.
The typical approach, which casts the conservation laws as a system of first-order hyperbolic equations, is characterised by the solution of a Riemann problem (propagation of a solution discontinuity) at the control volume faces to determine forces.
The resulting discretisation evaluates the force at a control volume face as a weighted average of the force evaluated at each side of the face.
%
Godunov-type methods were first applied to structural problems by Trangenstein and Colella, when they modelled the 1-D propagation of waves in elasto-plastic solids \citep{Trangenstein1991, Trangenstein1992, Trangenstein1994};
in their approach, the primitive conservation variables were the linear momentum vector and the deformation gradient (or displacement gradient) tensor.
Subsequently, the method has been extended to unstructured 3-D grids in a variety of forms, differing in terms of discretisations and primitive variables \citep{Miller1996, Tang1999, Berezovski2001, Howell2002, Berezovski2003, LeVeque2004, Kluth2008, Carre2009, Kluth2010, Lee2013, Maire2013, Sambasivan2013, Sijoy2015, Despres2015, Ndanou2015, Cheng2015, Loubere2016, Boscheri2016, Vilar2016, Haider2017, Georges2017, Hueze2017, Cheng2017a, Cheng2017b, Fridrich2017, Heuze2018, Sevilla2018a, Sevilla2018b, Haider2018}, for example, see Figure \ref{fig:godunovMeshes}.
Although cell-centred formulations are the most common form of Godunov-type method, vertex-centred \citep{Aguirre2014, Aguirre2015} and, recently, face-centred approaches \citep{Sevilla2018a, Sevilla2018b} have also been explored.

To-date, Godunov-type approaches have been used to model a wide range of physical mechanisms:
\begin{itemize}
	\item Linear elasticity \citep{Trangenstein1991, Trangenstein1992, Hueze2017, Sevilla2018b};
	\item Material nonlinearity \citep{Tang1999, Howell2002, Kluth2010, Maire2013, Sambasivan2013, Sijoy2015, Ndanou2015, Cheng2015, Loubere2016, Vilar2016, Boscheri2016, Haider2017, Georges2017, Hueze2017, Cheng2017a, Cheng2017b, Fridrich2017, Heuze2018, Haider2018};
	\item Fracture and cavitation \citep{Tang1999, Ndanou2015};
	\item Finite strains \citep{Maire2013, Sambasivan2013, Sijoy2015, Ndanou2015, Cheng2015, Georges2017, Boscheri2016, Loubere2016, Haider2017, Cheng2017a, Cheng2017b, Fridrich2017, Heuze2018, Haider2018};
	\item Material heterogeneity \citep{Berezovski2001, Berezovski2003};
	\item Wave propagation and impacts \citep{Trangenstein1991, Trangenstein1992, Trangenstein1994, Miller1996, Tang1999, Berezovski2001, Howell2002, Berezovski2003, LeVeque2004, Maire2013, Sambasivan2013, Despres2015, Sijoy2015, Ndanou2015, Cheng2015, Loubere2016, Boscheri2016, Georges2017, Hueze2017, Cheng2017a, Cheng2017b, Fridrich2017, Heuze2018}.
\end{itemize}

A distinctive characteristic of Godunov-type methods is the adoption of fully explicit solution algorithms, where the time increment size is restricted by the standard Courant-Friedrichs-Lewy constraint \citep{Courant1928}.

From an implementation and philosophy perspective, the implicit approaches stemming from \citet{Demirdzic1988} differ greatly from the explicit approaches deriving from \citet{Trangenstein1991}.
The most obvious difference is the need of the implicit approach to form and solve a linear system of equations, and the time increment size restriction for explicit methods.
In addition to this, the implicit methods have typically assumed a smooth variation of the primitive variable within each cell (or across cell faces); in contrast, the explicit approaches have aimed to directly approximate the propagation of solution discontinuities.
Many of these technical distinctions are discussed further in Section \ref{sec:compareFiniteVolumeVariants}.
It should be noted that \emph{explicit} cell-centred approaches that do \emph{not} adopt a Godunov approach are also possible, for example, as developed by Selim \etal \citep{Selim2016, Selim2017}.

\subsection{Vertex-centred approaches}
\begin{figure}[htb]
	\centering
	\subfigure[2-D unstructured quadrilateral and trianglular mesh from \citet{Fryer1991}]
	{
		\includegraphics[width=0.45\textwidth]{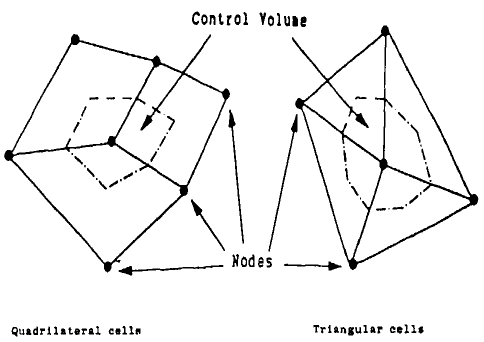}
	}
	\subfigure[3-D unstructured mesh from \citet{Taylor2003}]
	{
		\includegraphics[width=0.45\textwidth]{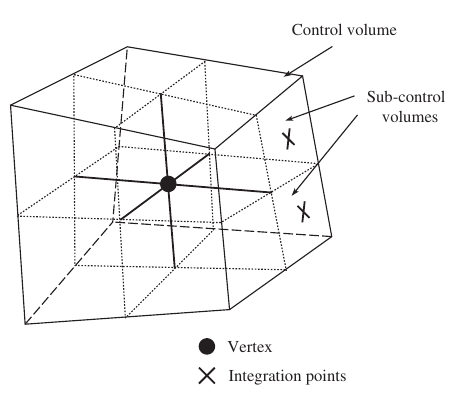}
	}
	\caption{Vertex-based finite volume formulation of \citet{Fryer1991} (left) showing the construction of control volumes around a mesh vertex for quadrilateral and triangular 2-D cells, and of \citet{Taylor2003} showing the extension to 3-D cells.}
	\label{fig:vertexBasedFryerBailey}
\end{figure}

\paragraph{Implicit vertex-centred methods}
\citet{Fryer1991} was the first to propose a vertex-centred finite volume method for solid mechanics.
The approach, initially termed a \emph{control volume-unstructured mesh} procedure, could analyse complex 2-D geometry using both quadrilateral and triangular cells (Fig. \ref{fig:vertexBasedFryerBailey}(a)).
The method of \citet{Fryer1991} followed closely the approach of \citet{Baliga1980} from a decade earlier, who developed a so-called \emph{control-volume-based finite element method} for convection-diffusion equations on unstructured triangular grids.

In addition to the designation \emph{vertex-based finite volume method}, the formulation is referred to by a number of other names, including \emph{vertex-centred finite volume method}, \emph{cell-vertex finite volume method}, \emph{control volume procedure}, \emph{control-volume finite element method}, \emph{control-volume-based finite-element method}, and \emph{element-based finite volume method}.

In contrast to cell-centred and staggered-grid approaches, vertex-centred approaches store the primary unknowns at the \emph{primary} mesh vertices and integrate the governing equation over secondary/\emph{dual} mesh control volumes surrounding each \emph{primary} mesh vertex.

The original 2-D approach of \citet{Fryer1991} was subsequently extended to three dimensions by \citet{Bailey1995} (Fig. \ref{fig:vertexBasedFryerBailey}(b)) and has since been applied to a wide range of physical phenomena:
\begin{itemize}
	\item Compressible and incompressible elasticity \citep{Bailey1995, Wheel1996a, Wheel1996b, Wheel1998, Wheel1999a, Costa2000, Fallah2005, Xia2007, Tsui2013, Wu2013};
	\item Elasto-plasticity, elasto-visco-plasticity, creep and material nonlinearity \citep{Taylor1995a, Ferguson1998, Taylor2003, Hambleton2011, Fallah2014a, Fallah2017b};
	\item Finite strains and geometric nonlinearity \citep{Fallah1999, Fallah2000a, Fallah2000b, Fallah2000c, Slone2002b, Teran2003, Limache2007};
	\item Beams, plates and shells  \citep{Wheel1997, Beveridge2009, Fallah2004};
	\item Multi-physics and fluid-solid interaction \citep{Taylor1996, Bailey1996, Bailey1997, Bailey1999a, Bailey1999b, Oldroyd1999, Wheel1999b, Slone2000, Slone2001, Slone2002, Slone2003, Slone2004, Slone2007, Cross2007, Lv2007, Croft2008, Xia2008, Tsui2013, Hejranfar2016};
	\item Casting \citep{Bailey1991, Cross1992, Cross1993, Fryer1993, Taylor1995b, Taylor1998, Cross1996, Bailey1996, Bailey1997, Bailey1999a, Bailey1999b, Bounds2000},
	extrusion and forging \citep{Williams2001, Williams2002, Williams2010};
	\item Welding \citep{Taylor2000a, Taylor2002a, Taylor2002b};
	\item Contact mechanics \citep{Taylor2000b};
	\item Functionally graded solids \citep{Gong2013, Gong2014} and wood drying \citep{Perre1995, Salinas2011};
	\item Micropolar/Cosserat elasticity (shear stress may not be symmetric) \citep{Wheel2008, Beveridge2013};
	\item Vibrations, acoustics and wave propagation \citep{Zhu2012, Fallah2013a, Xuan2014a, Xuan2014b}.
\end{itemize}
Aalternative solution methodologies have been examined, including mixed displacement-rotation approaches \citep{Wenke2003, Pan2010, Pan2011}, and addressing parallelisation on distributed memory supercomputers \citep{McManus2000, McManus2002}.
Also of note are the developments of \citet{Maitre2002} and \citet{Souhail2004}, who proposed a higher order extension of the vertex-based approach.

It is important to note the distinction between \emph{overlapping} and \emph{non-overlapping} vertex-centred methods:
the more common approach involves non-overlapping control volumes, where the control volumes about each primary mesh vertex do not overlap with the control volumes of neighbouring primary mesh vertices. In the less popular overlapping control volumes approaches, for example, as discussed by \citet{Onate1994}, primary mesh vertices partly share control volumes with primary mesh neighbouring vertices, resulting in the governing equations being integrated more than once over regions of the domain;
this overlapping local integration domain approach is more similar to the finite element method.
Additionally, \citet{Tsui2013} proposed a variation on the common approach to construct the control volumes:
the primary mesh cell centres are joined together to make the control volumes, rather than joining the cell-centres to the face-centres as per the classic vertex-based method.




\paragraph{Explicit vertex-centred methods}
Similar to cell-centred methods, both implicit and explicit solution algorithms have been developed, although explicit vertex-centred methods have seen less development.
Within the category of explicit vertex-centred methods, there exists a variety of sub-classes, including:
so-called dual-time-stepping explicit methods, where the solution is calculated in an explicit manner and a linear system is not directly formed \citep{Xia2007, Xia2008, Lv2007, Zhu2012, Tsui2013};
Godunov-type approaches \citep{Aguirre2014, Aguirre2015} similar to those seen in cell-centred approaches, as noted in Section \ref{sec:GodunovHistory};
and the so-called \emph{grid method}  \citep{Zhang1999, Zhang2002, Zhang2004, Liu2004, Liu2005, Gao2006}.
Regarding the grid-method, the originators of the method,
\citet{Zhang1999}, argue for the distinction between the grid method and the vertex-based finite volume method, however, the authors of the current article believe this distinction is unwarranted:
like other finite volume methods, the grid method starts from the governing momentum equation in strong integral form and approximates the forces over the boundaries of the control volumes; although there are minor differences in the techniques used to approximate the surface forces, for example, comparing \citet{Dormy1995} and \citet{Zhang2002}, the \emph{grid method} still remains a form of vertex-centred finite volume method.

\subsection{Staggered-grid approaches}
In staggered-grid approaches, as originally proposed for CFD by Harlow and Welch \citep{Harlow1965},
the components of the primary solution variable, for example, the $x$ and $y$ components of displacement, are stored at different locations.
In addition, different sets of control volumes are used when integrating the discretised governing equation in each Cartesian direction.
For example, the $x$-momentum component equation employs a different grid to the $y$-momentum component equation.
The primary motivation for such staggered-grid approaches is the avoidance of solution instabilities, namely the ``checker-boarding'' phenomenon, where high frequency variations appear in the solution variables that are unobservable to the discretisation.
Consequently, staggered-grid approaches do not need to explicitly included a stabilisation term, as is common in cell-centred methods.
The extension of staggered-grid approaches to general unstructured 3-D meshes is, however, not trivial; and consequently, this major limitation has resulted in declining popularity for such approaches in solid mechanics.

The first application of a staggered-grid formulation to solid mechanics was by Beale and Elias \citep{Beale1990a, Beale1990b}, who showed how the finite volume CFD code {PHOENICS} could be applied to stress analysis problems.
In their method, Beale and Elias used the analogy between the stream functions in creeping-fluid-flow and the Airy stress functions in solid mechanics, where the stress components were the primary unknowns.
This method was later generalised by Spalding, Bukhari, Qin, Hamill and co-workers \citep{Bukhari1990, Bukhari1991, Spalding1993, Spalding1997, Spalding1998a, Spalding1998b, Spalding2002, Spalding2004, Spalding2008}, where a displacement, rather than a force analogy, was adopted, having the major advantage of being more general in three dimensions.
Spalding noted that a CFD solution procedure designed for computing velocities is suitable for computing displacements if the convection terms are set to zero and the volume/dilatation stress term is introduced by inclusion of a pressure- and temperature-dependent source term.
As a consequence, the resulting solution method used an \emph{implicit} SIMPLE algorithm \citep{Patankar1972, Patankar1980}, still popular in modern CFD codes.
The employed primitive variables were displacement and pressure, as opposed to velocity and pressure in standard fluid flow analyses.

The use of mixed displacement-pressure approaches, however, is not a requirement of staggered-grid approaches.
Hattel, Hansen and collaborators \citep{Hattel1990, Hattel1992, Hattel1993a, Hattel1993b, Hattel1993c, Hattel1993d, Hattel1993e, Hattel1994, Hattel1995, Hattel1997, Pryds1997, Hattel1998, Hattel2001, Thorborg2001, Hattel2003, Thorborg2003} demonstrated a staggered-grid approach where the sole primary variable was displacement.
\citet{Hattel1990} initially proposed their staggered-grid approach for the analysis of thermally induced stresses in casting problems, and subsequently extended it for a variety of thermo-elasto-plasticity problems.
They termed their approach ``a control volume based finite difference method'', further indicating the close-relationship with finite difference methods.
The authors noted that their method resulted in an elegant formulation for non-constant material properties, a benefit of the staggered grid approach.


In addition to the displacement and displacement-pressure approaches, a number of alternative staggered grid formulations have also been proposed:
\citet{Spalding2002} proposed a staggered formulation where \emph{rotation} and displacement were the primitive variables; Spalding surmised that the rotation-based method may provide a more efficient solution algorithm in certain situations.
The approach is described in documentation from an early version of the {PHOENICS} software; however, issues with boundary conditions are noted and no further articles appear on the formulation. A similar idea was subsequently proposed by Wenke, Wheel, Pan and Qin \citep{Wenke2003, Pan2010, Pan2011} in the framework of vertex-based finite volume methods, where both displacements and rotations are considered the primitive variables.
More recently, \citet{Wang2007} put forward a staggered finite volume approach for analysis of shape memory alloys, defined on 2-D structured rectangular grids.
In contrast to the majority of staggered grid approaches, which employ implicit solution algorithms, \citet{Rajagopal2014} proposed a one-step explicit staggered grid finite volume approach to investigate the response of a layered viscoelastic plate.
\subsection{Other approaches}
Across the spectrum of finite volume methods for solid mechanics, not all procedures align with the previously discussed divisions.

\paragraph{Approaches for periodic heterogenous microstructures}
Although arguably not a fundamental class of finite volume method, there has been significant development related to specialised versions of the finite volume method for periodic microstructures.
Depending on the formulation, the related methods can be referred to by a number of names, including: the \emph{higher-order theory for functionally graded material} (HOTFGM) \citep{Aboudi1999}, the \emph{high-fidelity generalised method of cells} (HFGMC) \citep{Aboudi2001a, Aboudi2001b, Haj-Ali2009, Haj-Ali2012}, the \emph{finite volume direct averaging micromechanics} (FVDAM) theory \citep{Bansal2005, Bansal2006}, or some variant thereof.
Although there is some disagreement over the naming convention \citep{Haj-Ali2012, Bansal2006, Haj-Ali2009}, these methods have a common origin in the \emph{method of cells} and the \emph{generalised method of cells} developed by Aboudi and Paley \citep{Aboudi1982, Aboudi1991, Paley1992}.
Recent discussions around the development of these methods can be found in \citet{Haj-Ali2012}, \citet{Cavalcante2012c}, \citet{Gong2013}, and \citet{Cavalcante2016}, along with reviews of the high fidelity generalised method of cells approaches by \citet{Aboudi2004} and of microstructural analysis approaches by \citet{Pindera2009} and \citet{Charalambakis2006}. Reviews of the application of such methods can be found in \citet{Aboudi2007} and \citet{Aboudi2008}.
A brief overview of the discretisation used in HOTFGM/HFGMC/FVDAM approaches is given in Appendix B.

\paragraph{Meshless finite volume approaches}
Meshless methods have been proposed in order to overcome the drawbacks of mesh-based finite element and finite volume methods, particularly related to large deformations and cracking.
\citet{Atluri2002} were the first to propose a meshless finite volume formulation, based on the earlier generalised \emph{Meshless Local Petrov-Galerkin (MLPG) method} by \citet{Atluri1998}.
The approach has subsequently been extended and applied to a variety of problems in elasto-statics, elasto-dynamics and fracture mechanics \citep{Ching2001, Warlock2002, Qian2003, Raju2003, Atluri2004, Han2004, Batra2004, Han2005, Sladek2008, Moosavi2008, Moosavi2009, Moosavi2011a, Moosavi2011b, Hosseini2011, Soares2012, Moosavi2012a, Moosavi2012b, Moosavi2013}.
Based on the initial development of Atluri and co-workers, \citet{Moosavi2008} proposed a novel meshless form of the finite volume method for elasto-static analysis, that combined the finite volume concept with the meshless local Petrov-Galerkin approach with moving least squares interpolation.
Moosavi and co-workers have since extended the method to elasto-dynamics, beams, plates, shells and crack problems \citep{Moosavi2009, Moosavi2011a, Moosavi2011b, Moosavi2012b, Moosavi2013}, where the method has been named the \emph{orthogonal meshless finite volume method}.
Whereas the \citet{Atluri2002} approach uses overlapping control volumes, the meshless finite volume method developed by Ebrahimnejad \etal \citep{Ebrahimnejad2014, Ebrahimnejad2015, Ebrahimnejad2017} employs non-overlapping control volumes. The method has been applied to 2-D elasticity with adaptive mesh refinement, and was later extended to problems with material discontinuities \citep{Fallah2018a, Fallah2018c},  to free vibration analysis of laminated composite plates \citep{Davoudi-Kia2017, Davoudi-Kia2018}, and to fractures in orthotropic media \citep{Fallah2018d}.

\paragraph{Eulerian approaches}
There are many solid mechanics problems which can be equivalently considered from the fluid mechanics perspective, for example, in the analysis of extrusion and drawing.
Eulerian ``fluid'' approaches are often used to solve such problems, for example, \citep{Basic2001, Williams2002, Basic2005, Basic2008, Basic2009, Chen2007, Jafari2007, Lou2008, Bressan2010, Williams2010, Al-Athel2011, Wang2011, Wang2012, Bressan2013, Bressan2015, Zhang2016, DeBrauer2016, Martins2016, Bressan2017, DeBrauer2017, Martins2017}.
These methods are more closely related to CFD procedures and are not discussed further here.

\paragraph{Miscellaneous}
\citet{Teng1999} and \citet{Chen2001} developed a finite volume method to simulate the draping of woven fabrics, where the governing nonlinear equations were solved using a single-step full Newton-Raphson method.
Martin and Pascal \citep{Martin2011, Martin2012} proposed a novel \emph{discrete duality finite volume} method for solving linear elasticity problems on unstructured meshes; the main characteristic of the discretisation is the integration of the governing equations over two meshes: the given primal mesh and also over a dual mesh built from the primal one.
\citet{Pietro2011} proposed a novel discretisation scheme for linear elasticity with only one degree of freedom per control-volume face, corresponding to the normal component of the displacement.



\section{Comparing variants of the finite volume method for computational solid mechanics}
\label{sec:compareFiniteVolumeVariants}

%

The finite volume method, like other related numerical approaches, consists of the following main components:
\begin{itemize}
	\item[a)] Discretisation of space and time;
	\item[b)] Discretisation of the mathematical model equations;
	\item[c)] Solution algorithm.
\end{itemize}
To facilitate comparison between the major variants, three popular formulations will be considered here:
\begin{itemize}
	\item Implicit cell-centred approach originating from the work of \citet{Demirdzic1988};
	\item Implicit vertex-centred approach of \citet{Fryer1991} and \citet{Bailey1995};
	\item Explicit cell-centred approach emanating from \citet{Trangenstein1991}.
\end{itemize}
Analysis and insight into the similarities and differences between these variants is provided.

\subsection{Mathematical model for dynamic linear elasticity} \label{sec:mathModel}
\begin{figure}[htb]
	\centering
	\includegraphics[width=0.4\textwidth]{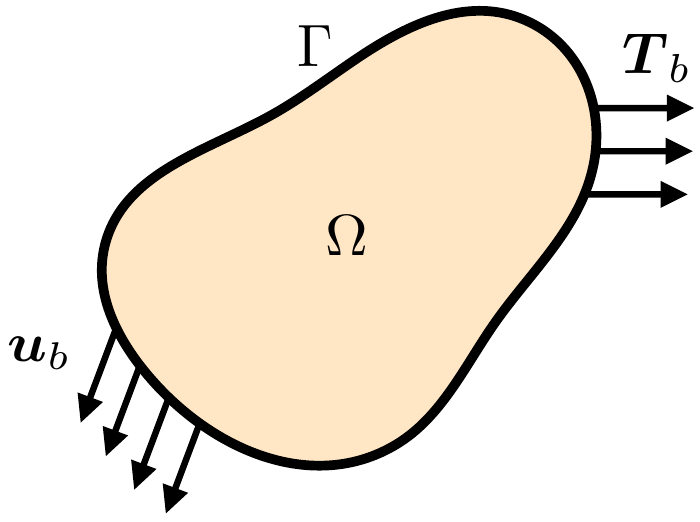}
	\caption{Generic solid body, with volume $\Omega$ and surface $\Gamma$, subjected to boundary displacements, $\boldsymbol{u}_b$, and boundary tractions, $\boldsymbol{T}_b$}
	\label{fig:body}
\end{figure}

To allow a clear comparison of the methods, the dynamic behaviour of a body with volume $\Omega$ and surface $\Gamma$ is analysed (Fig. \ref{fig:body}), where part of its boundary is subjected to a specified displacement, $\boldsymbol{u}_b$, and the remainder is subjected to a specified traction, $\boldsymbol{T}_b$.

Assuming the relationship between stress and strain to be described by Hooke's law, the governing conservation of linear momentum equation, a generalisation of Newton's second law, can be given in strong integral form as:
\begin{eqnarray} \label{eq:momentumIntegral}
	\underbrace
	{
	\int_\Omega \rho \frac{\partial^2 \boldsymbol{u}}{\partial t^2}\;\text{d}\Omega
	}_{\text{Inertia forces}}
	&=&
	\underbrace
	{
	\oint_\Gamma \boldsymbol{n} \cdot
	\left[
	\mu \boldsymbol{\nabla} \boldsymbol{u}
	+ \mu (\boldsymbol{\nabla} \boldsymbol{u})^T
	+ \lambda \; \text{tr}(\boldsymbol{\nabla} \boldsymbol{u}) \textbf{I}
	\right]
	\;\text{d}\Gamma
	}_{\text{Surface forces}}
	\;+\;
	\underbrace
	{
	\int_\Omega \rho \boldsymbol{f}_b \;\text{d}\Omega
	}_{\text{Body forces}}
\end{eqnarray}
Small deformations are assumed \ie no distinction is made between the initial and deformed configurations; $\rho$ is the initial density, $\boldsymbol{u}$ is the unknown total displacement vector, $\boldsymbol{n}$ is the outward-pointing surface unit normal vector, $\boldsymbol{\nabla}$ is the del operator, $\lambda$ is the first Lam\'{e} parameter, $\mu$ is the second Lam\'{e} parameter, synonymous with the shear modulus, $\textbf{I}$ is the second-order identity tensor, and $\boldsymbol{f}_b$ is a body force acceleration.
To avoid unnecessary complexity, the material properties ($\rho$, $\mu$ and $\lambda$) are assumed to be uniform and isotropic; however, this assumption is not required.
It should be noted that the finite volume method directly discretises this strong integral form of the governing equation, without requiring weighting functions, the weak form of the equation or the use of the Gauss divergence theorem.

\subsection{Implicit cell-centred approach} \label{sec:standardFiniteVolume}
In this sub-section, the implicit cell-centred approach stemming from \citet{Demirdzic1988} is described.

\paragraph{Discretisation of time}
For all described variants of the finite volume method, discretisation of the solution domain comprises time discretisation and space discretisation.
The total specified simulation time is divided into a finite number of time increments, $\Delta t$, and the discretised governing equations are solved in a time-marching manner.

\begin{figure}[htb]
   \centering
   \includegraphics[height=0.35\textwidth]{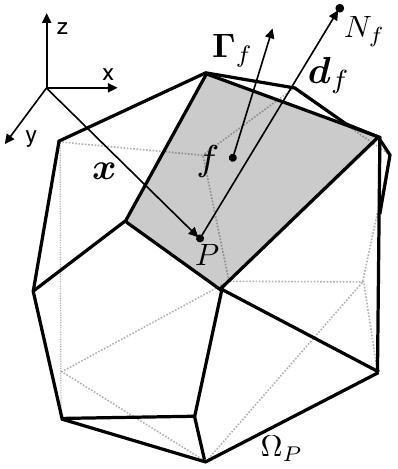}
   \caption{General convex polyhedral control volume (adapted from \citep{Cardiff2016a, Demirdzic1995, Jasak2007})}
   \label{fig:generalFiniteVolumeCell}
\end{figure}

\paragraph{Discretisation of space} \label{sec:standardFiniteVolumeDiscretisationOfSpace}
For the implicit cell-centred approach, the spatial domain is divided into a finite number of contiguous convex polyhedral cells bounded by polygonal faces that do not overlap and fill the space completely. 
A typical control volume is shown in Figure \ref{fig:generalFiniteVolumeCell}, with the computational node $P$ located at the cell centre/centroid, and the cell volume is $\Omega_P$;
$N_{f}$ is the centroid of a neighbouring control volume, which shares face $f$ with the current control volume;
$\boldsymbol{\Gamma}_{f}$ is the area vector of face $f$, vector $\boldsymbol{d}_{f}$ joins $P$ to $N_f$, and $\boldsymbol{x}$ is a positional vector.
No distinction is made between different cell volume shapes, as all general convex polyhedra are discretised in the same fashion.


\paragraph{Discretisation of the mathematical model equations} \label{sec:discretisationFV}
The governing conservation law (Equation \ref{eq:momentumIntegral}) is applied to each polyhedral cell in the mesh.
The discretisation of each of the three terms (inertia, surface forces, body forces) within the equation is now discussed in turn.

Considering first the spatial discretisation of the inertia term in Equation \ref{eq:momentumIntegral}.
The term can be approximated by making an assumption about the variation of the displacement vector $\boldsymbol{u}$ over the cell; typically in the implicit cell-centred approach a linear variation is assumed.
This linear variation can be expressed in terms of a truncated Taylor series expansion about the cell centre:
\begin{eqnarray} \label{eq:taylorSeries}
	\boldsymbol{u}(\boldsymbol{x}) =
	\boldsymbol{u}_P + (\boldsymbol{x} - \boldsymbol{x}_P) \cdot (\boldsymbol{\nabla} \boldsymbol{u})_P 
\end{eqnarray}
This expression says that the displacement $\boldsymbol{u}(\boldsymbol{x})$ at any point in the cell can be calculated using the cell-centre displacement $\boldsymbol{u}_P$ and the constant gradient of displacement within the cell $(\boldsymbol{\nabla} \boldsymbol{u})_P$.
This approximation is second-order order in space \ie as the mesh spacing is reduced, the error in this approximation reduces at second-order rate.
In principle, any other distribution could be used; for example, \citet{Demirdzic2016} extended the approach to use a fourth-order cubic distribution.

Using the approximation in Equation \ref{eq:taylorSeries}, the inertia term can be calculated using the midpoint rule as a function of the cell-centred values:
\begin{eqnarray}
	\int_\Omega \rho \frac{\partial^2 \boldsymbol{u}}{\partial t^2}\;\text{d}\Omega &\approx&
	\rho \left(\frac{\partial^2 \boldsymbol{u}}{\partial t^2}\right)_P \Omega_P
\end{eqnarray}
where the subscript $P$ on the density has been dropped as a uniform density field is assumed \ie $\rho_P = \rho$.

The acceleration $\partial^2 \boldsymbol{u}/\partial t^2$ may be discretised in time using any appropriate finite difference scheme; in the original cell-centred approach of \citet{Demirdzic1988}, the bounded first-order backward Euler method was used:
\begin{eqnarray} \label{eq:discretiseAcceleration}
	\left(\frac{\partial^2 \boldsymbol{u}}{\partial t^2}\right)_P &\approx&
	\frac{\boldsymbol{u}_P - 2\boldsymbol{u}_P^{[m-1]} + \boldsymbol{u}_P^{[m-2]}}{\Delta t^2}
\end{eqnarray}
where $m$ is the time-step counter;
the superscript on the unknown current time value of displacement has been dropped for brevity.
Here it is assumed that the time-step size $\Delta t$ is constant; however, the method can be generalised to variable time-step sizes.
There are numerous other temporal discretisations that may be used; for example, the unbounded second-order backward Euler scheme \citep{Jasak2000a}, the second-order trapezoidal rule \citep{Gonzalez2018}, or the second-order Newmark schemes that are popular with the finite element method.

The final discretised inertia term is:
\begin{eqnarray} \label{eq:inertiaTermDiscretisedCellCentred}
	\int_\Omega \rho \frac{\partial^2 \boldsymbol{u}}{\partial t^2}\;\text{d}\Omega
	&\approx&
	\rho \; \frac{\boldsymbol{u}_P - 2\boldsymbol{u}_P^{[m-1]} + \boldsymbol{u}_P^{[m-2]}}{\Delta t^2}  \; \Omega_P
\end{eqnarray}

In a similar fashion, the body force term (second term on the right-hand side of Equation \ref{eq:momentumIntegral}) is approximated by assuming $\boldsymbol{f}_b$ to vary over the cell according to Equation \ref{eq:taylorSeries}.
Consequently, the term can be approximated in terms of the cell-centre values using the midpoint rule as:
\begin{eqnarray}
	\int_\Omega \rho \boldsymbol{f}_b \;\text{d}\Omega &\approx& \rho \boldsymbol{f}_{b_P} \Omega_P
\end{eqnarray}

Towards the discretisation of the surface force term, the closed surface integral in Equation \ref{eq:momentumIntegral} is converted into a sum of surface integrals over each polygonal face:
\begin{eqnarray} \label{eq:momentumIntegralSurfaceSum}
	\oint_\Gamma \boldsymbol{n} \cdot
	\left[
	\mu \boldsymbol{\nabla} \boldsymbol{u}
	+ \mu (\boldsymbol{\nabla} \boldsymbol{u})^T
	+ \lambda \; \text{tr}(\boldsymbol{\nabla} \boldsymbol{u}) \textbf{I}
	\right]
	\;\text{d}\Gamma
	&=&
	\sum^{nFaces}_{f=1}
	\int_{\Gamma_f}
	\boldsymbol{n} \cdot \left[
	\mu \boldsymbol{\nabla} \boldsymbol{u}
	+ \mu (\boldsymbol{\nabla} \boldsymbol{u})^T
	+ \lambda \; \text{tr}(\boldsymbol{\nabla} \boldsymbol{u}) \textbf{I}
	\right]
	\;\text{d}\Gamma \notag \\
\end{eqnarray}
where $nFaces$ is the number of faces in the cell.
This form represents a force balance for the cell.
To approximate the stress term $[\mu \boldsymbol{\nabla} \boldsymbol{u} 	+ \mu (\boldsymbol{\nabla} \boldsymbol{u})^T + \lambda \; \text{tr}(\boldsymbol{\nabla} \boldsymbol{u}) \textbf{I}]$ at each cell face, the assumed variation of displacement in Equation \ref{eq:taylorSeries} is once again used;
accordingly, the stress can be calculated in terms of the displacement gradient values at the face centre (centroid):
\begin{eqnarray} \label{eq:momentumIntegralSurfaceSum2}
	\sum^{nFaces}_{f=1}
	\int_{\Gamma_f}
	\boldsymbol{n} \cdot \left[
	\mu \boldsymbol{\nabla} \boldsymbol{u}
	+ \mu (\boldsymbol{\nabla} \boldsymbol{u})^T
	+ \lambda \; \text{tr}(\boldsymbol{\nabla} \boldsymbol{u}) \textbf{I}
	\right]
	\;\text{d}\Gamma
	\approx
	\quad\quad\quad\quad\quad\quad\quad\quad
	\notag \\
	\sum^{nFaces}_{f=1}
	\boldsymbol{n}_{f} \cdot \left[
	\mu (\boldsymbol{\nabla} \boldsymbol{u})_{f}
	+ \mu (\boldsymbol{\nabla} \boldsymbol{u})^T_{f}
	+ \lambda \; \text{tr}\left[(\boldsymbol{\nabla} \boldsymbol{u})_{f} \right] \textbf{I}
	\right] |\boldsymbol{\Gamma}_{f}|
\end{eqnarray}
where subscript $f$ indicates a quantity at a face centre, for example, $(\boldsymbol{\nabla} \boldsymbol{u})_{f}$ is the gradient of displacement tensor at the centre of face $f$.
The approach used to approximate this face displacement gradient is one of the principal differences between variants of the finite volume method.
In addition, polygonal faces may not be flat and approaches have been examined to accommodate for this; for example, \citet{Tukovic2018a} proposed decomposing all faces into triangles before evaluating the forces.

As a consequence of the assumed variation (Equation \ref{eq:taylorSeries}), the gradient of displacement $\boldsymbol{\nabla} \boldsymbol{u}$ is constant within each cell and so the displacement gradient at a face, between two cells, is discontinuous;
to resolve this, the standard cell-centred approach expresses the displacement gradient at a face $(\boldsymbol{\nabla} \boldsymbol{u})_{f}$ as a weighted mean of the displacement gradient at the two adjacent cell-centres:
\begin{eqnarray} \label{eq:approxFaceGradient}
	(\boldsymbol{\nabla} \boldsymbol{u})_{f} &\approx&
	\gamma_{f} (\boldsymbol{\nabla} \boldsymbol{u})_P
	+ (1 - \gamma_{f}) (\boldsymbol{\nabla} \boldsymbol{u})_{N_{f}}
\end{eqnarray}
where subscript $P$ indicates a quantity at the current cell-centre, subscript $N_f$ indicates a quantity at the neighbour cell-centre adjacent to face $f$, and $0 < \gamma_{f} < 1$ is the interpolation weight.
Typically, the interpolation weight is calculated using an inverse distance method:
\begin{eqnarray} \label{eq:linearInterpFaceGradient}
	\gamma_{f} =
	\frac{|\boldsymbol{x}_{N_f} - \boldsymbol{x}_{f}|}
	{|(\boldsymbol{x}_{f} - \boldsymbol{x}_P) + (\boldsymbol{x}_{N_f} - \boldsymbol{x}_{f})|}
\end{eqnarray}
In addition to the standard cell-centred method for approximating the face centre displacement gradient given in Equation \ref{eq:approxFaceGradient}, a number of alternative methods have also been examined, for example, temporary elements with isoparametric formulations \citep{Fallah2008a} or using a compact stencil for the normal gradient as in the original discretisation of \citet{Demirdzic1988};
in this compact stencil form, the normal component is calculated using central differencing: 
\begin{align} \label{eq:approxFaceGradientCompact}
	(\boldsymbol{\nabla} \boldsymbol{u})_{f} &\approx
	\underbrace{
	\boldsymbol{n}_f \boldsymbol{n}_f \cdot 
	\left[
		|\boldsymbol{\Delta}_{f}| \frac{\boldsymbol{u}_{N_{f}} - \boldsymbol{u}_P}{|\boldsymbol{d}_{f}|}
		+ (\boldsymbol{\Gamma}_{f} - \boldsymbol{\Delta}_{f}) \cdot (\boldsymbol{\nabla} \boldsymbol{u})_{f} 
	\right]
	}_{\text{normal component}}
	\notag \\
	&\qquad + \underbrace{
	(\textbf{I} - \boldsymbol{n}_f \boldsymbol{n}_f) \cdot
	\left[
	\gamma_{f} (\boldsymbol{\nabla} \boldsymbol{u})_P
	+ (1 - \gamma_{f}) (\boldsymbol{\nabla} \boldsymbol{u})_{N_{f}}
	\right]
	}_{\text{tangential component}} \notag \\
\end{align}
where $\boldsymbol{n}_f = \boldsymbol{\Gamma}_{f}/|\boldsymbol{\Gamma}_{f}|$ are the face unit normals, and:
\begin{align}
	\boldsymbol{\Delta}_{f} &=
	\frac{\boldsymbol{d}_{f}}{\boldsymbol{d}_{f} \cdot \boldsymbol{n}_{f}} |\boldsymbol{\Gamma}_{f}|
\end{align}


At this point, it is worth noting the locally conservative nature of the discretisation: adjacent cells share integration points at the face centres, resulting in the force at cell faces being locally and hence globally conserved.
This is a characteristic shared by all finite volume methods.

To complete the discretisation of the surface force in terms of displacement $\boldsymbol{u}$, the cell-centred displacement gradients need to be expressed in terms of the cell-centred displacements.
For its accuracy on unstructured grids and ease of implementation, the least-squares method is the most popular:
\begin{eqnarray} \label{eq:leastSquaresGradient}
	(\boldsymbol{\nabla} \boldsymbol{u})_P &\approx&
	\left[  \sum_{f=1}^{nFaces} w_{f}^2 \boldsymbol{d}_{f} \boldsymbol{d}_{f} \right]^{-1}
	\cdot \sum_{f=1}^{nFaces} \left[w_{f}^2 \boldsymbol{d}_{f} \left( \boldsymbol{u}_{N_{f}} - \boldsymbol{u}_P \right) \right]
\end{eqnarray}
As shown previously in Figure \ref{fig:generalFiniteVolumeCell}, vector $\boldsymbol{d}_{f}$ joins cell-centre $P$ to the neighbour cell-centre $N_f$.
The scalar weighting function can be taken as unity ($w_{f} = 1$) \citep{Demirdzic1995} or as the inverse distance ($w_{f} = \nicefrac{1}{|\boldsymbol{d}_{f}|}$) \citep{Jasak2000a}.
Alternative gradient calculation methods, such as the Gauss divergence method or point Gauss divergence method \citep{Tukovic2018a} have also been proposed.
With respect to the Gauss divergence method of gradient calculation, this approach stems from the Lawrence Livermore Laboratory `hydrocodes' of the 1960s developed by Wilkins \citep{Wilkins1963, Wilkins1964}.
Zienkiewicz and O{\~{n}}ate \citep{Zienkiewicz1991, Onate1994} claimed this Gauss divergence cell-gradient method to be ``\emph{an early attempt to use FV concepts in CSM}'';
this link is, however, tenuous;
the Gauss divergence cell-gradient calculation is not a core postulate of the finite volume method and is in fact not required.
For descriptions of subsequent finite difference developments based on the original \citet{Wilkins1963} approach, the interested reader is referred to \citep{Wilkins1999} and \citep{Bessonov2009}.


Although the discretisation of the surface force in terms of displacement $\boldsymbol{u}$ is complete, the presented discretisation is unstable and known to suffer from so-called \emph{checker-boarding} errors.
Without an appropriate stabilisation term, oscillations in the displacement field, which are twice the period of the cell size, will go unnoticed.
These unstable oscillations are analogous to the spurious singular modes that appear in reduced-integration finite element discretisations, also known as zero-energy modes or \emph{hourglassing}.
Typically the implicit cell-centred approach adds the so-called Rhie-Chow stabilisation term to the discretised divergence of stress (Equation \ref{eq:momentumIntegralSurfaceSum2}), as introduced to solid mechanics by \citet{Demirdzic1995}:
\begin{align}
	\mathscrsfs{D}_{\text{Rhie-Chow}} & =
	\sum^{nFaces}_{f=1}
	\left\{
	K_{f} 
	\left[
		|\boldsymbol{\Delta}_{f}| \frac{\boldsymbol{u}_{N_{f}} - \boldsymbol{u}_P}{|\boldsymbol{d}_{f}|}
		+ (\boldsymbol{\Gamma}_{f} - \boldsymbol{\Delta}_{f}) \cdot (\boldsymbol{\nabla} \boldsymbol{u})_{f} 
	\right]
       - \boldsymbol{\Gamma}_{f} \cdot  \left[ K_{f}  (\boldsymbol{\nabla} \boldsymbol{u})_{f} \right]
       \right\} 
	\notag \\
	&=
	\sum^{nFaces}_{f=1}
	K_{f} \left[ |\boldsymbol{\Delta}_{f}| \frac{\boldsymbol{u}_{N_{f}} - \boldsymbol{u}_P}{|\boldsymbol{d}_{f}|}
	- \boldsymbol{\Delta}_{f} \cdot (\boldsymbol{\nabla} \boldsymbol{u})_{f} 
       \right] 
\end{align}
This third-order diffusion term corresponds to the difference between two ways of calculating the normal gradient of displacement at a face, resulting in an ability to `sense' high-frequency oscillations in $\boldsymbol{u}$.
The approach was first proposed by \citet{Rhie1983} in the context of cell-centred finite volume methods for incompressible fluid flow, and is commonly used in cell-centred finite volume fluid formulations.
The $K_{f}$ coefficient controls the magnitude of the smoothing effect and is typically taken as $K_{f} = \mu$ \citep{Demirdzic1995}, $K_{f} = \mu + \lambda$ \citep{Jasak2000a} or $K_{f} = 2\mu + \lambda$ \citep{Cardiff2016b}.
The third-order diffusion term also serves a purpose towards choice of implicit components within the segregated solution algorithm: this is discussed further below.
Alternative forms of diffusion/smoothing terms have been also been proposed, for example, the fourth-order Jameson-Schmidt-Turkel \citep{Jameson1981} term employed in Godunov-type approaches \citep{Aguirre2014, Aguirre2015}, which takes the form of a Laplacian of a Laplacian:
\begin{eqnarray}
	\mathscrsfs{D}_{\text{JST}} \;\; = \;\; \boldsymbol{\nabla}^2 \left[ K_{\text{JST}} \left( \boldsymbol{\nabla}^2 \boldsymbol{u} \right) \right]
\end{eqnarray}
The scalar coefficient $K_{\text{JST}}$ gives the correct dimension to the dissipation as well as controlling its magnitude.



The final discretised form of the governing momentum equation, employing the Rhie-Chow form of stabilisation and $K_f = 2\mu + \lambda$, is expressed as:
\begin{align} \label{eq:finalDiscretisedEqn}
	\rho \, \frac{\boldsymbol{u}_P - 2\boldsymbol{u}_P^{[m-1]} + \boldsymbol{u}_P^{[m-2]}}{\Delta t^2} \, \Omega_P
	&=
	\sum^{nFaces}_{f=1}
	\boldsymbol{n}_{f} \cdot \left[
	\mu (\boldsymbol{\nabla} \boldsymbol{u})_{f}
	+ \mu (\boldsymbol{\nabla} \boldsymbol{u})^T_{f}
	+ \lambda \; \text{tr}\left[(\boldsymbol{\nabla} \boldsymbol{u})_{f} \right] \textbf{I}
	\right] |\boldsymbol{\Gamma}_{f}| \notag \\
	& \quad \;+\;
	\sum^{nFaces}_{f=1}
	(2\mu + \lambda)
	\left[
		|\boldsymbol{\Delta}_{f}| \frac{\boldsymbol{u}_{N_{f}} - \boldsymbol{u}_P}{|\boldsymbol{d}_{f}|}
		- \boldsymbol{\Delta}_{f} \cdot (\boldsymbol{\nabla} \boldsymbol{u})_{f} 
	\right]
	 \notag \\
       & \quad \;+\; \rho \boldsymbol{f}_{b_P}  \Omega_P
\end{align}
where the face displacement gradients $(\boldsymbol{\nabla} \boldsymbol{u})_{f}$ are calculated using Equations \ref{eq:approxFaceGradient}, \ref{eq:linearInterpFaceGradient} and \ref{eq:leastSquaresGradient}.
The primitive unknown variables are the cell-centre displacement vectors $\boldsymbol{u}_P$ at time $t$ (time index $[m]$).

Boundary conditions are incorporated through appropriate modification of the surface force term discretisation at  faces coinciding with the boundary of the solution domain.
In the case of a displacement condition $\boldsymbol{u}_b$ (Dirichlet boundary condition), the face displacement gradients are calculated at the face, while in the case of a traction condition $\boldsymbol{T}_b$ (Neumann boundary condition), the specified traction directly replaces the surface stress expression.
Initial conditions, in the form of the displacement field at $t = 0$, $t = -\Delta t$, and $t = -2\Delta t$, must also be specified.

\paragraph{Solution algorithm}
In order to solve the discretised governing equation (Equation \ref{eq:finalDiscretisedEqn}) for the unknown displacement vector, 
the typical cell-centre approach uses a \emph{segregated} solution procedure, where the central-differencing component in Equation \ref{eq:finalDiscretisedEqn}, $(2\mu + \lambda) |\boldsymbol{\Delta}_{f}| (\boldsymbol{u}_{N_{f}} - \boldsymbol{u}_P)/|\boldsymbol{d}_{f}|$, and $\rho \left(\boldsymbol{u}_P/\Delta t^2 \right) \Omega_P$ within the inertia term are treated implicitly; all other terms are calculated explicitly using the latest available displacement field.
The purpose of the segregated approach is to temporarily decouple/segregate the three scalar components of the vector momentum equation so that they can be solved sequentially;
outer fixed-point/Gauss-Seidel/Picard iterations provide the necessary coupling, where the displacement gradient terms are explicitly updated each outer iteration using the latest available displacement field.
Employing this implicit-explicit split, the discretised equation for each cell can be written in the form of a linear algebraic equation:
\begin{eqnarray} \label{eq:algebraicEqn}
	a_P \; \boldsymbol{u}_P \; \; - \; \; \sum_{f=1}^{nFaces} a_{N_{f}} \boldsymbol{u}_{N_{f}} &=& \boldsymbol{b}_P
\end{eqnarray}
where
\begin{eqnarray}
	a_P &=& \frac{\rho \Omega_P}{\Delta t^2} + \sum_{f=1}^{nFaces} a_{N_{f}} \\
	a_{N_{f}} &=& (2\mu + \lambda) \frac{|\boldsymbol{\Delta}_{f}|}{|\boldsymbol{d}_{f}|} \\
	\boldsymbol{b}_P &=&
	\rho \, \frac{2\boldsymbol{u}_P^{[m-1]} - \boldsymbol{u}_P^{[m-2]}}{\Delta t^2}  \, \Omega_P \notag \\
	&& +
	\sum^{nFaces}_{f=1}
	\boldsymbol{n}_{f} \cdot \left\{
	\mu (\boldsymbol{\nabla} \boldsymbol{u})_{f}
	+ \mu (\boldsymbol{\nabla} \boldsymbol{u})^T_{f}
	+ \lambda \; \text{tr}\left[(\boldsymbol{\nabla} \boldsymbol{u})_{f} \right] \textbf{I}
	\right\} |\boldsymbol{\Gamma}_{f}| \notag \\
	&&-
	\sum^{nFaces}_{f=1}
	(2\mu + \lambda)
	\boldsymbol{\Delta}_{f}
	\cdot (\boldsymbol{\nabla} \boldsymbol{u})_{f} \;\;+\;\; \rho \boldsymbol{f}_{b_P}  \Omega_P
\end{eqnarray}
In contrast, typical implicit vertex-centred and implicit finite element solution algorithms treat the entire divergence of stress term implicitly within the linear system matrix, or a linearisation of it when nonlinearities are present.
This so-called \emph{block-coupled} approach has also been proposed for the cell-centred approach \citep{Das2012, Cardiff2016a}:
in this case, $a_P$ and $a_{N_{f}} $ in Equation \ref{eq:algebraicEqn} are second-order tensors.

The algebraic equations (Equation \ref{eq:algebraicEqn}) can be assembled for all $M$ cells in the domain into the form of three decoupled linear systems:
\begin{eqnarray} \label{eq:linearSystem}
	\left[ \boldsymbol{K} \right]
	\left[ \boldsymbol{U} \right] =
	\left[ \boldsymbol{F} \right]
\end{eqnarray}
where $\left[ \boldsymbol{K} \right]$ is a $M \times M$ sparse matrix with diagonal coefficients $a_P$ and off-diagonal coefficients $a_{N_{f}}$, $\left[ \boldsymbol{U} \right]$ is a vector of the unknown cell-centre displacement vectors, and $\left[ \boldsymbol{F} \right]$ is the source vector containing contributions from $\boldsymbol{b}_P$.
In finite element parlance, $\left[ \boldsymbol{K} \right]$ is the global stiffness matrix and $\left[ \boldsymbol{F} \right]$ is the global force vector.
The segregated cell-centred discretisation ensures that matrix $\left[ \boldsymbol{K} \right]$, has the following properties \citep{Demirdzic1995}:
\begin{itemize}
	\item It is sparse with the number of non-zero elements in each row equal to the number of nearest
neighbours cells (those sharing a face with the cell) plus one;
	\item It is symmetric;
	\item It is positive definite;
	\item It is diagonally dominant ($|a_P| \geq \sum_{f=1}^{nFaces} |a_{N_{f}}|$), which makes the linear system efficiently solved by a number of iterative methods, which retain the sparsity of matrix $\left[ \boldsymbol{K} \right]$: this results in significantly lower memory requirements than equivalent direct linear solvers, for example, see \citep{Cardiff2016a}. The most common iterative method used is the conjugate gradient method with incomplete Cholesky preconditioning \citep{Jacobs1980}.
\end{itemize}
It is worth noting that the segregated solution algorithm has been shown to suffer from slow convergence for slender geometry undergoing bending \citep{Cardiff2016a}; in that case, a block-coupled algorithm is significantly faster \citep{Cardiff2016a};
and $\left[ \boldsymbol{K} \right]$ is a $M \times M$ sparse matrix, where each coefficient is a second-order tensor.

For the segregated approach, the linear system (Equation \ref{eq:linearSystem}) need not be solved to a tight tolerance as coefficients and source terms are approximated from the previous outer iteration; instead, a reduction in the residuals of one order of magnitude is typically sufficient. Outer iterations are performed until the predefined solution tolerance has been achieved.
Under-relaxation of the displacement field and/or the linear system may improve convergence, depending on the boundary conditions and mesh.
Acceleration of the approach can be achieved through geometric multi-grid procedures \citep{Fainberg1996, Demirdzic1997c, Ivankovic1997a}, block-coupled algorithms \citep{Das2011a, Das2011b, Cardiff2016b, Gonzalez2018}, Aitken acceleration \citep{Tukovic2014, Tang2015, Gonzalez2018}, and/or parallelisation on distributed memory clusters \citep{Jasak2000a, Cardiff2014a}. 
A favourable characteristic of the solution procedure is the straightforward extension to nonlinearity:
nonlinear terms (material, geometric or boundary conditions) are resolved \emph{on the fly};
after each outer iteration the coefficients and the source terms are updated and the procedure continues as in the linear case, for example, compare the linear elasticity approach of \citet{Jasak2000a} with the finite strain elasto-plastic frictional contact procedure of \citet{Cardiff2016b}.
In addition to the displacement-based approach described above, alternative solution algorithms, where pressure and displacement are the primary variables, have been proposed by Bijelonja \etal \citep{Bijelonja2005a, Bijelonja2006, Bijelonja2017} and \citet{Fowler2003}; the benefit of such approaches is their ability to deal with incompressible and quasi-incompressible solids in a straightforward manner, while avoiding pressure instabilities.

\subsection{Implicit vertex-centred approach} \label{sec:vertexCentredDiscretisation}
\begin{figure}[htb]
   \centering
	\subfigure[Primary mesh]
	{
		\includegraphics[width=0.45\textwidth]{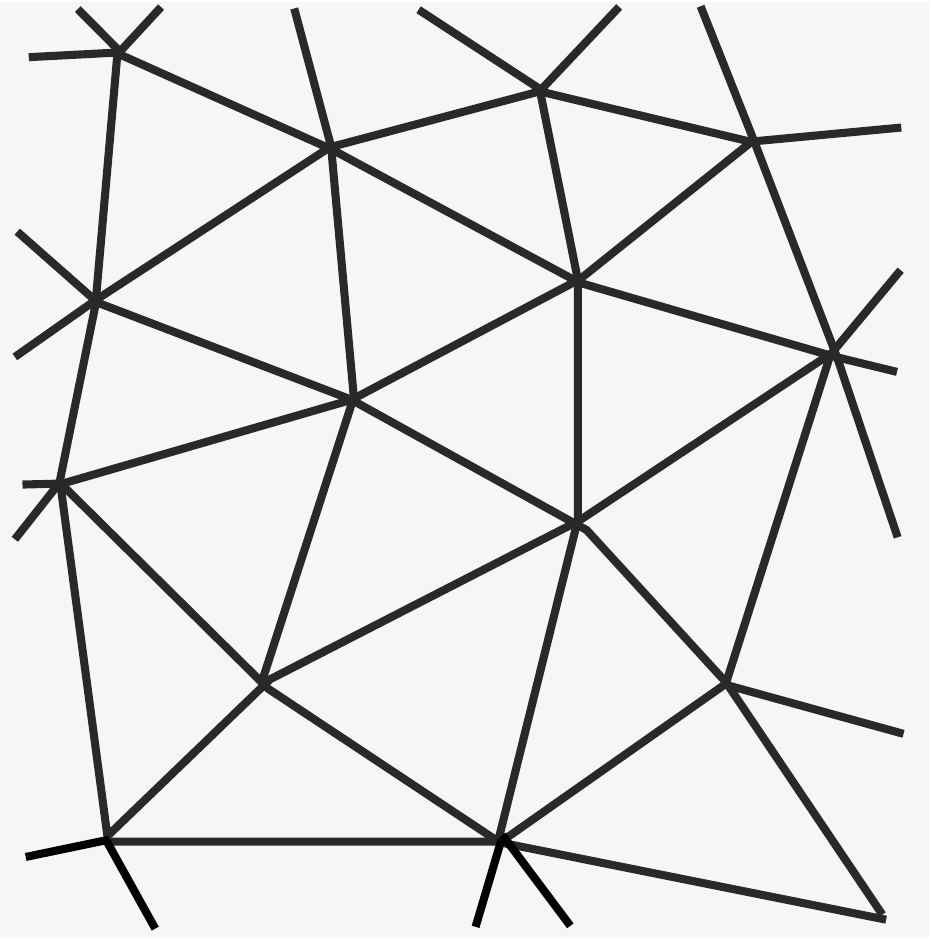}
	}
	\subfigure[secondary-grid used to integrate the governing equations]
	{
	    	\includegraphics[width=0.45\textwidth]{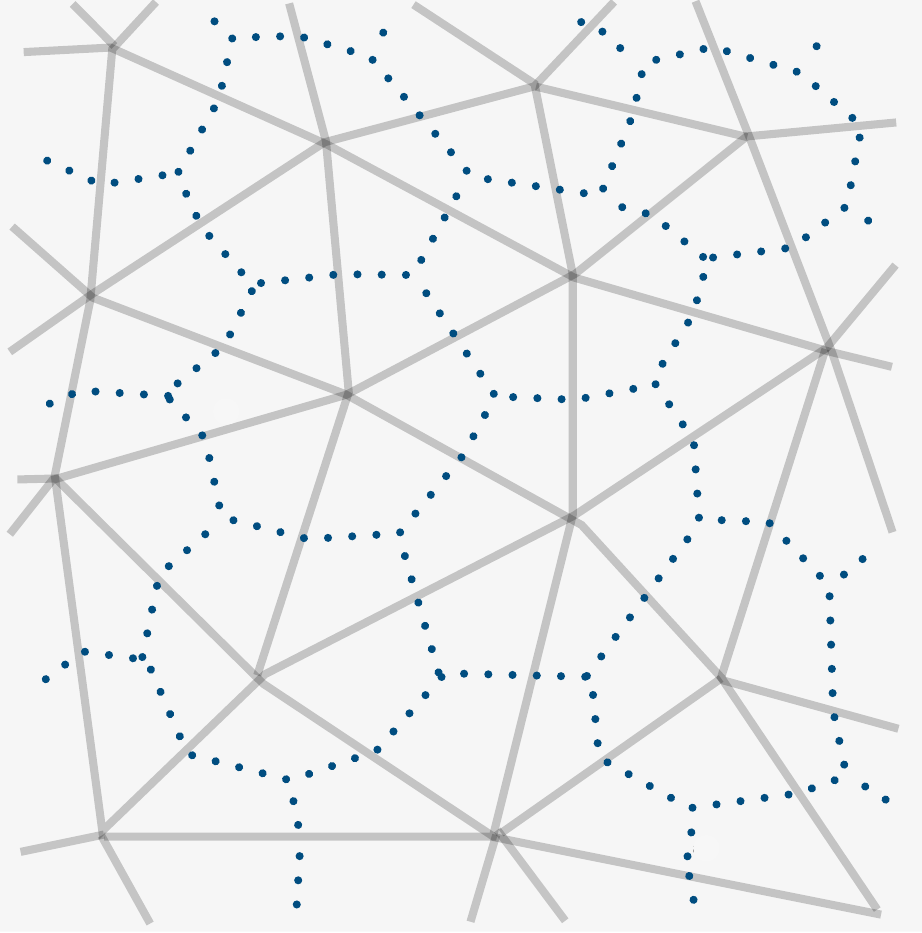}
	}
   \caption{2-D vertex-centred grid showing the (a) primary mesh, (b) secondary-grid used to integrate the governing equations. Figure adapted from \citet{Hassan2019}}
   \label{fig:vertexCentredGrids}
\end{figure}

\paragraph{Discretisation of space}
Like the other approaches, the vertex-centred approach divides the spatial domain into a finite number of contiguous cells that do not overlap and fill the space completely.
When compared with the \emph{typical} cell-centred approach, there are two key differences related to the mesh arrangement:
\begin{itemize}
	\item The vertex-centred approach integrates the governing equations over cells in a secondary-grid, 
	with cells that are typically constructed around the vertices/points in the primary grid (Figure \ref{fig:vertexCentredGrids}).
	The resulting secondary grid cells may not preserve convexity;
	\item The primitive unknowns are stored at the vertices/points of the primary grid, corresponding to the approximate, but not necessarily exact, centre of the secondary-grid cells.
\end{itemize}
However, in essence vertex-centred approaches can be viewed as a form of cell-centred method integrated over a secondary mesh.

In principal, the primary mesh can consist of arbitrary convex polyhedral cells; however, depending on the chosen discretisation (for example, if shape functions are used), the method may be limited to triangular/quadrilateral meshes in 2-D and tetrahedral/hexahedral meshes in 3-D; in fact, no vertex-centred solid mechanics examples using polyhedral meshes were found when preparing this article.

In addition to this form of the vertex-centred approach, a variant exists where the primary mesh cells around a vertex are used to perform the integration, for example, as discussed by \citet{Onate1994};
this produces overlapping regions of integration, where neighbouring vertices share part of their integrated volume.
Secondary-grid cells can also be created by joining the mesh cell centres together \citep{Tsui2013}, rather than joining the cell-centres to the face-centres as per the classic vertex-based method.

\begin{figure}[htb]
	\centering
	\includegraphics[width=0.4\textwidth]{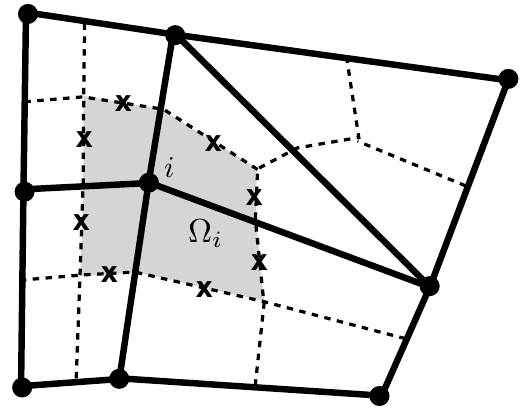}
	\caption{
	A typical control volume constructed around vertex $i$, with a volume $\Omega_i$.
	The solid lines show the primary grid and the dashed lines show the secondary-grid.
	When calculating the momentum balance for a secondary-grid cell (shaded), the displacement gradient is calculated at the boundary edge centres (boundary face centres in 3-D) of the secondary-grid cell (marked by ``x''). 
Figure adapted from \citet{Bailey1995}}
	\label{fig:vertex_centred_integration}
\end{figure}

\paragraph{Discretisation of the mathematical model equations} \label{sec:vertexCentredDiscretisation}
The vertex-centred approach starts from the strong integral form of the governing momentum equation (Equation \ref{eq:momentumIntegral}), which is integrated over a secondary-grid \citep{Fryer1991}.
The vertex-centred approach discretises each term of the governing equation over the cells in the secondary-grid.
An example secondary-grid cell is shown in Figure \ref{fig:vertex_centred_integration}, where each vertex $i$ in the primary grid is uniquely associated with a cell in the secondary-grid.

As with the cell-centred variant, the discretisation of each of the three terms (inertia, surface forces, body forces) in the governing conservation law (Equation \ref{eq:momentumIntegral}) will now be discussed in turn.

To approximate the volume integral temporal term, the displacement $\boldsymbol{u}$ within each secondary-grid cell is assumed to be constant \citep{Baliga1980, Slone2003, Bailey1995}.
Consequently, for a volume about vertex $i$, the term may be approximated in terms of the acceleration at vertex $i$ and the secondary-grid cell volume $\Omega_i$ about vertex $i$:
\begin{eqnarray}
	\int_\Omega \rho \frac{\partial^2 \boldsymbol{u}}{\partial t^2}\;\text{d}\Omega
	&\approx&
	\rho \left(\frac{\partial^2 \boldsymbol{u}}{\partial t^2}\right)_i \Omega_i
\end{eqnarray}

As with other finite volume approaches, the acceleration term $\partial^2 \boldsymbol{u}/\partial t^2$ may be discretised in time using any appropriate finite difference scheme.
For example, \citet{Slone2003} employed the Newmark finite difference scheme, which is popular in the finite element community.
For ease of comparison with the cell-centred approach (Equations \ref{eq:inertiaTermDiscretisedCellCentred}), 
a first-order Euler scheme is assumed here, giving the final discretised term as:
\begin{eqnarray} \label{eq:inertiaTermDiscretisedVertexCentred}
	\int_\Omega \rho \frac{\partial^2 \boldsymbol{u}}{\partial t^2}\;\text{d}\Omega
	&\approx&
	\rho \; \frac{\boldsymbol{u}_i - 2\boldsymbol{u}_i^{[m-1]} + \boldsymbol{u}_i^{[m-2]}}{\Delta t^2} \; \Omega_i
\end{eqnarray}
Comparing this expression with the equivalent cell-centred expression (Equation \ref{eq:inertiaTermDiscretisedCellCentred}), the only difference is the mesh location where the unknown displacement is stored.
A conceptual difference comes from the fact that vertex $i$ is not in general situated at the centre of the secondary-grid cell volume $\Omega_i$; consequently, in the approximation of the inertia term, the cell-centred approach allows a local linear displacement distribution, whereas the vertex-centred approach requires a constant displacement distribution.

The volume integral body force term is discretised in a similar manner to the temporal term, where the body force $\boldsymbol{f}_b$ is assumed to be constant within each secondary-grid cell and equal to the value at vertex $i$:
\begin{eqnarray}
	\int_\Omega \rho \boldsymbol{f}_b \;\text{d}\Omega &\approx& \rho \boldsymbol{f}_{b_i} \Omega_i
\end{eqnarray}
Once again the discretised term differs from the the equivalent cell-centred term in that a constant rather than a linear local variation is assumed.
In the case that vertex $i$ does lie at the centroid of the control volume, then the vertex-centred and cell-centred terms are the same.

To discretise the surface force term, the closed surface integral in Equation \ref{eq:momentumIntegral} is converted into a sum of surface integrals over the faces of the secondary-grid cell, taking the same form as Equation \ref{eq:momentumIntegralSurfaceSum}.
Approximation of the stress term at each face requires an assumption about the local displacement distribution.
In contrast to the standard cell-centred approach, \emph{most} vertex-centred approaches explicitly use shape functions to describe the local displacement field variation.
Following the standard finite element notation, these shape (or interpolation) functions refer to the displacement within each cell of the primary mesh:
\begin{eqnarray} \label{eq:shapeFunctions}
	\boldsymbol{u}(\boldsymbol{x})
	\quad = \quad \sum_{i=1}^{nVertices} N_i (\boldsymbol{x}) \boldsymbol{u}_i 
	\quad = \quad \left[ \boldsymbol{N} \right] \left[ \boldsymbol{u} \right]
\end{eqnarray}
where $nVertices$ is the number of vertices in the primary mesh cell of interest, $N_i  (\boldsymbol{x})$ is the shape function associated with vertex $i$ within the cell, $\boldsymbol{u}_i$ is the displacement at vertex $i$, and $\left[ \boldsymbol{N} \right]$ and $\left[ \boldsymbol{u} \right]$ represent a vector of all shape functions and nodal displacements within the cell respectively.
It should be emphasised that these shape functions refer to the vertices and cells of the primary mesh, as opposed to of the secondary mesh cells.

As is standard in the conventional continuous Bubnov-Galerkin finite element method, the shape functions, which are defined in the reference domain and mapped to the physical domain, approximate the displacement field as a continuous piecewise distribution.
As the shape functions describe the displacement distribution between the vertices on the primary mesh, this allows convenient calculation of the displacement gradients at the secondary-grid cell boundaries when applying the momentum balance (Figure \ref{fig:vertex_centred_integration}).
Although shape functions are a key characteristic of the vertex-centred approaches presented in literature, they are \emph{not} essential.
As shown by \citet{Tsui2013}, it is possible to develop a vertex-centred finite volume approach that does not directly use shape functions;
similar to many cell-centred approaches, \citet{Tsui2013} applied the Gauss divergence theorem to calculate secondary-grid cell displacement gradients; these gradients were then interpolated to secondary-grid cell faces.
In the cell-centred approach, the truncated Taylor series expansion about the cell centres (Equation \ref{eq:taylorSeries}) fulfills the role of shape functions. As shape functions are specific to the cell geometry (\eg triangle, quadrilateral, tetrahedral), shape function based approaches are limited in their choice of element shape, whereas all convex polyhedral meshes are valid for cell-centred approaches.

To approximate the stress at the faces of a secondary-grid cell, the gradient of displacement must be calculated.
The use of shape functions allows this spatial gradient of displacement to be conveniently calculated as:
\begin{eqnarray} \label{eq:shapeFunctionsDeriv}
	\boldsymbol{\nabla} \boldsymbol{u}
	\quad \approx \quad \sum_{i=1}^{nVertices} \boldsymbol{\nabla} N_i (\boldsymbol{x}) \boldsymbol{u}_i 
	\quad = \quad \left[ \boldsymbol{B} \right] \left[ \boldsymbol{u} \right]
\end{eqnarray}
where $\left[ \boldsymbol{B} \right]$ is a vector of the shape functions gradients.
Substituting Equation \ref{eq:shapeFunctionsDeriv} into Equation \ref{eq:momentumIntegralSurfaceSum} allows the surface force term to be expressed in terms of the unknown displacements at the primary mesh vertices:
\begin{eqnarray}
	\sum^{nFaces}_{f=1}
	\int_{\Gamma_i}
	\boldsymbol{n} \cdot \left[
	\mu \boldsymbol{\nabla} \boldsymbol{u}
	+ \mu (\boldsymbol{\nabla} \boldsymbol{u})^T
	+ \lambda \; \text{tr}(\boldsymbol{\nabla} \boldsymbol{u}) \textbf{I}
	\right]
	\;\text{d}\Gamma
	\approx
	\quad\quad\quad\quad\quad\quad\quad\quad
	\notag \\
	\sum_{f=1}^{nFaces} 
	\sum_{v=1}^{nVertices} 
	\boldsymbol{n}_f \cdot \left[
	\mu \boldsymbol{\nabla} N_{vf} \boldsymbol{u}_v 
	+ \mu \boldsymbol{u}_v  \boldsymbol{\nabla} N_{vf}
	+ \lambda \; \text{tr}(\boldsymbol{\nabla} N_{vf} \boldsymbol{u}_v ) \textbf{I}
	\right]
\end{eqnarray}
where $nVertices$ refers to the number of vertices in the primary-grid cell in which the secondary-grid face $f$ is situated;
$N_{vf}$ refers to the corresponding shape function for vertex $v$ within the primary-gird cell.
Consequently, the surface force term for a secondary-grid cell surrounding vertex $v$ will be a function of the displacement at vertex $v$ as well as the displacement at all vertices that share a primary-grid cell with vertex $v$.
For example, for the secondary-grid cell shown in Figure \ref{fig:vertex_centred_integration}, the surface force term will be a function of the displacements at all vertices shown except for the vertex in the top-right.

The final discretised form of the governing momentum equation for the vertex-centred approach is expressed for a secondary-grid cell as:
\begin{eqnarray}
	\rho \, \frac{\boldsymbol{u}_i - 2\boldsymbol{u}_i^{[m-1]} + \boldsymbol{u}_i^{[m-2]}}{\Delta t^2} \, \Omega_i
	\quad\quad =
	\quad\quad\quad\quad\quad\quad\quad\quad\quad\quad\quad\quad\quad\quad
	\notag \\
	\sum_{f=1}^{nFaces} 
	\sum_{v=1}^{nVertices} 
	\boldsymbol{n}_f \cdot \left[
	\mu \boldsymbol{\nabla} N_{vf} \boldsymbol{u}_v
	+ \mu \boldsymbol{u}_v  \boldsymbol{\nabla} N_{vf}
	+ \lambda \; \text{tr}(\boldsymbol{\nabla} N_{vf} \boldsymbol{u}_v ) \textbf{I}
	\right]
	\quad+\quad 
	\rho \boldsymbol{f}_{b_v} \Omega_v
\end{eqnarray}
As with the other finite volume approaches, boundary conditions are incorporated through appropriate modification of the surface force term, and initial conditions must be specified for dynamic cases.

\paragraph{Solution algorithm}
The discretised equation for each secondary-grid cell can be written in the form of an algebraic equation:
\begin{eqnarray} \label{eq:algebraicEqnVertexCentred}
	\boldsymbol{A}_i \cdot \boldsymbol{u}_i \; \; - \; \; \sum_{v=1}^{nVertices} \boldsymbol{A}_{v} \cdot \boldsymbol{u}_{v} &=& \boldsymbol{b}_i
\end{eqnarray}
where $\boldsymbol{u}_i$ is the displacement at the primary-grid vertex associated with the secondary-grid cell, and $nVertices$ indicates all vertices which share a primary-grid cell with vertex $i$.
$\boldsymbol{A}_i$ and $\boldsymbol{A}_{v}$ are the corresponding block coefficients (second-order tensors in 3-D).

In earlier publications \citep{Fryer1991}, the vertex-centred approach followed a similar solution algorithm to the standard cell-centred approach, where the surface force term was partitioned into \emph{implicit} and \emph{explicit} components and a segregated solution algorithm was employed.
In later publications \citep{Bailey1995}, a block-coupled solution procedure is used, where the entire surface force term is discretised implicitly.
As noted previously, similar coupled solution algorithms were later proposed for the cell-centred variant by \citet{Das2011a} and \citet{Cardiff2016a}.

Following the block-coupled approach, the algebraic equations (Equation \ref{eq:algebraicEqnVertexCentred}) can be assembled for all $M$ secondary-grid cells (primary-grid vertices) in the domain to form a system of linear equations:
\begin{eqnarray} \label{eq:linearSystem}
	\left[ \boldsymbol{K} \right]
	\left[ \boldsymbol{U} \right] =
	\left[ \boldsymbol{F} \right]
\end{eqnarray}
where the $M \times M$ global stiffness matrix $\left[ \boldsymbol{K} \right]$ is sparse with block diagonal coefficients $\boldsymbol{A}_i$ and off-diagonal block coefficients $\boldsymbol{A}_{v}$, $\left[ \boldsymbol{U} \right]$ is a vector of the unknown primary-grid vertex displacement vectors, and the global force source vector $\left[ \boldsymbol{F} \right]$ contains contributions from $\boldsymbol{b}_i$.

This linear system can then be solved using direct or iterative linear solvers, where iterative conjugate gradient solvers with diagonal/Jacobi scaling have been favoured in literature, for example, \citep{Fryer1991, Bailey1995}.
Parallelisation on distributed memory supercomputers has been addressed by McManus \etal \citep{McManus2000, McManus2002}, and mixed displacement-rotation approaches have also been employed \citep{Wenke2003, Pan2010, Pan2011}.

\subsection{Explicit cell-centred Godunov-type approach}
\label{sec:explicitGodunov}
Godunov-type procedures are specialised approaches for analysis of problems that involve wave propagation, shocks and solution discontinuities.
Features that distinguish this variant of finite volume approach from the others are the use of acoustic Riemann solvers, solution reconstruction procedures, slope limiters for gradient calculations, occasionally nodal integration, as well as representing the governing equations as a coupled system of first-order equations;
in addition, explicit time marching solution algorithms are employed.




\paragraph{Discretisation of space}
The majority of Godunov-type approaches are spatially discretised using cell-centred approaches, for example, \citep{Kluth2008, Kluth2010, Lee2013, Maire2013, Sambasivan2013, Despres2015, Ndanou2015, Cheng2015, Loubere2016, Boscheri2016, Haider2017, Georges2017, Hueze2017, Cheng2017a, Cheng2017b, Fridrich2017, Heuze2018, Haider2018}; consequently, this section will focus on such approaches;
however, Godunov-type approaches have also used vertex-centred \citep{Aguirre2014, Aguirre2015}, staggered grid \citep{Sijoy2015, Loubere2016} and face-centred \citep{Sevilla2018a, Sevilla2018b} formulations.
To-date, many of the developments have been limited to one and two dimensions, but in general the procedures can be extended to three dimensions.
Like the \emph{implicit} cell-centred approach described in Section \ref{sec:standardFiniteVolumeDiscretisationOfSpace}, the solution spatial domain is divided into a finite number of contiguous convex polyhedral cells bounded by polygonal faces that do not overlap and fill the space completely.


\paragraph{Discretisation of the mathematical model equations} \label{sec:GodunovDiscretisation}
In a notable deviation from the other finite volume variants, Godunov-type approaches portray the second-order governing momentum equation (Equation \ref{eq:momentumIntegral}) as a coupled system of first-order equations:
\begin{eqnarray}
	\label{eq:rateOfStressGoverningEqn1}
	\int_\Omega \rho \frac{\partial \boldsymbol{v}}{\partial t}\;\text{d}\Omega
	&=&
	\oint_\Gamma \boldsymbol{n} \cdot \boldsymbol{\sigma}	\;\text{d}\boldsymbol{\Gamma}
	\;+\;
	\int_\Omega \rho \boldsymbol{f}_b \;\text{d}\Omega \label{eq:rateOfVelocityGoverningEqn}
	\\ 
	\int_\Omega \frac{\partial \boldsymbol{F}}{\partial t}\;\text{d}\Omega
	&=&
	\oint_\Gamma	\boldsymbol{v} \boldsymbol{n}	\;\text{d}\boldsymbol{\Gamma}
	\label{eq:rateOfStressGoverningEqn}
\end{eqnarray}
where $\boldsymbol{v}$ is the velocity vector and the deformation gradient, $\boldsymbol{F}$, defines the local deformation as:
\begin{eqnarray} \label{eq:deformationGradient}
	\boldsymbol{F} = \textbf{I} + (\boldsymbol{\nabla} \boldsymbol{u})^T
\end{eqnarray}
Equation \ref{eq:rateOfStressGoverningEqn} is obtained by taking the time-derivative of Equation \ref{eq:deformationGradient} and  employing the Gauss divergence theorem.

To ensure the compatibility conditions are satisfied, the deformation gradient at the initial time should be curl free:
\begin{eqnarray}
	\boldsymbol{\nabla} \times \boldsymbol{F}^{[m = 0]} = \boldsymbol{0} \notag \\
\end{eqnarray}
In addition, the discrete evolution of $\boldsymbol{F}$ (or $\boldsymbol{\nabla} \boldsymbol{u}$) should not allow curl errors to escalate.
For further discussion of this point, see \citep{Trangenstein1991, Trangenstein1994, Lee2013, Despres2015, Haider2017}.
It is also possible to employ the evolution of the displacement gradient $\boldsymbol{\nabla} \boldsymbol{u}$ directly rather than the deformation gradient, as shown by \citet{Trangenstein1991}; however, this is less popular in literature.



The governing system of coupled first-order equations (Equations \ref{eq:rateOfVelocityGoverningEqn} and \ref{eq:rateOfStressGoverningEqn}) can be expressed concisely as:
\begin{eqnarray} \label{eq:GodunovConcise}
\int_\Omega \frac{\partial \boldsymbol{\mathcal{U}}}{\partial t} \;\text{d}\Omega
= \oint_\Gamma \boldsymbol{\mathcal{F}}_n \;\text{d}\boldsymbol{\Gamma}
+ \int_\Omega \boldsymbol{\mathcal{S}} \;\text{d}\Omega
\end{eqnarray}
where $\boldsymbol{\mathcal{U}}$ is the primary unknown vector, $\boldsymbol{\mathcal{F}}_n$ is the flux vector, and $\boldsymbol{\mathcal{S}}$ is the source vector, given as:
\begin{eqnarray}
	\boldsymbol{\mathcal{U}} = \begin{pmatrix} \boldsymbol{v} \\ \boldsymbol{F} \end{pmatrix}, \quad\quad
	\boldsymbol{\mathcal{F}}_n =
		\begin{pmatrix} \nicefrac{\boldsymbol{t}}{\rho} \\ \boldsymbol{v}\boldsymbol{n} \end{pmatrix}, \quad\quad
	\boldsymbol{\mathcal{S}} = \begin{pmatrix} \boldsymbol{f}_b \\ \boldsymbol{0} \end{pmatrix}
\end{eqnarray}
and the traction vector, $\boldsymbol{t}$, gives the stress on a plane as $\boldsymbol{t} = \boldsymbol{n} \cdot \boldsymbol{\sigma}$.
To close the system, we give the constitutive relation (Hooke's law) in terms of the deformation gradient:
\begin{eqnarray} \label{eq:HookesLawF}
	\boldsymbol{\sigma} =
	\mu \boldsymbol{F}^T
	+ \mu \boldsymbol{F}
	+ \lambda \; \text{tr} \left( \boldsymbol{F} \right) \textbf{I}
	 - (2\mu + 3 \lambda)\textbf{I}
\end{eqnarray}

We next describe the discretisation of the three terms in Equation \ref{eq:GodunovConcise}: the time derivative term, the diffusion term and the body force term.

Assuming the primary unknowns ($\boldsymbol{v}$ and $\boldsymbol{F}$) vary linearly within each cell according to Equation \ref{eq:taylorSeries}, the volume integrals can be expressed in terms of the cell-centre values (subscript $P$) and the surface integral becomes a sum over the face-centre values (subscript $f$);
consequently, Equation \ref{eq:GodunovConcise} becomes:
\begin{eqnarray} \label{eq:GodunovSemiDiscrete}
	\frac{\partial \boldsymbol{\mathcal{U}}_P}{\partial t} \Omega_P
	= \sum_{f=1}^{nFaces} \boldsymbol{\mathcal{F}}_{n_f} |\boldsymbol{\Gamma}_f|
	\; \; + \; \; \boldsymbol{\mathcal{S}}_P \Omega_P
\end{eqnarray}
Similar to the other approaches, discretisation of the time-rate term can be achieved using a variety of finite difference methods.
Here we assume a first-order \emph{forward} Euler discretisation to allow straight-forward comparison with the other methods:
\begin{eqnarray} \label{eq:fluxOverFaces}
	\frac{\boldsymbol{\mathcal{U}}_P - \boldsymbol{\mathcal{U}}_P^{[m-1]}}{\Delta t} \Omega_P
	= \sum_{f=1}^{nFaces}  \boldsymbol{\mathcal{F}}_{n_f}^{[m-1]} |\boldsymbol{\Gamma}_f|
	\; \; + \; \; \boldsymbol{\mathcal{S}}_P \Omega_P
\end{eqnarray}
where, as before, $m$ is the time-step counter, which for brevity has been dropped on the unknown current time value.
Of course, as an explicit time-marching solution algorithm will be employed, the time step size is limited by the Courant-Friedrichs-Lewy constraint \citep{Courant1928};
this condition is necessary for stability but is not sufficient.
It is also necessary that the time integrator avoids the creation of new local extrema, known as the Total Variation Diminishing property or the local maximum principle;
the first-order forward Euler approach is one such method that obeys this condition.

Alternatively, the flux calculation in Equation \ref{eq:fluxOverFaces}, which sums over the cell faces, can  be expressed as a sum over the cell points/vertices according to:
\begin{eqnarray} \label{eq:fluxOverPoints}
	\frac{\boldsymbol{\mathcal{U}}_P - \boldsymbol{\mathcal{U}}_P^{[m-1]}}{\Delta t} \Omega_P
	= \sum_{v=1}^{nVertices}  \boldsymbol{\mathcal{F}}_{n_v}^{[m-1]} |\boldsymbol{C}_v|
	\; \; + \; \; \boldsymbol{\mathcal{S}}_P \Omega_P
\end{eqnarray}
where the calculation of point/vertex/nodal area vector, $\boldsymbol{C}_v$, is given in \citet{Carre2009} and \citet{Kluth2010} or equivalently in \citet{Maire2013}.

Up to this point, the presented discretisation coincides with the cell-centred approach described in Section \ref{sec:standardFiniteVolume}, apart from the introduction of a system of first-order conservation equation and the use of the forward Euler method as opposed to the backward Euler method;
however, in the discretisation of the face flux, $\boldsymbol{\mathcal{F}}_n$, Godunov-type methods deviate from the other approaches.
As a result of the assumed piecewise linear distribution (Equation \ref{eq:taylorSeries}), there is a discontinuity at the cell internal faces.
A distinguishing characteristic of Godunov-type methods is the acknowledgement of this discontinuity in the solution field, known as a Riemann problem, and the development of appropriate methods (Riemann solvers) to deal with the propagation of this discontinuity.
Accordingly, the face flux, $\boldsymbol{\mathcal{F}}_n$, is defined as a function of the solution variable at either side of the interface:
\begin{eqnarray}
	\boldsymbol{\mathcal{F}}_{n_f} = f(\boldsymbol{\mathcal{U}}_{Pf}, \boldsymbol{\mathcal{U}}_{Nf})
\end{eqnarray}
The solution at either side of the interface ($\boldsymbol{\mathcal{U}}_{P_f}$ and $\boldsymbol{\mathcal{U}}_{N_f}$) is determined via extrapolation from the adjacent cell centres:
\begin{eqnarray} \label{eq:GodunovReconstruction}
\boldsymbol{\mathcal{U}}_{P_f} = \boldsymbol{\mathcal{U}}_{P} + \boldsymbol{G}_P \cdot (\boldsymbol{x}_f - \boldsymbol{x}_P) \notag \\
\boldsymbol{\mathcal{U}}_{N_f} = \boldsymbol{\mathcal{U}}_{N} + \boldsymbol{G}_N \cdot (\boldsymbol{x}_f - \boldsymbol{x}_N)
\end{eqnarray}
where $\boldsymbol{G}_P$ is the gradient of $\boldsymbol{\mathcal{U}}$ within cell $P$, and $\boldsymbol{G}_N$ is the gradient within cell $N$.
A critical component of Godunov-type methods is the definition of this solution \emph{reconstruction}, such that the local maximum principle is preserved \ie the value of $\boldsymbol{\mathcal{U}}$ at face $f$ should not be greater than the value at adjacent cell centres $P$ and $N$.
There are a number of ways to define such discrete gradients;
here, as a typical example, the Monotone Upstream Scheme for Conservation Law (MUSCL) scheme is described \citep{Lee2013}.
The MUSCL approach consists of two steps: first, the gradient is predicted based on local neighbouring values, then this gradient is corrected/limited to respect the local maximum principle.

To predict the gradient, the second-order least-squares method described in Section \ref{sec:standardFiniteVolume} can be used;
this approach does not prohibit overshoots and undershoots at the faces, and hence does not satisfy the local maximum principle.
To remedy this, a so-called \emph{slope limiter} is used to restrict the value of the gradient, $\boldsymbol{G}$.
The slope limiter is included through modification of the reconstruction expression (Equation \ref{eq:GodunovReconstruction}):
\begin{eqnarray} \label{eq:GodunovReconstructionLimited}
\boldsymbol{\mathcal{U}}_{P_f} = \boldsymbol{\mathcal{U}}_{P} + \phi_P \boldsymbol{G}_P \cdot (\boldsymbol{x}_f - \boldsymbol{x}_P) \notag \\
\boldsymbol{\mathcal{U}}_{N_f} = \boldsymbol{\mathcal{U}}_{N} + \phi_N \boldsymbol{G}_N \cdot (\boldsymbol{x}_f - \boldsymbol{x}_N)
\end{eqnarray}
where $0 \leq \phi \leq 1$ is a scalar slope limiter.
When $\phi = 1$, no limiting is applied, whereas when $\phi = 0$, full limiting is applied and the value near the face is assumed equal to the cell-centre value.
Choosing the value of $\phi$ is a balance between stability (lower value of $\phi$) and accuracy (higher value of $\phi$).
For the MUSCL procedure, $\phi$ is determined as described by \citet{Lee2013}:
\begin{enumerate}
\item Find the smallest and largest values among the current and adjacent cells:
\begin{eqnarray}
	\boldsymbol{\mathcal{U}}^{\text{min}} = \text{min}(\boldsymbol{\mathcal{U}}_P, \boldsymbol{\mathcal{U}}_{Ni}),
	\quad\quad
	\boldsymbol{\mathcal{U}}^{\text{max}} = \text{max}(\boldsymbol{\mathcal{U}}_P, \boldsymbol{\mathcal{U}}_{Ni}) 
\end{eqnarray}
where $\boldsymbol{\mathcal{U}}_{Ni}$ represents the values at neighbour cells, which share an internal face with cell $P$.
\item Calculate the un-restricted reconstructed value $\boldsymbol{\mathcal{U}}_{Pf}$ with $\phi_P = 1$ at each internal face within cell $P$.
\item Find the maximum allowable value of $\phi_{Pf}$ for each face in cell $P$:
\[
	\phi_{Pf} =
	\begin{cases}
		\text{min}\left(1, \frac{\boldsymbol{\mathcal{U}}^{\text{max}} - \boldsymbol{\mathcal{U}}_P}{\boldsymbol{\mathcal{U}}_{P_f} - \boldsymbol{\mathcal{U}}_P} \right), & \text{if} \quad \boldsymbol{\mathcal{U}}_{P_f} - \boldsymbol{\mathcal{U}}_{P} > 0  \\
		\text{min}\left(1, \frac{\boldsymbol{\mathcal{U}}^{\text{min}} - \boldsymbol{\mathcal{U}}_P}{\boldsymbol{\mathcal{U}}_{P_f} - \boldsymbol{\mathcal{U}}_P} \right), & \text{if} \quad \boldsymbol{\mathcal{U}}_{P_f} - \boldsymbol{\mathcal{U}}_{P} < 0  \\
		1, & \text{if} \quad \boldsymbol{\mathcal{U}}_{P_f} - \boldsymbol{\mathcal{U}}_{P} = 0 
	\end{cases}
\]
\item Select $\phi_P = \text{min}_f(\phi_{Pf})$
\item Calculate the reconstructed value at each internal face using $\phi_P$ determined in step 4 and Equation \ref{eq:GodunovReconstructionLimited}.
\end{enumerate}


To complete the discretisation, all that remains is to express the flux vector, $\boldsymbol{\mathcal{F}}_{n_f}$, in terms of the solution vector, $\boldsymbol{\mathcal{U}}$.
To achieve this, the Rankine-Hugoniot jump conditions \citep{LeVeque2004} are employed.
These jump conditions describe the relationship between the states on both sides of a shock wave.
The jump conditions corresponding to Equations \ref{eq:rateOfStressGoverningEqn1} and \ref{eq:rateOfStressGoverningEqn} are \citep{Lee2013}:
\begin{align} 
	\label{eq:RankineHugoniotJump1}
	U \dbr{\boldsymbol{v}} &= -\left(\nicefrac{1}{\rho}\right) \boldsymbol{n} \cdot \dbr{\boldsymbol{\sigma}}  \\
	\label{eq:RankineHugoniotJump2}
	U\dbr{\boldsymbol{F}} &= -\dbr{\boldsymbol{v}} \boldsymbol{n} 
\end{align}
where the operator $\dbr{\cdot}$ represents the jump across the shock, for example, $\dbr{\boldsymbol{\mathcal{U}}_f} = \boldsymbol{\mathcal{U}}_{N_f} - \boldsymbol{\mathcal{U}}_{P_f}$.
The wave speed is indicated by $U$.
Using equations \ref{eq:RankineHugoniotJump1} and \ref{eq:RankineHugoniotJump2} and assuming a constant wave speed, the face flux can be expressed as the sum of an average flux and a stabilisation flux \citep{Lee2013, Haider2017}:
\begin{align}
	\boldsymbol{\mathcal{F}}_{n_f} &= \boldsymbol{\mathcal{F}}_{n_f}^{\text{average}} + \boldsymbol{\mathcal{F}}_{n_f}^{\text{stab}}
\end{align}
where the average flux is:
\begin{align}
	\boldsymbol{\mathcal{F}}_{n_f}^{\text{average}} &=
	\frac{1}{2} \left[ \boldsymbol{\mathcal{F}}_n(\boldsymbol{\mathcal{U}}_{P_f})
	+ \boldsymbol{\mathcal{F}}_n(\boldsymbol{\mathcal{U}}_{N_f}) \right]
\end{align}
and the so-called upwinding stabilisation term is:
\begin{align}
	\boldsymbol{\mathcal{F}}_{n_f}^{\text{stab}} &=
	\begin{bmatrix}
	\frac{\rho}{2}
	\left[c_s \textbf{I} + (c_p - c_s) \boldsymbol{n}_f \boldsymbol{n}_f \right]
	\cdot \left( \boldsymbol{v}_{N_f} - \boldsymbol{v}_{P_f} \right)
	 \\
	\frac{1}{2} 
	\left[\frac{1}{c_s}\textbf{I} + (\frac{1}{c_p} - \frac{1}{c_s})\boldsymbol{n}_f \boldsymbol{n}_f \right]
	\cdot \left( \boldsymbol{n}_f \cdot \left[ \boldsymbol{\sigma}_{N_f} - \boldsymbol{\sigma}_{P_f} \right] \right)
	\end{bmatrix}
\end{align}
The volumetric and shear wave speeds are indicated by $c_p$ and $c_s$.
Jameson-Schmidt-Turkel stabilisation has been used as an alternative to this upwinding stabilisation, and in principal Rhie-Chow stabilisation could also be used.

The final discretised governing equations, given in terms of unknowns $\boldsymbol{v}^{[m]}_P$ and $\boldsymbol{F}^{[m]}_P$ at time-step $m$, are expressed as:
\begin{eqnarray}
	\frac{\boldsymbol{v}_P^{[m]} - \boldsymbol{v}_P}{\Delta t} \Omega_P
	&=&
	\sum_{f=1}^{nFaces}
	\frac{1}{2\rho} \boldsymbol{n}_f \cdot \left[ 
	\mu \boldsymbol{F}_{P_f}^T
	+ \mu \boldsymbol{F}_{P_f}
	+ \lambda \; \text{tr} \left( \boldsymbol{F}_{P_f} \right) \textbf{I}
	 - (2\mu + 3 \lambda)\textbf{I}
	\right] |\Gamma_f| \notag \\
	&& \;+\; \sum_{f=1}^{nFaces}
	\frac{1}{2\rho} \boldsymbol{n}_f \cdot  \left[
	\mu \boldsymbol{F}_{N_f}^T
	+ \mu \boldsymbol{F}_{N_f}
	+ \lambda \; \text{tr} \left( \boldsymbol{F}_{N_f} \right) \textbf{I}
	 - (2\mu + 3 \lambda)\textbf{I}	
	\right] |\Gamma_f| \notag \\
	&& \;+\; \sum_{f=1}^{nFaces}
	\frac{1}{2}
	\left[c_s \textbf{I} + (c_p - c_s) \boldsymbol{n}_f \boldsymbol{n}_f \right] \cdot \left( \boldsymbol{v}_{N_f} - \boldsymbol{v}_{P_f} \right)
	|\Gamma_f| \notag \\
	&& \;+ \; \boldsymbol{f}_{b_P} \; \Omega_P \label{eq:discretisedGoverningEquation1}
\end{eqnarray}
\begin{eqnarray}
	\frac{\boldsymbol{F}_P^{[m]} - \boldsymbol{F}_P}{\Delta t} \Omega_P
	&=&
	\sum_{f=1}^{nFaces}
	\frac{1}{2} \left(
	\boldsymbol{v}_{P_f} \boldsymbol{n}_f + \boldsymbol{v}_{N_f} \boldsymbol{n}_f
	\right) |\Gamma_f| \notag \\
	&& \;+\; \sum_{f=1}^{nFaces}
	\frac{1}{2} 
	\left[\frac{1}{c_s}\textbf{I} + (\frac{1}{c_p} - \frac{1}{c_s})\boldsymbol{n}_f \boldsymbol{n}_f \right]
	\cdot \left[ \boldsymbol{n}_f \cdot \left( \boldsymbol{\sigma}_{N_f} - \boldsymbol{\sigma}_{P_f} \right) \right]
	|\Gamma_f|
	\label{eq:discretisedGoverningEquation2}
\end{eqnarray}
where the $m - 1$ time index on $\boldsymbol{F}$ and $\boldsymbol{v}$ has been omitted for brevity \ie $\boldsymbol{F} \equiv \boldsymbol{F}^{[m-1]}$ and $\boldsymbol{v} \equiv \boldsymbol{v}^{[m-1]}$.
The stress either side of a face is calculated according to Equation \ref{eq:HookesLawF}.

%

\paragraph{Solution algorithm}
The solution of the governing discretised equations (Equations \ref{eq:discretisedGoverningEquation1} and \ref{eq:discretisedGoverningEquation2}) proceeds in an explicit manner as follows:
\begin{enumerate}
	\item Increase the total time by $\Delta t = \alpha_{\text{CFL}} \frac{h_{\text{min}}}{c_p}$, according to the Courant-Friedrichs-Lewy \citep{Courant1928} condition, where $h_{\text{min}}$ is the shortest length within the mesh and $\alpha_{\text{CFL}}$ is typically chosen to be less than $\nicefrac{1}{2}$
	\item Store the old-time solution values: $\boldsymbol{v}^{[m - 1]} \leftarrow \boldsymbol{v}^{[m]}$ and $\boldsymbol{F}^{[m - 1]} \leftarrow \boldsymbol{F}^{[m - 1]}$
	\item Calculate the reconstructed solution values, $\boldsymbol{v}_{N_f}$, $\boldsymbol{v}_{P_f}$, $\boldsymbol{F}_{N_f}$ and $\boldsymbol{F}_{P_f}$, at each cell face according to Equation \ref{eq:GodunovReconstructionLimited}
	\item Solve Equation \ref{eq:discretisedGoverningEquation1} for $\boldsymbol{v}^{[m]}$
	\item Solve Equation \ref{eq:discretisedGoverningEquation2} for $\boldsymbol{F}^{[m]}$
	\item Repeat steps 1 to 4 until the end time has been reached
\end{enumerate}

Like all explicit methods, this approach does not require the solution of an implicit linear system of equations and hence each time-step can be evaluated more rapidly than in the previously discussed implicit methods;
however, the time-step size limit essentially restricts explicit methods to hyperbolic-style problems.
Finally, the explicit nature of the algorithm allows straight-forward and efficient parallelisation on distributed memory supercomputer.

\subsection{Discussion} \label{sec:finiteVolumeVariantsDiscussion}
The field of finite volume solid mechanics comprises more approaches than those presented above, however, the selected three variants capture the primary differences in the approaches.
It can be seen that the distinction between the variants can be narrowed down to four components:
\begin{enumerate}
	\item Control volume construction;
	\item Face gradient calculation;
	\item Stabilisation approach;
	\item Solution methodology.
\end{enumerate}
Each of these components is briefly reviewed below, followed by the natural definition of a unified approach linking all variants.

\begin{figure}[htb]
	\centering
	\includegraphics[width=0.65\textwidth]{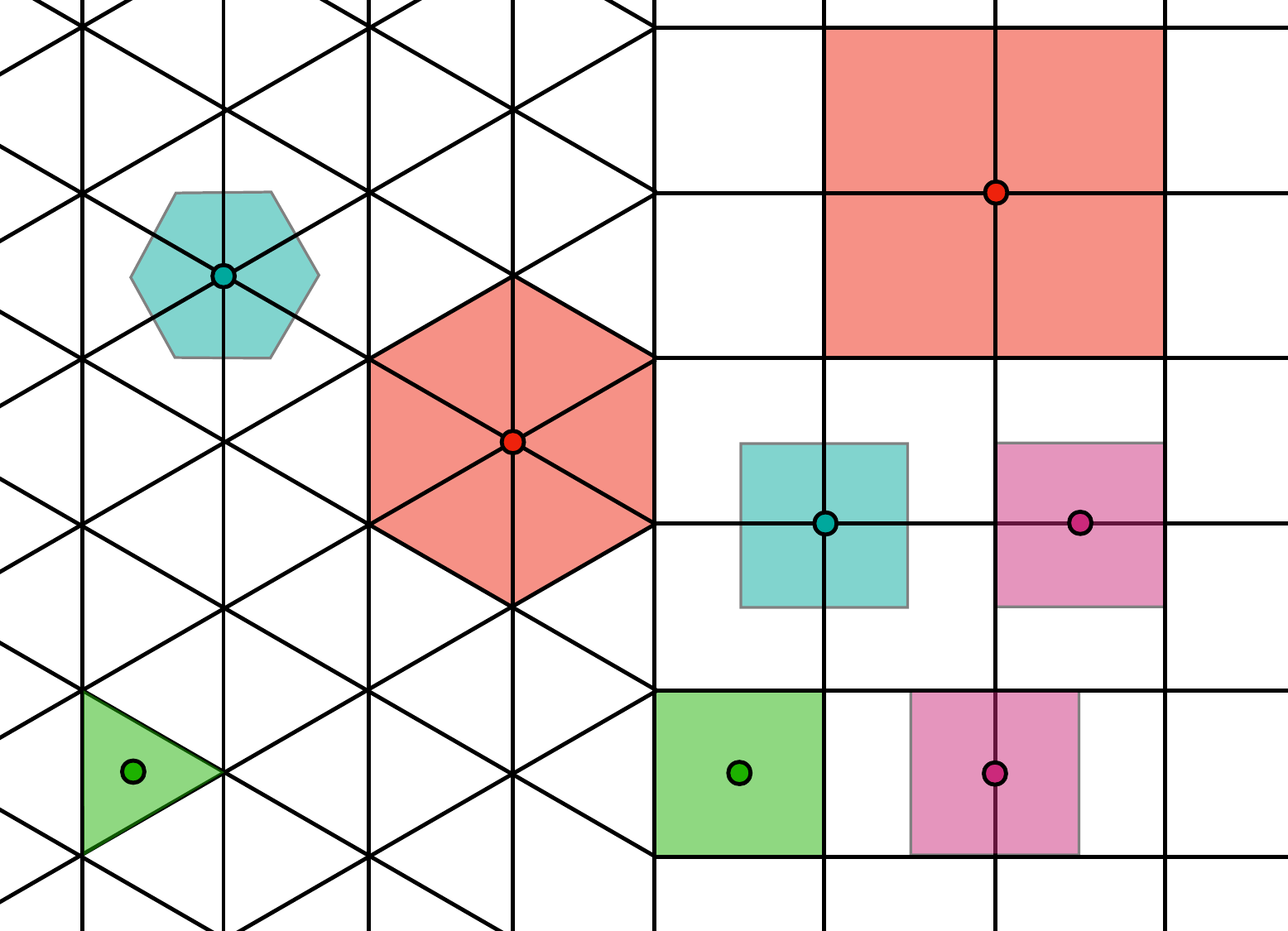}
	\caption{Illustration of the four ways to construct control volumes on a primary mesh:
	1) cell-centred (green); 2) vertex-centred with non-overlapping control volumes (blue); 3) vertex-centred with overlapping control volumes (red); and 4) staggered-grid (purple). The nodal locations are indicated by filled circles.}
	\label{fig:controlVolumeConstruction}
\end{figure}

\paragraph{Control volume construction}
There are predominantly four ways to construct control volumes (Figure \ref{fig:controlVolumeConstruction}): 1) cell-centred, 2) vertex-centred with non-overlapping volumes; 3) vertex-centred with overlapping control volumes, 4) and a staggered-grid.
In terms of the relative merits of the different approaches, a number of observations can be made:
\begin{enumerate}
	\item The staggered-grid is restricted to structured quadrilateral/hexahedral meshes in 2-D/3-D;
	\item The vertex-centred approach with non-overlapping control volumes requires the construction and storage of a second mesh;
	\item The vertex-centred approaches have nodes on the boundary, whereas the cell-centre approach typically does not;
	\item The cell-centred approach allows convenient approximation of the face normal gradients but requires interpolation to approximate the tangential gradients;
	\item In contrast, the vertex-centred approach with overlapping volumes allows convenient approximation of the tangential gradients but requires interpolation to approximate the normal gradients.
\end{enumerate}

\paragraph{Face gradient calculation} \label{sec:faceGradientDiscussion}
Once the control volume is constructed, all approaches must approximate the traction at each control volume face.
To achieve this, the gradient of displacement at the face $(\boldsymbol{\nabla} \boldsymbol{u})_f$ is approximated in terms of the nodal displacement values.
Neglecting any additional stabilisation terms (addressed in the next section), the most common approaches are:
\begin{enumerate}


\item Interpolate the gradient from the adjacent cell-centres:
\begin{align}
	(\boldsymbol{\nabla} \boldsymbol{u})_f & =
	\frac{
	|\boldsymbol{d}_{Pf}| (\boldsymbol{\nabla} \boldsymbol{u})_P
	+ |\boldsymbol{d}_{Nf}| (\boldsymbol{\nabla} \boldsymbol{u})_N
	}{|\boldsymbol{d}_f|}
\end{align}
where $\boldsymbol{d}_f$ is the vector from the centre of cell $P$ to the centre of cell $N$, $\boldsymbol{d}_{Pf}$ is the vector from the centre of cell $P$ to the centre of face $f$, and $\boldsymbol{d}_{Nf}$ is the vector from the centre of cell $N$ to the centre of face $f$. These interpolation weights may differ depending on the specific approach.

\item Calculate the normal gradient using central-differencing, and interpolate the tangential gradient from the adjacent cell-centres:
\begin{align}
	(\boldsymbol{\nabla} \boldsymbol{u})_f & =
	\boldsymbol{n}_f 
	\frac{\boldsymbol{u}_N - \boldsymbol{u}_P}{|\boldsymbol{d}_f|}
	\;+\; \left( \textbf{I} - \boldsymbol{n}_f  \boldsymbol{n}_f  \right) \cdot
	\,
	\frac{
	|\boldsymbol{d}_{Pf}| (\boldsymbol{\nabla} \boldsymbol{u})_P
	+ |\boldsymbol{d}_{Nf}| (\boldsymbol{\nabla} \boldsymbol{u})_N
	}{|\boldsymbol{d}_f|}
\end{align}

\item Calculate the normal gradient using central-differencing, and calculate the tangential gradient using the point or edge values.
In the case of a staggered-grid, these point or edge values correspond to nodes and no interpolation is necessary.
For the cell-centred approach, these point or edge values are interpolated from the cell-centred nodes, and an appropriate face tangential gradient calculation method used;
for example, the face-Gauss method \citep{Cardiff2016a, Tukovic2018a} is a generalisation of the approach used in the original \citet{Demirdzic1988} approach:
\begin{align}
	(\boldsymbol{\nabla} \boldsymbol{u})_f & =
	\boldsymbol{n}_f 
	\frac{\boldsymbol{u}_N - \boldsymbol{u}_P}{|\boldsymbol{d}|_f}
	\;+\;  \sum_{e=1}^{nEdges} \boldsymbol{m}_e \, \boldsymbol{u}_e \, L_e
\end{align}
where $\boldsymbol{m}_e$ is the outward-facing bi-normal at edge $e$, $L_e$ is the length of edge $e$, and $\boldsymbol{u}_e$ is the displacement at the centre of edge $e$, which has been interpolated from adjacent cell-centres.

\item Extrapolate the gradient from each adjacent cell-centre, and take the average:
\begin{align}
	(\boldsymbol{\nabla} \boldsymbol{u})_f & =
	\frac{
	(\boldsymbol{\nabla} \boldsymbol{u})_P
	+ \boldsymbol{d}_{Pf} \cdot \boldsymbol{\nabla}(\boldsymbol{\nabla} \boldsymbol{u})_P
	+ 
	(\boldsymbol{\nabla} \boldsymbol{u})_N 
	+ \boldsymbol{d}_{Nf} \cdot \boldsymbol{\nabla} (\boldsymbol{\nabla} \boldsymbol{u})_N
	}{2}
\end{align}

\item For each of the previous face gradient calculation methods, a limiter can be applied to preserve the local maximum principle, as discussed in Section \ref{sec:GodunovDiscretisation}.

\item Use shape functions to evaluate the face gradient.
As an example, taking a vertex-based method and a 2-D triangular grid with linear shape functions, the face gradient is given in Voigt notation as:

\scriptsize
\begin{align}
	\begin{bmatrix}
	(\boldsymbol{\nabla} \boldsymbol{u})_f 
	\end{bmatrix}
	\quad
	=
	\quad
	\begin{bmatrix}
	\partial u_x/\partial x \\ \partial u_y/\partial y \\
	\partial u_x/\partial y  \\ \partial u_y/\partial x
	\end{bmatrix}
	\quad
	=
	\quad
	\underbrace
	{
	\frac{1}{x_{13} y_{23} - x_{23} y_{13}}
	\begin{bmatrix}
	y_{23} & 0 & y_{31} & 0 & y_{12} & 0 \\
	0 & x_{32} & 0 & x_{13} & 0 & x_{21} \\
	x_{32} & 0 & x_{13} & 0 &  x_{21} & 0 \\
	0 & y_{23} & 0 & y_{31} & 0 & y_{12}
	\end{bmatrix}
	}_{\text{Derivative of shape functions}}
	\cdot
	\begin{bmatrix}
	u_{x_1} \\ u_{y_1} \\ 
	u_{x_2} \\ u_{y_2} \\ 
	u_{x_3} \\ u_{y_3}
	\end{bmatrix}
\end{align} \normalsize
where suffixes $1$, $2$ and $3$ refer to the three nodes in the triangle, and the relative coordinates are given as $x_{ab} = x_a - x_b$.
Here we are referring to the vertex-centred method, however, shape functions could in principle be used with any control volume construction.

\item Use a higher-order approach, for example, the fourth-order approach proposed by \citet{Demirdzic2016}.

\end{enumerate}

\begin{figure}[htb]
	\centering
	\subfigure[Mesh over which the equations are integrated]
	{
		\includegraphics[height=0.47\textwidth]{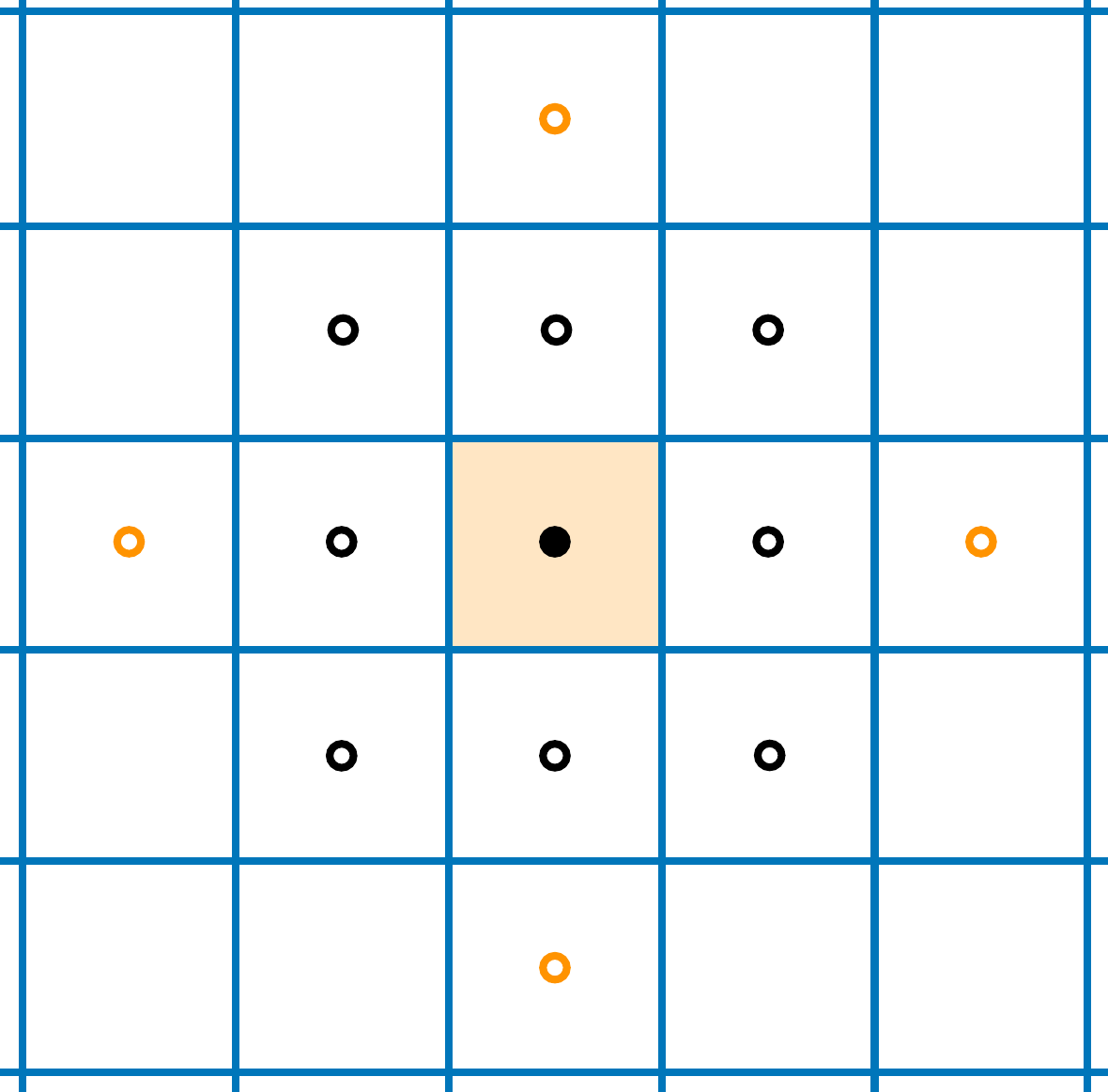}
	}
	\subfigure[Mesh over which the equations are integrated, as well as the primary mesh where the shape functions are defined (in gray)]
	{
		\includegraphics[height=0.47\textwidth]{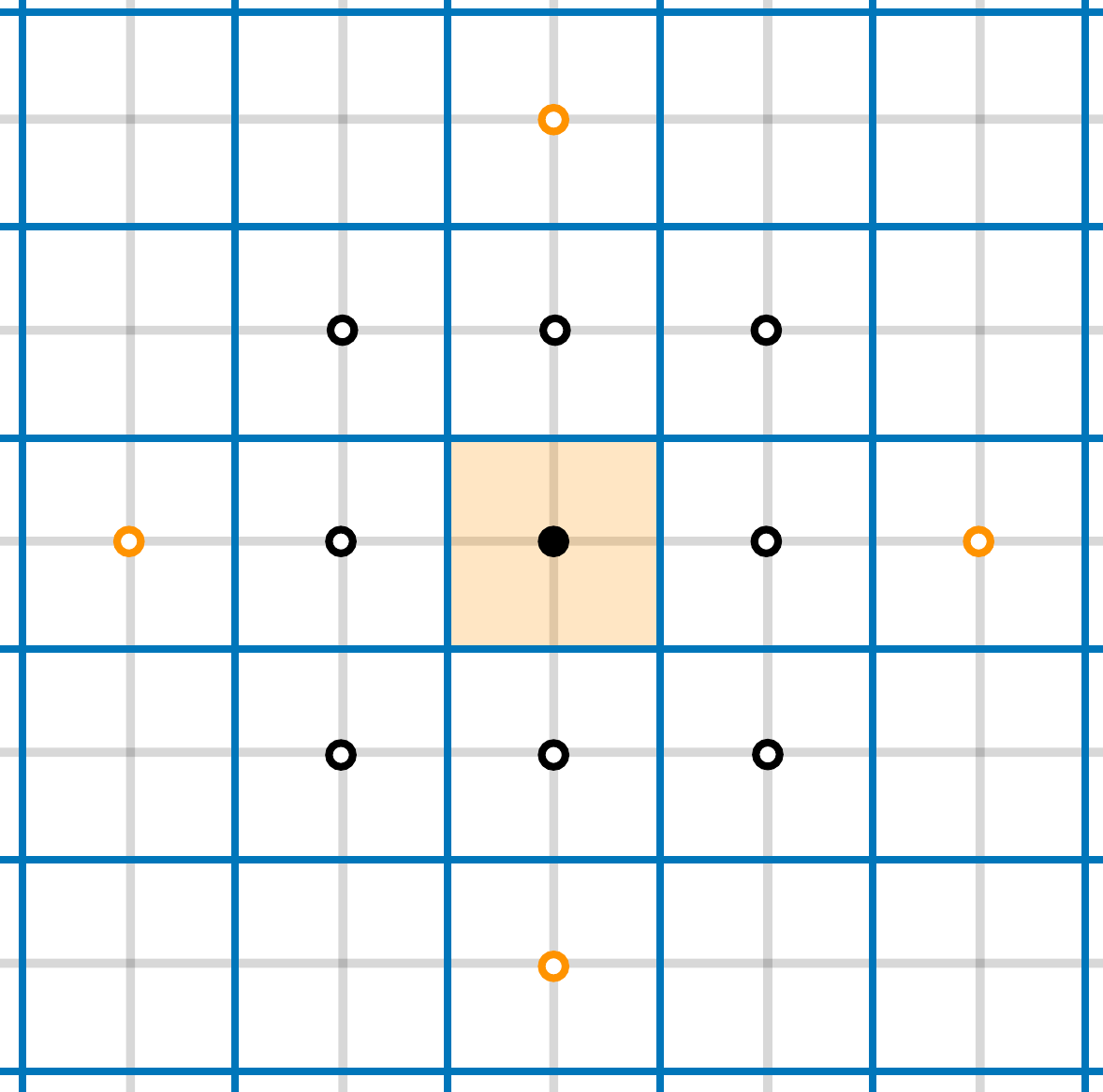}
	}
	\caption{Local integration domain (shaded orange) and computational stencil (white-filled dots) for a node (black-filled dot) in a 2-D quadrilateral mesh.
	 All methods include close neighbours (white-filled dots with black border) within their stencil, while method 1 additionally includes distant neighbours (white-filled dots with orange border).}
	\label{fig:finiteVolume2DQuadGrid}
\end{figure}

\begin{figure}[htb]
	\centering
	\subfigure[Method 1: Interpolate gradients from adjacent cell-centres]
	{
		\includegraphics[height=0.47\textwidth]{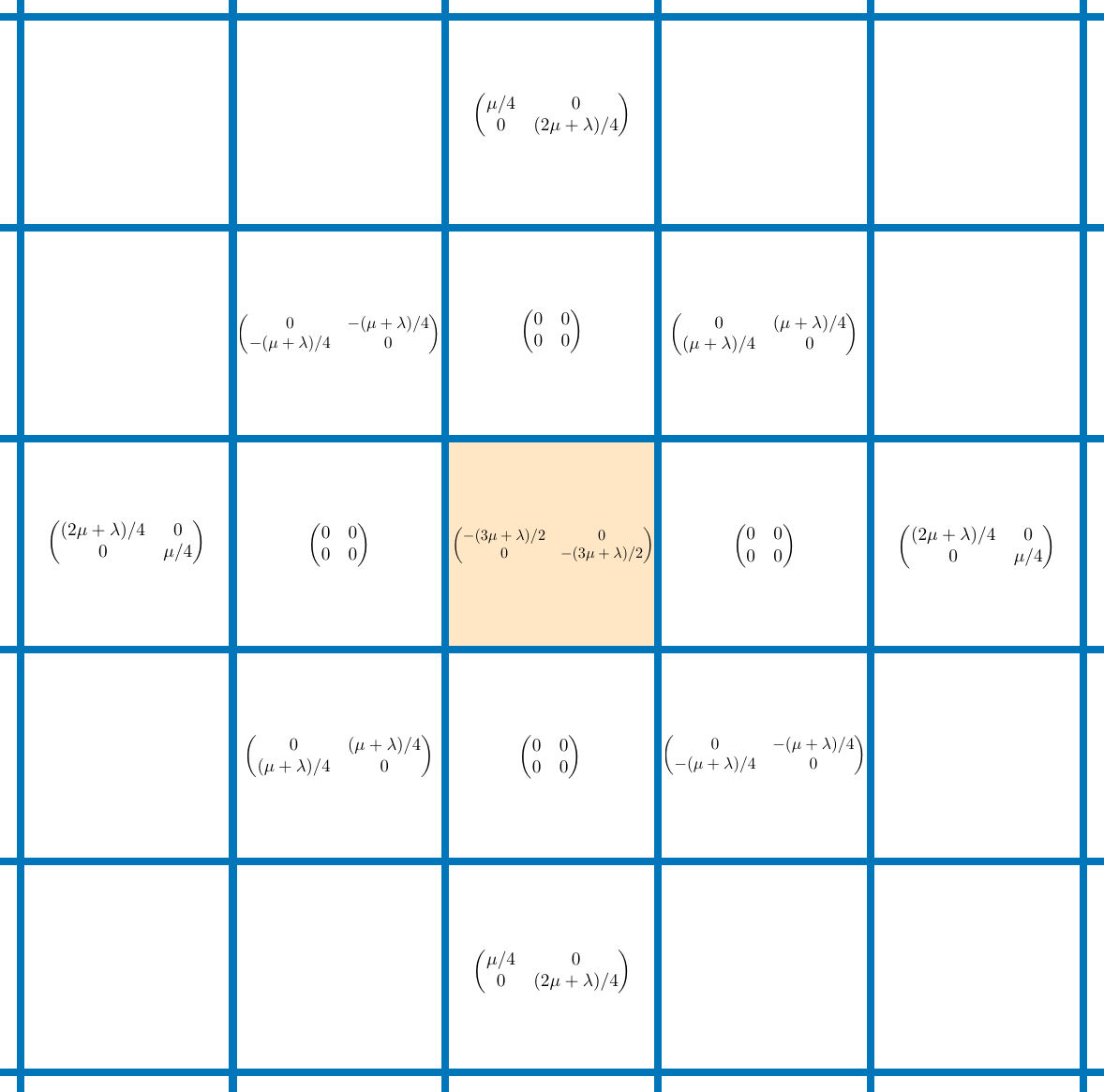}
	}
	\subfigure[Method 2: Calculate the normal gradient at the face using central-differencing, and interpolate the tangential gradient from the adjacent cell-centres]
	{
	    	\includegraphics[height=0.47\textwidth]{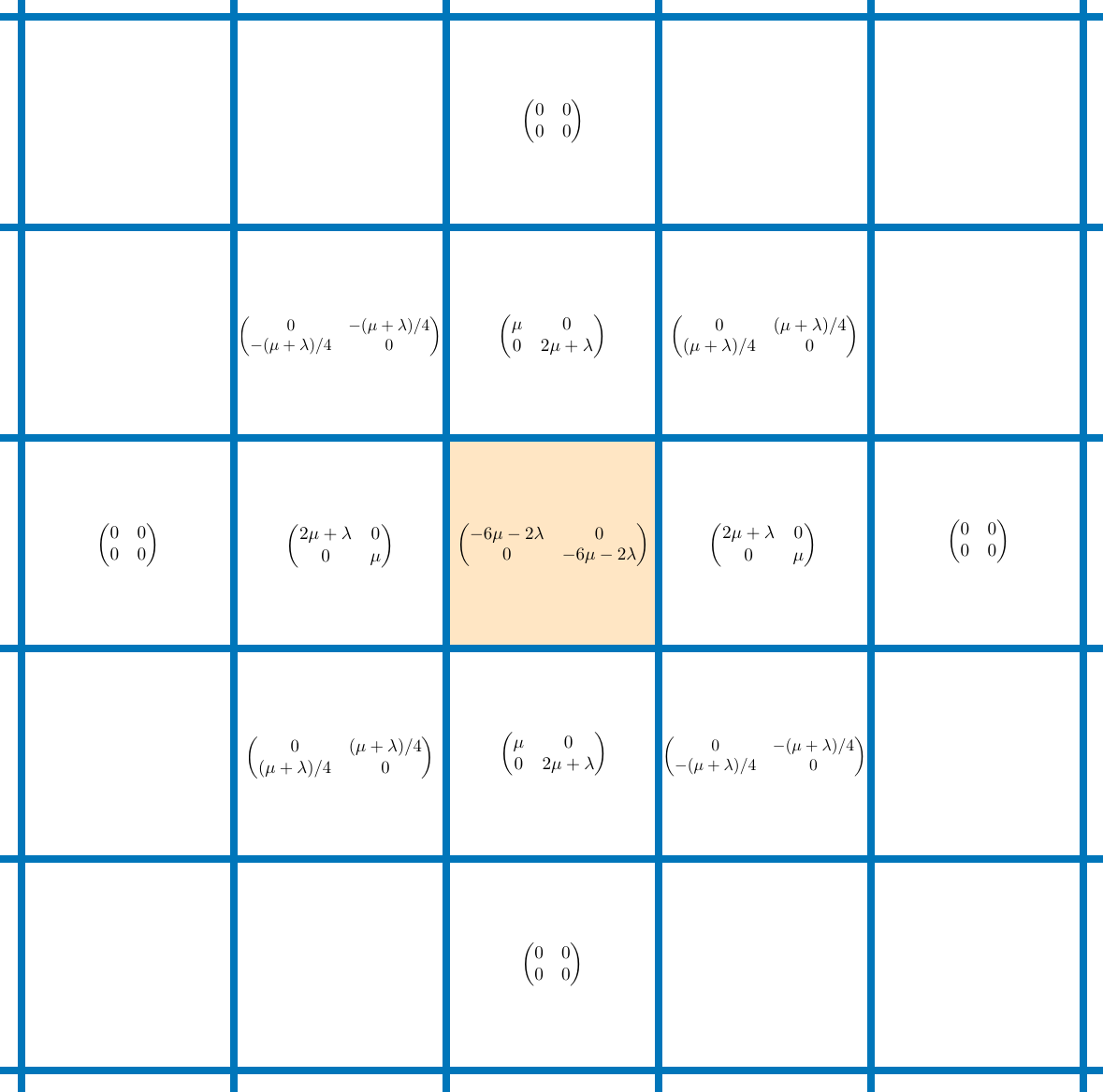}
	}
	\subfigure[Method 3: Calculate the normal gradient at the face using central-differencing, and calculate the tangential gradient at the face using the point or edge values]
	{
		\includegraphics[height=0.47\textwidth]{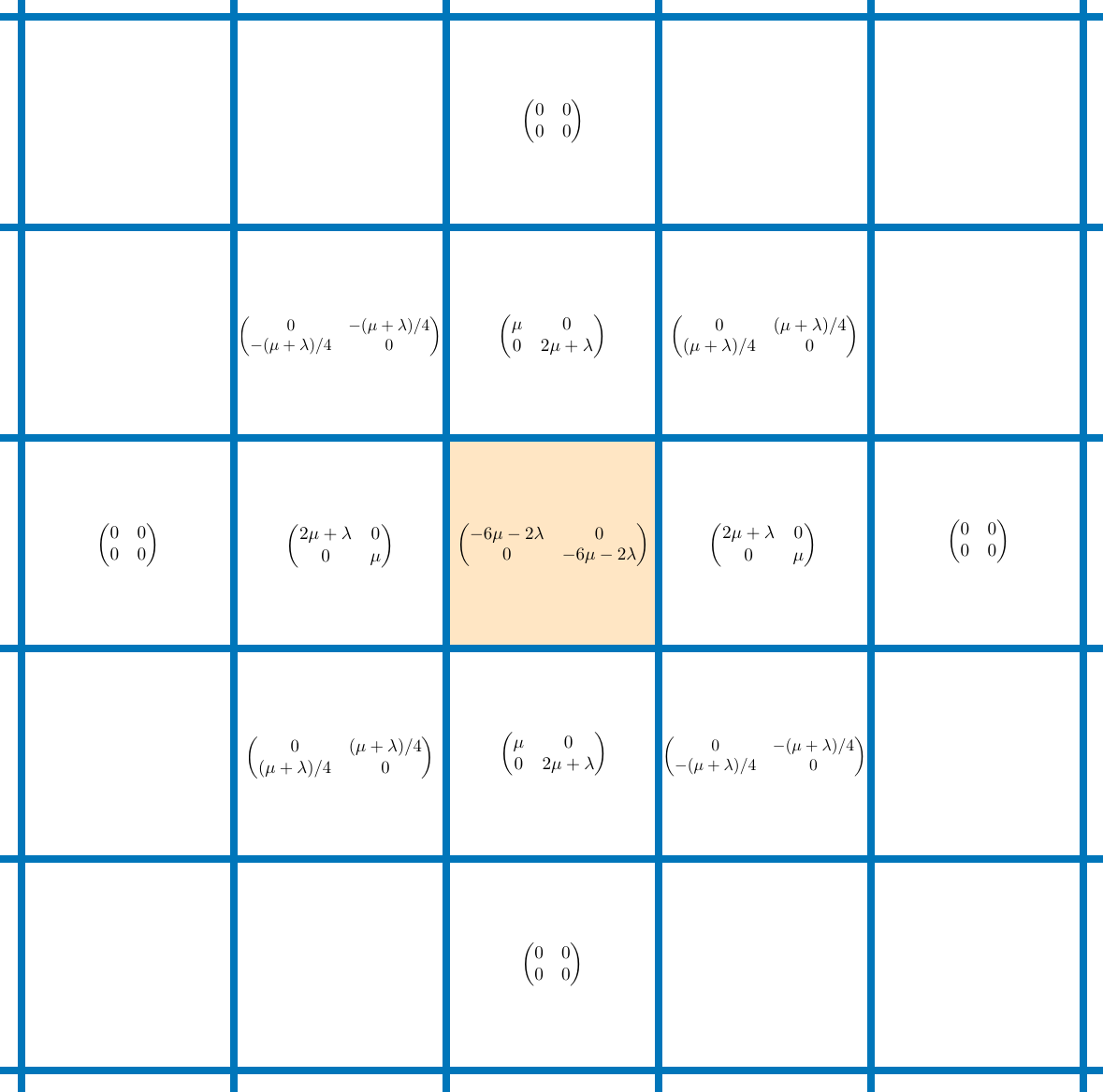}
	}
	\subfigure[Method 6: Use shape functions with reduced integration]
	{
		\includegraphics[height=0.47\textwidth]{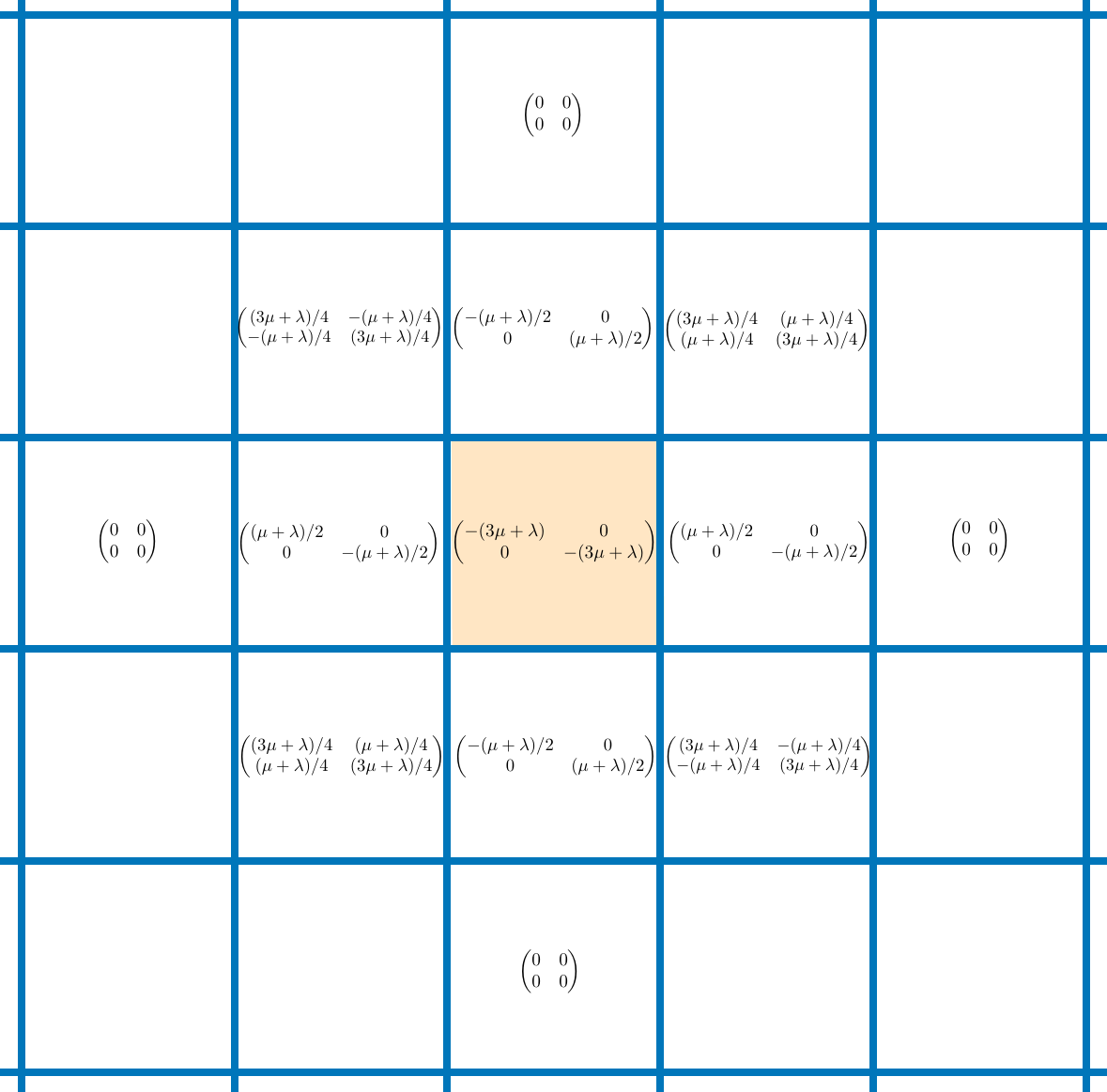}
	}
	\caption{A comparison of the block coefficients in the stiffness matrix for the centre node/cell, where different face gradient calculation methods are employed.
	A unit thickness is assumed and stabilisation terms have been disregarded.}
	\label{fig:faceGradientMethodCoeffs}
\end{figure}

To provide additional insight, let us consider the face gradient calculations on a simple grid.
Taking a uniform quadrilateral mesh (Figure \ref{fig:finiteVolume2DQuadGrid}), we will calculate the gradient at a face using the most popular of the methods above.
The discretised governing equation for the cell in Figure \ref{fig:finiteVolume2DQuadGrid} can be written in the form of a tensor algebraic equation:
\begin{align}
	\sum_{i=1}^{i=13}\boldsymbol{A}_i \cdot \boldsymbol{u}_i \quad = \quad \boldsymbol{0}
\end{align}
which corresponds to a block row in the resulting stiffness matrix, assuming a block-coupled solution methodology is used.
It is possible to present the coefficients $\boldsymbol{A}_i$ of this equation graphically to allow direct comparison between the different face gradient calculation approaches:
Figure \ref{fig:faceGradientMethodCoeffs} compares the resulting coefficients from using methods 1, 2, 3 and 6,
neglecting any stabilisation terms.
Noting that boundary conditions have not been considered, a number of observations can be made:
\begin{itemize}
	\item Methods 2, 3 and 6 employ the same compact computational stencil, whereas method 1 (interpolated gradient) uses a larger stencil, including second face-neighbours;
	\item All coefficients are symmetric;
	\item All methods show geometric symmetries, for example, the top-right coefficient is equal to the bottom-left coefficient; this is a consequence of the momentum/force being conserved between nodes;
	\item In methods 1, 2, and 3, the normal forces are calculated entirely from the central cell displacement as well as the left, right, top and bottom cell displacements (far-left/right/top/bottom cells in the case of method 1). Similarly, the shear forces are calculated entirely from the top-left, top-right, bottom-left and bottom-right cell displacements;
	\item Methods 2 (interpolated tangential gradient) and 3 (tangential gradient calculated at the face) are equivalent for this simple grid. For meshes including skewness and non-orthogonality, this may not be the case;
	\item Method 1 differs from methods 2 and 3 in only one way: the left, right, top and bottom cell coefficients have been moved to a more distant neighbour and scaled in magnitude; in addition, the central coefficient has been scaled in magnitude.
	This corresponds to the normal component being calculated using a larger stencil (interpolated gradient) in method 1 vs methods 2 and 3 (central differencing at the face);
	\item Method 6 (shape functions) requires the force to be integrated over 8 faces compared with 4 faces for methods 1, 2 and 3; this comes from the way in which the mesh is constructed from a primary mesh;
	\item Half the components of the corner coefficients are zero for methods 1, 2 and 3, while they are all non-zero for method 6.

\end{itemize}

\begin{figure}[htb]
	\centering
	\includegraphics[height=0.47\textwidth]{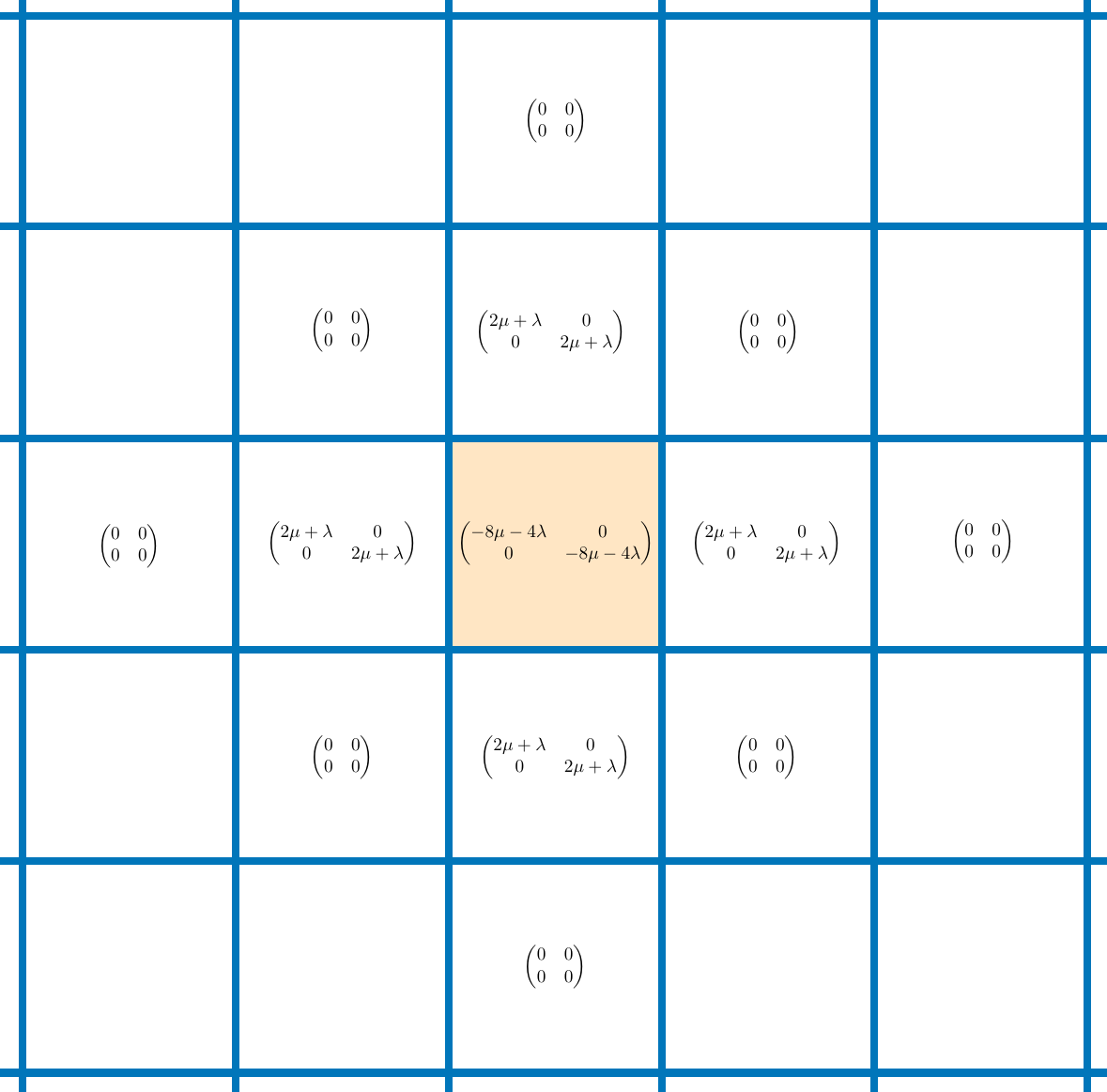}
	\caption{Stiffness matrix coefficients for methods 2 and 3 when a segregated solution algorithm is used; the inter-component coupling terms are included via the source vector in a deferred correction manner}
	\label{fig:faceGradientMethodCoeffsSegregated}
\end{figure}

It should be noted that the coefficients that appear in the linear system matrix depend on the chosen solution algorithm; for block-coupled approaches, the coefficients will be as shown, whereas for segregated approaches, the coefficients will be more sparse and will not contain inter-component couplings (the off-diagonal components of the coefficients will be zero).
For example, the matrix coefficients corresponding to method 2 (and equivalently method 3) are given for the segregated approach in Figure \ref{fig:faceGradientMethodCoeffsSegregated}.
In contrast to the coefficients given in Figure \ref{fig:faceGradientMethodCoeffs}(b), the segregated coefficients (Figure \ref{fig:faceGradientMethodCoeffsSegregated}) differ in the following ways:
\begin{itemize}
	\item As noted, all inter-component coupling is zero \ie all block off-diagonal coefficients are zero; this allows the two scalar displacement equations to be solved separately; 
	\item The matrix has greater sparsity that the other approaches, leading to reduced memory requirements;
	\item The coefficients produce a matrix which is \emph{weakly} diagonally dominant: the magnitude of the centre node coefficient is equal to the sum of the magnitudes of the other nodes.
	Such a weakly diagonally dominant system (which becomes strongly diagonally dominant with the inclusion of Dirichlet/essential boundary conditions) promotes the convergence of iterative linear solvers.
\end{itemize}


\paragraph{Stabilisation approach} \label{sec:stabilisationApproach}
For the discretisations that require stabilisation, three distinct forms of stabilisation term can be found;
these can be expressed as a stabilisation traction $\boldsymbol{t}^{\text{stab}}$ applied at the control volume face:
\begin{equation} \label{eq:stabilisationTractions}
\boldsymbol{t}^{\text{stab}} =
	\begin{cases}
	 \alpha^{\text{stab}} K_f 
	 \left(
	\dfrac{\boldsymbol{u}_N - \boldsymbol{u}_P}{|\boldsymbol{d}_{f}|}
	\;- \; \boldsymbol{n}_f \cdot
	\dfrac{
	|\boldsymbol{d}_{Nf}| \boldsymbol{\nabla}\boldsymbol{u}_N
	+ |\boldsymbol{d}_{Pf}| \boldsymbol{\nabla}\boldsymbol{u}_P}
	{|\boldsymbol{d}_{f}|}
	\right) \, |\boldsymbol{\Gamma}_f|
	& \text{Rhie-Chow} 
	\vspace{3mm} \\
	  -\alpha^{\text{stab}} \rho \, c_p
	  \, |\boldsymbol{d}_{f}|^2 \left(
	 \boldsymbol{\nabla}^2 \boldsymbol{v}_N - \boldsymbol{\nabla}^2 \boldsymbol{v}_P
	\right) \, |\boldsymbol{\Gamma}_f|
	& \text{Jameson-Schmidt-Turkel}
	\vspace{4mm} \\
	 \rho
	 \left[ c_s \textbf{I} + (c_p - c_s) \boldsymbol{n}_f\boldsymbol{n}_f \right] \cdot
         \dfrac{
         \boldsymbol{v}_{N} + \boldsymbol{d}_{Nf} \cdot \boldsymbol{\nabla} \boldsymbol{v}_N
         - \boldsymbol{v}_{P} - \boldsymbol{d}_{Pf} \cdot \boldsymbol{\nabla} \boldsymbol{v}_P
         }{2}
         \,|\boldsymbol{\Gamma}_f|
      & \text{Godunov upwinding}
	\end{cases}
\end{equation}
where $\alpha^{\text{stab}}$ is a user-defined scaling factor, 
the pressure wave speed of sound is $c_p = \sqrt{\frac{2\mu + \lambda}{\rho}}$,
and the shear wave speed of sound is $c_s = \sqrt{\frac{\mu}{\rho}}$.

In Equation \ref{eq:stabilisationTractions}, Jameson-Schmidt-Turkel and Godunov-upwinding terms are given in a form that is only suitable for dynamic problems;
however, it is possible to define similar stabilisation terms for quasi-static analyses:
\begin{equation} \label{eq:stabilisationTractions2}
\boldsymbol{t}^{\text{stab}} =
	\begin{cases}
	 \alpha^{\text{stab}} K_f 
	 \left(
		 \dfrac{\boldsymbol{u}_N - \boldsymbol{u}_P}{|\boldsymbol{d}_{f}|}
		\;-\; \boldsymbol{n}_f \cdot
		\dfrac{ |\boldsymbol{d}_{Nf}|\boldsymbol{\nabla}\boldsymbol{u}_N 
		+ |\boldsymbol{d}_{Pf}| \boldsymbol{\nabla}\boldsymbol{u}_P}
		{|\boldsymbol{d}_{f}|}
	\right) \, |\boldsymbol{\Gamma}_f|
	& \text{Rhie-Chow} 
	\vspace{3mm} \\
	  -\alpha^{\text{stab}}
	  K_f 
	  \, |\boldsymbol{d}_{f}| \left(
		 \boldsymbol{\nabla}^2 \boldsymbol{u}_N
	   - \boldsymbol{\nabla}^2 \boldsymbol{u}_P
	\right) \, |\boldsymbol{\Gamma}_f|
	& \text{Jameson-Schmidt-Turkel}
	\vspace{4mm} \\
	 \alpha^{\text{stab}}
	\dfrac{ K_f}{2}
     \dfrac{
	     \boldsymbol{u}_{N}
	     + \boldsymbol{d}_{Nf} \cdot \boldsymbol{\nabla} \boldsymbol{u}_N
	     -
	     \boldsymbol{u}_{P}
	     - \boldsymbol{d}_{Pf} \cdot \boldsymbol{\nabla} \boldsymbol{u}_P
	     }{|\boldsymbol{d}_{f}|}
      \,|\boldsymbol{\Gamma}_f|
      & \text{Godunov upwinding}
	\end{cases}
\end{equation}
where the Rhie-Chow term is given for comparison.
A user-defined scaling factor $\alpha^{\text{stab}}$ is added to the Godunov-type term in \ref{eq:stabilisationTractions2}, as this form of the term does not have the same physical significance as the dynamic term in Equation \ref{eq:stabilisationTractions}.
After some algebraic manipulation, where we include the $(\nicefrac{1}{2})$ factor in $\alpha^{\text{stab}}$ and note that $\boldsymbol{d}_{Nf} = -|\boldsymbol{d}_{Nf}|\boldsymbol{n}_{f}$,
the Godunov upwinding-type stabilisation term is seen to be identical to Rhie-Chow stabilisation.
This shows that even in its original dynamic form (Equation \ref{eq:stabilisationTractions}), it is in fact just a scaled version of Rhie-Chow stabilisation, and is equivalent for a specific choice of scaling parameters.

\begin{figure}[htb]
	\centering
	\subfigure[Rhie-Chow with a scale factor $\alpha^{\text{stab}} = 1$]
	{
		\includegraphics[height=0.47\textwidth]{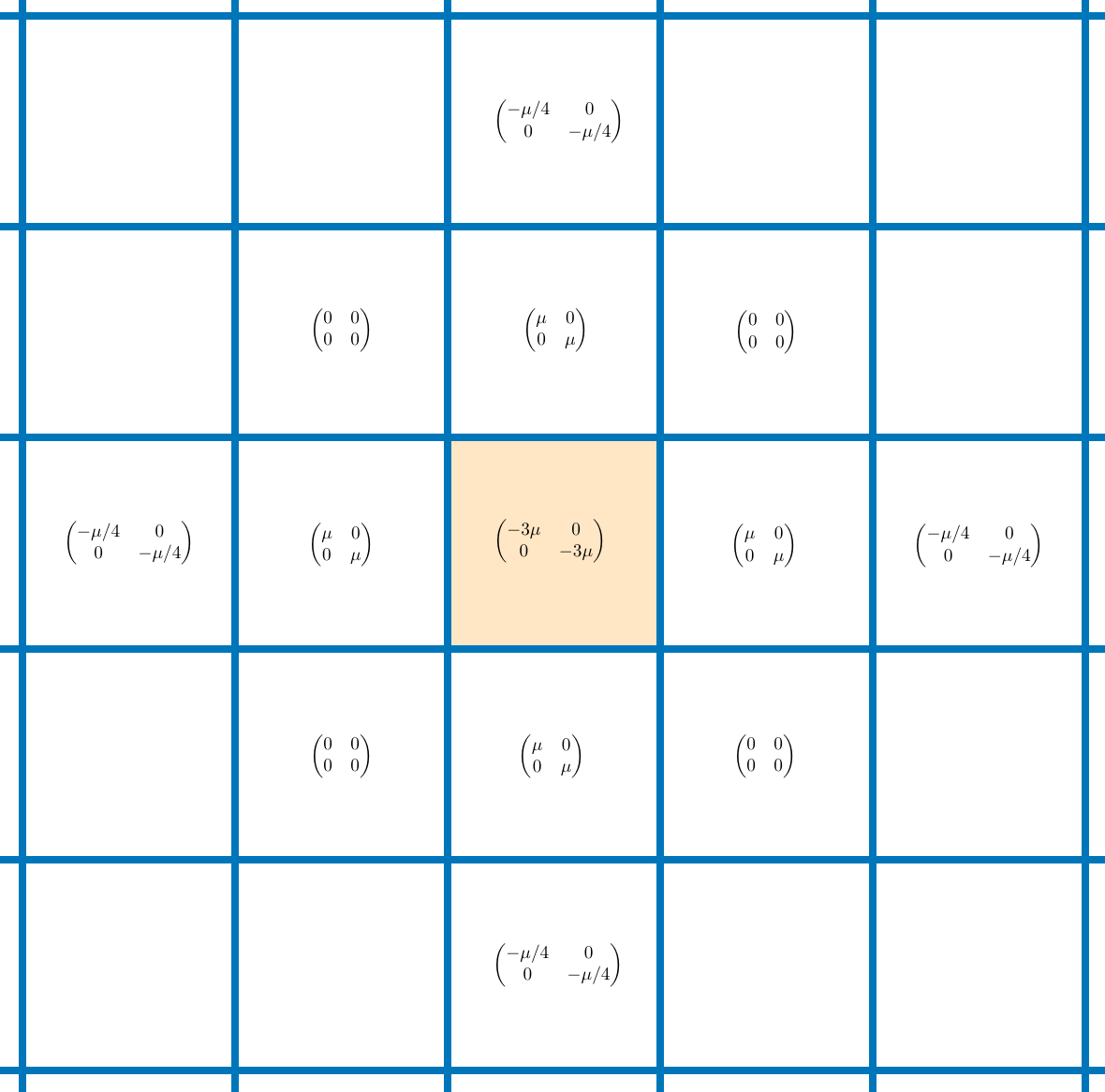}
	}
	\subfigure[Jameson-Schmidt-Turkel with a scale factor $\alpha^{\text{stab}} = (1/4)$]
	{
	    	\includegraphics[height=0.47\textwidth]{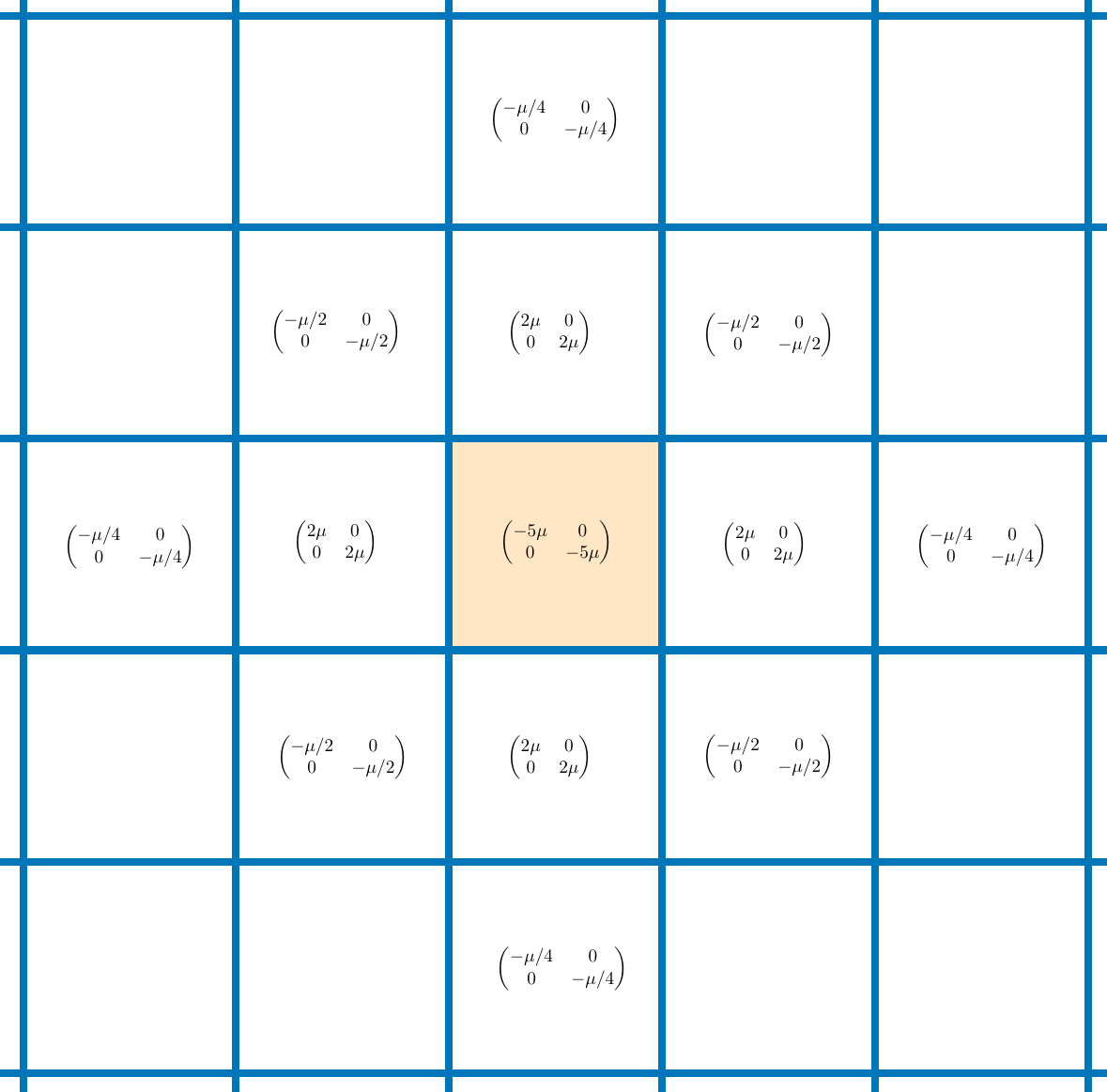}
	}
	\caption{A comparison of the block coefficients in the stiffness matrix for the centre node/cell for two styles of stabilisation term: Rhie-Chow and Jameson-Schmidt-Turkel}
	\label{fig:stabilisationCoeffs}
\end{figure}

\begin{figure}[htb]
	\centering
	\subfigure[Face gradient calculation method 1 with Rhie-Chow stabilisation, to be compared with Figure \ref{fig:faceGradientMethodCoeffs}(a)]
	{
		\includegraphics[height=0.47\textwidth]{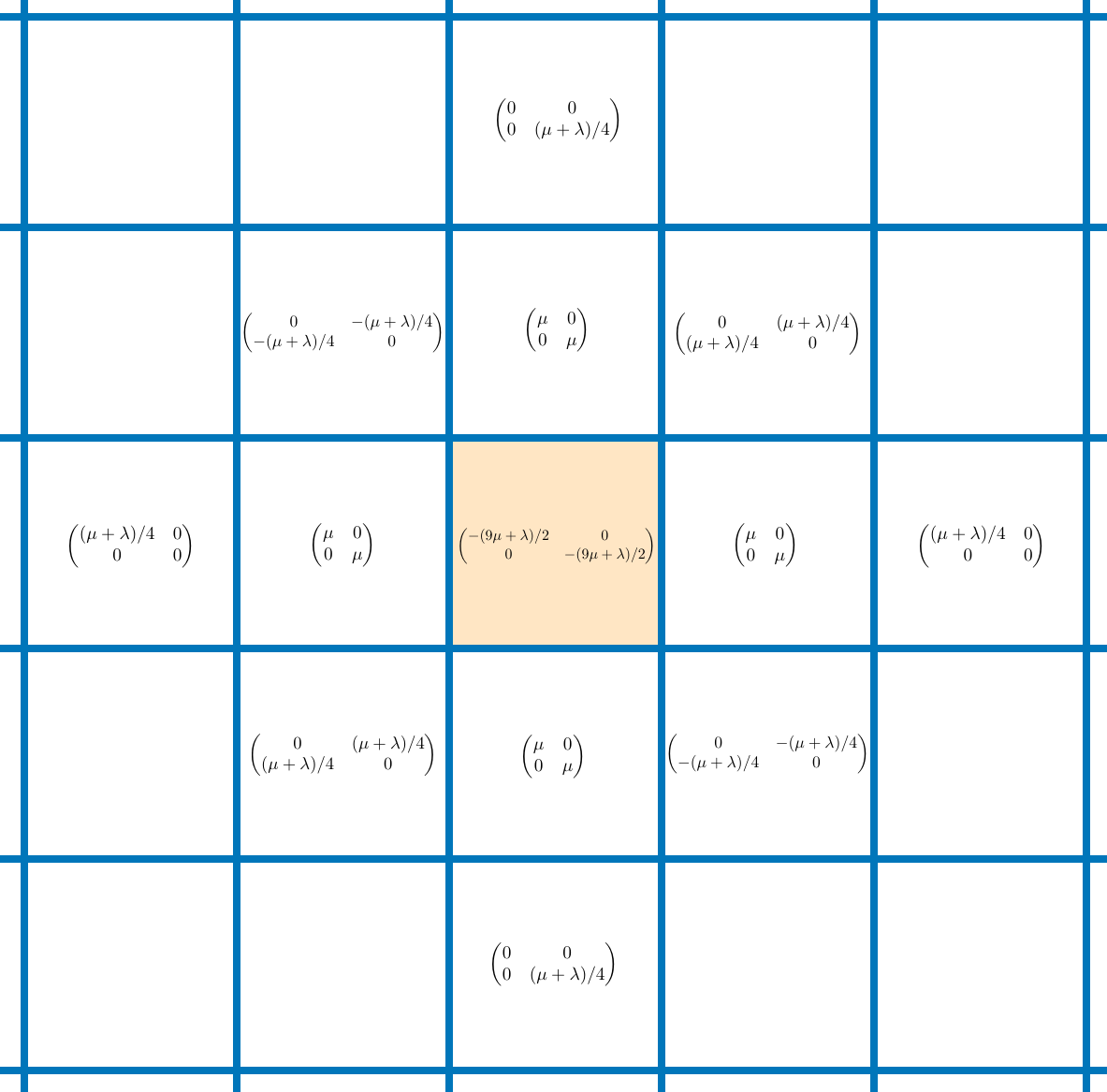}
	}
	\subfigure[Face gradient calculation method 2 with Jameson-Schmidt-Turkel stabilisation, to be compared with Figure \ref{fig:faceGradientMethodCoeffs}(b)]
	{
	    	\includegraphics[height=0.47\textwidth]{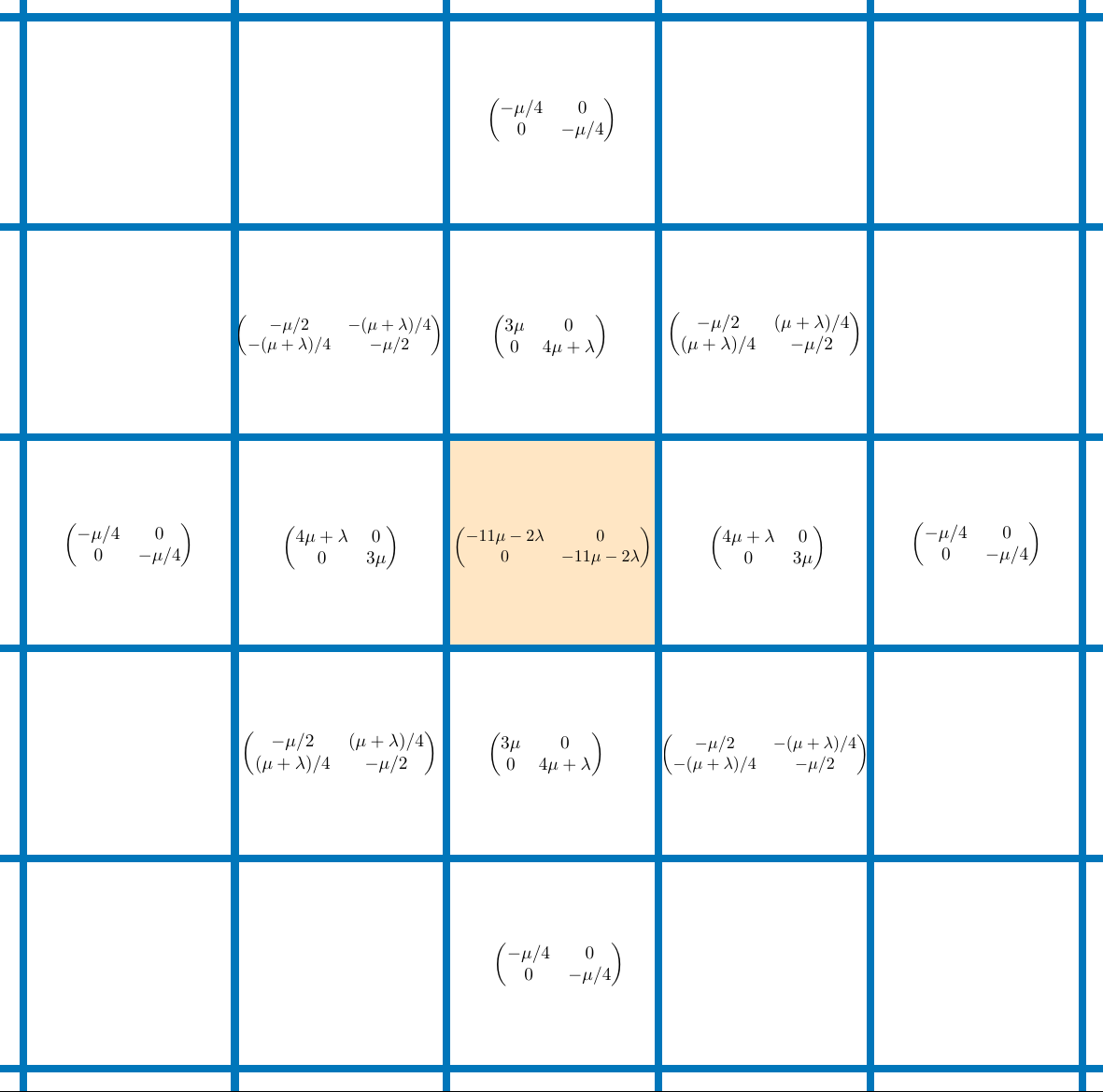}
	}
	\caption{The effect on the stiffness matrix coefficients by including stabilisation terms}
	\label{fig:faceGradientMethodCoeffsStabilisation}
\end{figure}

Using the 2-D square grid (Figure \ref{fig:finiteVolume2DQuadGrid}) as before, the computational stencil and coefficients resulting from the Rhie-Chow and Jameson-Schmidt-Turkel stabilisation terms can be graphed (Figure \ref{fig:stabilisationCoeffs}).
Scale factors of $\alpha^{\text{stab}} = (1/4)$ for the Jameson-Schmidt-Turkel term and $\alpha^{\text{stab}} = 1$ for the Rhie-Chow term are chosen so that the magnitude of the terms is similar.
For this grid, both approaches produce similar coefficients, however, the Jameson-Schmidt-Turkel approach differs in that it includes additional coupling in the corner coefficients.
It can also be seen that both approaches require a large computational stencil, in that second face-neighbours are needed.
When either of these stabilisation approaches is combined with one of the face gradient calculations methods discussed above, the stencil of coefficients are summed.
For example, the  computational stencil and coefficients from face gradient calculation method 1 with Rhie-Chow stabilisation are shown in Figure \ref{fig:faceGradientMethodCoeffsStabilisation}(a), and from face gradient calculation method 2 with Jameson-Schmidt-Turkel stabilisation in Figure \ref{fig:faceGradientMethodCoeffsStabilisation}(b).

\begin{figure}[htb]
	\centering
	\includegraphics[height=0.47\textwidth]{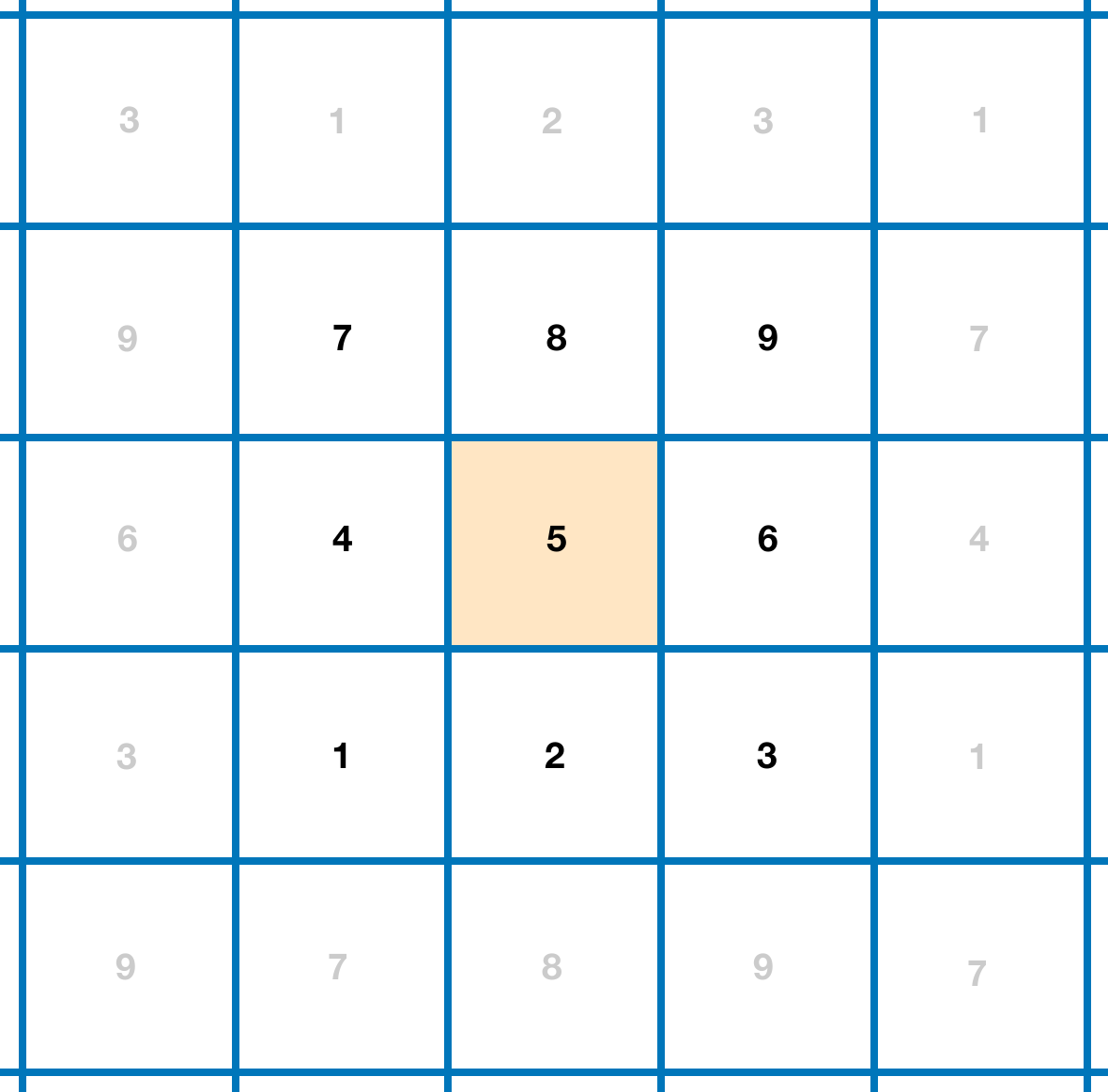}
	\caption{A periodic patch of finite volume cells, containing a central cell (cell 5) and all neighbour cells within its computational stencil.
	The periodic/cyclic neighbours are indicated as cells in grey numbering.}
	\label{fig:periodicFiniteVolume}
\end{figure}

Within literature, a number of authors have described stabilisation techniques, however, stability analysis of finite volume discretisations for solid mechanics is not common.
Of course, the stability of a formulation quickly becomes apparent in use, however, we can also take inspiration from finite element approaches \citep{Belytschko2014} and analyse the stiffness matrix eigenvalues.
Taking a single unconstrained finite element, the number of zero eigenvalues of the stiffness matrix indicates the number of zero energy modes.
For a stable formulation, the number of zero valued eigenvalues is equal to the number of rigid degrees of freedom; in 3-D, there are three rigid translations and three rigid rotations, whereas in 2-D there are two rigid translations and one rigid rotation.
For an unstable formulation, there will be additional zero valued eigenvalues, where the corresponding eigenvector indicates the unstable mode.
For the finite volume method, an equivalent analysis of an individual cell/element is less obvious.
One such approach is to consider a periodic patch of finite volume cells, containing a central cell and all neighbour cells within its computational stencil (Figure \ref{fig:periodicFiniteVolume}), and analyse the eigenvalues of its (block-coupled) global stiffness matrix.
In this case, as there are 9 cells, each with 2 degrees of freedom, the global stiffness matrix is $18 \times 18$.
As each cells contains all other eight cells in their stencil (due to the periodic conditions), the stiffness matrix is fully dense (no zero block entries).
In this case, using face gradient calculation method 2, the stiffness matrix contains 2 zero eigenvalues, corresponding to the two rigid translation directions; the periodic conditions prohibit rigid rotation;
this indicates that the discretisation is stable in this configuration.
In addition to analysing a periodic patch of \emph{internal} cells, it may also be necessary to examine a patch of cells adjacent to a boundary.
The discretisation at boundary faces is typically different than at internal faces, and so boundaries may quell or excite spatial instabilities.

\paragraph{Solution methodology}
Like other popular numerical methods, the finite volume method can employ implicit or explicit solution algorithms.
The relative merits of implicit vs explicit approaches are independent of the finite volume method;
the interested reader can find numerous textbooks addressing this topic,
for example, \citep{Bathe1996, Belytschko2014, Schafer2006}.

\begin{figure}[htb]
	\centering
	\includegraphics[width=0.95\textwidth]{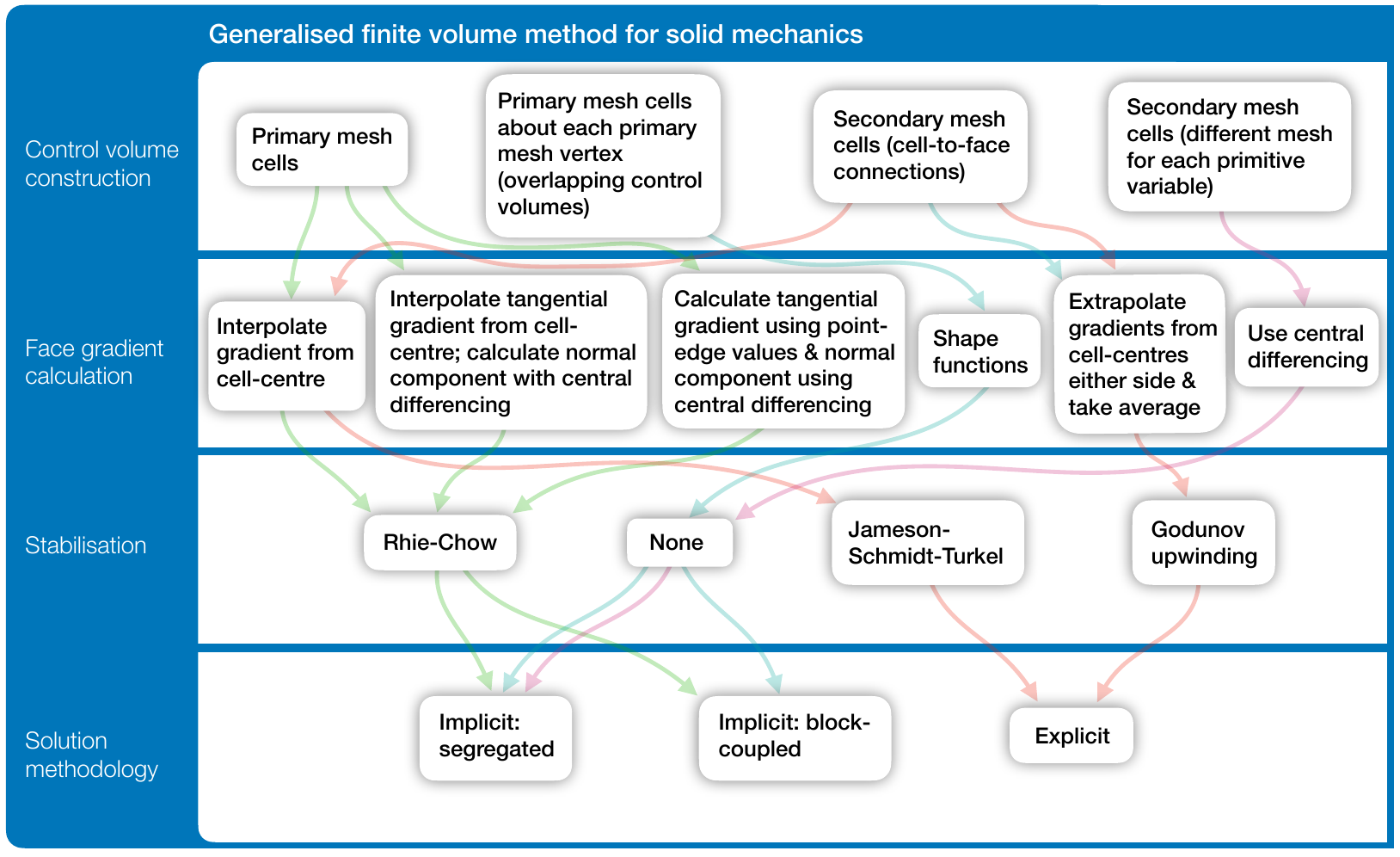}
	\caption{Generalisation of the finite volume method for solid mechanics, allowing the creation of any scheme through appropriate selection of four key components: 1) control volume construction method, 2) face gradient calculation approach, 3) stabilisation technique, and 4) solution methodology.
	Common approaches are indicated by the coloured lines: implicit cell-centred (green), implicit vertex-centred (blue), explicit Godunov-type (red), and staggered-grid (purple).}
	\label{fig:generalisedFiniteVolumeMethod}
\end{figure}

\paragraph{A generalised finite volume method for solid mechanics}
Based on these four components, it is possible to describe a generalised approach encompassing all individual variants, as shown in Figure \ref{fig:generalisedFiniteVolumeMethod}.
Common approaches are indicated by the coloured lines: implicit cell-centred (green), implicit vertex-centred (blue), explicit Godunov-type (red), and staggered-grid (purple).
This figure allows the relationship between the variants to be concisely expressed.
What is also apparent from this figure is that there are avenues that have yet to be fully explored.

\section{Comparing of the finite volume method for computational solid mechanics with the finite element method}
\label{sec:compareFVandFE}


Within this section, the finite volume method for solids mechanics is compared with the ``standard'' 
continuous Bubnov-Galerkin finite element method, for example, as described by \citet{Bathe1996}, \citet{Zienkiewicz2000}, and \citet{Belytschko2014}.
%
Following the approach taken in the previous section, the finite element method will be compared to the finite volume method in terms of: (a) discretisation of space and time; (b) discretisation of the mathematical model equations; and (c) solution algorithm.


\subsection{Discretisation of time and space}
\begin{figure}[htb]
	\centering
	\includegraphics[width=0.8\textwidth]{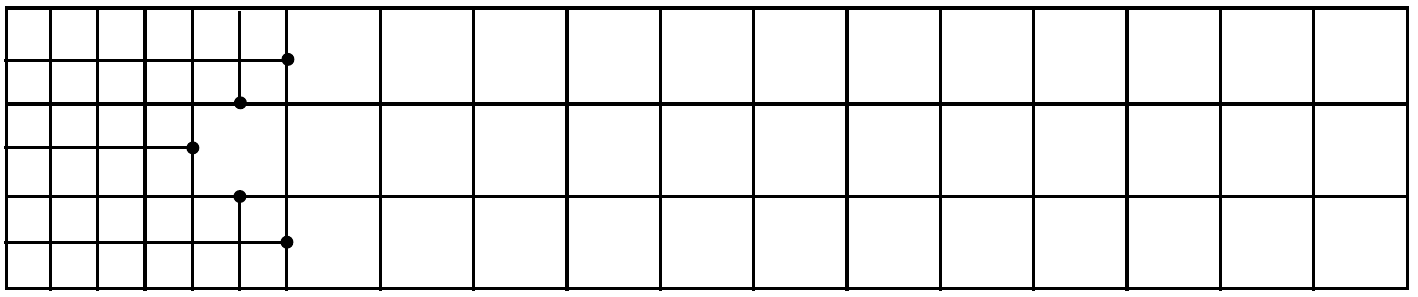}
	\caption{2-D quadrilateral mesh of a beam, where `hanging nodes' are shown as black dots.
	}
	\label{fig:hangingNodes}
\end{figure}

Like the finite volume method, the finite element method follows the standard time-marching temporal discretisation approach.
Similarly, the solution domain space is divided into a finite number of convex cells (or elements) that do not overlap and fill the space completely.
With regard to the mesh, the following differences can be noted between the finite element method and each of the finite volume variants:
\begin{itemize}
	\item The cell-centred approach and vertex-centred approaches which \emph{do not} use shape functions are applicable to general convex polyhedral in 3-D and general polygons in 2-D, whereas the standard finite element method is limited to standard element shapes, such as hexahedra/tetrahedra in 3-D and quadrilaterals/triangles in 2-D.
	A consequence of this it that `hanging nodes' (Fig. \ref{fig:hangingNodes}) are common in finite volume analyses but not directly possible with the standard finite element method;
	\item The vertex-centred (and hypothetically cell-centred) approaches which \emph{do} use shape function are limited to the same types of meshes used by the standard finite element method;
\end{itemize}


\subsection{Discretisation of the mathematical model equations} \label{sec:discretisationFE}
In Section \ref{sec:mathModel}, the conservation of linear momentum in strong integral form (Equation \ref{eq:momentumIntegral}) was taken as the starting point for the finite volume discretisation.
In contrast, the finite element method requires the weak form of the governing equation, and hence begins with the strong differential form:
\begin{eqnarray} \label{eq:momentumDifferential}
	\rho \frac{\partial^2 \boldsymbol{u}}{\partial t^2}
	&=&
	\boldsymbol{\nabla} \cdot
	\left[
	\mu \boldsymbol{\nabla} \boldsymbol{u}
	+ \mu (\boldsymbol{\nabla} \boldsymbol{u})^T
	+ \lambda \; \text{tr}(\boldsymbol{\nabla} \boldsymbol{u}) \textbf{I}
	\right]
	\;+\;
	\rho \boldsymbol{f}_b
\end{eqnarray}
This form is then multiplied by an arbitrary continuous weighting function $\boldsymbol{\omega}$ and integrated over the material volume to give the conservation of momentum in weak form:
\begin{eqnarray} \label{eq:momentumWeakIntegral}
	\int_\Omega
	\boldsymbol{\omega} \cdot
	\left\{
		\rho \frac{\partial^2 \boldsymbol{u}}{\partial t^2}
		\;-\;
		\boldsymbol{\nabla} \cdot
		\left[
		\mu \boldsymbol{\nabla} \boldsymbol{u}
		+ \mu (\boldsymbol{\nabla} \boldsymbol{u})^T
		+ \lambda \; \text{tr}(\boldsymbol{\nabla} \boldsymbol{u}) \textbf{I}
		\right]
		\;-\;
		\rho \boldsymbol{f}_b
	\right\}
	\;\text{d}\Omega
	&=& \boldsymbol{0}
\end{eqnarray}
Unsurprisingly a finite volume method is recovered if the weighting functions $\boldsymbol{\omega}$ are taken as unity within the control volumes and zero elsewhere.

To derive the finite element method, the weak form (Equation \ref{eq:momentumWeakIntegral}) is rearranged using integration by parts combined with the Gauss divergence theorem:
\begin{eqnarray} \label{eq:virtualWork}
	\int_\Omega
	\rho \boldsymbol{\omega} \cdot \frac{\partial^2 \boldsymbol{u}}{\partial t^2}
	\;\text{d}\Omega
	\;+\;
	\int_\Omega
	\boldsymbol{\nabla} \boldsymbol{\omega} :
	\overbrace{
	\left[
	\mu \boldsymbol{\nabla} \boldsymbol{u}
	+ \mu (\boldsymbol{\nabla} \boldsymbol{u})^T
	+ \lambda \; \text{tr}(\boldsymbol{\nabla} \boldsymbol{u}) \textbf{I}
	\right]
	}^{\boldsymbol{\sigma}}
	\;\text{d}\Omega
	&=& \notag \\
	\oint_{\Gamma_T}
	\boldsymbol{\omega}
	\cdot
	\boldsymbol{T}_\Gamma
	\;\text{d}\Gamma
	\;+\;
	\int_\Omega
	\rho \boldsymbol{\omega} \cdot \boldsymbol{f}_b
	\;\text{d}\Omega
\end{eqnarray}
where $\Gamma_T$ is the region of the domain boundary where tractions $\boldsymbol{T}_\Gamma$ are applied.
In the second term on the left-hand side of Equation \ref{eq:virtualWork}, the differential operator now acts on the weighting function $\boldsymbol{\omega}$, in contrast to Equation \ref{eq:momentumIntegral} where it acts on the stress tensor $\boldsymbol{\sigma}$.
It is assumed that the weighting functions, also known as test functions, satisfy the following requirements \citep{Belytschko2014}:
\begin{itemize}
	\item They are \emph{not} functions of time;
	\item They are $C^0$ continuous;
	\item They vanish on displacement boundaries.
\end{itemize}
Interpreting the weighting functions $\boldsymbol{\omega}$ as \emph{virtual displacements} allows Equation \ref{eq:virtualWork} to be viewed as the principle of virtual work.
In this way, the original problem of trying to find a displacement field which satisfies the conservation of linear momentum can be reinterpreted as trying to find a displacement field which minimises the total energy.

The finite element method then assumes the displacement $\boldsymbol{u}$ within each mesh element to vary according to \emph{shape functions}, as presented previously in Equation \ref{eq:shapeFunctions}.
The standard \emph{Bubnov-Galerkin} form of the finite element method is achieved by assuming that the weighting functions $\boldsymbol{\omega}$ in Equation \ref{eq:virtualWork} are approximated using the same shape functions as the displacement field:
\begin{align} \label{eq:shapeFunctionsWeighting}
	\boldsymbol{\omega}
	 &= \sum_{A=1}^{nVertices} N_A  \boldsymbol{\bar{\omega}}_A 
\end{align}
where $\boldsymbol{\bar{\omega}}_A$ represents the discrete weighting function value at vertex $A$;
uppercase letter $A$ is used as an index for vertices to avoid confusion with index notation $i$, $j$ and $k$.

Combining Equations \ref{eq:shapeFunctions} (displacement shape functions), \ref{eq:shapeFunctionsWeighting} and \ref{eq:virtualWork}, the finite element equations are expressed as \citep{Laursen2002}:
\scriptsize
\begin{align} \label{eq:virtualWorkDiscrete}
	&\underbrace{
	\int_\Omega
	\rho
	\sum_{A=1}^{nVertices} N_A  \boldsymbol{\bar{\omega}}_A 
	\; \cdot \;
	  \sum_{B=1}^{nVertices} N_B  \frac{\partial^2 \boldsymbol{u}_B}{\partial t^2} 
	\;\text{d}\Omega}_{\text{Inertial force term}}
	\notag \\
	&\quad\;+\;
	\underbrace{
	\int_\Omega
	\sum_{A=1}^{nVertices} \boldsymbol{\nabla} N_A  \boldsymbol{\bar{\omega}}_A 
	\; : \;
	\left[
	\mu 
	\sum_{B=1}^{nVertices} \boldsymbol{\nabla} N_B  \boldsymbol{u}_B 
	+ \mu 
	 \sum_{B=1}^{nVertices} \boldsymbol{u}_B \boldsymbol{\nabla} N_B 
	+ \lambda \; \text{tr}\left( \sum_{B=1}^{nVertices} \boldsymbol{\nabla} N_B  \boldsymbol{u}_B \right) \textbf{I}
	\right]
	\;\text{d}\Omega
	}_{\text{Internal force term}}
	\notag \\
	&\quad\;=\; 
	\underbrace{
	\oint_{\Gamma_T}
	\sum_{A=1}^{nVertices} N_A  \boldsymbol{\bar{\omega}}_A 
	\; \cdot \;
	\boldsymbol{T}_\Gamma
	\;\text{d}\Gamma 
	}_{\text{External (surface) force term}}
	\;+\;
	\underbrace{
	\int_\Omega
	\sum_{A=1}^{nVertices} N_A  \boldsymbol{\bar{\omega}}_A 
	\;\cdot \;
	\rho \boldsymbol{f}_b
	\;\text{d}\Omega
	}_{\text{External (body) force term}}
\end{align}
\normalsize

To proceed, we will switch to index notation, not referring to spatial directions $i$ and $j$ or nodes $A$ and $B$, but instead to the degrees of freedom \citep{Laursen2002}.
For this, the index for a global degree of freedom $P$ can be given uniquely in terms of the spatial index and node index:
\begin{align}
	P = f\left( i, A\right)
\end{align}
 
Examining Equation \ref{eq:virtualWorkDiscrete} term by term and using both index notation and degree of freedom notation, the inertial force term can be expanded as:
\begin{align}
	&\int_\Omega
	\rho
	\sum_{A=1}^{nVertices} N_A  \boldsymbol{\bar{\omega}}_A 
	 \; \cdot
	 \sum_{B=1}^{nVertices} N_B  \frac{\partial^2 \boldsymbol{u}_B}{\partial t^2} 
	\;\text{d}\Omega \notag \\
	&\quad =
	\sum_{A=1}^{nVertices}
	\sum_{i=1}^{nDims}
	\int_\Omega
	\rho
	N_A  \bar{\omega}_{iA}
	\sum_{B=1}^{nVertices} N_B  \frac{\partial^2 U_{iB}}{\partial t^2} 
	\;\text{d}\Omega \notag \\
	&\quad =
	\sum_{A=1}^{nVertices}
	\sum_{i=1}^{nDims}
	\bar{\omega}_{iA}
	\sum_{B=1}^{nVertices} 
	\sum_{j=1}^{nDims}
	\int_\Omega
	\rho
	N_A  
	\delta_{ij}
	N_B  \frac{\partial^2 U_{iB}}{\partial t^2}
	\;\text{d}\Omega 
	\notag \\
	&\quad =
	\sum_{P=1}^{nDoF}
	\bar{\omega}_{P}
	\sum_{Q=1}^{nDoF} 
	M_{PQ}
	\frac{\partial^2 U_{Q}}{\partial t^2}
	\;\text{d}\Omega 
\end{align}
where $nDims$ is the number of spatial dimensions, $\delta_{ij}$ is the Kronecker delta, and $nDoF = nDims \times nVertices$ is the number of global degrees of freedom.
Matrix $M_{PQ}$ is known as the mass matrix and is given as:
\begin{align} \label{eq:massMatrix}
	M_{PQ} =	\int_\Omega \rho N_A  \delta_{ij} N_B  \;\text{d}\Omega 
\end{align}
The volume integral in Equation \label{eq:massMatrix} can be calculated (or approximated) using Gaussian/numerical quadrature:
\begin{align} \label{eq:numericalQuadratureMassMatrix}
	M_{PQ} \approx 
	\sum_{p=1}^{nQuadPoints} w_p\, \rho(\xi_p)\, N_A(\xi_p)\,  \delta_{ij}\, N_B(\xi_p) 
\end{align}
where $nQuadPoints$ is the number of quadrature points, $w$ is the quadrature weight, and $\xi$ is the quadrature location within the element. For example, using one-point quadrature, the quadrature locations are situated at the centroids of the elements; in that case, the density and the derivative of the shape functions are evaluated at element centroids.

Unlike finite volume approaches, the resulting mass matrix will not in general be diagonal; however, diagonalisation or \emph{lumping} of the mass matrix is common, although often ad hoc \citep{Belytschko2014}; for example, the row-sum technique calculates the diagonal elements as the sum of the coefficients for that row.

The internal force term in Equation \ref{eq:virtualWorkDiscrete} can be simplified as:
\scriptsize
\begin{align}
	&\int_\Omega
	\sum_{A=1}^{nVertices} \boldsymbol{\nabla} N_A  \boldsymbol{\bar{\omega}}_A 
	\; : \;
	\left[
	\mu 
	 \sum_{B=1}^{nVertices} \boldsymbol{\nabla} N_B  \boldsymbol{u}_B 
	+ \mu
	\sum_{B=1}^{nVertices} \boldsymbol{u}_B \boldsymbol{\nabla} N_B 
	+ \lambda \; \text{tr}\left( \sum_{B=1}^{nVertices} \boldsymbol{\nabla} N_B  \boldsymbol{u}_B \right) \textbf{I}
	\right]
	\;\text{d}\Omega \notag \\
	&\quad=
	\int_\Omega
	\sum_{A=1}^{nVertices} \sum_{i=1}^{nDims}  \sum_{j=1}^{nDims} 
	N_{A,j}  \bar{\omega}_{iA}
	\left(
	\mu N_{B,j} U_{iB} + \mu  N_{B,i} U_{jB} + \lambda N_{B,k}  U_{kB} \delta_{ij}
	\right)
	\;\text{d}\Omega \notag \\
	&\quad=
	\sum_{P=1}^{nDoF} \bar{\omega}_{P}  K_{PQ} U_{Q} 
\end{align} \normalsize
where matrix $K_{PQ}$ is known as the stiffness matrix:
\begin{align}
	K_{PQ} &=
	\int_\Omega
	\sum_{j=1}^{nDims} 
	N_{A,j}  \left[ \mu N_{B,j} + \mu  N_{B,i} + \lambda N_{B,i} \delta_{ij} \right] 
	\;\text{d}\Omega
\end{align}
The integral can once again be evaluated using Gaussian/numerical quadrature:
\begin{align} \label{eq:numericalQuadratureStiffnessMatrix}
	K_{PQ} \approx 
	\sum_{p=1}^{nQuadPoints}
	\sum_{j=1}^{nDims} 
	w_p\,
	N_{A,j}(\xi_p)\,  \left[ \mu(\xi_p)\, N_{B,j}(\xi_p)\, + \mu(\xi_p)\,  N_{B,i}(\xi_p)\, + \lambda(\xi_p)\, N_{B,i}(\xi_p)\, \delta_{ij} \right] 
\end{align}

Evaluating this integral exactly (using sufficient quadrature points) would naively appear to be the best approach;
however, formulations that use \emph{full} integration tend to suffer from locking \citep{Belytschko2014}, which can be described as an overly stiff behaviour in bending.
Consequently, \emph{reduced} integration is often favoured, where the local field is under-integrated.
This reduced integration has the benefit of relieving this locking phenomena as well as reducing the  computational time, due to the lower number of integration points.

A downside of reduced integration is the introduction of spatial instabilities into the discretisation.
In essence, the element is capable of deforming in certain \emph{modes} which offer no resistance.
These spurious singular or zero energy modes produce an accordion-like deformation pattern known as hourglassing.
Similar to the finite volume method, a stabilisation term is included in the formulation to suppress such spatial instabilities.
For the finite element method, this hourglass stabilisation is incorporated through the inclusion of a stabilisation stiffness within the element stiffness matrix, or a viscous stabilisation term for dynamic problems \citep{Belytschko2014}.
For static problems a variety of stabilisation methods have been proposed, but typically a \emph{stabilisation stiffness} term $K_{PQ}^{\text{stab}}$ is added to the element stiffness matrix of the form \citep{Belytschko2014, Abaqus}:

\begin{align}
	K_{PQ}^{\text{stab}} = 
	\alpha_{\text{stab}} \, \mu \, N_{A,j}^{\text{stab}} \, N_{B,j}^{\text{stab}} \, \Omega
\end{align}
where $N_{A/B,j}^{\text{stab}}$ represents the gradient interpolators used to define the hourglass deformation modes.
The $\alpha_{\text{stab}}$ is a scaling factor typically set between 0.005 and 0.1 \citep{Belytschko2014, Abaqus}, depending on the form of $N_{A,j}^{\text{stab}}$.
For volumetric spatial instabilities, $\mu$ is replaced by the bulk modulus, $\kappa = (2/3)\mu + \lambda$, and $N_{A,j}^{\text{stab}}$ represents the gradient interpolators for the volumetric/pressure hourglass mode.


Finally, the two external force terms in Equation \ref{eq:virtualWorkDiscrete} are expressed as:
\begin{align}
	&\oint_{\Gamma_T}
	\sum_{A=1}^{nVertices} N_A  \boldsymbol{\bar{\omega}}_A 
	\; \cdot \;
	\boldsymbol{T}_\Gamma
	\;\text{d}\Gamma 
	\;+\;
	\int_\Omega
	\sum_{A=1}^{nVertices} N_A  \boldsymbol{\bar{\omega}}_A 
	\;\cdot \;
	 \rho \boldsymbol{f}_b
	\;\text{d}\Omega
	\; = \;
	\sum_{P=1}^{nDoF} \bar{\omega}_{P} F_P
\end{align}
where vector $ F_P$ is known as the global force vector:
\begin{align}
	F_P =
	&\oint_{\Gamma_T}
	\sum_{A=1}^{nVertices} N_A
	T_{i\Gamma}
	\;\text{d}\Gamma
	\;+\;
	\int_\Omega
	\sum_{A=1}^{nVertices} \rho N_A  f_{ib}
	\;\text{d}\Omega
\end{align}
and once again the integrals are evaluated using quadrature.

Equation \ref{eq:virtualWorkDiscrete} can now be expressed as:
\begin{align} \label{eq:virtualWorkDiscreteAlmostFinal}
	\boldsymbol{\bar{\omega}}^T
	\left(
	\boldsymbol{M} \frac{\partial^2 \boldsymbol{U}}{\partial t^2} 
	+ \boldsymbol{K} \boldsymbol{U}
	- \boldsymbol{F}
	\right)
	= \boldsymbol{0}
\end{align}
where the global vectors of dimension $nDoF \times 1$ are:
\begin{align}
	\boldsymbol{\bar{\omega}} = [\bar{\omega}_P], \quad\quad
	\boldsymbol{U} = [U_P], \quad\quad
	\boldsymbol{F} = [F_P]
\end{align}
and global matrices of dimension $nDoF \times nDoF$ are:  
\begin{align}
	\boldsymbol{M} = [M_P], \quad\quad
	\boldsymbol{K} = [K_P]
\end{align}

As Equation \ref{eq:virtualWorkDiscreteAlmostFinal} is satisfied for all values of $\boldsymbol{\bar{\omega}}$, this requires that the bracketed term is equal to zero.
This gives the semi-discrete form (discrete in space, not in time) of the finite element equations:
\begin{align}
	\boldsymbol{M} \frac{\partial^2 \boldsymbol{U}}{\partial t^2} 
	+ \boldsymbol{K} \boldsymbol{U}
	&=
	\boldsymbol{F}
\end{align}

To complete the discretisation, the acceleration term $\frac{\partial^2 \boldsymbol{U}}{\partial t^2}$ is discretised using a finite difference scheme.
Similar to the finite volume approaches, many finite difference schemes can be used, where Newmark schemes are popular.
The simple Euler backward scheme is given here for comparative purposes:
\begin{align} \label{eq:virtualWorkDiscreteFinal}
	\boldsymbol{M} \;
	 \frac{\boldsymbol{U} - 2 \boldsymbol{U}^{[m-1]} + \boldsymbol{U}^{[m-2]}}{\Delta t^2}
	+ \boldsymbol{K} \boldsymbol{U}
	&=
	\boldsymbol{F}
\end{align}

Boundary tractions/forces have already been included in $\boldsymbol{F}$ via the external (surface) force term in Equation \ref{eq:virtualWorkDiscrete}.
For the incorporation of displacement conditions, this signifies that some of the degrees of freedom in $\boldsymbol{U}$ are known;
these equations can hence be disregarded.
As in the case with the finite volume methods, initial conditions, in the form of the displacement field at $t = 0$, $t = -\Delta t$, and $t = -2\Delta t$, must also be specified.

\subsection{Solution algorithm}
Once known degrees of freedom have been incorporated, Equation \ref{eq:virtualWorkDiscreteFinal} represents a system of $nDoF - nKnownDof$  linear algebraic equations, where $nKnownDof$ is the number of known degrees of freedom.
Similar to the analogous matrices in the finite volume method,
the mass matrix can be seen to be a function of the density and element geometry;
the stiffness matrix is a function of the mechanical properties and element geometry;
and the force vector $\boldsymbol{F}$ contains surface and body force contributions as well as inertial terms and non-zero known degree of freedom contributions.
Like the finite volume approaches, this system of algebraic equations can be solved using either an implicit or explicit time marching procedure.

For implicit approaches, Equation \ref{eq:virtualWorkDiscreteFinal} can be rearranged and solved for $\boldsymbol{U}$:
\begin{align} \label{eq:finiteElementSolution}
	\left( \frac{1}{\Delta t^2} \boldsymbol{M} + \boldsymbol{K} \right) \boldsymbol{U}
	&=
	\boldsymbol{F}
	+ \boldsymbol{M} \; \frac{2 \boldsymbol{U}^{[m-1]} - \boldsymbol{U}^{[m-2]}}{\Delta t^2} 
\end{align}
Assuming a stable discretisation, the matrix of the linear system $\left( \frac{1}{\Delta t^2} \boldsymbol{M} + \boldsymbol{K} \right)$ (or $\boldsymbol{K}$ for quasi-static) has the following properties:
\begin{itemize}
	\item It is sparse with the number of non-zero elements in each row equal to the number of vertices which share an element with the current vertex, plus one;
	\item It is symmetric;
	\item It is positive definite;
	\item In general, unlike for the cell-centred finite volume approaches where a segregated algorithm is used, it is not diagonally dominant; however, finite volume discretisations that employed block-coupled solution methodologies will produce a similar non-diagonally dominant matrix.
\end{itemize}
To solve Equation \ref{eq:finiteElementSolution}, typically a direct linear solver is employed, for example, Gaussian elimination, {LU} decomposition or multi-frontal methods; however, iterative methods such as preconditioned conjugate gradient are also possible \citep{Bathe1996}.


\subsection{Discussion}

\paragraph{Overview}
Given the close relationship between the finite volume and finite element methods, it is not surprising that both have been compared previously:
some notable dissections include those from O{\~{n}}ate, Zienkiewicz, Idelsohn \citep{Onate1992, Onate1993, Onate1994, Idelsohn1994, Onate1998, Onate2000, Zienkiewicz1991, Zienkiewicz1995}, \citet{Lahrmann1992}, \citet{Perre1995}, \citet{Harrild1997}, \citet{Zarrabi1999}, \citet{Fang2002}, \citet{Yamamoto2002}, Jacquemet and Henriquez
\citep{Jacquemet2005a, Jacquemet2005b}, \citet{Vaz2009}, \citet{Filippini2014}, and recently \citet{Demirdzic2020}.
The general consensus is that both methods share the same data structure and general strategy to assemble the corresponding characteristic matrices, with the main conceptual difference being in the local integration domain and local integration method.
Although a number of authors have claimed superior accuracy of one method over the other \citep{Lahrmann1992, Harrild1997, Vaz2009}, the majority of authors have found both methods to produce similar predictions with no significant differences \citep{Perre1995, Fang2002, Yamamoto2002, Jacquemet2005b}.

One generally accepted appeal of the finite volume method is the ease with which it can be followed and implemented: the finite volume method is simply based on balancing forces acting on a volume, requiring no knowledge of advanced mathematical frameworks.
Comparing the derivations in Section \ref{sec:standardFiniteVolumeDiscretisationOfSpace} and \ref{sec:discretisationFE}, the finite element approach would appear to be more mathematically involved;
however, proponents of the finite element method claim that the principle of virtual work and energy minimisation techniques are equally interpretable as the balance equation form.

\paragraph{Weak vs strong forms of the conservation law and the implications}
As discussed, finite volume and finite element methods differ is their philosophy:
where the finite volume method deals with the strong balance form of the governing law, the finite element method deals with the equivalent weak virtual energy form.
One consequence of this is the manner in which local conservation is enforced.
As finite volume approaches discretise surface integrals at the (typically non-overlapping) control volume boundaries, strong local conservation is achieved:
forces are equal and opposite at cell boundaries. As a consequence on this local conservation, global conservation within the domain is automatically achieved.
In contrast, finite element methods discretise the surface force term as a volume integral using locally overlapping integration domains: this results in local conservation in an average sense, rather than directly for each element.
Global conservation is ensured in an average sense, assuming there are sufficient numbers of elements.
An additional consequence of these contrasting approaches is how each method treats Neumann/natural boundary conditions:
considering traction conditions, finite volume methods satisfy these conditions exactly regardless of the mesh density; whereas finite element approaches satisfy them in an approximate sense, and as the mesh is refined strong enforcement is approached.

As noted by a number of authors, for example, \citep{Idelsohn1994, Onate1994, Bailey1995}, by choosing unity weighting functions in the weak form of the governing equation (Equation \ref{eq:virtualWork}), it is possible to recover the finite volume method.
\citet{Spalding2008} alludes to this point by referring to finite volume approaches as \emph{unity-weighting function methods} and to finite element approaches as \emph{non-unity weighting function methods}.

\paragraph{Geometric flexibility} 
As \emph{standard} finite element methods use shape functions to define a continuous displacement distribution, this limits their application to meshes containing \emph{standard} element shapes.
Concretely, only standard shapes such as triangles and quadrilaterals are allowed in 2-D, and tetrahedra and hexahedra in 3-D.
Similarly, finite volume approaches that explicitly use shape functions are constrained in the same way.
More commonly finite volume methods describe the local solution distribution using a truncated Taylor expansion, resulting in discontinuous jumps at the cell interfaces.
An outcome of this is the ability of finite volume approaches to deal with general convex polyhedral meshes.
In this way, finite volume approaches are more flexible in terms of mesh generation and dynamic remeshing.
In addition to hanging nodes, typical finite volume discretisations can deal with overset/\emph{chimera} and immersed boundary meshes in a straight-forward manner.

\paragraph{Mass matrix and stiffness matrix properties}
A desirable property of final volume methods for explicit implementation and parallelisation is that the mass matrix is automatically diagonal.
For finite element methods, the lack of a diagonal consistent mass matrix results in the use of so-called \emph{mass lumping} approaches, which are often ad-hoc, as noted previously.
The stiffness matrix, however, bares many similarities between the methods.

\begin{figure}[htb]
	\centering
	\subfigure[Local integration domain around a vertex in the finite element method]
	{
		\includegraphics[height=0.4\textwidth]{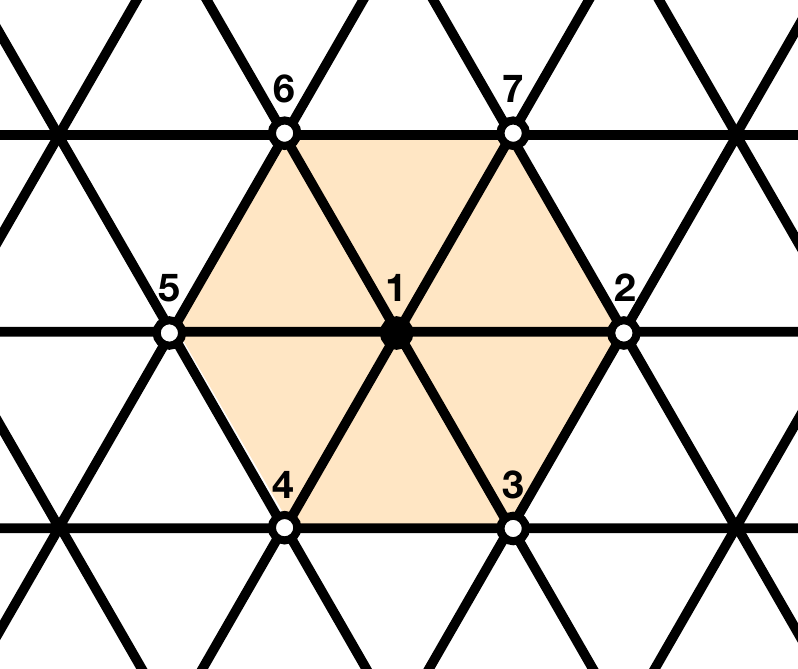}
	}
	\subfigure[Local integration domain for a cell in the finite volume method]
	{
	    	\includegraphics[height=0.4\textwidth]{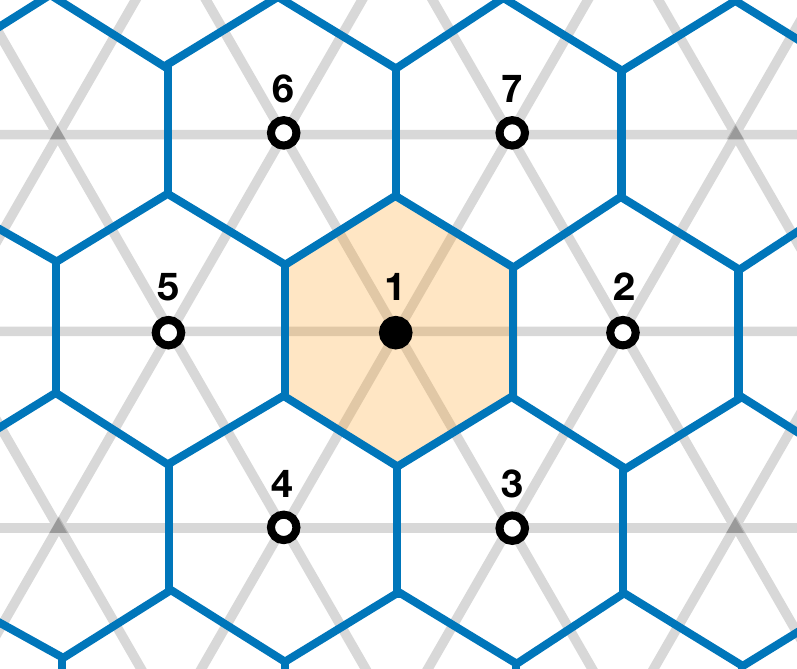}
	}
	\caption{Uniform 2-D triangular grid showing the local integration domains and computational stencils for the finite element and finite volume methods. The centre node is shown as a black dot and all neighbouring nodes are shown as white-filled dots, where nodes are vertices in the finite element method, and are cell-centres in the finite volume method.}
	\label{fig:uniformGrid}
\end{figure}

To allow direct comparison, let us consider finite volume and finite element discretisations on a uniform equilateral-triangular grid (Figure \ref{fig:uniformGrid}).
We will consider a typical domain of integration around a node, where a node represents a vertex in the finite element method and a cell-centre in the finite volume method.
For ease of comparison, the dual mesh is used for the finite volume method.
In the finite element case (Figure \ref{fig:uniformGrid}(a)), the local integration domain (shaded in orange) consists of all elements adjacent to the centre node (black-filled dot).
Consequently, the computational stencil for the centre vertex consists of the six neighbouring vertices (indicated by white-filled dots), which share an element with the centre node.
In the finite volume case (Figure \ref{fig:uniformGrid}(b)), the local integration domain (shaded in orange) is the cell itself containing the node at its centre (black-filled dot);
the computational stencil includes all neighbouring cell-centres (white-filled dots).
In this case, as both methods share the same computational stencil, the resulting stiffness matrices will share the same structure and sparsity, \emph{assuming} that the same solution strategy is used, for example, block-coupled or segregated.
For a triangular grid, \emph{either} the nodal locations \emph{or} the orientation of the local integration domains can be aligned between the methods, but not both;
as such, the comparison of coefficients from both methods is not entirely direct, however, it does still provide insight.
Here it was chosen to align the nodal locations rather than the local integration domains:
the finite volume cell can be seen to be rotated 90\degree\ relative to the finite element cell, as well as being one third the area.

It is important to mention that the size and shape of the computational stencil in the finite volume method depend on the approach used to calculate the displacement gradients at the cell faces.
Here face gradient calculation method 2 in Section \ref{sec:faceGradientDiscussion} has been employed.
Alternative face gradient calculation methods can result in the \emph{second cell-neighbours} being included in the stencil;
in contrast, a characteristic of the finite element method is that the integration domain is fixed regardless of local discretisation.

An additional observation is that the local integration domains overlap in the finite element method, but do not in the typical finite volume method \ie integration domains for neighbouring finite element nodes overlap.
An overlapping version of the vertex-centred finite volume method has been considered by \citet{Onate1994} but has not received significant attention.

\begin{figure}[htb]
	\centering
	\subfigure[Finite element method \emph{consistent} mass matrix]
	{
		\includegraphics[height=0.4\textwidth]{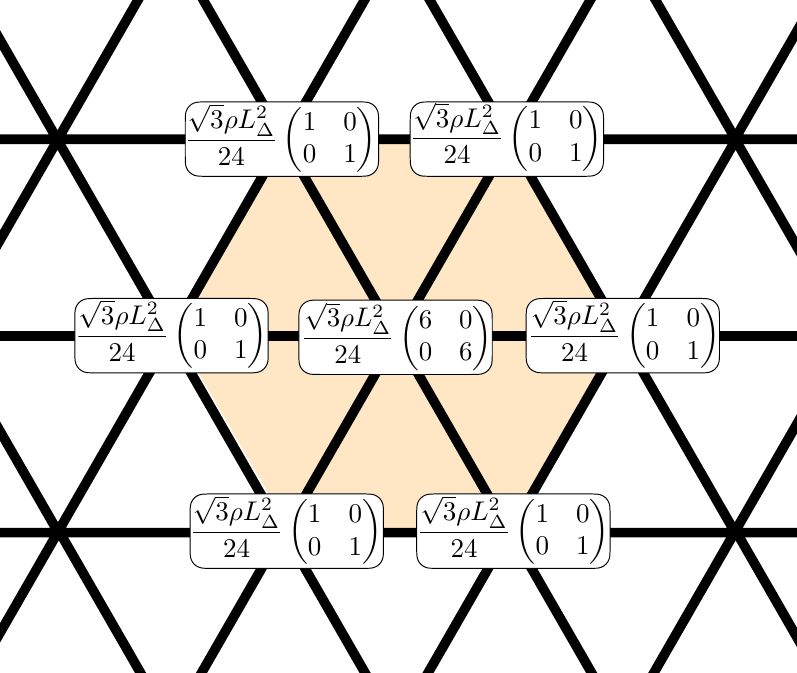}
	}
	\subfigure[Finite element method \emph{lumped} mass matrix]
		{
	    	\includegraphics[height=0.4\textwidth]{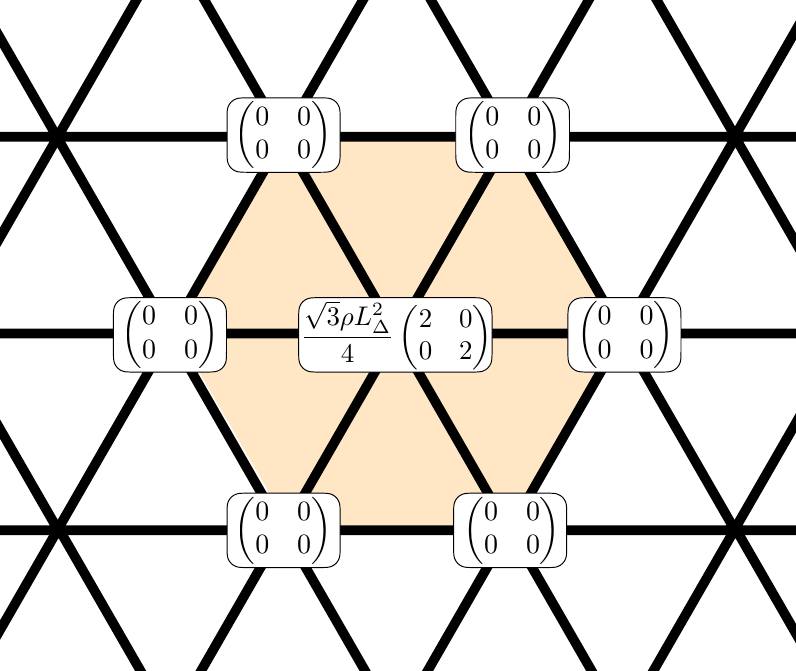}
	}
	\subfigure[Finite volume method]
	{
	    	\includegraphics[height=0.4\textwidth]{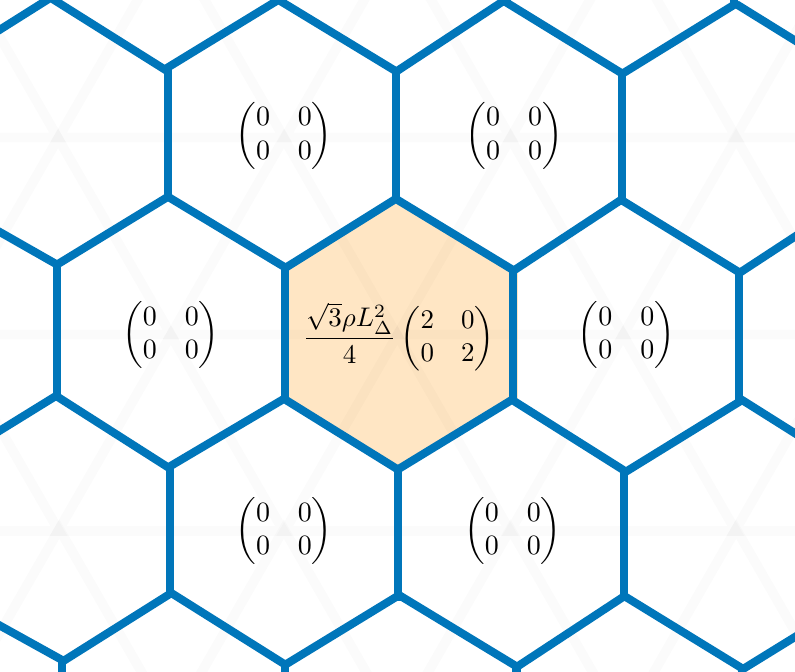}
	}
	\caption{Block coefficients for a row in the mass matrix corresponding to the centre node, where $L_\Delta$ is the triangle side-length, a unit thickness is assumed, and $\frac{\sqrt{3}}{4} L_\Delta^2$ is the area of a triangular element. The area of the hexagonal finite volume cell is equal to $\frac{\sqrt{3}}{2} L_\Delta^2$.}
	\label{fig:massMatrices}
\end{figure}

\begin{figure}[htb]
	\centering
	\subfigure[Finite element method stiffness matrix]
	{
		\includegraphics[height=0.4\textwidth]{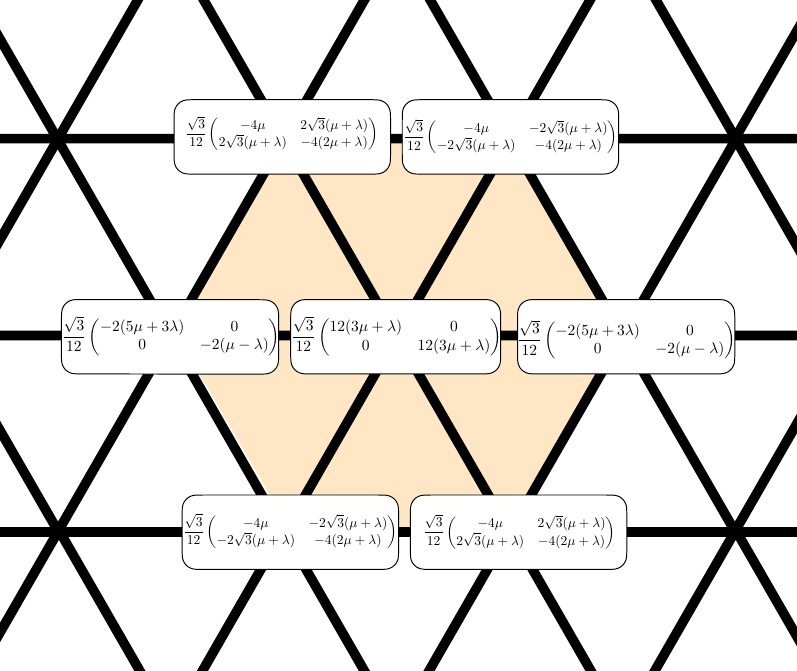}
	}\\
	\subfigure[Finite volume method stiffness matrix with a \emph{coupled} solution algorithm]
		{
	    	\includegraphics[height=0.4\textwidth]{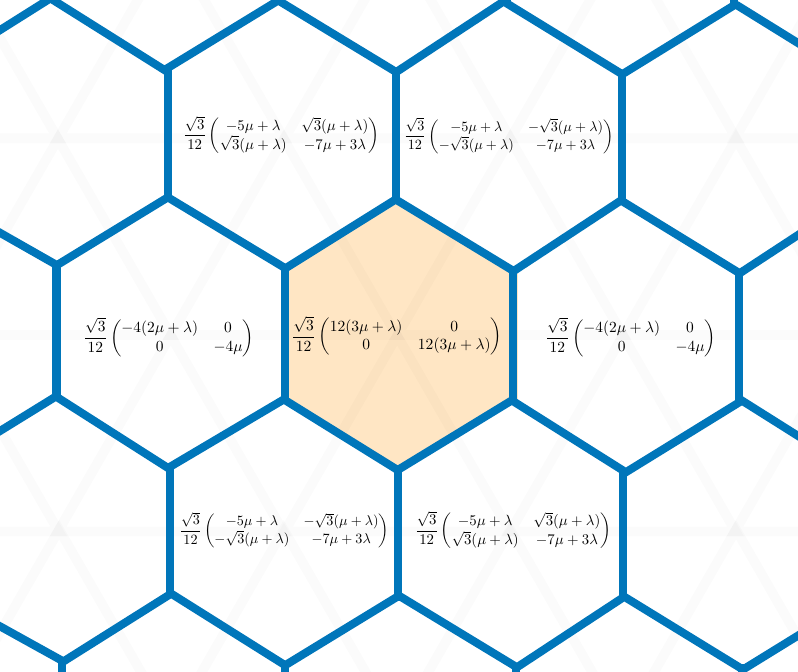}
	}
	\subfigure[Finite volume method stiffness matrix with a \emph{segregated} solution algorithm]
	{
	    	\includegraphics[height=0.4\textwidth]{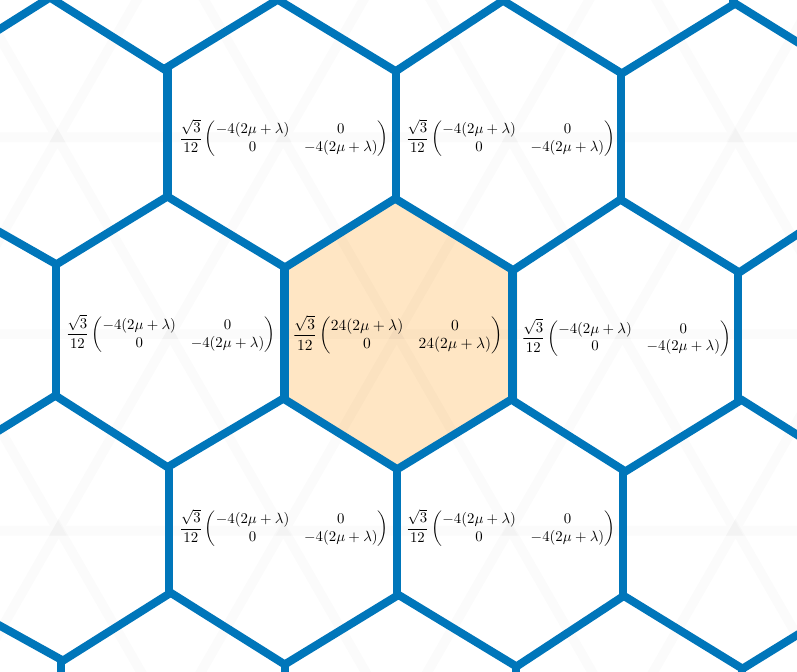}
	}
	\caption{Block coefficients for a row in the stiffness matrix corresponding to the centre node, where $L_\Delta$ is the triangle side-length and a unit thickness is assumed}
	\label{fig:stiffnessMatrices}
\end{figure}

Referring to the node numbering given in Figure \ref{fig:uniformGrid}, the block row for the centre node in the global mass matrix is given for the finite element method \citep{Belytschko2014}, consistent and lumped, as well as the finite volume method in Figure \ref{fig:massMatrices}.
The finite element mass-lumped approach can be seen to coincide with the finite volume approach; from this, we can see that the finite element lumped approach is essentially assigning the mass of the hexagon surrounding a node to the node itself, corresponding with the integration domain of the finite volume approach.
Considering next the stiffness matrix.
Once again taking the same uniform 2-D triangular grid given in Figure \ref{fig:uniformGrid}, the block row for the centre node in the global stiffness matrix for the finite element method \citep{Belytschko2014} is compared with the equivalent row from the finite volume method in Figure \ref{fig:stiffnessMatrices}.
For the finite volume method, the matrix using a block-coupled solution algorithm is given, as well as for the segregated approach.
Examining the general structure of the coefficients for all three methods (Figure \ref{fig:stiffnessMatrices}(a), (b) and (c)), the following observations can be made:
\begin{itemize}
	\item All three methods employ the same computational stencil;
	\item The coefficients are symmetric in all three cases, albeit the effect of boundary conditions has not been considered here;
	\item As expected, all three methods show geometric symmetries, for example, the top-right coefficient is equal to the bottom-left coefficient;
	\item The finite element and coupled finite volume approaches show the same sparsity structure \ie there are zeros in the same locations;
	in addition, the coefficients have the same signs and show similar magnitudes, but do not have the same value, apart from for the central node;
	\item The segregated finite volume approach differs from the others by removing all inter-component coupling from the coefficients, as discussed in Section \ref{sec:faceGradientDiscussion};
	\item The segregated finite volume approach produces a row of coefficients which is weakly diagonally dominant, also discussed in Section \ref{sec:faceGradientDiscussion};
	\item The momentum/force \emph{between nodes} is conserved for all three methods (based on the symmetry of the coefficients), however, conservation of momentum/force \emph{between elements/cells} is \emph{not} directly enforced for the finite element method, whereas it \emph{is} for the finite volume method.
\end{itemize}

\begin{figure}[htb]
	\centering
	\subfigure[Example of the \emph{hourglassing} spatial instabilities in reduced integration finite element formulations (taken from \citep{Belytschko2014})]
	{
		\includegraphics[width=0.35\textwidth]{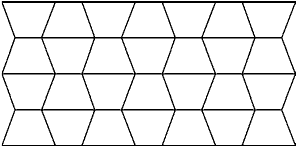}
	}
	\quad\quad
	\subfigure[Example of the \emph{checkboarding} spatial instabilities (hydrostatic pressure field) in cell-centred finite volume formulations (taken from \citep{Cardiff2016b})]
	{
		\includegraphics[width=0.45\textwidth]{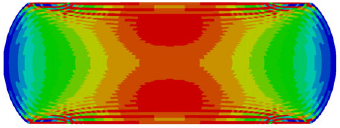}
	}
	\caption{The equivalence of spatial instabilities in finite element and finite volume formulations}
	\label{fig:hourglassingCheckerboarding}
\end{figure}

\paragraph{Discretisation stabilisation}
As described in Section \ref{sec:discretisationFE}, finite element formulations which under-integrate the local domain require the inclusion of a stabilisation term, in the same way one is required in many finite volume formulations.
In the finite element method, these spurious singular, zero energy modes produce an accordion-like deformation pattern known as hourglassing (Figure \ref{fig:hourglassingCheckerboarding}(a)), named due to its visual similarity to an hourglass timing device.
These spatial instabilities are equivalent to checkerboarding instabilities that occur in finite volume formulations (Figure \ref{fig:hourglassingCheckerboarding}(b));
these checkerboarding pressure-type instabilities are also possible in finite element formulations.
The origin of these spatial instabilities in both finite element and finite volume formulations is a \emph{rank deficiency} is their respective stiffness matrices.
Essentially, a stable discretisation should not support any deformation modes which do not offer resistance, apart from rigid body translations and rotations.
More precisely, rank deficiency refers to the fact that the discretised stiffness matrix has zero-valued eigenvalues which are not related to rigid body motions (as discussed in Section \ref{sec:stabilisationApproach}).
For example, considering a 3-D element/cell, the stiffness matrix for a \emph{full rank} formulation has six zero valeed eigenvalues: three rigid translations and three rigid rotations;
for a rank deficient formulation, additional zero valued eigenvalues are present corresponding to zero energy instability modes.

To address these spatial instabilities, both finite element and finite volume approaches add a stabilisation term to the discretised governing equations.
In both cases, there are two primary constraints on the form of this stabilisation term \citep{Belytschko2014}:
\begin{itemize}
	\item It should not significantly affect the accuracy of the discretisation, and
	\item It should not affect the linear completeness of the discretisation;
		in other terms, the discretisation should still be able capable of describing a linear solution field exactly after the inclusion of the stabilisation term.
\end{itemize}

For both methods, a variety of techniques exist and the stabilisation \emph{force} often takes a similar form:
\begin{align}
\begin{aligned}[l]
	\boldsymbol{f}^{\text{stab}} = \qquad \\
	\\
	f^{\text{stab}}_P = \qquad
\end{aligned}
\begin{aligned}[l]
	 &\alpha^{\text{stab}} C^{\text{stab}} \Omega \,
	 \left(
	 \boldsymbol{\nabla}^2 \boldsymbol{u} 
	 - \boldsymbol{\nabla} \cdot \boldsymbol{\nabla} \boldsymbol{u} 
	 \right) \\\
	 &-\alpha^{\text{stab}} C^{\text{stab}} \Omega \,
	 \boldsymbol{\nabla}^2 \left( \boldsymbol{\nabla}^2 \boldsymbol{u} \right)  \\
	 &\alpha^{\text{stab}} C^{\text{stab}} \Omega \,
	 N_{A,j} \, N_{B,j} U_Q
\end{aligned}
\begin{aligned}[l]
	&\qquad \text{Rhie-Chow stabilisation} \\
	&\qquad \text{Jameson-Schmidt-Turkel stabilisation} \\
	&\qquad \text{Hourglass stabilisation}
\end{aligned}
\end{align}
where $C^{\text{stab}}$ is some measure of mechanical property that gives an appropriate dimension to the dissipation, such as $\mu$, $2\mu + \lambda$, $\mu + \lambda$, or $\kappa$, for quasi-static analyses, or a function of the speed of sound for dynamic analyses.
The $\alpha^{\text{stab}}$ factor allows the magnitude of the stabilisation to be scaled.

\paragraph{Performance of the methods}
The differences in local integration domain and local integration method between the finite volume and finite element methods has consequences for the robustness, accuracy and efficiency of the resulting methods.
Particularly for nonlinear analysis, the ideal discretisation is unclear due to a variety of numerical challenges which are yet to be fully resolved, including \citep{Hassan2019paper}:
(1) spurious hourglassing and pressure checker-boarding, (2) bending difficulties, (3) shear and volumetric locking, (4) high frequency noise in the vicinity of shocks, (5) lower order of convergence for strains and stresses in comparison with displacements, and (6) sensitivity to mesh distortions. 
It is these challenges that finite volume discretisations can potentially solve in a novel way.

Apart from this, a major motivation for finite volume solid mechanics schemes is the challenge of multi-physics problems.
As long as the finite volume method is prominent in the world of computational fluid mechanics, there will be a demand for straight-forward finite volume solid mechanics implementations. These solid mechanics implementations can share the same computational framework, discretisation and solution methodologies as their fluid counterparts and can be integrated seamlessly into the code base.

Regarding accuracy and overly stiff behaviour, finite volume methods have not shown the same locking behaviour typical in fully integrated finite element methods.
Given their similarity with reduced integrated finite elements, it is perhaps not surprising that this is the case; however, like finite elements, the absolute accuracy depends on the details of the formulation.
Concerning order of accuracy, an attractive property of most finite volume methods is that the error in the strain and stress reduces at a second-order rate, like the displacement \citep{Vaz2009, Tukovic2018a}; this may not be the case for many finite element schemes, where the error in the strain and stress reduce at a rate closer to first order.

As models becomes larger and the availability of supercomputers and cloud computing increases, code parallelisation is becoming critical.
Due to the widespread use of iterative linear solvers, finite volume methods (fluid and solid) commonly exploit hundreds or thousands of CPU cores, for example, the OpenFOAM software.
In contrast, as direct linear solvers have often been the chosen solution approach for finite element schemes, parallel efficiency is inherently limited (relative to iterative solvers) and the use of large numbers of CPU cores has been less common; however, there are a number of projects focussed on the application of finite element methods to supercomputers using iterative solvers, for example, {ParaFEM} \citep{ParaFEM}.
Apart from the choice between direct and iterative linear solver, cell-finite volume approaches possess a convenient advantage over vertex-centred methods (such as the finite element method) when it comes to domain decomposition parallelisation:
each node (cell-centre) is uniquely located on one CPU core domain.
As a result, duplication of nodal data at processor-to-processor boundaries is not required.
In contrast, in the finite element method, nodes (vertices) that lie on a processor-to-processor boundary are present on at least two processor domains, requiring the use of \emph{ghost elements} or similar data duplication techniques.
A final point related to parallelisation is the debate around open-source vs commercial software, which equally affects both finite volume and finite element methods.
For many, the current \emph{per-CPU-core} pricing of some well known commercial codes may in fact be the greatest practical obstacle to parallel efficiency.


\paragraph{Higher order discretisations}
Apart from the size of the finite element solid mechanics community, arguably the next greatest advantage of the finite element method over the finite volume method is its straight-forward extension to higher-order discretisation.
A key feature of the finite element method is that an element is completely characterised by the coordinates and degrees of freedom associated with its nodes/vertices \citep{Laursen2002}.
This compact computational stencil allows uncomplicated inclusion of higher-order local distributions.
As noted previously, higher-order finite volume approaches have been developed for solid mechanics \citep{Demirdzic2016}, however, as the order increases, so does the size of the computational stencil;
this introduces significant challenges for unstructured polyhedral grids, and as such, higher-order schemes are not as common as in the finite element field.

For linear analyses, the power of higher-order schemes is undeniable, notwithstanding challenges with locking behaviour, however, for nonlinear analysis such higher order schemes are not commonly used in practice.
For example, considering an elasto-plastic analysis, as noted by \citet{Belytschko2014}, the stress may have discontinuous derivatives at the surface separating elastic and plastic material.
In this case, the errors in Gaussian quadrature of an element that contains an elastic-plastic interface are likely to be large; higher-order quadrature is not a solution as it often leads to stiff behaviour or locking.
Additionally, finite element users are typically recommended to use lower order formulations for problems that include either contact, large strains or plasticity \citep{Abaqus}.

\begin{figure}[htb]
	\centering
	\subfigure[Discontinuous field representation within the finite volume method and discontinuous Galerkin finite element method]
	{
		\includegraphics[width=0.45\textwidth]{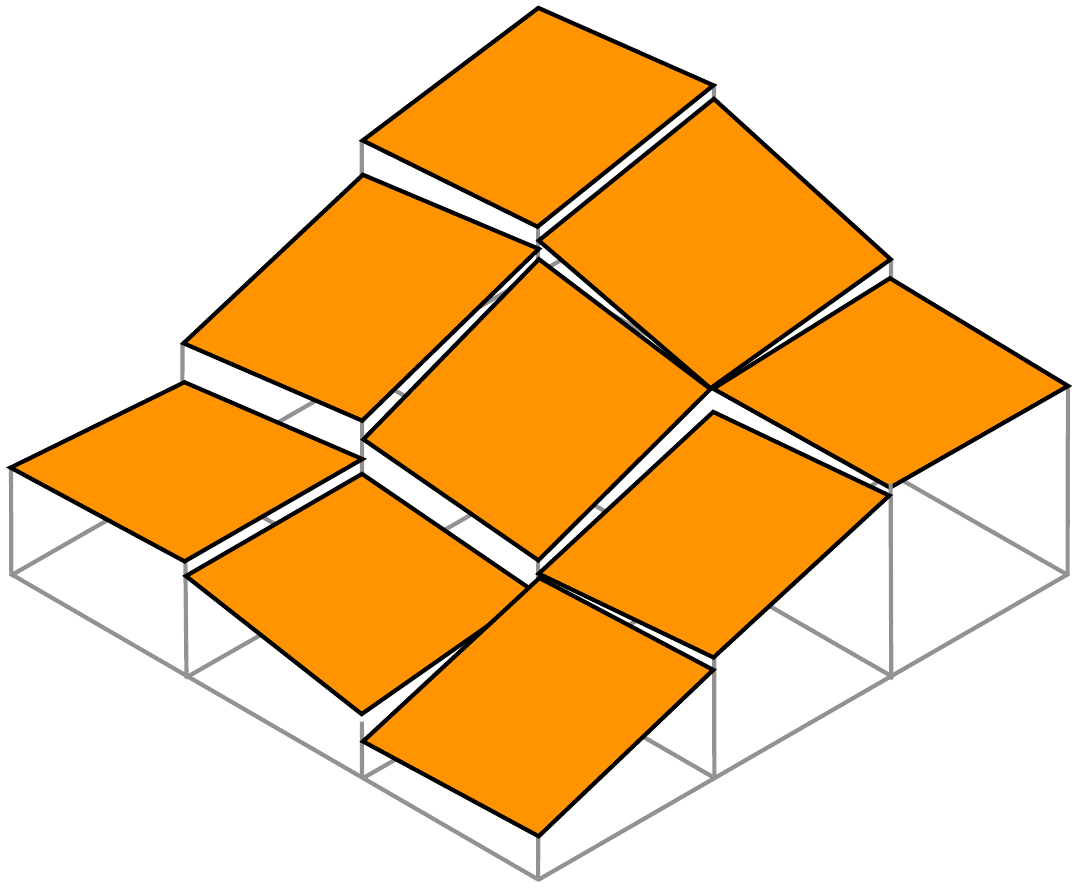}
	}
	\subfigure[Continuous field representation within the finite element method]
	{
	    	\includegraphics[width=0.45\textwidth]{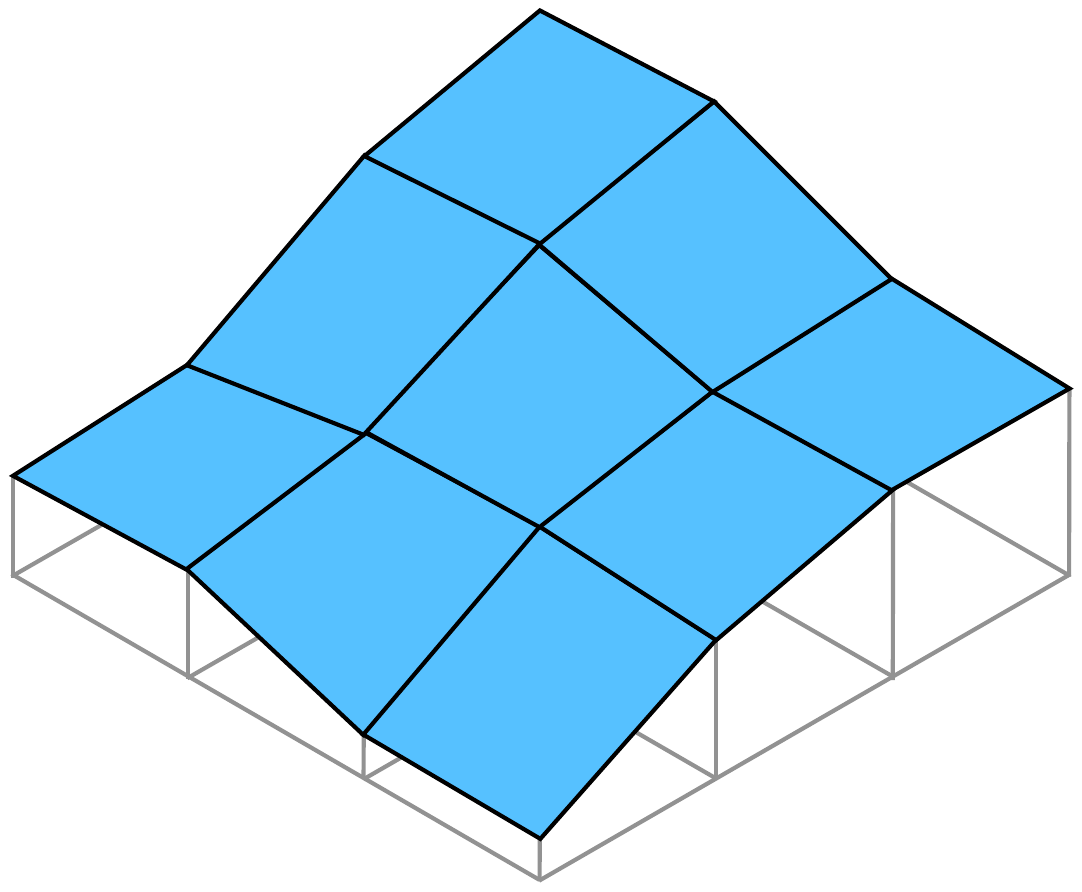}
	}
	\caption{A comparison between the representation of the displacement field in the finite volume method (and discontinuous Galerkin finite element method) and the finite element method (adapted from \citet{Lee2013})}
	\label{fig:FEvsFVdisplacementDistribution}
\end{figure}

\paragraph{Discontinuous Galerkin methods}
The standard \emph{continuous} Bubnov-Galerkin finite element method, as presented above, assumes a continuous distribution of the displacement between elements.
In contrast, as shown in Figure \ref{fig:FEvsFVdisplacementDistribution}, the finite volume method \emph{typically} assumes discontinuous jumps in the displacement field at the interfaces between cells (although it does not have to).
There is, however, a class of finite element approaches known as discontinuous Galerkin methods, introduced by \citet{Reed1973}, which assume similar jumps in the solution field across element boundaries.
In this way, the local integration method adopted by the discontinuous Galerkin method combines features of the finite volume and finite element methods.
In particular, discontinuous Galerkin schemes bare the following desirable finite volume properties \citep{Cockburn2003}:
\begin{itemize}
	\item They produce mass matrices that are block-diagonal;
	\item They easily handle irregular meshes with hanging nodes;
	\item They are locally conservative, which is a property that is a critical property for computational fluid dynamics applications.
\end{itemize}
In addition, a potential advantage of discontinuous Galerkin schemes over finite volume schemes is their ease of extension to higher orders, and the mixing of lower and higher order elements.

Discontinuous Galerkin methods are, however, typically implemented using explicit solution algorithms. The reason they are less suitable for implicit implementations (and analysis of quasi-static type problems) is they possess large numbers of globally-coupled degrees of freedom.
To overcome this disadvantage, the so-called \emph{hybridisable} discontinuous Galerkin (HDG) method was introduced by \citet{Cockburn2009}.
The key characteristic of the hybridisable form was defining \emph{traces} of field variables as single values at cell interfaces, allowing a significant reduction in the number of global unknowns.
Detailed analysis of the HDG method is outside the scope of the current article and readers are referred to recent articles on the application of HDG to solid mechanics \citep{Fu2015, Qiu2017, Sevilla2018HDG, Hesthaven2007}.

\section{Applications of the finite volume method for computational solid mechanics}
Some of the main areas where finite volume methods have been applied to solid mechanics are summarised below, including example images.
\subsection{Fluid-solid interaction}
\begin{figure}[htb]
	\centering
	\subfigure[Wing and fluid domain meshes for the prediction of flutter \citep{Slone2004}]
	{
		\includegraphics[height=0.27\textwidth]{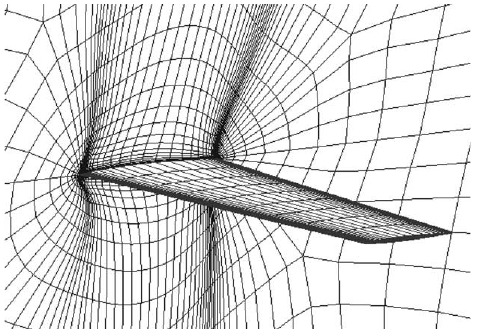}
	}
	\subfigure[Fluid velocity magnitude and solid displacement magnitude distribution for channel flow over an elastic thick plate \citep{Tukovic2018a}]
	{
		\includegraphics[height=0.27\textwidth]{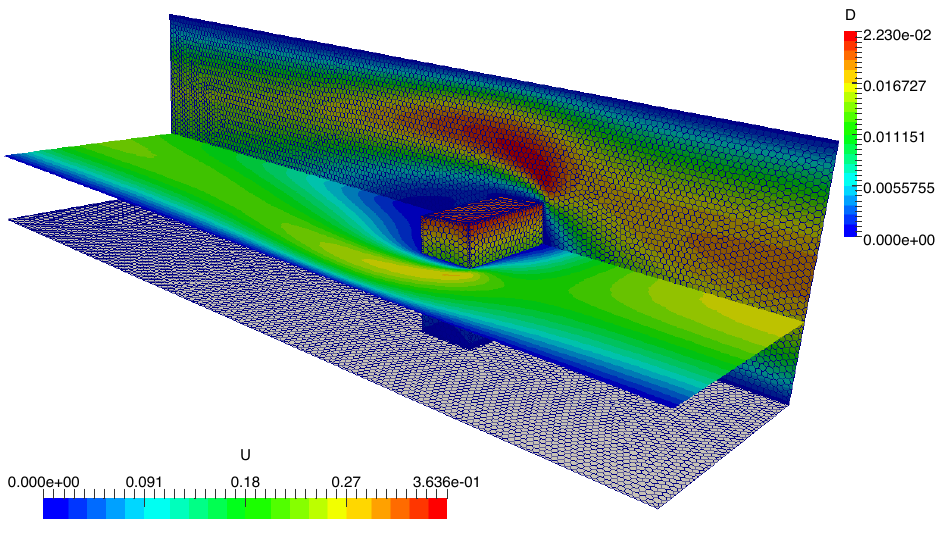}
	}
	\subfigure[Pressure (top) and velocity (bottom) distribution around a red blood cell \citep{Wiedemair2012}]
	{
		\includegraphics[height=0.35\textwidth]{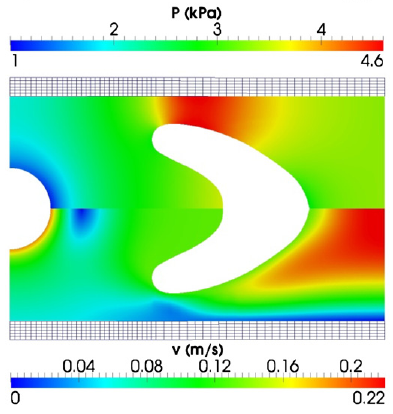}
	}
	\subfigure[Pressure field for flow over a rigid square with a thin flexible plate attached \citep{Tsui2013}]
	{
		\includegraphics[height=0.35\textwidth]{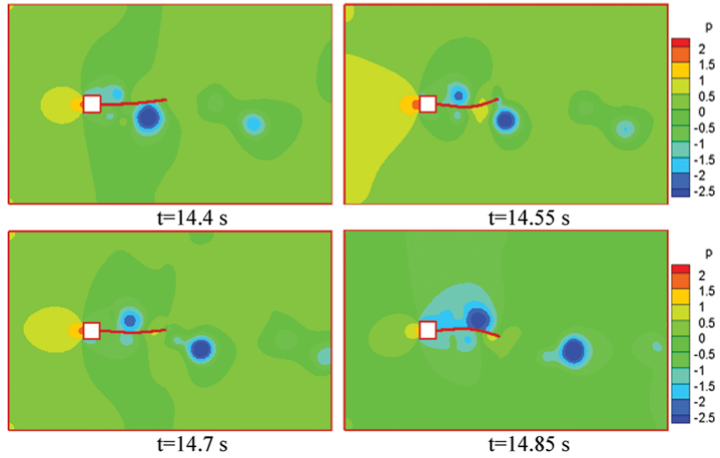}
	}
	\caption{Fluid-solid interaction examples}
	\label{fig:examplesFSI}
\end{figure}

Example cases are shown in Figure \ref{fig:examplesFSI}, and examples references include: \citep{Henry1993a, Henry1993b, Demirdzic1995, Ivankovic1997a, Slone1997, Demirdzic1998, Ivankovic1998, Greenshields1999b, Oldroyd1999, Tang1999, Greenshields2000, Schafer2000, Slone2000b, Slone2000c, Schafer2001b, Slone2001, Slone2001b, Slone2001c, Dzaferovic2002, Ivankovic2002b, Karac2002, Slone2002, Schafer2002, Slone2002c, Slone2002d, Torlak2002a, Karac2003a, Karac2003b, Karac2003c, Cross2004, Karac2004, Slone2004, Kovacevic2004b, Giannopapa2004a, Giannopapa2004b, Slone2004b, Slone2004c, Zhang2004, Greenshields2005, Shaw2005, Cross2006, Giannopapa2006, Papadakis2006, Slone2007, Cross2007, Giannopapa2007, Jasak2007,  Kovacevic2007, Lv2007, Lv2007, Stosic2007, Tukovic2007b, Croft2008, Papadakis2008, Giannopapa2008, Xia2008, Kanyanta2009a, Karac2009a, Karac2009b, Safari2009, Slone2009, Das2010, Kelly2010, Jagad2011, Kelly2012, Wiedemair2012, Das2013, Tsui2013, Habchi2013, Tukovic2014, Hejranfar2016, Sekutkovski2016, Selim2017, Cardiff2017b, Oliveira2017a, Oliveira2017b, Cardiff2018a, Gonzalez2018, Tukovic2018a, Tukovic2018b};

\subsection{Fracture and adhesive joints}
\begin{figure}[htb]
	\centering
	\subfigure[Crack path predictions through a material interface \citep{Carolan2013a}]
	{
		\includegraphics[height=0.27\textwidth]{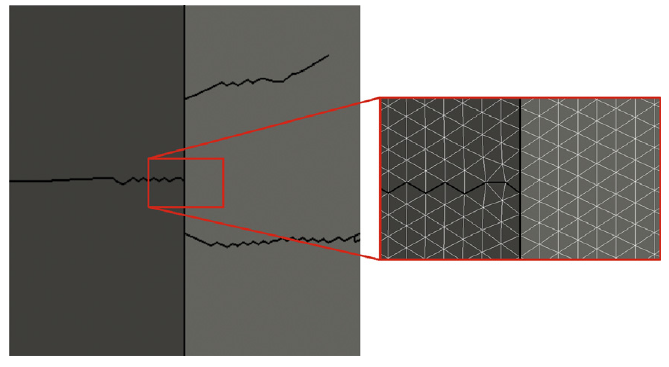}
	}
	\subfigure[Dynamic crack path predictions in PMMA \citep{Murphy2006}]
	{
		\includegraphics[height=0.27\textwidth]{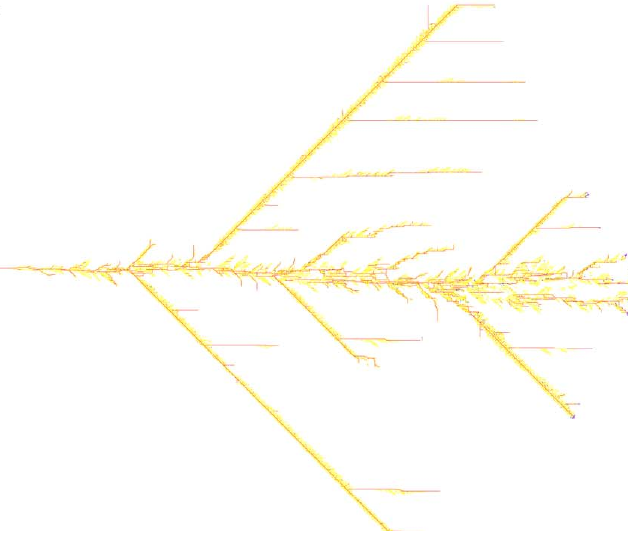}
	}
	\subfigure[Crack path predictions during hydraulic fracturing \citep{Manchanda2017}]
	{
		\includegraphics[height=0.3\textwidth]{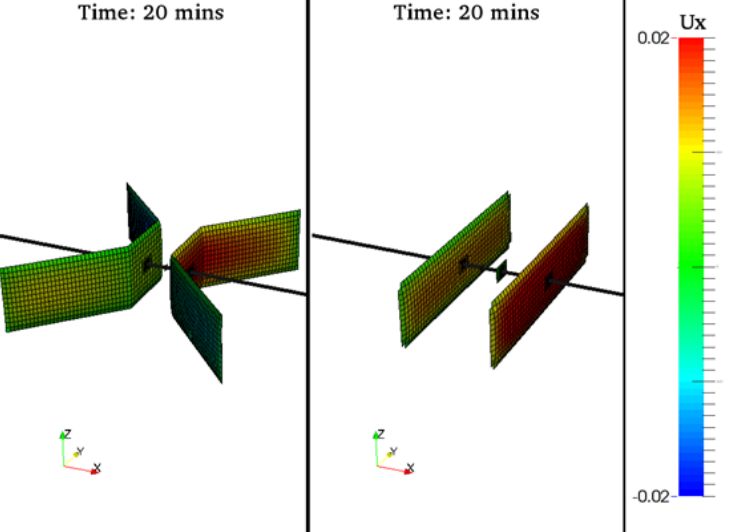}
	}
	\subfigure[Predicting rapid crack propagation in gas-pressurised pipes \citep{Ivankovic2002a}]
	{
		\includegraphics[height=0.3\textwidth]{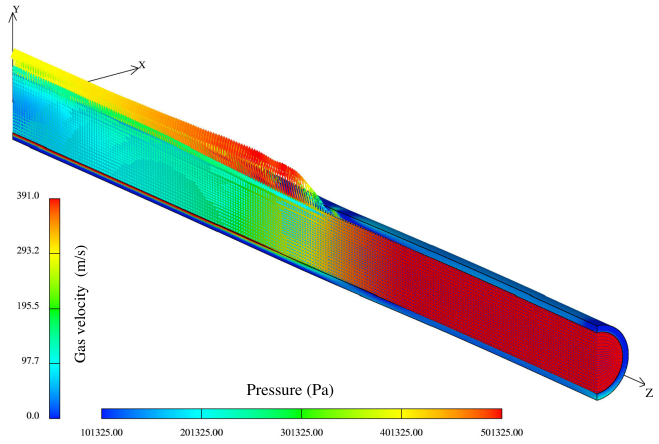}
	}
	\caption{Fracture and adhesive joint examples}
	\label{fig:examplesFracture}
\end{figure}

Example cases are shown in Figure \ref{fig:examplesFracture}, and examples references include:
 \citep{Leevers1993, Ivankovic1994, Demirdzic1996b, Ivankovic1997a, Ivankovic1998, Ivankovic1999, Murphy1999, Stylianou1999, Greenshields2000, Pandya2000a, Pandya2000b, Pandya2000c, Ivankovic2001a, Ivankovic2002a, Ivankovic2002c, Stylianou2002a, Stylianou2002b, Karac2003b, Ivankovic2004b, Rager2005, Murphy2005, Murphy2006, Tropsa2006, Murphy2007, Karac2009a, Karac2009b, Tukovic2010, Karac2011, McAuliffe2011, McAuliffe2012a, McAuliffe2012b, Cooper2012, Georgiou2003a, Georgiou2003b, Ivankovic2004a, Cooper2008, Cooper2010, Tabakovic2010, Carolan2011, Petrovic2011, Carolan2012b, Carolan2013b, Alveen2014, McNamara2014, Alveen2015, Manchanda2015, McNamara2015, McNamara2015b, Lee2017, Manchanda2017, Fallah2018d, ShitingYi2018, Sabbagh-Yazdi2018};

\subsection{Microstructure analysis}
\begin{figure}[htb]
	\centering
	\subfigure[Stress predictions for a hexagonal array of cylindrical inclusions in a boron/aluminum composite \citep{Cavalcante2016}]
	{
		\includegraphics[height=0.45\textwidth]{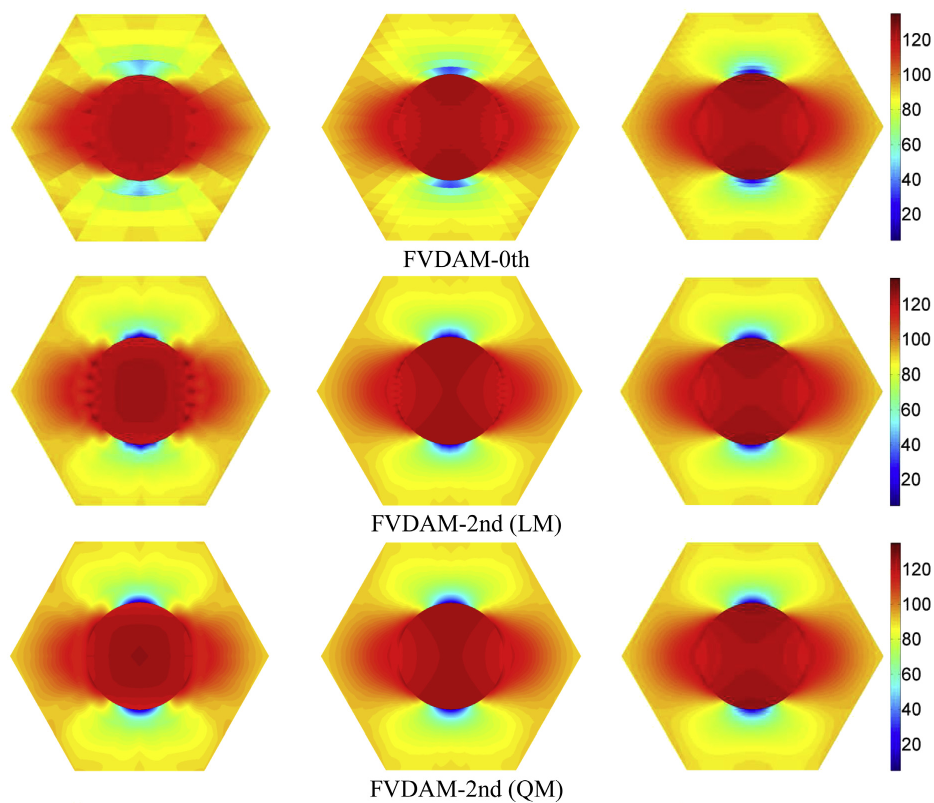}
	}
	\subfigure[Effective plastic strain distributions for a boron/aluminium unidirectional composite with an applied average transverse shear strain \citep{Bansal2006}]
	{
		\includegraphics[height=0.45\textwidth]{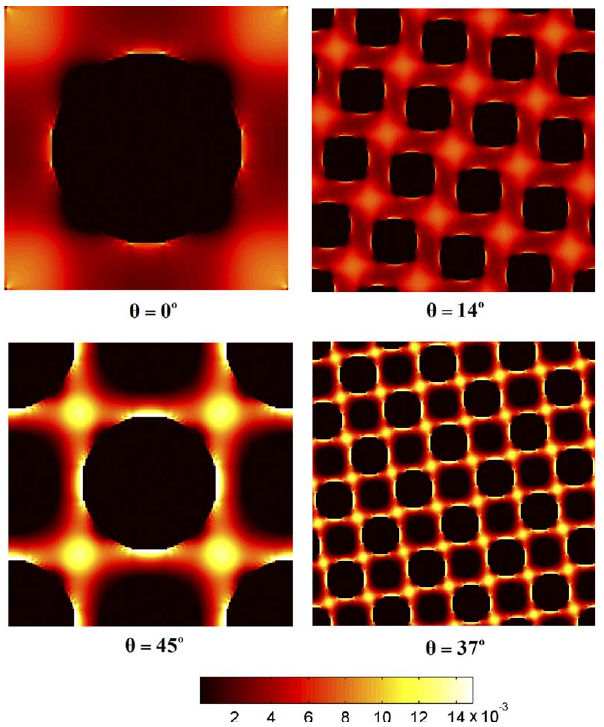}
	}
	\subfigure[Transverse stress distribution in the repeating unit cell of a carbon/aluminium composite loaded by a longitudinal uniaxial stress \citep{Aboudi2011}]
	{
		\includegraphics[height=0.21\textwidth]{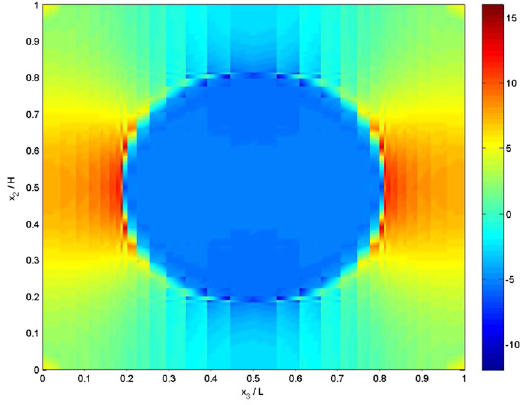}
	}
	\subfigure[Micro-structural analysis of advanced ceramics \citep{McNamara2014}]
	{
		\includegraphics[height=0.21\textwidth]{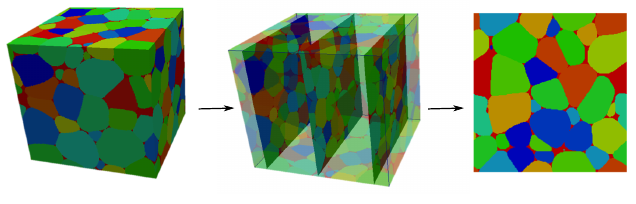}
	}
	\caption{Microstructure analysis examples}
	\label{fig:examplesMicrostructure}
\end{figure}

Example cases are shown in Figure \ref{fig:examplesMicrostructure}, and examples references include:
\citep{Aboudi1991, Pindera1991,
Paley1992,
Aboudi1994,
Aboudi1999,
Aboudi2001a, Aboudi2001b,
Aboudi2002a, Aboudi2002b, Zhong2002,
Aboudi2003, Bansal2003,
Arnold2004, Aboudi2004, Bansal2004, Bednarcyk2004, Pindera2004, Zhong2004,
Aboudi2005a, Aboudi2005b, Bansal2005,
Bansal2006, Cavalcante2006, Charalambakis2006, Pindera2006, 
Aboudi2007, Bruck2007, Cavalcante2007a, Cavalcante2007b, Drago2007, Gattu2007, Pindera2007, Ryvkin2007,
Aboudi2008, Bednarcyk2008, Cavalcante2008, Gattu2008, Paulino2008,
Cavalcante2009, Gao2009, Haj-Ali2009, Khatam2009a, Khatam2009b, Khatam2009c, Pindera2009,
Aboudi2010, Bednarcyk2010, Haj-Ali2010, Khatam2010,
Aboudi2011, Cavalcante2011a, Cavalcante2011b, Chareonsuk2011, Khatam2011,
Carolan2012a, Cavalcante2012a, Cavalcante2012b, Cavalcante2012c, Cavalcante2012d, Haj-Ali2012, Khatam2012, Leonard2012, 
Cavalcante2013,
Alveen2014, Cavalcante2014a, Cavalcante2014b, Cardiff2014g, Leonard2014, McNamara2014, Tu2014,
Alveen2015, Carolan2015, McNamara2015, McNamara2015b,
Cavalcante2016, Tu2016,
Chen2017,
Chen2018a, Chen2018b, Ye2018};

\subsection{Metal forming and casting}
\begin{figure}[htb]
	\centering
	\subfigure[Shape and effective stress predictions in a cast part \citep{Teskeredzic2015b}]
	{
		\includegraphics[height=0.27\textwidth]{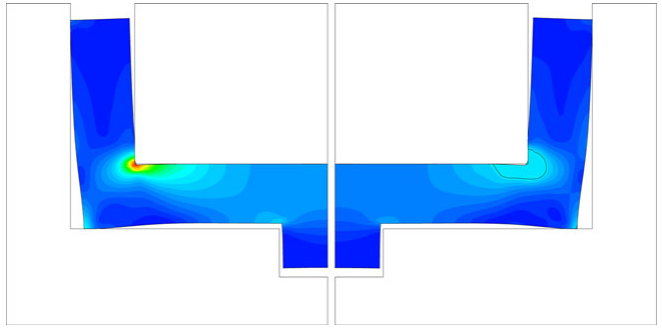}
	}
	\subfigure[Casting of a tunnel section showing the \emph{crack criterion} distribution \citep{Hattel2003}]
	{
		\includegraphics[height=0.27\textwidth]{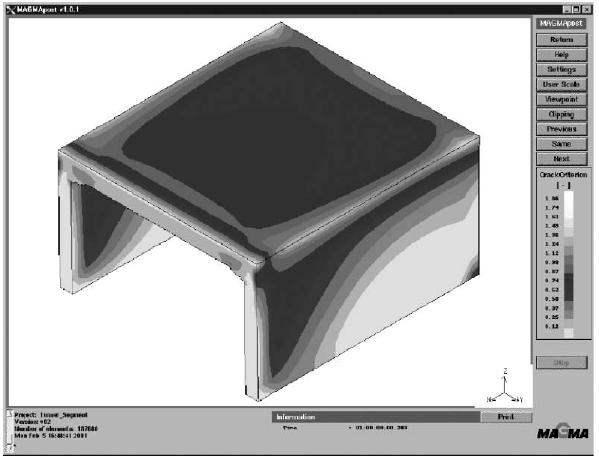}
	}
	\subfigure[Temperature distribution during extrusion \citep{Williams2010}]
	{
		\includegraphics[height=0.35\textwidth]{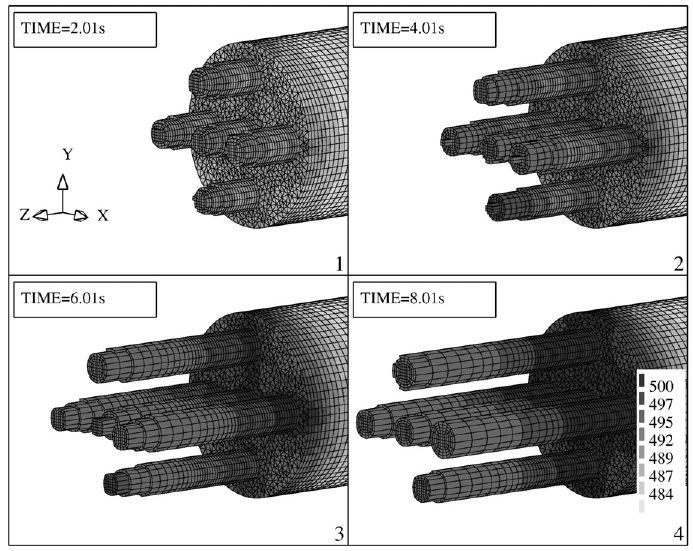}
	}
	\subfigure[Predicted geometry, equivalent plastic strain and hydrostatic pressure distributions in flat wire rolling \citep{Cardiff2016b}]
	{
		\includegraphics[height=0.35\textwidth]{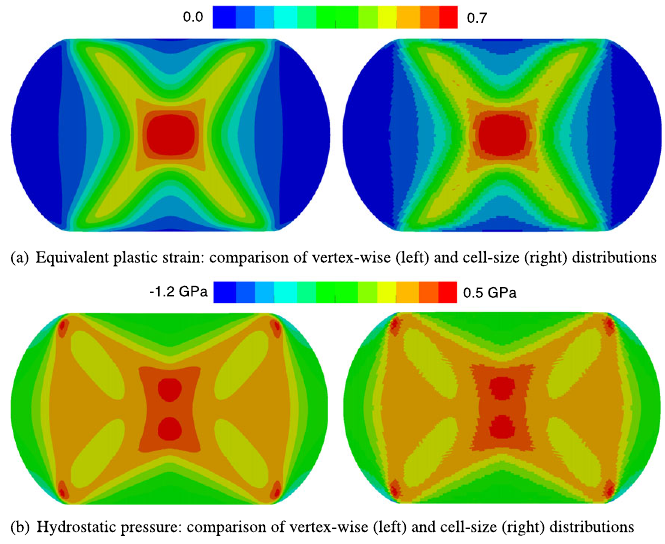}
	}
	\caption{Metal forming and casting examples}
	\label{fig:examplesMetals}
\end{figure}

Example cases are shown in Figure \ref{fig:examplesMetals}, and examples references include:
\citep{Taylor1995b, Pryds1997, Bijelonja2000, Thorborg2001, Williams2001, Basic2002, Williams2002, Teskeredzic2002, Cross2003, Demirdzic2003, Thorborg2003, Williams2003, Basic2005, Sato2006, Cross2007, Jafari2007, Lou2008, Martins2010, Williams2010, Wang2011, Kalkan2011, Mohan2011, Wang2012, Cardiff2014f, Teskeredzic2015a, Teskeredzic2015b, Martins2016, Cardiff2016b, Zhang2016, Bressan2015, Martins2016, Cardiff2017a, Martins2017, Bressan2017, Cardiff2018b, Clancy2018};

\subsection{Biomechanics}
\begin{figure}[htb]
	\centering
	\subfigure[Pressure distribution in a loaded stent-like structure \citep{Haider2018}]
	{
		\includegraphics[height=0.55\textwidth]{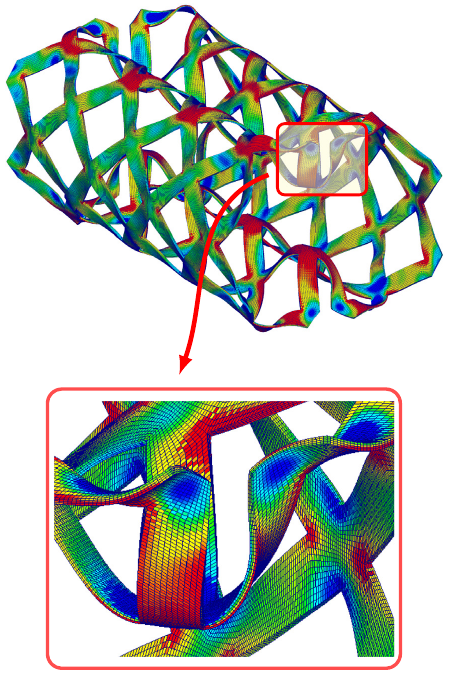}
	}
	\subfigure[Von Mises stress distribution in the natural hip joint during stance \citep{Cardiff2014b}]
	{
		\includegraphics[height=0.55\textwidth]{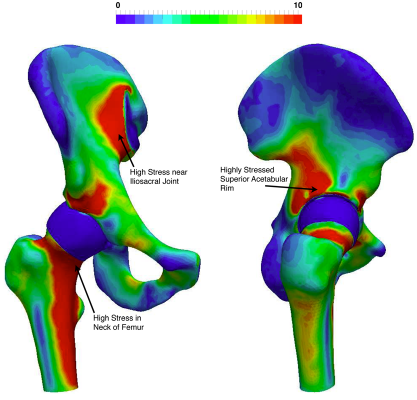}
	}
	\subfigure[Simulation of isometric contraction in the triceps and biceps muscles \citep{Teran2003}]
	{
		\includegraphics[height=0.4\textwidth]{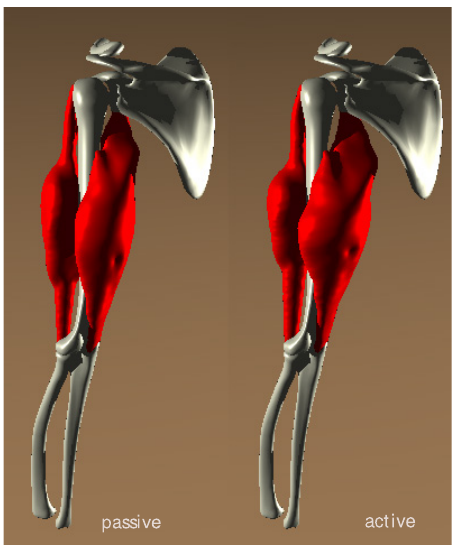}
	}
	\subfigure[Mesh of an idealised carotid artery wall \citep{Kanyanta2009c}]
	{
		\includegraphics[height=0.4\textwidth]{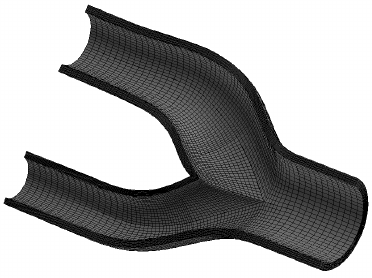}
	}
	\caption{Biomechanics examples}
	\label{fig:examplesBiomechanics}
\end{figure}

Example cases are shown in Figure \ref{fig:examplesBiomechanics}, and examples references include:
\citep{Henry1993a, Harrild1997, Grossman2002, Anthony2003, Teran2003, Alakija2005, Jacquemet2005a, Papadakis2006, Quinn2007a, Quinn2007b, Wang2007, Croft2008, Safari2008, Kanyanta2009a, Kanyanta2009b, Kanyanta2009c, Kelly2009, Safari2009, Cardiff2010, Cardiff2011a, Cardiff2011b, Cardiff2011c, Quinn2011, Cardiff2012b, Wiedemair2012, Cardiff2014b, Cardiff2014c, Parsa2014, Safari2015, Fitzgerald2016, Safari2016, Fitzgerald2017, Muralidharan2017, Oliveira2017a, Oliveira2017b, Haider2018}.

\subsection{Screw compressors}
\begin{figure}[htb]
	\centering
	\subfigure[Displacement and temperature distribution of a rotor in an oil-free compressor \citep{Kovacevic2004b}]
	{
		\includegraphics[height=0.4\textwidth]{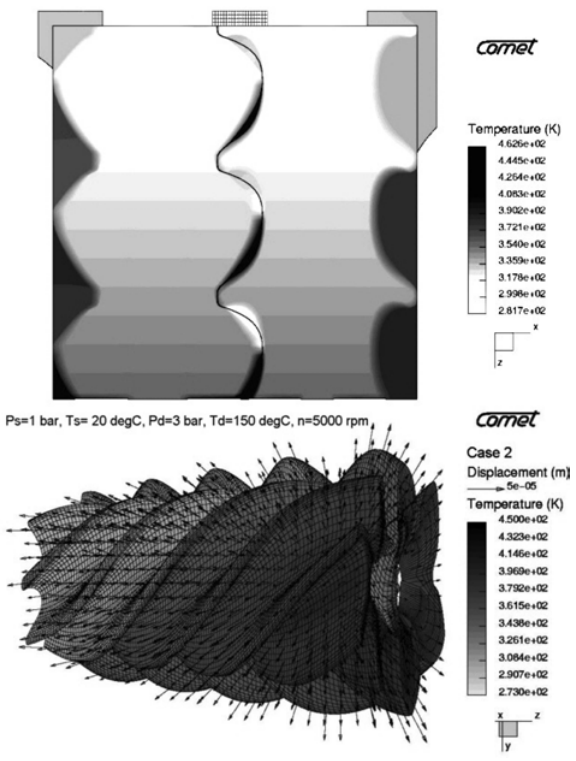}
	}
	\subfigure[Mesh of screw compressor \citep{Kovacevic2007}]
	{
		\includegraphics[height=0.4\textwidth]{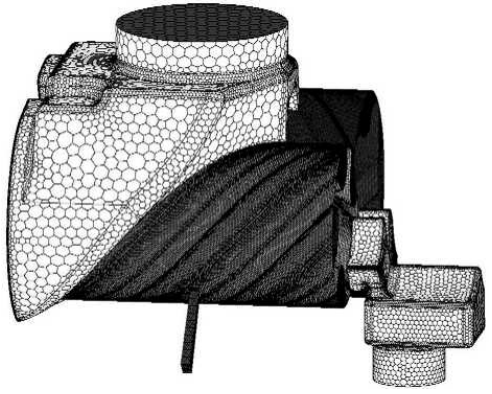}
	}
	\caption{Screw compressor examples}
	\label{fig:examplesScrewCompressor}
\end{figure}

Example cases are shown in Figure \ref{fig:examplesScrewCompressor}, and examples references include:
\citep{Kovacevic2002a, Kovacevic2002b, Kovacevic2002c, Kovacevic2002d, Kovacevic2003, Kovacevic2004a, Kovacevic2004b,  Stosic2005, Kovacevic2006, Kovacevic2007, Stosic2007, Kovacevic2011, Smith2014};

\subsection{Geomechanics and poroelasticity}
\begin{figure}[htb]
	\centering
	\subfigure[Elastic waves in a heterogenous medium \citep{Zhang2002}]
	{
		\includegraphics[height=0.4\textwidth]{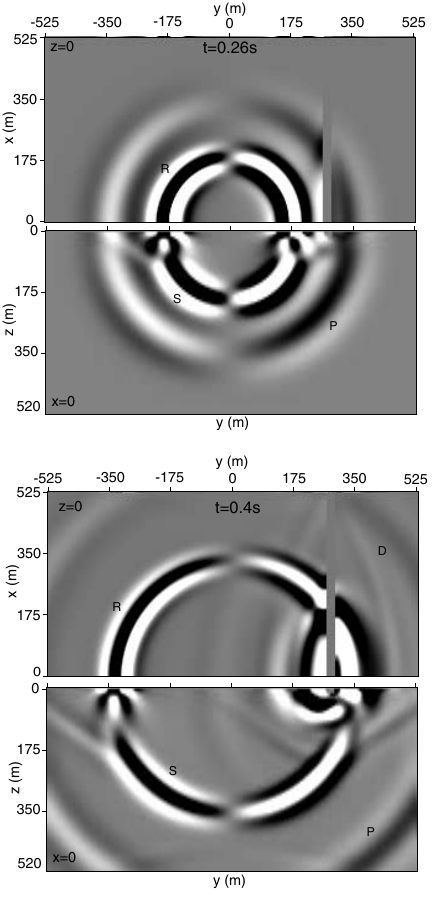}
	}
	\quad\quad
	\subfigure[Pore-pressure predictions near a hydraulic fracture \citep{Lee2015}]
	{
		\includegraphics[height=0.35\textwidth]{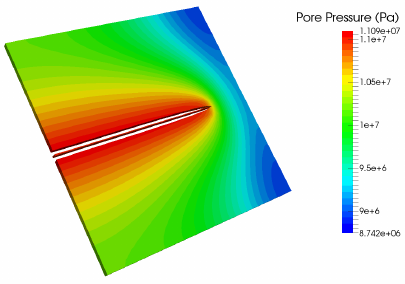}
	}
	\caption{Geomechanics and poroelasticity examples}
	\label{fig:examplesGeomechanics}
\end{figure}

Example cases are shown in Figure \ref{fig:examplesGeomechanics}, and examples references include:
\citep{Dormy1995, Zhang1999, Zhang2002, Zhang2004, Liu2004, Liu2005, Shaw2005, Gao2006, Bijelonja2011b, Tang2015, Cardiff2015a, Cardiff2015b, Lee2015, Manchanda2017}.


\section{Softwares employing the finite volume method for solid mechanics}
A number of software have, or previously have, implemented versions of the finite volume method for solid mechanics;
these software, in alphabetical order, include:
\begin{itemize}
	\item {COMET / STAR-CD} (commercial software) \citep{Demirdzic1994a, Demirdzic1995, Ivankovic1997a, Demirdzic1997a, Demirdzic1997c, Dzaferovic2000, Teskeredzic2002, Kovacevic2002a, Kovacevic2002b, Kovacevic2002c, Kovacevic2002d, Kovacevic2003, Kovacevic2004a, Kovacevic2004b, Basic2005, Bijelonja2005a, Bijelonja2005b, Demirdzic2005, Bijelonja2006, Kovacevic2006, Kovacevic2007, Stosic2005, Stosic2007, Bijelonja2011a, Smith2014, Teskeredzic2015a, Teskeredzic2015b, Bijelonja2017};

	\item {FOAM / OpenFOAM} (open-source software) \citep{Weller1998, Greenshields1999b, Greenshields2000, Jasak2000a, Jasak2000b, Karac2003b, Giannopapa2004a, Giannopapa2004b, Greenshields2005, Giannopapa2006, Jasak2007, Tukovic2007b, Giannopapa2008, Kanyanta2009a, Karac2009a, Karac2009b, Kelly2010, Tabakovic2010, Petrovic2011, Cardiff2012a, Carolan2012b, Kelly2012, Wiedemair2012, Carolan2013a, Carolan2013b, Habchi2013, Alveen2014, Cardiff2014a, Cardiff2014b, Cardiff2014c, Tukovic2014, Cardiff2015a, Carolan2015, Tang2015, McNamara2014, McNamara2015, Cardiff2016a, Cardiff2016b, Elsafti2016, Sekutkovski2016, Cardiff2017b, Haider2017, Manchanda2017, Cardiff2018a, Haider2018, Holzmann2018, Tukovic2018a, Tukovic2018b};

	\item GTEA (in-house software) \citep{Zhu2012};

	\item {MulPhys} (in-house software) \citep{Limache2007};

	\item {NASIR} (in-house software)
	\citep{Sabbagh-Yazdi2008a, Sabbagh-Yazdi2008b, Alkhamis2008, Sabbagh-Yazdi2009, Sabbagh-Yazdi2011a, Sabbagh-Yazdi2011b, Sabbagh-Yazdi2011c, Sabbagh-Yazdi2012a, Sabbagh-Yazdi2012b, Sabbagh-Yazdi2018};

	\item {PHOENICS} (commercial software) \citep{Beale1990a, Beale1990b, Bukhari1990, Bukhari1991, Spalding1993, Spalding1997, Spalding1998a, Spalding1998b, Spalding2002, Spalding2004, Spalding2008};

	\item {PHYSICA} (in-house software) \citep{Bailey1997, Bailey1999a, Bailey1999b, Bailey1999c, Taylor2003};

	\item {TRANSPORE} (in-house software) \citep{Perre1995}.
\end{itemize}


\section{Conclusions and challenges}
The various finite volume approaches to solid mechanics can be seen to share many similarities, and in fact can all be described in a general unified manner (Section \ref{sec:finiteVolumeVariantsDiscussion}).
When comparing the finite volume approach with the finite element approach, the likenesses are clear:
both approaches adopt the same general strategy to discretise space into cells/elements,
both use similar data storage structures,
and both follow similar approaches to assemble their corresponding characteristic matrices.
The main differences lie in how the local integration domains are constructed, how the local integration is performed, and their fundamental philosophy:
finite volume approaches are rooted in balance laws, where the governing equation is enforced by summing forces/fluxes acting on a control volume;
in contrast, finite element methods adopt a more mathematical approach, based on variational methods, where the weak form of the governing equation is imposed in a volumetrically-averaged sense.

To end this article, we state three main challenges we see for the development of finite volume solid mechanics, such that its strengths and weaknesses may be rigorously explored in the context of solid mechanics, and its merits, relative to other similar approaches, may become clear to the computational solid mechanics community at large:
\vspace{-0.015\textwidth}
\begin{enumerate}
	\item Awareness:\\
		Three decades after the first contributions to the field, there is still a general lack of awareness around the capabilities of the finite volume method in the sphere of solid mechanics;
		consequently, in the worst case, prestigious journals inadvertently assign inappropriate reviewers,
		reviewers who, having only limited knowledge of the area, accept poor or reject good articles, for example, see \citep{Demirdzic2015},
		and authors are unaware of major developments in the field when surveying the literature.

	\item Benchmarking:\\
		As should be clear from this review, numerous differing finite volume formulations are possible;
		although some variants have been developed specifically for specialised applications, the comparison of the differing approaches on standard benchmarks is rare.
		Without such comparisons, in terms of efficiency, accuracy and robustness, it will not be possible to determine which approaches are optimal for certain classes of problems.
		Furthermore, given the trends in modern computing, the suitability of proposed approaches for execution on large-scale distributed memory clusters ($>$1000s CPU cores) should be further explored.
		Simulations using hundreds of millions of cells are commonplace at CFD conferences: clearly finite volume-based solid mechanics procedures have great potential.
		Similarly, finite volume variants should be benchmarked against alternative approaches, such as the finite element method, to determine relative merits, not just on academic standard cases but on complex industrial cases to test robustness.
		Processes such as \emph{round robin} benchmarking series, for example \citep{Hitchings1987}, may offer one solution.

	\item Code dissemination:\\
		Where possible, code for published procedures should be shared for academic scrutiny:
		such distribution has the potential to:
		(a) accelerate academic progress, as others learn from and build on methods, as well as aiding in the discovery and resolution of errors;
		(b) facilitate ease of understanding and ease of implementation;
		(c) allow direct comparison of methods;
		(d) provide insight into the algorithm intricacies that may not be clear from academic articles;
		and (e) allow faster integration into commercial software and industrial use.
\end{enumerate}



%

\appendix
\section{Appendix A: Table of most cited articles related to the finite volume method for solid mechanics}

Table \ref{table:mostCitedArticles} lists the most cited articles related to the finite volume method for solid mechanics;
the references have been listed in order of decreasing number of citations, and only articles with greater than fifty citations have been included, according to Google Scholar citations on 25$^{th}$ August 2018.
As noted in the body of the article, care should be taken when interpreting the data, as the number of citations may not be directly proportional to impact on the field; for example, \citet{Weller1998} has by far the greatest number of citations; however, a significant percentage of its received citations are related to its computational fluid mechanics developments, rather than its solid mechanics contributions.
\begin{table}[htb]
  \centering
 	\ra{1.3}
  	\begin{tabular}{@{}lll@{}}
	\toprule
	Reference & Formulation-type & Number of citations \\
	\midrule
	\citet{Weller1998} & cell-centred & 1988 \\
	\citet{Demirdzic1995} & cell-centred & 408 \\
	\citet{Idelsohn1994} & vertex-centred and cell-centred & 221 \\
	\citet{Demirdzic1994a} & cell-centred & 189 \\
	\citet{Jasak2000a} & cell-centred & 185 \\
	\citet{Pindera2009} & HFGMC/HOTFGM/FVDAM & 185 \\
	\citet{Demirdzic1993} & cell-centred & 172 \\
	\citet{Bailey1995} & vertex-centred & 143 \\
	\citet{Onate1994} & vertex-centred and cell-centred & 128 \\
	\citet{Fryer1991} & vertex-centred & 103 \\
	\citet{Slone2002} & vertex-centred & 101 \\
	\citet{Slone2003} & vertex-centred & 97 \\
	\citet{Bansal2006} & HFGMC/HOTFGM/FVDAM & 94 \\
	\citet{Voller2009} & vertex-centred & 91 \\
	\citet{Bijelonja2006} & cell-centred & 87 \\
	\citet{Taylor1995a} & vertex-centred & 86 \\
	\citet{Cavalcante2007a} & HFGMC/HOTFGM/FVDAM & 80 \\
	\citet{Taylor2003} & vertex-centred & 70 \\
	\citet{Khatam2009a} & HFGMC/HOTFGM/FVDAM & 70 \\
	\citet{Fallah2000a} & vertex-centred & 64 \\
	\citet{Ivankovic1994} & cell-centred & 61 \\
	\citet{Karac2011} & cell-centred & 59 \\
	\citet{Taylor2002a} & vertex-centred & 58 \\
	\citet{Bailey1999b} & vertex-centred & 57 \\
	\citet{Slone2004} & vertex-centred & 56 \\
	\citet{Murphy2005} & cell-centred & 55 \\
	\citet{Tukovic2007a} & cell-centred & 54 \\
	\citet{Demirdzic1997c} & cell-centred & 52 \\
	\citet{Bailey1996} & vertex-centred & 51 \\
	\bottomrule
	\end{tabular}
   \caption{Most cited articles related to the finite volume method for solid mechanics from Google Scholar citations on 25$^{th}$ August 2018. Approaches for heterogeneous periodic microstructures are indicated by HFGMC/HOTFGM/FVDAM.}
   \label{table:mostCitedArticles}
\end{table}

\section{Appendix B: Overview of the discretisation used in HOTFGM/HFGMC/FVDAM approaches}

There are a variety of related methods with finite volume attributes which have been designed for the analysis of heterogenous microstructures.
The related methods include the \emph{higher-order theory for functionally graded material} (HOTFGM) \citep{Aboudi1999}, the \emph{high-fidelity generalised method of cells} (HFGMC) \citep{Aboudi2001a, Aboudi2001b, Haj-Ali2009, Haj-Ali2012}, and the \emph{finite volume direct averaging micromechanics} (FVDAM) theory \citep{Bansal2005, Bansal2006}.
A brief summary of the methods is given here, and readers are referred to \citet{Aboudi1999}, \citet{Bansal2006} and \citet{Cavalcante2012c} for further details.

\begin{figure}[htb]
   \centering
   \includegraphics[height=0.5\textwidth]{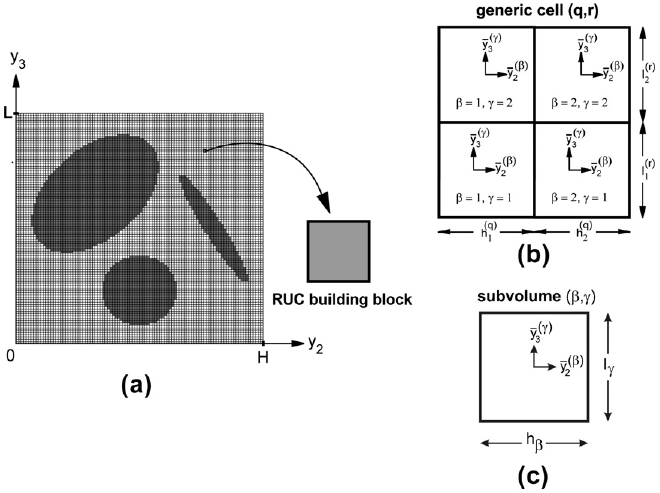}
   \caption{(a) Unit cell of a periodic material with microstructure discretised into rectangular building blocks. (b) Two-level discretisation employed by HFGMC into $(q, r)$ generic cells further subdivided into four $(b, c)$ subcells. (c) Single-level discretisation employed in the FVDAM theory into stand-alone $(b, c)$ sub-volumes. Figure taken from \citet{Cavalcante2012c}}
   \label{fig:cavalcante2012mesh}
\end{figure}

The HOTFGM and HFGMC approaches start by spatially discretising the solution domain into rectangular so-called \emph{generic cells}, which are further split into a second discretisation level containing four rectangular sub-cells (Figure \ref{fig:cavalcante2012mesh}); for brevity and clarity, the description here has been limited to two dimensions; however, the approaches have been extended to three dimensions, as described in \citet{Aboudi1999}.
As a consequence of the assumed orthogonal Cartesian mesh, curved interfaces between material phases are approximated in a castellated \emph{staircase} manner, as shown in Figure \ref{fig:cavalcante2012mesh}; this limitation was later removed by the FVDAM approach with extension to unstructured quadrilateral meshes.
By considering the unit cell of a periodic material, 
the displacement field can be decomposed into average and fluctuating components, $\boldsymbol{u} = \bar{\boldsymbol{u}} + \boldsymbol{u}'$, where the average displacement is determined from the specified macroscopic average strains, $\bar{\boldsymbol{u}} = \bar{\boldsymbol{\epsilon}} \boldsymbol{x}$.
Within each sub-cell, the fluctuating displacement is then assumed to vary quadratically as a function of the local coordinates, $\bar{y}_2$ and $\bar{y}_3$:
\begin{eqnarray}
	\boldsymbol{u}' (\bar{y}_2, \bar{y}_3) =
	\boldsymbol{W}_{00}
	+ \bar{y}_2 \boldsymbol{W}_{10}
	+ \bar{y}_3 \boldsymbol{W}_{01}
	+ \frac{1}{2}\left( 3 \bar{y}_2^2 - \frac{h^2}{4} \right) \boldsymbol{W}_{20}
	+ \frac{1}{2}\left( 3 \bar{y}_3^2 - \frac{l^2}{4} \right) \boldsymbol{W}_{02}
\end{eqnarray}
where $h$ and $l$ are the width and height respectively of the sub-cell;
$\boldsymbol{W}_{00}$, $\boldsymbol{W}_{10}$, $\boldsymbol{W}_{01}$, $\boldsymbol{W}_{20}$, and $\boldsymbol{W}_{02}$ are unknown vector displacement coefficients, each with three components; the $\boldsymbol{W}_{00}$ component corresponds to the unknown displacement at the centre of the sub-cell, while the remaining coefficients correspond to higher-order displacement contributions within the sub-cell.
Accordingly, there are $5 \times 3 = 15$ unknown displacement coefficients within each sub-cell and hence $4 \times 15 = 60$ within each generic cell; in three dimensions, there are 168 unknown quantities.
For brevity here, the $(\gamma)$ and $(\beta)$ superscripts indicating the sub-cell have been dropped \ie $\bar{y}_2 = \bar{y}_2^{(\beta)}$, $\bar{y}_3 = \bar{y}_3^{(\gamma)}$, \etc 
It is also worth pointing out that although the approach has been developed for periodic microstructures, the method can also be used for general structural stress analysis by assuming the average displacement $\bar{\boldsymbol{u}}$ to be zero.

To determine the unknown displacement coefficients, the $0^{th}$, $1^{st}$ and $2^{nd}$ \emph{moments} of momentum conservation are applied to each sub-cell, in addition to the enforcement of traction and displacement continuity between sub-cells and generic cells, and inclusion of boundary conditions.
A characteristic of the method, which is not possessed by the other finite volume variants, is the enforcement of these so-called \emph{moments} of the governing equation.
To achieve this, the governing equation (Equation \ref{eq:momentumIntegral}), where temporal and body force terms have been neglected, is written in terms of a so-called stress moment $\boldsymbol{S}$:
\begin{eqnarray}
	\oint_\Gamma \boldsymbol{n} \cdot
	\boldsymbol{S}
	\;\text{d}\Gamma
	&=& \boldsymbol{0}
\end{eqnarray}
with the stress moment is defined as:
\begin{eqnarray}
	\boldsymbol{S} = \frac{1}{h l} \int_{-\frac{h}{2}}^{\frac{h}{2}} \int_{-\frac{l}{2}}^{\frac{l}{2}}
	\left[ \boldsymbol{\sigma} \, \bar{y}_2^m \, \bar{y}_3^n \right]
	\; \text{d} \bar{y}_2 \; \text{d} \bar{y}_3
\end{eqnarray}
The exponents $m$ and $n$ indicate the order of the equation; for example, when $m = n = 0$, the relation reduces to conservation of force; when $m = 1$ and $n =1 $ the relation represents conservation of angular momentum; while for $m > 1$ and $n > 1$ the relation represents conservation of higher stress moments.
Note: $m$ is not related to the time-step counter in Equation \ref{eq:discretiseAcceleration}.

In this way, it is possible to assemble a system of $60 N$ algebraic equations of the standard form, $[\boldsymbol{K}][\boldsymbol{U}] = [\boldsymbol{F}]$, where $N$ is the number of generic cells in the solution domain, $[\boldsymbol{U}]$ is a vector of unknown displacement coefficients $\boldsymbol{W}$, the global stiffness matrix $[\boldsymbol{K}]$ is a function of the sub-cell dimensions and mechanical properties, and the global force vector $[\boldsymbol{F}]$ contains contributions from boundary condition and nonlinear material stresses.
The linear system is inverted to give the displacements distributions within the sub-cells.

The HOTFGM and HFGMC approaches described above provide the basis for the subsequent FVDAM approach;
the FVDAM approach differs from the HOTFGM and HFGMC methods in a number of ways:
\vspace{-0.02\textwidth}
\begin{enumerate}
	\item[(a)] The two-level spatial domain decomposition (generic cells and sub-cells) of the HOTFGM/HFGMC methods is replaced by one-level of discretisation/cells;
	\item[(b)] The displacement coefficients within each cell $\boldsymbol{W}$ are expressed in terms of surface-averaged displacements \ie displacement averaged at each cell surface;
	\item[(c)] Higher order moments of the equilibrium equation are not used;
	\item[(d)] In the parametric form of the FVDAM, the use of parametric mapping with a parent/reference cell allows the use of an unstructured mesh (similar to the finite element method), instead of the orthogonal Cartesian mesh of the HOTFGM/HFGMC approaches (see Figure \ref{fig:HFGMCmesh});
	\item[(e)] In the assembled system of algebraic equations $[\boldsymbol{K}][\boldsymbol{U}] = [\boldsymbol{F}]$, the solution vector $[\boldsymbol{U}]$ contains cell surface-averaged displacements, as opposed to sub-cell displacement coefficients.
\end{enumerate}

Further technicals details of the HOTFGM, HFGMC and FVDAM methods can be found in \citet{Aboudi1999, Aboudi2001a, Aboudi2001b, Haj-Ali2009, Haj-Ali2012, Bansal2005, Bansal2006, Cavalcante2012c} and \citet{Cavalcante2016}.

\section*{References}


\Urlmuskip=0mu plus 1mu\relax

\bibliographystyle{unsrtnat}
\bibliography{Bibliography}

\begin{thebibliography}{611}
\providecommand{\natexlab}[1]{#1}
\providecommand{\url}[1]{\texttt{#1}}
\expandafter\ifx\csname urlstyle\endcsname\relax
  \providecommand{\doi}[1]{doi: #1}\else
  \providecommand{\doi}{doi: \begingroup \urlstyle{rm}\Url}\fi

\bibitem[Runchal(2017)]{Runchal2017}
A.~K. Runchal.
\newblock Tributes to an exceptional life: {D. Brian Spalding}, {9 January 1923
  - 27 November 2016}.
\newblock \emph{Available at:
  \url{https://www.astfe.org/doc/Brian_Spalding-Tributes_to_an_Exceptional_Life.pdf}},
  2017.

\bibitem[Demird{\v{z}}i\'{c} et~al.(1988)Demird{\v{z}}i\'{c}, Martinovi\'{c},
  and Ivankovi\'{c}]{Demirdzic1988}
I.~Demird{\v{z}}i\'{c}, D.~Martinovi\'{c}, and A.~Ivankovi\'{c}.
\newblock Numerical simulation of thermal deformation in welded workpiece.
\newblock \emph{Zavarivanje}, 31:\penalty0 209--219, 1988.
\newblock In Croatian. English translation available at
  \url{https://www.researchgate.net/profile/Alojz_Ivankovic/publication/296148474_Numerical_simulation_of_thermal_deformation_in_welded_workpiece/links/5d07642ba6fdcc39f12219eb/Numerical-simulation-of-thermal-deformation-in-welded-workpiece.pdf}.

\bibitem[Demird{\v{z}}i\'{c} and Martinovi\'{c}(1993)]{Demirdzic1993}
I.~Demird{\v{z}}i\'{c} and D.~Martinovi\'{c}.
\newblock Finite volume method for thermo-elasto-plastic stress analysis.
\newblock \emph{Computer Methods in Applied Mechanics and Engineering},
  109:\penalty0 331--349, 1993.

\bibitem[Bailey et~al.(September, 1991)Bailey, Fryer, Cross, and
  Chow]{Bailey1991}
C.~Bailey, Y.D. Fryer, M.~Cross, and P.~Chow.
\newblock Predicting the deformation of castings in moulds using a control
  volume approach on unstructured meshes.
\newblock In M.~Cross, J.~F.~T. Pittman, and R.~D. Wood., editors,
  \emph{Mathematical modelling for materials processing : based on the
  proceedings of a conference on Mathematical Modelling of Materials
  Processing, organized by the Institute of Mathematics and its Applications},
  University of Bristol, September, 1991.

\bibitem[Fryer et~al.(1991)Fryer, Bailey, Cross, and Lai]{Fryer1991}
Y.~D. Fryer, C.~Bailey, M.~Cross, and C.-H. Lai.
\newblock A control volume procedure for solving the elastic stress-strain
  equations on an unstructured mesh.
\newblock \emph{Applied Mathematical Modelling}, 15:\penalty0 639--645, 1991.

\bibitem[Zienkiewicz and O{\~{n}}ate(1991)]{Zienkiewicz1991}
O.~C. Zienkiewicz and E.~O{\~{n}}ate.
\newblock Finite volumes vs finite elements: Is there really a choice?
\newblock \emph{Nonlinear Computational Mechanics. State of the Art}, pages
  240--254, 1991.

\bibitem[Cross et~al.(1992)Cross, Bailey, Chow, and Peircleous]{Cross1992}
M.~Cross, C.~Bailey, P.~Chow, and K.~Peircleous.
\newblock Towards an integrated control volume unstructured mesh code for the
  simulation of all macroscopic processes involved in shape casting.
\newblock In R.~D.~Wood Chenot and O.~C. Zienkiewicz, editors, \emph{Numerical
  Methods in Industrial Processes NUMIFORM 92}, pages 787--792, Rotterdam,
  1992. Belkema.

\bibitem[O{\~{n}}ate et~al.(1992)O{\~{n}}ate, Cervera, and
  Zienkiewicz]{Onate1992}
E.~O{\~{n}}ate, M.~Cervera, and O.~C. Zienkiewicz.
\newblock A study of the finite volume format for structural mechanics.
\newblock Technical report, {Internal Report Publication No. 15, CIMNE,
  Barcelona}, 1992.

\bibitem[Cross(1993)]{Cross1993}
M.~Cross.
\newblock Development of novel computational technique for the next generation
  of software tools for casting simulation.
\newblock In L.~Katgerman T.S~Piwonka, V.R.~Voler, editor, \emph{Modelling of
  Casting, Welding and Advanced Solidification Processes VI, TMS}, pages
  115--126, 1993.

\bibitem[Fryer et~al.(1993)Fryer, Bailey, Cross, and Chow]{Fryer1993}
Y.~D. Fryer, C.~Bailey, M.~Cross, and P.~Chow.
\newblock Predicting micro-porosity in shape casting using an integrated
  control volume unstructured mesh framework.
\newblock In V.~R.~Voler T.S.~Piwonka and L.~Katgerman, editors,
  \emph{Modelling of Casting, Welding and Advanced Solidification Processes
  {VI}}, pages 143--152, 1993.

\bibitem[Idelsohn and O{\~{n}}ate(1994)]{Idelsohn1994}
S.~R. Idelsohn and E.~O{\~{n}}ate.
\newblock Finite volumes and finite elements: Two `good friends'.
\newblock \emph{International Journal for Numerical Methods in Engineering},
  37:\penalty0 3323--3341, 1994.

\bibitem[O{\~{n}}ate et~al.(1994)O{\~{n}}ate, Cervera, and
  Zienkiewicz]{Onate1994}
E.~O{\~{n}}ate, M.~Cervera, and O.~C. Zienkiewicz.
\newblock A finite volume format for structural mechanics.
\newblock \emph{International Journal for Numerical Methods in Engineering},
  37:\penalty0 181--201, 1994.

\bibitem[Beale and Elias(1990{\natexlab{a}})]{Beale1990a}
S.~B. Beale and S.~R. Elias.
\newblock Numerical solution of two-dimensional elasticity problems by means of
  a {SIMPLE}-based finite-difference scheme.
\newblock Technical report, Institute for Mechanical Engineering, National
  Research Council, Ottawa, Ont. TR-LT-020 (NRC No. 32090), 1990{\natexlab{a}}.

\bibitem[Beale and Elias(1990{\natexlab{b}})]{Beale1990b}
S.~B. Beale and S.~R. Elias.
\newblock Stress distribution in a plate subject to uniaxial loading.
\newblock \emph{{PHOENICS} Journal of Computational Fluid Dynamics}, 3\penalty0
  (3):\penalty0 255--287, 1990{\natexlab{b}}.

\bibitem[Bukhari et~al.(1990)Bukhari, Qin, and Spalding]{Bukhari1990}
K.~M. Bukhari, H.~Q. Qin, and D.~B. Spalding.
\newblock Progress report (to {Rolls-Royce Ltd}) on the calculation of thermal
  stresses in bodies of revolution.
\newblock Technical report, {CHAM Ltd}, 1990.

\bibitem[Hattel and Hansen(Oct. 1990)]{Hattel1990}
J.~H. Hattel and P.~N. Hansen.
\newblock {FDM} solutions of the thermoelastic equations using a staggered
  grid.
\newblock In \emph{{Danish-German-Polish Workshop on Application of Computer
  Methods in Practice}}, Warsaw, Poland, Oct. 1990.

\bibitem[Bukhari et~al.(1991)Bukhari, Hamill, Qin, and Spalding]{Bukhari1991}
K.~M. Bukhari, I.~S. Hamill, H.~Q. Qin, and D.~B. Spalding.
\newblock Stress-analysis simulations in {PHOENICS}.
\newblock Technical report, {CHAM Ltd}, 1991.

\bibitem[Hattel(1992)]{Hattel1992}
J.~H. Hattel.
\newblock Analysis of thermal induced stresses in die molds.
\newblock In J.~T.~Cross Lacaze, M.~Rappaz, G.~Sciama, and I.~Svensson,
  editors, \emph{Examples of European Expertise in the Computer Simulation of
  Solidification and Casting}. B--18, Ecole de Mines, Nancy, France, 1992.

\bibitem[Hattel et~al.(1993{\natexlab{a}})Hattel, Hansen, and
  Hansen]{Hattel1993a}
J.~H. Hattel, P.~N. Hansen, and L.~F. Hansen.
\newblock Analysis of thermal induced stresses in die casting using a novel
  control volume technique.
\newblock In V.~R.~Voller T.S.~Piwonka and L.~Katgerman, editors,
  \emph{Modelling of Casting, Welding and Advanced Solidification Processes
  {VI, TMS}}, pages 585--592, Palm Coast, Florida, USA, 1993{\natexlab{a}}.

\bibitem[Sevilla et~al.(2018{\natexlab{a}})Sevilla, Giacomini, and
  Huerta]{Sevilla2018a}
R.~Sevilla, M.~Giacomini, and A.~Huerta.
\newblock A face-centred finite volume method for second-order elliptic
  problems.
\newblock \emph{International Journal for Numerical Methods in Engineering},
  0:\penalty0 1--29, 2018{\natexlab{a}}.

\bibitem[Sevilla et~al.(2018{\natexlab{b}})Sevilla, Giacomini, and
  Huerta]{Sevilla2018b}
R.~Sevilla, M.~Giacomini, and A.~Huerta.
\newblock A locking-free face-centred finite volume {(FCFV)} method for linear
  elasticity.
\newblock \emph{preprint}, 2018{\natexlab{b}}.
\newblock arXiv:1806.07500v1 [math.NA], available at:
  \url{https://arxiv.org/pdf/1806.07500}.

\bibitem[Ebrahimnejad et~al.(2014)Ebrahimnejad, Fallah, and
  Khoei]{Ebrahimnejad2014}
M.~Ebrahimnejad, N.~Fallah, and A.~R. Khoei.
\newblock New approximation functions in the meshless finite volume method for
  {2D} elasticity problems.
\newblock \emph{Engineering Analysis with Boundary Elements}, 46:\penalty0
  10--22, 2014.

\bibitem[Fallah and Delzendeh(2018)]{Fallah2018a}
N.~Fallah and M.~Delzendeh.
\newblock Free vibration analysis of laminated composite plates using meshless
  finite volume method.
\newblock \emph{Engineering Analysis with Boundary Elements}, 88:\penalty0
  132--144, 2018.

\bibitem[Bailey and Cross(1995)]{Bailey1995}
C.~Bailey and M.~Cross.
\newblock A finite volume procedure to solve elastic solid mechanics problems
  in three dimensions on an unstructured mesh.
\newblock \emph{International Journal for Numerical Methods in Engineering},
  38:\penalty0 1757--1776, 1995.

\bibitem[Cardiff et~al.(2016{\natexlab{a}})Cardiff, Tukovi\'{c}, Jasak, and
  Ivankovi\'{c}]{Cardiff2016a}
P.~Cardiff, Tukovi\'{c}, H.~Jasak, and A.~Ivankovi\'{c}.
\newblock A block-coupled finite volume methodology for linear elasticity and
  unstructured meshes.
\newblock \emph{Computers \& Structures}, 175:\penalty0 100--122,
  2016{\natexlab{a}}.
\newblock \doi{10.1016/j.compstruc.2016.07.004}.

\bibitem[Kluth and Despr\'{e}s(2010)]{Kluth2010}
G.~Kluth and B.~Despr\'{e}s.
\newblock Discretization of hyperelasticity on unstructured mesh with a
  cell-centered {Lagrangian} scheme.
\newblock \emph{Journal of Computational Physics}, 229:\penalty0 9092--9118,
  2010.

\bibitem[Lee et~al.(2013)Lee, Gil, and Bonet]{Lee2013}
C.~H. Lee, A.~J. Gil, and J.~Bonet.
\newblock Development of a cell centred upwind finite volume algorithm for a
  new conservation law formulation in structural dynamics.
\newblock \emph{Computers \& Structures}, 118:\penalty0 13--38, 2013.

\bibitem[Haider et~al.(2017)Haider, Lee, Gil, and Bonet]{Haider2017}
J.~Haider, C.~H. Lee, A.~J. Gil, and J.~Bonet.
\newblock A first order hyperbolic framework for large strain computational
  solid dynamics: An upwind cell centred total {Lagrangian} scheme.
\newblock \emph{International Journal for Numerical Methods in Engineering},
  109:\penalty0 407--456, 2017.

\bibitem[Demird{\v{z}}i\'{c} and Muzaferija(1995)]{Demirdzic1995}
I.~Demird{\v{z}}i\'{c} and S.~Muzaferija.
\newblock Numerical method for coupled fluid flow, heat transfer and stress
  analysis using unstructured moving meshes with cells of arbitrary topology.
\newblock \emph{Computer Methods in Applied Mechanics and Engineering},
  125:\penalty0 235--255, 1995.

\bibitem[Cardiff et~al.(2016{\natexlab{b}})Cardiff, Tukovi\'{c}, De~Jaeger,
  Clancy, and Ivankovi\'{c}]{Cardiff2016b}
P.~Cardiff, Tukovi\'{c}, P.~De~Jaeger, M.~Clancy, and A.~Ivankovi\'{c}.
\newblock A {Lagrangian} cell-centred finite volume method for metal forming
  simulation.
\newblock \emph{International journal for numerical methods in engineering},
  109\penalty0 (13):\penalty0 1777--1803, 2016{\natexlab{b}}.
\newblock \doi{10.1002/nme.5345}.

\bibitem[Aguirre et~al.(2014)Aguirre, Gil, Bonet, and {n}o]{Aguirre2014}
M.~Aguirre, A.~J. Gil, J.~Bonet, and A.~A.~Carre\ {n}o.
\newblock A vertex centred finite volume {Jameson-Schmidt-Turkel (JST)}
  algorithm for a mixed conservation formulation in solid dynamics.
\newblock \emph{Journal of Computational Physics}, 259:\penalty0 672 -- 699,
  2014.
\newblock ISSN 0021-9991.
\newblock \doi{10.1016/j.jcp.2013.12.012}.
\newblock URL
  \url{http://www.sciencedirect.com/science/article/pii/S0021999113008115}.

\bibitem[LeVeque(2004)]{LeVeque2004}
R.~LeVeque.
\newblock \emph{Finite Volume Methods for Hyperbolic Problems}.
\newblock Cambridge University Press, 2004.

\bibitem[Baliga and Atabaki(2009)]{Baliga2009}
B.~R. Baliga and N.~Atabaki.
\newblock \emph{Control-Volume-Based Finite-Difference and Finite-Element
  Methods}, chapter~6, pages 191--224.
\newblock Wiley-Blackwell, 2009.
\newblock ISBN 9780470172599.
\newblock \doi{10.1002/9780470172599.ch6}.
\newblock URL
  \url{https://onlinelibrary.wiley.com/doi/abs/10.1002/9780470172599.ch6}.

\bibitem[Courant et~al.(1928)Courant, Friedrichs, and Lewy]{Courant1928}
R.~Courant, K.~Friedrichs, and H.~Lewy.
\newblock On the partial difference equations of mathematical physics.
\newblock \emph{Mathematische Annalem}, 100:\penalty0 32--74, 1928.

\bibitem[Demird{\v{z}}i\'{c}(2015)]{Demirdzic2015}
I.~Demird{\v{z}}i\'{c}.
\newblock On the discretization of the diffusion term in finite-volume
  continuum mechanics.
\newblock \emph{Numerical Heat Transfer, Part B: Fundamentals: An International
  Journal of Computation and Methodology}, 68:1:\penalty0 1--10, 2015.
\newblock \doi{10.1080/10407790.2014.985992}.

\bibitem[Maneeratana(2000)]{Maneeratana2000b}
K.~Maneeratana.
\newblock \emph{Development of the finite volume method for non-linear
  structural applications}.
\newblock PhD thesis, Imperial College London, 2000.

\bibitem[{Vaz Jr.} et~al.(2009){Vaz Jr.}, {n}oz Rojas, and Filippini]{Vaz2009}
M.~{Vaz Jr.}, P.~A.~Mu\ {n}oz Rojas, and G.~Filippini.
\newblock On the accuracy of nodal stress computation in plane elasticity using
  finite volumes and finite elements.
\newblock \emph{Computers \& Structures}, 87:\penalty0 1044--1057, 2009.

\bibitem[Cavalcante et~al.(2012)Cavalcante, Pindera, and
  Khatam]{Cavalcante2012c}
M.~A.~A. Cavalcante, M.-J. Pindera, and H.~Khatam.
\newblock Finite-volume micromechanics of periodic materials: Past, present and
  future.
\newblock \emph{Composites Part B: Engineering}, 43:\penalty0 2521--2543, 2012.

\bibitem[Van~Eck and Waltman(2007)]{VanEck2007}
N.~J. Van~Eck and L.~Waltman.
\newblock {VOS}: a new method for visualizing similarities between objects.
\newblock In \emph{Advances in Data Analysis: Proceedings of the 30$^{th}$
  Annual Conference of the German Classification Society}, pages 299--306.
  Springer, 2007.

\bibitem[Trangenstein and Colella(1991)]{Trangenstein1991}
J.~A. Trangenstein and P.~Colella.
\newblock A higher-order {Godunov} method for modeling finite deformation in
  elastic-plastic solids.
\newblock \emph{Communications on Pure and Applied Mathematics}, 44:\penalty0
  41--100, 1991.

\bibitem[Demird{\v{z}}i\'{c} and Muzaferija(1994)]{Demirdzic1994a}
I.~Demird{\v{z}}i\'{c} and S.~Muzaferija.
\newblock Finite volume method for stress analysis in complex domains.
\newblock \emph{International Journal for Numerical Methods in Engineering},
  37:\penalty0 3751--3766, 1994.

\bibitem[Dioh et~al.(1994)Dioh, Ivankovi\'{c}, Leevers, and Williams]{Dioh1994}
N.~Dioh, A.~Ivankovi\'{c}, P.~Leevers, and J.~G. Williams.
\newblock The high strain rate behaviour of polymers.
\newblock \emph{Journal de Physique IV}, 04 (C8), 1994.
\newblock \doi{10.1051/jp4:1994818}.
\newblock C8-119-C8-124.

\bibitem[Rente and Oliveira(2000)]{Rente2000}
C.~J. Rente and P.~J. Oliveira.
\newblock Extension of a finite volume method in solid stress analysis to cater
  for non-linear elasto-plastic effects.
\newblock In J.L. Tassoulas, editor, \emph{{Proceedings of EM 14$^{th}$
  Engineering Mechanics Conference}}, University of Texas at Austin, USA, 2000.

\bibitem[Teskered{\v{z}}i\'{c} et~al.(2002)Teskered{\v{z}}i\'{c},
  Demird{\v{z}}i\'{c}, and Muzaferija]{Teskeredzic2002}
A.~Teskered{\v{z}}i\'{c}, I.~Demird{\v{z}}i\'{c}, and S.~Muzaferija.
\newblock Numerical method for heat transfer, fluid flow, and stress analysis
  in phase-change problems.
\newblock \emph{Numerical Heat Transfer, Part B: Fundamentals}, 42:\penalty0
  437--459, 2002.

\bibitem[Teskered{\v{z}}i\'{c}(2004)]{Teskeredzic2004}
A.~Teskered{\v{z}}i\'{c}.
\newblock \emph{Application of the Finite Volume Method to Casting Problems}.
\newblock PhD thesis, University of Sarajevo, 2004.

\bibitem[Ba\v{s}i\'{c} et~al.(2005)Ba\v{s}i\'{c}, Demird{\v{z}}i\'{c}, and
  Muzaferija]{Basic2005}
H.~Ba\v{s}i\'{c}, I.~Demird{\v{z}}i\'{c}, and S.~Muzaferija.
\newblock Finite volume method for simulation of extrusion processes.
\newblock \emph{Internation Journal for Numerical Methods in Engineering},
  62:\penalty0 475--494, 2005.

\bibitem[Martins et~al.(2010)Martins, Bressan, Button, and
  Ivankovi\'{c}]{Martins2010}
M.~M. Martins, J.~D. Bressan, S.~T. Button, and A.~Ivankovi\'{c}.
\newblock Extrusion process by finite volume method using {OpenFOAM} software.
\newblock In Y.~Chastel Chinesta and M.~El Mansori, editors,
  \emph{International Conference on Advances in Material and Processing
  Technologies - {AMPT2010}}, pages 1461--1466. American Institute of Physics,
  Paris, France, 2010.

\bibitem[Bressan et~al.(2010)Bressan, Martins, and {Vaz Jr.}]{Bressan2010}
J.~D. Bressan, M.~M. Martins, and M.~{Vaz Jr.}
\newblock Stress evolution and thermal shock computation using the finite
  volume method.
\newblock \emph{Journal of Thermal Stresses}, 33:\penalty0 533--558, 2010.

\bibitem[Leonard et~al.(2012)Leonard, Murphy, Kara\v{c}, and
  Ivankovi\'{c}]{Leonard2012}
M.~Leonard, N.~Murphy, A.~Kara\v{c}, and A.~Ivankovi\'{c}.
\newblock A numerical investigation of spherical void growth in an
  elastic-plastic continuum.
\newblock \emph{Computational Materials Science}, 64:\penalty0 38--40, 2012.

\bibitem[Bressan et~al.(2015)Bressan, Martins, and Button]{Bressan2015}
J.~D. Bressan, M.~M. Martins, and S.~T. Button.
\newblock Analysis of aluminium hot extrusion by finite volume method.
\newblock \emph{Materials Today: Proceedings}, 2\penalty0 (10, Part
  {A}):\penalty0 4740--4747, 2015.

\bibitem[Martins et~al.(2016)Martins, Bressan, and Button]{Martins2016}
M.~M. Martins, J.~D. Bressan, and S.~T. Button.
\newblock Analysis of aluminum extrusion in a $90^o$ die by finite volume
  method.
\newblock \emph{Advanced Materials Research}, 1135:\penalty0 153--160, 2016.

\bibitem[Martins et~al.(2017)Martins, Bressan, and Button]{Martins2017}
M.~M. Martins, J.~D. Bressan, and S.~T. Button.
\newblock Finite volume analysis with the maccormack method applied to metal
  flow in forward extrusion.
\newblock \emph{Universal Journal of Mechanical Engineering}, 5:\penalty0 1--8,
  2017.

\bibitem[Cardiff et~al.(2017{\natexlab{a}})Cardiff, Tang, Tukovic, Jasak,
  Ivankovic, and Jaeger]{Cardiff2017a}
P.~Cardiff, T.~Tang, Z.~Tukovic, H.~Jasak, A.~Ivankovic, and P.~De Jaeger.
\newblock {An Eulerian-inspired Lagrangian finite volume method for wire
  drawing simulations}.
\newblock In \emph{{IUTAM} Symposium on Multi-scale Fatigue, Fracture and
  Damage of Materials in Harsh Environments}, Galway, Ireland,
  2017{\natexlab{a}}. National University of Ireland Galway.

\bibitem[Bressan et~al.(2017)Bressan, Martins, and Button]{Bressan2017}
J.~D. Bressan, M.~M. Martins, and S.~T. Button.
\newblock Analysis of metal extrusion by the finite volume method.
\newblock \emph{Procedia Engineering}, 207:\penalty0 425--430, 2017.

\bibitem[Zarrabi and Basu(1999)]{Zarrabi1999}
K.~Zarrabi and A.~Basu.
\newblock An axisymmetric finite volume formulation for creep analysis.
\newblock \emph{Journal of the Mechanical Behavior of Materials}, 10:\penalty0
  325--340, 1999.

\bibitem[Zarrabi and Basu(2000)]{Zarrabi2000}
K.~Zarrabi and A.~Basu.
\newblock A finite volume element formulation for solution of elastic
  axisymmetric pressurized components.
\newblock \emph{International Journal of Pressure Vessels and Piping},
  77:\penalty0 479--484, 2000.

\bibitem[D{\v{z}}aferovi\'{c} et~al.(2000)D{\v{z}}aferovi\'{c}, Ivankovi\'{c},
  and Demird{\v{z}}i\'{c}]{Dzaferovic2000}
E.~D{\v{z}}aferovi\'{c}, A.~Ivankovi\'{c}, and I.~Demird{\v{z}}i\'{c}.
\newblock Finite volume modelling of linear viscoelastic deformation.
\newblock In \emph{Proceedings of 3$^{rd}$ Congress of Croatian Society of
  Mechanics}, Dubrovnik, Croatia, 2000.

\bibitem[Demird{\v{z}}i\'{c} et~al.(2005)Demird{\v{z}}i\'{c},
  D{\v{z}}aferovi\'{c}, and Ivankovi\'{c}]{Demirdzic2005}
I.~Demird{\v{z}}i\'{c}, E.~D{\v{z}}aferovi\'{c}, and A.~Ivankovi\'{c}.
\newblock Finite-volume approach to thermoviscoelasticity.
\newblock \emph{Numerical Heat Transfer, Part B: Fundamentals}, 47:\penalty0
  213--237, 2005.

\bibitem[Das et~al.(2012)Das, Mathur, and Murthy]{Das2012}
S.~Das, S.~R. Mathur, and J.~Y. Murthy.
\newblock {Finite-volume method for creep analysis of thin RF MEMS devices
  using the theory of plates}.
\newblock \emph{Numerical Heat Transfer, Part B: Fundamentals}, 61:\penalty0
  71--90, 2012.

\bibitem[Safari et~al.(2016)Safari, Tukovi\'{c}, Cardiff, Walter, Casey, and
  Ivankovi\'{c}]{Safari2016}
A.~Safari, {\v{Z}}.~Tukovi\'{c}, P.~Cardiff, M.~Walter, E.~Casey, and
  A.~Ivankovi\'{c}.
\newblock Interfacial separation of a mature biofilm from a glass surface-a
  combined experimental and cohesive zone modelling approach.
\newblock \emph{Journal of the mechanical behavior of biomedical materials},
  54:\penalty0 205--218, 2016.

\bibitem[Demird{\v{z}}i\'{c} et~al.(1994)Demird{\v{z}}i\'{c}, Ivankovi\'{c},
  Martinovi\'{c}, and Muzaferija]{Demirdzic1994b}
I.~Demird{\v{z}}i\'{c}, A.~Ivankovi\'{c}, D.~Martinovi\'{c}, and S.~Muzaferija.
\newblock Numerical method for solving linear and non-linear solid body
  problems.
\newblock In \emph{Proceedings of 1$^{st}$ Congress of Croatian Society of
  Mechanics}, Pula, Croatia, 1994.

\bibitem[Demird{\v{z}}i\'{c}(1996)]{Demirdzic1996a}
I.~Demird{\v{z}}i\'{c}.
\newblock Finite volumes in solid mechanics.
\newblock In \emph{{NAFEMS ASME} Seminar - Alternative Strategies in
  Computational Mechanics, London}, February 1996.

\bibitem[Osman et~al.(2011)Osman, Ahmad, and Arshad]{Osman2011}
H.~Osman, S.~Ahmad, and K.~A. Arshad.
\newblock A one-dimensional simulation of an electrofusion welding process.
\newblock \emph{International Conference on Modeling, Simulation and Applied
  Optimization ({ICMSAO})}, 4:\penalty0 1--5, 2011.

\bibitem[Junior et~al.(2013)Junior, Xavier, de~Castro, Campos, Palmeira, and
  Habibe]{Junior2013}
H.~G.~D. Junior, C.~R. Xavier, J.~A. de~Castro, M.~F. Campos, A.~A. Palmeira,
  and A.~F. Habibe.
\newblock Mathematical modeling and experimental investigation of the stress
  evolution at the steel welding.
\newblock \emph{{Cadernos UniFOA}}, 8:\penalty0 59--66, 2013.

\bibitem[Bibin and Ramarajan(2013)]{Bibin2013}
K.~S. Bibin and A.~Ramarajan.
\newblock Unstructured finite volume approach for {3-D} unsteady
  thermo-structural analysis using bi-conjugate gradient stabilized method.
\newblock \emph{International Journal of Innovative Research in Science
  Engineering and Technology}, 2:\penalty0 1389--1400, 2013.

\bibitem[Martinovi{\'c} and Horman(1998)]{Martinovic1998}
D.~Martinovi{\'c} and I.~Horman.
\newblock Numerical simulation of clay bricks drying process.
\newblock In \emph{{II} Meunarodni nauno-struni skup Proizvodnja i prerada
  nemetalnih mineralnih sirovina i njihova primjena u industriji}, Zenica,
  Bosnia and Herzegovina, 1998.
\newblock In Bosnian.

\bibitem[Demird{\v{z}}i\'{c} et~al.(2000)Demird{\v{z}}i\'{c}, Horman, and
  Martinovi\'{c}]{Demirdzic2000}
I.~Demird{\v{z}}i\'{c}, I.~Horman, and D.~Martinovi\'{c}.
\newblock Finite volume analysis of stress and deformation in
  hygro-thermo-elastic orthotropic body.
\newblock \emph{Computer Methods in Applied Mechanics and Engineering},
  190:\penalty0 1221--1232, 2000.

\bibitem[Martinovi\'{c} and Horman(2000)]{Martinovic2000}
D.~Martinovi\'{c} and I.~Horman.
\newblock Drying induced stresses in clay bricks an hygro-thermo-elastic model.
\newblock In \emph{International scientific and expert symposium Nonmetal
  Inorganic Materials}, pages 249--257, Zenica, Bosnia-Herzegovina, 2000.

\bibitem[Martinovi\'{c} et~al.(2001)Martinovi\'{c}, Horman, and
  Demird{\v{z}}i\'{c}]{Martinovic2001}
D.~Martinovi\'{c}, I.~Horman, and I.~Demird{\v{z}}i\'{c}.
\newblock Numerical and experimental analysis of wood drying process.
\newblock \emph{Wood Science and Technology}, 35:\penalty0 143--156, 2001.

\bibitem[Martinovi\'{c}(2002)]{Martinovic2002}
D.~Martinovi\'{c}.
\newblock \emph{A numerical method for analysis of thermo-deformational
  processes during the welding}.
\newblock PhD thesis, University of Sarajevo, 2002.
\newblock In Bosnian.

\bibitem[Horman et~al.(2008)Horman, Martinovi\'{c}, and
  Hajdarevi\'{c}]{Horman2008}
I.~Horman, D.~Martinovi\'{c}, and S.~Hajdarevi\'{c}.
\newblock Numerical analysis of a phenomena in the wood caused by heat.
\newblock In \emph{Moisture or External Load, in: Proceedings of International
  scientific conference: Challanges in forestry and wood technology in the},
  volume~21, pages 31--34, University of Zagreb, Faculty of Forestry, Zageb,
  2008.

\bibitem[Martinovi\'{c} et~al.(2008)Martinovi\'{c}, Horman, and
  Hajdarevi\'{c}]{Martinovic2008}
D.~Martinovi\'{c}, I.~Horman, and S.~Hajdarevi\'{c}.
\newblock Stress distribution in wooden corner joints.
\newblock \emph{Strojarstvo}, 50:\penalty0 193--204, 2008.

\bibitem[Horman et~al.(2009)Horman, Martinovi\'{c}, and
  Hajdarevi\'{c}]{Horman2009}
I.~Horman, D.~Martinovi\'{c}, and S.~Hajdarevi\'{c}.
\newblock Finite volume method for analysis of stress and strain in wood.
\newblock \emph{Drvna industrija}, 60:\penalty0 27--32, 2009.

\bibitem[Horman et~al.(2010)Horman, Hajdarevi\'{c}, Martinovi\'{c}, and
  Vukas]{Horman2010a}
I.~Horman, S.~Hajdarevi\'{c}, D.~Martinovi\'{c}, and N.~Vukas.
\newblock Numerical analysis of stress and strain in a wooden chair.
\newblock \emph{Drvna industrija}, 61:\penalty0 151--158, 2010.

\bibitem[Horman et~al.(2012)Horman, Martinovi\'{c}, Bijelonja, and
  Hajdarevi\'{c}]{Horman2012}
I.~Horman, D.~Martinovi\'{c}, I.~Bijelonja, and S.~Hajdarevi\'{c}.
\newblock Wood subjected to hygro-thermal and/or mechanical loads.
\newblock In \emph{Finite Volume Method - Powerful Means of Engineering
  Design}. Finite Volume Method - Powerful Means of Engineering Design, R.
  Petrova (ed) 5, 2012.

\bibitem[Fu(2018)]{Fu2018}
R.~Fu.
\newblock \emph{Thermo-Mechanical Coupling for Ablation}.
\newblock PhD thesis, University of Kentucky, 2018.

\bibitem[Holzmann et~al.(2018)Holzmann, Ludwig, and Raninger]{Holzmann2018}
T.~Holzmann, A.~Ludwig, and P.~Raninger.
\newblock Yield strength prediction in {3D} during local heat treatment of
  structural {A356} alloy components in combination with thermal-stress
  analysis.
\newblock \emph{Lambotte G., Lee J., Allanore A., Wagstaff S. (eds) Materials
  Processing Fundamentals. The Minerals, Metals \& Materials Series}, 2018.

\bibitem[Tang et~al.(2015)Tang, Hededal, and Cardiff]{Tang2015}
T.~Tang, O.~Hededal, and P.~Cardiff.
\newblock On finite volume method implementation of poro-elasto-plasticity soil
  model.
\newblock \emph{International journal for numerical and analytical methods in
  geomechanics}, 39:\penalty0 1410--1430, 2015.
\newblock \doi{10.1002/nag.2361}.

\bibitem[Bryant et~al.(2015)Bryant, Hwang, and Sharma]{Bryant2015}
E.~C. Bryant, J.~Hwang, and M.~M. Sharma.
\newblock Arbitrary fracture propagation in heterogeneous poroelastic
  formations using a finite volume-based cohesive zone model.
\newblock In \emph{SPE Hydraulic Fracturing Technology Conference}, The
  Woodlands, Texas, February, 2015.

\bibitem[Cardiff et~al.(2015{\natexlab{a}})Cardiff, Manchanda, Bryant, Lee,
  Ivankovi\'{c}, and Sharma]{Cardiff2015a}
P.~Cardiff, R.~Manchanda, E.~C. Bryant, D.~Lee, A.~Ivankovi\'{c}, and M.~M.
  Sharma.
\newblock Simulation of fractures in {OpenFOAM}: From adhesive joints to
  hydraulic fractures.
\newblock In \emph{10$^{th}$ {OpenFOAM} Workshop}, University of Michigan, Ann
  Arbor, MI, USA, 2015{\natexlab{a}}.

\bibitem[Cardiff et~al.(2015{\natexlab{b}})Cardiff, Manchanda, Bryant,
  Ivankovi\'{c}, and Sharma]{Cardiff2015b}
P.~Cardiff, R.~Manchanda, E.~C. Bryant, A.~Ivankovi\'{c}, and M.~M. Sharma.
\newblock Finite volume method for the simulation of hydraulic fractures.
\newblock In \emph{Joint Symposium of Irish Mechanics Society \& Irish Society
  for Scientific \& Engineering Computation Advances in Mechanics}, University
  College Dublin, Dublin, Ireland, 2015{\natexlab{b}}.

\bibitem[Lee et~al.(2015)Lee, Cardiff, Bryant, Manchanda, Wang, and
  M.~M.~Sharma]{Lee2015}
D.~Lee, P.~Cardiff, E.~C. Bryant, R.~Manchanda, H.~Wang, and A~M.~M.~Sharma.
\newblock New model for hydraulic fracture growth in unconsolidated sands with
  plasticity and leak-off.
\newblock \emph{{SPE} Annual Technical Conference and Exhibition}, pages
  28--30, September 2015.
\newblock 18174818-MS.

\bibitem[Elsafti and Oumeraci(2016)]{Elsafti2016}
H.~Elsafti and H.~Oumeraci.
\newblock A numerical hydro-geotechnical model for marine gravity structures.
\newblock \emph{Computers and Geotechnics}, 79:\penalty0 105 -- 129, 2016.
\newblock \doi{10.1016/j.compgeo.2016.05.025}.

\bibitem[Manchanda et~al.(2017)Manchanda, Bryant, Bhardwaj, Cardiff, and
  Sharma]{Manchanda2017}
R.~Manchanda, E.~C. Bryant, P.~Bhardwaj, P.~Cardiff, and M.~M. Sharma.
\newblock Strategies for effective stimulation of multiple perforation clusters
  in horizontal wells.
\newblock \emph{SPE}, preprint:\penalty0 28--30, December 2017.
\newblock \doi{10.2118/179126-PA}.

\bibitem[Asadollahi(2017)]{Asadollahi2017}
M.~Asadollahi.
\newblock Finite volume method for poroelasticity.
\newblock Master's thesis, Delft University of Technology, 2017.

\bibitem[Fainberg and Leister(1996)]{Fainberg1996}
J.~Fainberg and H.~J. Leister.
\newblock Finite volume multigrid solver for thermo-elastic stress analysis in
  anisotropic materials.
\newblock \emph{Computer Methods in Applied Mechanics and Engineering},
  137:\penalty0 167--174, 1996.

\bibitem[Horman(1999)]{Horman1999}
I.~Horman.
\newblock \emph{Finite volume method for analysis of timber drying}.
\newblock PhD thesis, University of Sarajevo, 1999.
\newblock In Bosnian.

\bibitem[Cardiff et~al.(2014{\natexlab{a}})Cardiff, Kara\v{c}, and
  Ivankovi\'{c}]{Cardiff2014a}
P.~Cardiff, A.~Kara\v{c}, and A.~Ivankovi\'{c}.
\newblock A large strain finite volume method for orthotropic bodies with
  general material orientations.
\newblock \emph{Computer Methods in Applied Mechanics and Engineering},
  268:\penalty0 318--335, 2014{\natexlab{a}}.
\newblock \doi{10.1016/j.cma.2013.09.008}.

\bibitem[Golubovi\'{c} et~al.(2017)Golubovi\'{c}, Demird{\v{z}}i\'{c}, and
  Muzaferija]{Golubovic2017a}
A.~Golubovi\'{c}, I.~Demird{\v{z}}i\'{c}, and S.~Muzaferija.
\newblock Finite volume analysis of laminated composite plates.
\newblock \emph{International Journal for Numerical Methods in Engineering},
  109\penalty0 (11):\penalty0 1607--1620, 2017.
\newblock ISSN 1097-0207.
\newblock \doi{10.1002/nme.5347}.
\newblock URL \url{10.1002/nme.5347}.

\bibitem[Tukovi\'{c} et~al.(2012)Tukovi\'{c}, Ivankovi\'{c}, and
  Kara\v{c}]{Tukovic2012}
{\v{Z}}.~Tukovi\'{c}, A.~Ivankovi\'{c}, and A.~Kara\v{c}.
\newblock Finite-volume stress analysis in multi-material linear elastic body.
\newblock \emph{International Journal for Numerical Methods in Engineering},
  93:\penalty0 400--419, 2012.

\bibitem[Carolan et~al.(2013{\natexlab{a}})Carolan, Tukovi\'{c}, Murphy, and
  Ivankovi\'{c}]{Carolan2013a}
D.~Carolan, {\v{Z}}.~Tukovi\'{c}, N.~Murphy, and A.~Ivankovi\'{c}.
\newblock Arbitrary crack propagation in multi-phase materials using the finite
  volume method.
\newblock \emph{Computational Materials Science}, 69:\penalty0 153--159,
  2013{\natexlab{a}}.

\bibitem[Cardiff(2012)]{Cardiff2012d}
P.~Cardiff.
\newblock \emph{{Development of the Finite Volume Method for Hip Joint Stress
  Analysis}}.
\newblock PhD thesis, University College Dublin, 2012.
\newblock URL
  \url{https://www.researchgate.net/publication/262772501_Development_of_the_Finite_Volume_Method_for_Hip_Joint_Stress_Analysis}.

\bibitem[Greenshields et~al.(1999{\natexlab{a}})Greenshields, Weller, and
  Ivankovi\'{c}]{Greenshields1999a}
C.~J. Greenshields, H.~G. Weller, and A.~Ivankovi\'{c}.
\newblock The finite volume formulation for fluid structure interaction.
\newblock In R.~Vilsmeirer Haenel and F.~Benkhaldoun, editors, \emph{Finite
  Volumes for Complex Applications, II - Problems and Perspectives,}, pages
  467--474. Hermes Science, 1999{\natexlab{a}}.

\bibitem[Greenshields et~al.(1999{\natexlab{b}})Greenshields, Weller, and
  Ivankovi\'{c}]{Greenshields1999b}
C.~J. Greenshields, H.~G. Weller, and A.~Ivankovi\'{c}.
\newblock The finite volume method for coupled fluid flow and stress analysis.
\newblock \emph{Computer Modeling and Simulation in Engineering}, 4:\penalty0
  213--218, 1999{\natexlab{b}}.

\bibitem[Bijelonja et~al.(2005)Bijelonja, Demird{\v{z}}i\'{c}, and
  Muzaferija]{Bijelonja2005a}
I.~Bijelonja, I.~Demird{\v{z}}i\'{c}, and S.~Muzaferija.
\newblock A finite volume method for large strain analysis of incompressible
  hyperelastic materials.
\newblock \emph{International Journal for Numerical Methods in Engineering},
  64:\penalty0 1594--1609, 2005.

\bibitem[Greenshields and Weller(2005)]{Greenshields2005}
C.~J. Greenshields and H.~G. Weller.
\newblock A unified formulation for continuum mechanics applied to
  fluid-structure interaction in flexible tubes.
\newblock \emph{International Journal for Numerical Methods in Engineering},
  64:\penalty0 1575--1593, 2005.

\bibitem[Bijelonja et~al.(2006)Bijelonja, Demird{\v{z}}i\'{c}, and
  Muzaferija]{Bijelonja2006}
I.~Bijelonja, I.~Demird{\v{z}}i\'{c}, and S.~Muzaferija.
\newblock A finite volume method for incompressible linear elasticity.
\newblock \emph{Computer Methods in Applied Mechanics and Engineering},
  195:\penalty0 6378--6390, 2006.

\bibitem[Giannopapa and Papadakis(2006)]{Giannopapa2006}
C.~G. Giannopapa and G.~Papadakis.
\newblock New formulations of the dynamic equations of elastic solids suitable
  for a unified methodology for fluid-structure interaction problems.
\newblock \emph{Computer Methods in Applied Mechanics and Engineering}, 2006.

\bibitem[Giannopapa and Papadakis(2008)]{Giannopapa2008}
C.~G. Giannopapa and G.~Papadakis.
\newblock Linear stability analysis and application of a new solution method of
  the elastodynamic equations suitable for a unified
  fluid-structure-interaction approach.
\newblock \emph{{ASME} Journal Pressure Vessels Technology}, 130:\penalty0
  31303--1, 2008.

\bibitem[Bijelonja et~al.(2017)Bijelonja, Demird{\v{z}}i\'{c}, and
  Muzaferija]{Bijelonja2017}
I.~Bijelonja, I.~Demird{\v{z}}i\'{c}, and S.~Muzaferija.
\newblock Mixed finite volume method for linear thermoelasticity at all
  {Poisson}'s ratios.
\newblock \emph{Numerical Heat Transfer, Part A: Applications}, 72:\penalty0
  215--235, 2017.

\bibitem[Jasak and Weller(2000{\natexlab{a}})]{Jasak2000b}
H.~Jasak and H.~G. Weller.
\newblock Finite volume methodology for contact problems of linear elastic
  solids.
\newblock In \emph{Proceedings of 3$^{rd}$ Congress of Croatian Society of
  Mechanics}, pages 253--260, Dubrovnik, Croatia, 2000{\natexlab{a}}.

\bibitem[Cardiff et~al.(2011{\natexlab{a}})Cardiff, Kara\v{c}, Flavin,
  FitzPatrick, and Ivankovi\'{c}]{Cardiff2011c}
P.~Cardiff, A.~Kara\v{c}, R.~Flavin, D.~FitzPatrick, and A.~Ivankovi\'{c}.
\newblock Contact stress analysis in {OpenFOAM} - application to hip joint
  bones.
\newblock In \emph{{OpenFOAM} Workshop, Penn State University}, Penn State, PA,
  USA, 2011{\natexlab{a}}.

\bibitem[Cardiff et~al.(2011{\natexlab{b}})Cardiff, Kara\v{c}, and
  Ivankovi\'{c}]{Cardiff2011d}
P.~Cardiff, A.~Kara\v{c}, and A.~Ivankovi\'{c}.
\newblock Development of a finite volume methodology for linear elastic contact
  problems.
\newblock In \emph{21$^{st}$ International Workshop on Computational Mechanics
  of Materials, IWCMM, Limerick}, Limerick, Ireland, 2011{\natexlab{b}}.

\bibitem[Cardiff et~al.(2012{\natexlab{a}})Cardiff, Kara\v{c}, and
  Ivankovi\'{c}]{Cardiff2012a}
P.~Cardiff, A.~Kara\v{c}, and A.~Ivankovi\'{c}.
\newblock Development of a finite volume contact solver based on the penalty
  method.
\newblock \emph{Computational Materials Science}, 64:\penalty0 283 -- 284,
  2012{\natexlab{a}}.

\bibitem[Cardiff et~al.(2014{\natexlab{b}})Cardiff, Kara\v{c}, FitzPatrick,
  Flavin, and Ivankovi\'{c}]{Cardiff2014b}
P.~Cardiff, A.~Kara\v{c}, D.~FitzPatrick, R.~Flavin, and A.~Ivankovi\'{c}.
\newblock Development of a hip joint model for finite volume simulations.
\newblock \emph{Journal of Biomechanical Engineering}, 136:\penalty0 1--8,
  2014{\natexlab{b}}.
\newblock \doi{10.1115/1.4025776}.

\bibitem[Maneeratana and Ivankovi\'{c}(1999{\natexlab{a}})]{Maneeratana1999a}
K.~Maneeratana and A.~Ivankovi\'{c}.
\newblock Finite volume method for structural applications involving material
  and geometrical non-linearities.
\newblock In \emph{European Council of Computational Mechanics Technische
  Universit{\"{u}}nchen, Proceedings of the European Conference on
  Computational Mechanics ({ECCM}'99)}, pages 874--875, 1999{\natexlab{a}}.

\bibitem[Maneeratana and Ivankovi\'{c}(1999{\natexlab{b}})]{Maneeratana1999b}
K.~Maneeratana and A.~Ivankovi\'{c}.
\newblock Finite volume method for geometrically non-linear stress analysis
  applications.
\newblock In \emph{Proceedings of the Seventh Annual Conference of the
  Association for Computational Mechanics in Engineering ({ACME'1999})},
  1999{\natexlab{b}}.

\bibitem[Maneeratana and Ivankovi\'{c}(1999{\natexlab{c}})]{Maneeratana1999c}
K.~Maneeratana and A.~Ivankovi\'{c}.
\newblock Finite volume method for large deformation with linear hypoelastic
  materials.
\newblock In R.~Vilsmeier and F.~Benkhaldoun, editors, \emph{{Second
  International Symposium on Finite Volumes for Complex Applications (FVCA II)
  for Complex Applications II: Problems and Perspectives}}, pages 459--466,
  University Duisburg, Germany, 1999{\natexlab{c}}. Hermes Science Publication.

\bibitem[Maneeratana and Ivankovi\'{c}(2000)]{Maneeratana2000a}
K.~Maneeratana and A.~Ivankovi\'{c}.
\newblock Modelling of high strain rate behaviour of a series 7108 aluminium
  alloy.
\newblock In \emph{Proceedings of the 14$^{th}$ Conference of the Mechanical
  Engineering Network of Thailand}, pages 227--233, 2000.

\bibitem[Ba\v{s}i\'{c}(2002)]{Basic2002}
H.~Ba\v{s}i\'{c}.
\newblock \emph{Application of the finite volume method to the analysis of
  plastic metal flow in extrusion technologies}.
\newblock PhD thesis, University of Sarajevo, 2002.
\newblock In Bosnian.

\bibitem[Bijelonja(2002)]{Bijelonja2002}
I.~Bijelonja.
\newblock \emph{Finite volume method for incremental analysis of small and
  large thermo-elasto-plastic deformations}.
\newblock PhD thesis, University of Sarajevo, 2002.
\newblock In Bosnian.

\bibitem[Tukovi\'{c} and Jasak(2007{\natexlab{a}})]{Tukovic2007a}
{\v{Z}}.~Tukovi\'{c} and H.~Jasak.
\newblock Updated {Lagrangian} finite volume solver for large deformation
  dynamic response of elastic body.
\newblock \emph{Transactions of {FAMENA}}, 31:\penalty0 55--70,
  2007{\natexlab{a}}.

\bibitem[Tukovi\'{c} and Jasak(2007{\natexlab{b}})]{Tukovic2007b}
{\v{Z}}.~Tukovi\'{c} and H.~Jasak.
\newblock {FVM} for fluid-structure interaction with large structural
  displacements.
\newblock In \emph{2$^{nd}$ {OpenFOAM} Workshop}, Zagreb, Croatia,
  2007{\natexlab{b}}.

\bibitem[Cardiff et~al.(2012{\natexlab{b}})Cardiff, Kara\v{c}, Tukovi\'{c}, and
  Ivankovi\'{c}]{Cardiff2012c}
P.~Cardiff, A.~Kara\v{c}, {\v{Z}}.~Tukovi\'{c}, and A.~Ivankovi\'{c}.
\newblock Development of a finite volume based structural solver for large
  rotation of non-orthogonal meshes.
\newblock In \emph{7$^{th}$ {OpenFOAM} Workshop}, Darmstadt, Germany,
  2012{\natexlab{b}}.

\bibitem[Cardiff et~al.(2013)Cardiff, Kara\v{c}, Tukovi\'{c}, and
  Ivankovi\'{c}]{Cardiff2013}
P.~Cardiff, A.~Kara\v{c}, Z.~Tukovi\'{c}, and A.~Ivankovi\'{c}.
\newblock An open-source finite method for computational solid mechanics.
\newblock In \emph{Joint Symposium of Irish Mechanics Society and Irish Society
  for Scientific and Engineering Computation}, University College Dublin,
  Dublin, Ireland, 2013.

\bibitem[Cardiff et~al.(2014{\natexlab{c}})Cardiff, Tukovi\'{c}, Kara\v{c}, and
  Ivankovi\'{c}]{Cardiff2014d}
P.~Cardiff, {\v{Z}.}~Tukovi\'{c}, A.~Kara\v{c}, and A.~Ivankovi\'{c}.
\newblock Nonlinear solid mechanics in {OpenFOAM}.
\newblock In \emph{9$^{th}$ {OpenFOAM} Workshop}, University of Zagreb,
  Croatia, 2014{\natexlab{c}}.

\bibitem[Liu et~al.(2018)Liu, Ming, and Zhang]{Liu2018}
Q.~Liu, P.-J. Ming, and W.-P. Zhang.
\newblock Research on the nonlinear finite volume numerical method for the
  large rotating of disk.
\newblock \emph{Journal of Harbin Engineering University}, 39:\penalty0
  1012--1018, 2018.

\bibitem[Ivankovi\'{c} et~al.(1993)Ivankovi\'{c}, Demird{\v{z}}i\'{c},
  Williams, and Leevers]{Ivankovic1993}
A.~Ivankovi\'{c}, I.~Demird{\v{z}}i\'{c}, J.~G. Williams, and P.~S. Leevers.
\newblock A new numerical method for analysing dynamic fracture problems.
\newblock In \emph{{ESIS} Symposium on Impact and Dynamic Fracture of Polymers
  and Composites}, Potro Cervo, Sardinia, Italy, 1993.

\bibitem[Ivankovi\'{c} et~al.(1994)Ivankovi\'{c}, Demird{\v{z}}i\'{c},
  Williams, and Leevers]{Ivankovic1994}
A.~Ivankovi\'{c}, I.~Demird{\v{z}}i\'{c}, J.~G. Williams, and P.~S. Leevers.
\newblock Application of the finite volume method to the analysis of dynamic
  fracture problems.
\newblock \emph{International Journal of Fracture}, 66:\penalty0 357--371,
  1994.

\bibitem[Ivankovi\'{c} et~al.(1997{\natexlab{a}})Ivankovi\'{c}, Muzaferija, and
  Demird{\v{z}}i\'{c}]{Ivankovic1997a}
A.~Ivankovi\'{c}, S.~Muzaferija, and I.~Demird{\v{z}}i\'{c}.
\newblock Finite volume method and multigrid acceleration in modelling of rapid
  crack propagation in full-scale pipe test.
\newblock \emph{Computational Mechanics}, 20:\penalty0 46--52,
  1997{\natexlab{a}}.

\bibitem[Ivankovi\'{c} and Venizelos(1998)]{Ivankovic1998}
A.~Ivankovi\'{c} and G.~P. Venizelos.
\newblock Rapid crack propagation in plastic pipe: Predicting full-scale
  critical pressure from {S4} test results.
\newblock \emph{Engineering Fracture Mechanics}, 59:\penalty0 607--622, 1998.

\bibitem[Ivankovi\'{c}(1999)]{Ivankovic1999}
A.~Ivankovi\'{c}.
\newblock Finite volume modelling of dynamic fracture problems.
\newblock \emph{Computer Modelling and Simulation in Engineering}, 4:\penalty0
  227--235, 1999.

\bibitem[Ivankovi\'{c} and Hillmansen(2001)]{Ivankovic2001a}
A.~Ivankovi\'{c} and S.~Hillmansen.
\newblock Evoluton of dynamic fractures in {PMMA}.
\newblock \emph{Plastics, Rubber and Composites}, 30:\penalty0 88--93, 2001.

\bibitem[Stylianou and Ivankovi\'{c}(2002{\natexlab{a}})]{Stylianou2002a}
V.~Stylianou and A.~Ivankovi\'{c}.
\newblock Finite volume analysis of dynamic fracture phenomena {I}: A node
  release methodology.
\newblock \emph{International Journal of Fracture}, 113:\penalty0 107--123,
  2002{\natexlab{a}}.

\bibitem[Stylianou and Ivankovi\'{c}(2002{\natexlab{b}})]{Stylianou2002b}
V.~Stylianou and A.~Ivankovi\'{c}.
\newblock Finite volume analysis of dynamic fracture phenomena {II}: A cohesive
  zone type methodology.
\newblock \emph{International Journal of Fracture}, 113:\penalty0 125--151,
  2002{\natexlab{b}}.

\bibitem[Ivankovi\'{c} and K.~C.~Pandya(2004)]{Ivankovic2004a}
A.~Ivankovi\'{c} and J.~G.~Williams K.~C.~Pandya.
\newblock Crack growth predictions in polyethylene using measured
  traction-separation curves.
\newblock \emph{Engineering Fracture Mechanics}, pages 657--668, 2004.

\bibitem[Rager et~al.(2005)Rager, Williams, and Ivankovi\'{c}]{Rager2005}
A.~Rager, J.~G. Williams, and A.~Ivankovi\'{c}.
\newblock Numerical analysis of the three point bend impact test for polymers.
\newblock \emph{International Journal of Fracture}, 135\penalty0
  (1-4):\penalty0 199--215, 2005.

\bibitem[Murphy and Ivankovi\'{c}(2005)]{Murphy2005}
N.~Murphy and A.~Ivankovi\'{c}.
\newblock The prediction of dynamic fracture evolution in {PMMA} using a
  cohesive zone model.
\newblock \emph{Engineering Fracture Mechanics}, 72:\penalty0 861--875, 2005.

\bibitem[Murphy et~al.(2006)Murphy, Ali, and Ivankovi\'{c}]{Murphy2006}
N.~Murphy, M.~Ali, and A.~Ivankovi\'{c}.
\newblock Dynamic crack bifurcation in {PMMA}.
\newblock \emph{Engineering Fracture Mechanics}, 73\penalty0 (16):\penalty0
  2569--2587, 2006.

\bibitem[Trop\v{s}a et~al.(2006)Trop\v{s}a, Georgiou, Ivankovi\'{c}, Kinloch,
  and Williams]{Tropsa2006}
V.~Trop\v{s}a, I.~Georgiou, A.~Ivankovi\'{c}, A.~J. Kinloch, and J.~G.
  Williams.
\newblock {OpenFOAM} in non-linear stress analysis: modelling of adhesive
  joints.
\newblock In \emph{1$^{st}$ {OpenFOAM} Workshop, Zagreb, Croatia}, 2006.

\bibitem[Tukovi\'{c}(2010)]{Tukovic2010}
{\v{Z}}.~Tukovi\'{c}.
\newblock {Arbitrary crack propagation model in {OpenFOAM}, }.
\newblock Technical report, Faculty of Mechanical Engineering and Naval
  Architecture, University of Zagreb, in association with the School of
  Mechanical and Materials Engineering, University College Dublin, 2010.

\bibitem[Kara\v{c} et~al.(2011)Kara\v{c}, Blackman, Cooper, Kinloch, Sanchez,
  Teo, and Ivankovi\'{c}]{Karac2011}
A.~Kara\v{c}, B.~R.~K. Blackman, V.~Cooper, A.~J. Kinloch, S.~Rodriguez
  Sanchez, W.~S. Teo, and A.~Ivankovi\'{c}.
\newblock Modelling the fracture behaviour of adhesively-bonded joints as a
  function of test rate.
\newblock \emph{Engineering Fracture Mechanics}, 78:\penalty0 973--989, 2011.

\bibitem[Ivankovi\'{c} et~al.(1997{\natexlab{b}})Ivankovi\'{c}, Trop\v{s}a, and
  Williams]{Ivankovic1997b}
A.~Ivankovi\'{c}, V.~Trop\v{s}a, and J.~G. Williams.
\newblock Finite volume modelling of residual stresses in cast plastic slabs.
\newblock In \emph{Fifth international conference on residual stresses
  {ICRS-5}}, pages 392--399, Linkoping Sweden, 1997{\natexlab{b}}.

\bibitem[Trop\v{s}a et~al.(2000)Trop\v{s}a, Ivankovic, and
  Williams]{Tropsa2000}
V.~Trop\v{s}a, A.~Ivankovic, and J.~G. Williams.
\newblock Predicting residual stresses due to solidification in cast plastic
  plates.
\newblock \emph{Plastics, Rubber and Composites}, 29\penalty0 (9):\penalty0
  468--474, 2000.
\newblock \doi{10.1179/146580100101541319}.
\newblock URL \url{https://doi.org/10.1179/146580100101541319}.

\bibitem[Sato et~al.(2006)Sato, Ohnaka, and Iwane]{Sato2006}
A.~Sato, I.~Ohnaka, and J.~Iwane.
\newblock Stress analysis by finite volume method for prediction of porosity
  and deformation defects of spheroidal graphite castings.
\newblock \emph{Journal of Japan Foundry Engineering Society}, 78:\penalty0
  231--237, 2006.
\newblock In Japanese.

\bibitem[Teskered{\v{z}}i\'{c} et~al.(2015{\natexlab{a}})Teskered{\v{z}}i\'{c},
  Demird{\v{z}}i\'{c}, and Muzaferija]{Teskeredzic2015a}
A.~Teskered{\v{z}}i\'{c}, I.~Demird{\v{z}}i\'{c}, and S.~Muzaferija.
\newblock Numerical method for calculation of complete casting process - part
  {I}: Theory.
\newblock \emph{Numerical Heat Transfer, Part B: Fundamentals: An International
  Journal of Computation and Methodology}, 68:\penalty0 295--316,
  2015{\natexlab{a}}.

\bibitem[Teskered{\v{z}}i\'{c} et~al.(2015{\natexlab{b}})Teskered{\v{z}}i\'{c},
  Demird{\v{z}}i\'{c}, and Muzaferija]{Teskeredzic2015b}
A.~Teskered{\v{z}}i\'{c}, I.~Demird{\v{z}}i\'{c}, and S.~Muzaferija.
\newblock Numerical method for calculation of complete casting process - part
  {II}: Validation and application.
\newblock \emph{Numerical Heat Transfer, Part B: Fundamentals: An International
  Journal of Computation and Methodology}, 68:\penalty0 317--335,
  2015{\natexlab{b}}.

\bibitem[Henry and Collins(1993{\natexlab{a}})]{Henry1993a}
F.~S. Henry and M.~W. Collins.
\newblock A novel predictive model with compliance for arterial flows.
\newblock \emph{{BED-Vol. 26}, Advances in Bioengineering, ASME},
  1993{\natexlab{a}}.

\bibitem[Henry and Collins(1993{\natexlab{b}})]{Henry1993b}
F.~S. Henry and M.~W. Collins.
\newblock Prediction of transient wall movement of an incompressible elastic
  tube using finite volume procedure.
\newblock In \emph{{Proceedings of BIOMED93}}, Bath, UK, 1993{\natexlab{b}}.

\bibitem[Greenshields et~al.(2000)Greenshields, Vanizelos, and
  Ivankovi\'{c}]{Greenshields2000}
C.~J. Greenshields, G.~P. Vanizelos, and A.~Ivankovi\'{c}.
\newblock A fluid-structure model for fast brittle fracture in plastic pipes.
\newblock \emph{Journal of Fluids and Structures}, 14:\penalty0 221--234, 2000.

\bibitem[Ivankovi\'{c} et~al.(2001)Ivankovi\'{c}, Kara\v{c}, Dendrinos, and
  Parker]{Ivankovic2001b}
A.~Ivankovi\'{c}, A.~Kara\v{c}, E.~Dendrinos, and K.~Parker.
\newblock Blood flow in deformable arteries: The finite volume method for
  fluid-structure interaction problem.
\newblock In \emph{Proceedings of the 9$^{th}$ {ACME} Conference on
  Computational Mechanics in Engineering}, Birmingham, {UK}, 2001.

\bibitem[Sch{\"{a}}fer and Teschauer(2001)]{Schafer2001b}
M.~Sch{\"{a}}fer and I.~Teschauer.
\newblock Numerical simulation of coupled fluid-solid problems.
\newblock \emph{Computer Methods in Applied Mechanics and Engineering},
  190:\penalty0 3645--3667, 2001.

\bibitem[Sch{\"{a}}fer et~al.(2002)Sch{\"{a}}fer, Teschauer, Kadinski, and
  Selder]{Schafer2002}
M.~Sch{\"{a}}fer, I.~Teschauer, L.~Kadinski, and M.~Selder.
\newblock A numerical approach for the solution of coupled fluid-solid and
  thermal stress problems in crystal growth processes.
\newblock \emph{Computational Materials Science}, 24:\penalty0 409--419, 2002.

\bibitem[Ivankovi\'{c} et~al.(2002{\natexlab{a}})Ivankovi\'{c}, Kara\v{c},
  Dendrinos, and Parker]{Ivankovic2002b}
A.~Ivankovi\'{c}, A.~Kara\v{c}, E.~Dendrinos, and K.~Parker.
\newblock Towards early diagnosis of artherosclerosis: The finite volume method
  for fluid-structure interaction.
\newblock \emph{Biorheology}, 39:\penalty0 401--407, 2002{\natexlab{a}}.

\bibitem[Torlak and Muzaferija(2002)]{Torlak2002a}
M.~Torlak and S.~Muzaferija.
\newblock Finite volume approach to computation of elastic plates and their
  interaction with fluid flows.
\newblock In \emph{{Finite Volumes for Complex Applications III - Problems and
  Perspectives}}. Kogan Page Science, 2002.

\bibitem[Torlak et~al.(2002)Torlak, Muzaferija, and Peri\'{c}]{Torlak2002b}
M.~Torlak, S.~Muzaferija, and M.~Peri\'{c}.
\newblock Application of a finite volume method to the computation of
  interaction between thin linearly elastic structures and incompressible fluid
  flows.
\newblock In \emph{{VDI-Berichte 1862, VDI Tagung
  Fluid-Struktur-Wechselwirkung, Wiesloch}}, 2002.

\bibitem[Kova\v{c}evi\'{c} et~al.(2004{\natexlab{a}})Kova\v{c}evi\'{c},
  Sto\v{s}i\'{c}, Smith, and Numerical]{Kovacevic2004b}
A.~Kova\v{c}evi\'{c}, N.~Sto\v{s}i\'{c}, I.~K. Smith, and A.~Numerical.
\newblock A numerical study of fluid -solid interaction in screw compressors.
\newblock \emph{International Journal on Computer Application in Technology},
  21:\penalty0 148--158, 2004{\natexlab{a}}.
\newblock \doi{10.1504/IJCAT.2004.006651}.

\bibitem[Sto\v{s}i\'{c} et~al.(2005)Sto\v{s}i\'{c}, Smith, and
  Kova\v{c}evi\'{c}]{Stosic2005}
N.~Sto\v{s}i\'{c}, I.~Smith, and A.~Kova\v{c}evi\'{c}.
\newblock \emph{Screw Compressors}.
\newblock Springer, Berlin, 2005.

\bibitem[Shaw and Stone(2005)]{Shaw2005}
G.~Shaw and T.~Stone.
\newblock Finite volume methods for coupled stress/fluid flow in commercial
  reservoir simulations.
\newblock In \emph{{SPE} Reservoir Simulation Symposium, Houston, Texas}, 2005.
\newblock {SPE} 93430.

\bibitem[Torlak(2006)]{Torlak2006}
M.~Torlak.
\newblock \emph{A Finite-Volume Method for Coupled Numerical Analysis of
  Incompressible Fluid Flow and Linear Deformation of Elastic Structures}.
\newblock PhD thesis, Technischen Universitaet Hamburg-Harburg, 2006.

\bibitem[Kova\v{c}evi\'{c} et~al.(2006)Kova\v{c}evi\'{c}, Sto\v{s}i\'{c}, and
  Smith]{Kovacevic2006}
A.~Kova\v{c}evi\'{c}, N.~Sto\v{s}i\'{c}, and I.~K. Smith.
\newblock Numerical simulation of combined screw compressor-expander machines
  for use in high pressure refrigeration systems.
\newblock \emph{Simulation Modelling Practice and Theory}, 14:\penalty0
  1143--1154, 2006.

\bibitem[Kova\v{c}evi\'{c} et~al.(2007)Kova\v{c}evi\'{c}, Sto\v{s}i\'{c},
  Muji\'{c}, and Smith]{Kovacevic2007}
A.~Kova\v{c}evi\'{c}, N.~Sto\v{s}i\'{c}, E.~Muji\'{c}, and I.~K. Smith.
\newblock {CFD} integrated design of screw compressors.
\newblock \emph{Journal of Engineering Applications of Computational Fluid
  Mechanics}, 1:\penalty0 96--108, 2007.

\bibitem[Sto\v{s}i\'{c} et~al.(2007)Sto\v{s}i\'{c}, Smith, and
  Kova\v{c}evi\'{c}]{Stosic2007}
N.~Sto\v{s}i\'{c}, I.~Smith, and A.~Kova\v{c}evi\'{c}.
\newblock Screw compressors.
\newblock \emph{Three Dimensional Computational Fluid Dynamics and Solid Fluid
  Interaction, Springer, Berlin}, 2007.

\bibitem[Papadakis(2008)]{Papadakis2008}
G.~Papadakis.
\newblock A novel pressure-velocity formulation and solution method for
  fluid-structure interaction problems.
\newblock \emph{Journal of Computational Physics}, 227:\penalty0 3383--3404,
  2008.

\bibitem[Jasak et~al.(2007)Jasak, Jemcov, and Tukovi\'{c}]{Jasak2007}
H.~Jasak, A.~Jemcov, and {\v{Z}}.~Tukovi\'{c}.
\newblock {OpenFOAM}: A {C++} library for complex physics simulations.
\newblock In \emph{International Workshop on Coupled Methods in Numerical
  Dynamics}, Dubrovnik, Croatia, 2007.

\bibitem[Kanyanta et~al.(2009{\natexlab{a}})Kanyanta, Ivankovi\'{c}, and
  Kara\v{c}]{Kanyanta2009a}
V.~Kanyanta, A.~Ivankovi\'{c}, and A.~Kara\v{c}.
\newblock Validation of a fluid-structure interaction numerical model for
  predicting flow transients in arteries.
\newblock \emph{Journal of Biomechanics}, 42:\penalty0 1705--1712,
  2009{\natexlab{a}}.

\bibitem[Jagad et~al.(2011)Jagad, Puranik, and Date]{Jagad2011}
P.~Jagad, B.~P. Puranik, and A.~W. Date.
\newblock A finite volume procedure on unstructured meshes for fluid-structure
  interaction problems.
\newblock \emph{World Academy of Science, Engineering and Technology
  International Journal of Mechanical, Aerospace, Industrial Mechatronic and
  Manufacturing Engineering}, 5:\penalty0 1406--1412, 2011.

\bibitem[Wiedemair et~al.(2012)Wiedemair, Tukovi\'{c}, Jasak, Poulikakos, and
  Kurtcuoglu]{Wiedemair2012}
W.~Wiedemair, {\v{Z}}.~Tukovi\'{c}, H.~Jasak, D.~Poulikakos, and V.~Kurtcuoglu.
\newblock On ultrasound-induced microbubble oscillation in acapillary blood
  vessel and its implications for the blood-brain barrier.
\newblock \emph{Physics in Medicine \& Biology}, 57:\penalty0 1019--1045, 2012.

\bibitem[Habchi et~al.(2013)Habchi, Russeil, Bougeard, Harion, Lemenand,
  Ghanem, Valle, and Peerhossaini]{Habchi2013}
C.~Habchi, S.~Russeil, D.~Bougeard, J.-L. Harion, T.~Lemenand, A.~Ghanem,
  D.~Della Valle, and H.~Peerhossaini.
\newblock Partitioned solver for strongly coupled fluid-structure interaction.
\newblock \emph{Computers \& Fluids}, 71:\penalty0 306--319, 2013.

\bibitem[Tukovi\'{c} et~al.(2014)Tukovi\'{c}, Cardiff, Ivankovi\'{c}, and
  Kara\v{c}]{Tukovic2014}
{\v{Z}}.~Tukovi\'{c}, P.~Cardiff, A.~Ivankovi\'{c}, and A.~Kara\v{c}.
\newblock {OpenFOAM} library for fluid structure interaction.
\newblock In \emph{9$^{th}$ {OpenFOAM} Workshop, Zagreb, Croatia}, Zagreb,
  Croatia, 2014.

\bibitem[Smith et~al.(2014)Smith, Sto\v{s}i\'{c}, and
  Kova\v{c}evi\'{c}]{Smith2014}
I.~Smith, N.~Sto\v{s}i\'{c}, and A.~Kova\v{c}evi\'{c}.
\newblock \emph{Power Recovery From Low Grade Heat Sources by the Use of Screw
  Expanders}.
\newblock Chandos Publishing, 2014.

\bibitem[\v{S}ekutkovski et~al.(2016)\v{S}ekutkovski, Kosti\'{c},
  Simonovi\'{c}, Cardiff, and Jazarevi\'{c}]{Sekutkovski2016}
B.~\v{S}ekutkovski, I.~Kosti\'{c}, A.~Simonovi\'{c}, P.~Cardiff, and
  V.~Jazarevi\'{c}.
\newblock Three-dimensional fluid-structure interaction simulation with a
  hybrid {RANS-LES} turbulence model for applications in transonic flow domain.
\newblock \emph{Aerospace Science and Technology}, 49:\penalty0 1 -- 16, 2016.

\bibitem[Jagad(2016)]{Jagad2016}
P.~Jagad.
\newblock A numerical procedure for elastic solids.
\newblock \emph{{GIT}-Journal of Engineering and Technology}, 9:\penalty0
  113--124, 2016.

\bibitem[Cardiff et~al.(2017{\natexlab{b}})Cardiff, Kara\v{c}, Jaeger, Jasak,
  Nagy, Ivankovi\'{c}, and Tukovi\'{c}]{Cardiff2017b}
P.~Cardiff, A.~Kara\v{c}, P.~De Jaeger, H.~Jasak, J.~Nagy, A.~Ivankovi\'{c},
  and {\v{Z}}.~Tukovi\'{c}.
\newblock Towards the development of an extendable solid mechanics and
  fluid-solid interactions toolbox for {OpenFOAM}.
\newblock In \emph{12th {OpenFOAM} Workshop}, volume~12, University of Exeter,
  UK, 2017{\natexlab{b}}.

\bibitem[Jagad et~al.(2017)Jagad, Puranik, and Date]{Jagad2017}
P.~Jagad, B.P. Puranik, and A.W. Date.
\newblock A numerical analysis of fluid-structure interaction problem with a
  flow channel embedded in a structural material.
\newblock \emph{Proceedings of the Indian National Science Academy},
  83:\penalty0 655--667, 2017.

\bibitem[Jagad et~al.(April 2018)Jagad, Puranik, and Date]{Jagad2018}
P.~Jagad, B.P. Puranik, and A.W. Date.
\newblock A finite volume procedure for fluid flow, heat transfer and
  solid-body stress analysis.
\newblock \emph{International Journal for Computational Methods in Engineering
  Science and Mechanics}, April 2018.
\newblock \doi{10.1080/15502287.2018.1434839}.

\bibitem[Demird{\v{z}}i\'{c} and Ivankovi\'{c}(1997)]{Demirdzic1997b}
I.~Demird{\v{z}}i\'{c} and A.~Ivankovi\'{c}.
\newblock Finite volume approach to modelling of plates.
\newblock In \emph{Proceedings of 2$^{nd}$ Congress of Croatian Society of
  Mechanics}, pages 101--108, Brac, Croatia, 1997.

\bibitem[Fallah(2004)]{Fallah2004}
N.~Fallah.
\newblock A cell vertex and cell centred finite volume method for plate bending
  analysis.
\newblock \emph{Computer Methods in Applied Mechanics and Engineering},
  193:\penalty0 3457--3470, 2004.

\bibitem[Fallah(2006{\natexlab{a}})]{Fallah2006a}
N.~Fallah.
\newblock A finite volume method for plate buckling analysis.
\newblock In C.~A. Mota~Soares et.al, editor, \emph{{III} European Conference
  on Computational Mechanics, Solids, Structures and Coupled Problems in
  Engineering,}, pages 5--8, Lisbon, Portugal, 2006{\natexlab{a}}.

\bibitem[Fallah et~al.(2006)Fallah, Hatami, and Displacement]{Fallah2006b}
N.~Fallah, F.~Hatami, and A.~Displacement.
\newblock Formulation based on finite volume method for analysis of
  {Timoshenko} beam.
\newblock In \emph{Proceedings of the 7$^{th}$ international conference on
  civil engineering}, pages 8--10, Tehran, Iran, May 2006.

\bibitem[Fallah(2006{\natexlab{b}})]{Fallah2006c}
N.~Fallah.
\newblock On the use of shape functions in the cell centered finite volume
  formulation for plate bending analysis based on {Mindlin-Reissner} plate
  theory.
\newblock \emph{Computers \& Structures}, 84:\penalty0 1664--1672,
  2006{\natexlab{b}}.

\bibitem[Fallah and Hatami(2006)]{Fallah2006d}
N.~Fallah and F.~Hatami.
\newblock Extension of the finite volume method for instability analysis of
  columns with shear effects.
\newblock In \emph{Proceedings of the Eighth International Conference on
  Computational Structures Technology}, Stirlingshire, Scotland, 2006.
  Civil-Comp Press.
\newblock Paper 192.

\bibitem[Hatami et~al.(2006)Hatami, Fallah, and Pourzeynali]{Hatami2006}
F.~Hatami, N.~Fallah, and S.~Pourzeynali.
\newblock Application of the finite volume method for shell analysis: A
  membrane study.
\newblock In B.~H.~V. Topping, G.~Montero, and R.~Montenegro, editors,
  \emph{Proceedings of the Eighth International Conference on Computational
  Structures Technology}. Civil-Comp Press, Stirlingshire, UK, 2006.
\newblock \doi{10.4203/ccp.83.160}.
\newblock Paper 160.

\bibitem[Isi\'{c} et~al.(2007{\natexlab{a}})Isi\'{c}, Dole\v{c}ek, and
  Karabegovi\'{c}]{Isic2007a}
S.~Isi\'{c}, V.~Dole\v{c}ek, and I.~Karabegovi\'{c}.
\newblock A comparison between finite element and finite volume methods on the
  problem of stability of timoshenko beam.
\newblock In \emph{The 12th International Conference on Problems of Material
  Engineering}, Jasna, Slovakia, 2007{\natexlab{a}}. Mechanics and Design.

\bibitem[Isi\'{c} et~al.(2007{\natexlab{b}})Isi\'{c}, Dole\v{c}ek, and
  Karabegovi\'{c}]{Isic2007b}
S.~Isi\'{c}, V.~Dole\v{c}ek, and I.~Karabegovi\'{c}.
\newblock Numerical and experimental analysis of postbuckling behaviour of
  prismatic beam under displacement dependent loading.
\newblock In \emph{Proceedings of the First Serbian Congress on Theoretical and
  Applied Mechanics}, Kopaonik, Serbia, 2007{\natexlab{b}}.

\bibitem[Isi\'{c} et~al.(2007{\natexlab{c}})Isi\'{c}, Dole\v{c}ek, and
  Karabegovi\'{c}]{Isic2007c}
S.~Isi\'{c}, V.~Dole\v{c}ek, and I.~Karabegovi\'{c}.
\newblock A comparison of finite element and finite volume method on stability
  analysis of rectangular plate.
\newblock In V.~Dole\v{c}ek Karabegovi\'{c} and M.~Jurkovi\'{c}, editors,
  \emph{I}. 6$^{th}$ International Scientific Conference on Production
  Engineering, Development and Modernization of Production, RIM Bihac, BiH,
  2007{\natexlab{c}}.

\bibitem[Isi\'{c}(2008)]{Isic2008}
S.~Isi\'{c}.
\newblock \emph{Numerical and experimental analysis of nonlinear stability
  phenomena in elastic systems}.
\newblock PhD thesis, University of Bihac, Bosnia and Herzegovina, 2008.
\newblock In Bosnian.

\bibitem[Fallah(2013)]{Fallah2013a}
N.~Fallah.
\newblock Finite volume method for determining the natural characteristics of
  structures.
\newblock \emph{Journal of Engineering Science and Technology}, 8:\penalty0
  93--106, 2013.

\bibitem[Fallah and Parayandeh-Shahrestany(2014)]{Fallah2014a}
N.~Fallah and A.~Parayandeh-Shahrestany.
\newblock A novel finite volume based formulation for the elasto-plastic
  analysis of plates.
\newblock \emph{Thin-Walled Structures}, 77:\penalty0 153--164, 2014.

\bibitem[Fallah and Ebrahimnejad(2014)]{Fallah2014b}
N.~Fallah and M.~Ebrahimnejad.
\newblock Finite volume analysis of adaptive beams with piezoelectric sensors
  and actuators.
\newblock \emph{Applied Mathematical Modelling}, 38:\penalty0 727--737, 2014.

\bibitem[Jing et~al.(2016)Jing, Ming, Zhang, Fu, and Cao]{Jing2016}
L.-L. Jing, P.-J. Ming, W.-P. Zhang, L.-R. Fu, and Y.-P. Cao.
\newblock Static and free vibration analysis of functionally graded beams by
  combination {Timoshenko} theory and finite volume method.
\newblock \emph{Composite Structures}, 138:\penalty0 192--213, 2016.

\bibitem[Fallah and Ghanbari(2017)]{Fallah2017a}
N.~Fallah and A.~Ghanbari.
\newblock A displacement finite volume formulation for the static and dynamic
  analysis of shear deformable circular curved beams.
\newblock \emph{Scientia Iranica}, pages~--, 2017.
\newblock ISSN 1026-3098.
\newblock \doi{10.24200/sci.2017.4259}.
\newblock URL \url{http://scientiairanica.sharif.edu/article_4259.html}.

\bibitem[Fallah et~al.(2017)Fallah, Parayandeh-Shahrestany, and
  Golkoubi]{Fallah2017b}
N.~Fallah, A.~Parayandeh-Shahrestany, and H.~Golkoubi.
\newblock A. finite volume formulation for the elasto-plastic analysis of
  rectangular {Mindlin-Reissner} plates, a non-layered approach.
\newblock \emph{Civil Engineering Infrastructures Journal}, 50:\penalty0
  293--310, 2017.

\bibitem[Mohebi et~al.(2017)Mohebi, Kaboudan, and Yazdanpanah]{Mohebi2017}
B.~Mohebi, A.~R. Kaboudan, and O.~Yazdanpanah.
\newblock Damage detection in beam-like structures using finite volume method.
\newblock \emph{Journal of Rehabilitation in Civil Engineering}, 5:\penalty0
  77--92, 2017.

\bibitem[Fallah and Ghanbary(2018)]{Fallah2018b}
N.~Fallah and A.~Ghanbary.
\newblock A displacement finite volume formulation for the static and dynamic
  analysis of shear deformable curved beams.
\newblock \emph{Scientia Iranica}, 25(3):\penalty0 999--1014, 2018.

\bibitem[Amraei and Fallah(2018)]{Amraei2018}
A.~Amraei and N.~Fallah.
\newblock A cell-centered finite volume formulation for the calculation of
  stress intensity factors in {Mindlin-Reissner} cracked plates.
\newblock \emph{Civil Engineering Journal}, 3:\penalty0 1366--1385, 2018.

\bibitem[Tukovi\'{c} et~al.(2019)Tukovi\'{c}, De~Jaeger, Cardiff, and
  Ivankovi\'{c}]{Tukovic2019}
{\v{Z}}.~Tukovi\'{c}, P.~De~Jaeger, P.~Cardiff, and A.~Ivankovi\'{c}.
\newblock A finite volume solver for geometrically exact {Simo-Reissner} beams.
\newblock In \emph{{ECCOMAS MSF 2019}, Thematic conference}, Sarajevo, Bosnian
  and Herzegovina, 2019.

\bibitem[Das et~al.(2011{\natexlab{a}})Das, Mathur, and Murthy]{Das2011a}
S.~Das, S.~R. Mathur, and J.~Y. Murthy.
\newblock {An unstructured finite-volume method for structure-electrostatic
  interactions in MEMS}.
\newblock \emph{Numerical Heat Transfer, Part B: Fundamentals}, 60:\penalty0
  425--451, 2011{\natexlab{a}}.

\bibitem[Das et~al.(2011{\natexlab{b}})Das, Koslowski, Mathur, and
  Murthy]{Das2011b}
S.~Das, M.~Koslowski, S.~R. Mathur, and J.~Y. Murthy.
\newblock {Finite volume method for simulation of creep in RF MEMS devices}.
\newblock In \emph{{ASME 2011 International Mechanical Engineering Congress and
  Exposition: Nano and Micro Materials, Devices and Systems; Microsystems
  Integration}}, volume~11, Denver, Colorado, USA, 2011{\natexlab{b}}.

\bibitem[Dioh et~al.(1995{\natexlab{a}})Dioh, Ivankovi\'{c}, and
  Demird{\v{z}}i\'{c}]{Dioh1995a}
N.~N. Dioh, A.~Ivankovi\'{c}, and I.~Demird{\v{z}}i\'{c}.
\newblock Dynamic thermo elastic plastic deformation of solids using finite
  volume technique.
\newblock In E.~O{\~{n}}ate D.R.J.~Owen, editor, \emph{Computational
  Plasticity, Fundamentals and Applications,}, pages 1947--1957. Pineridge
  Press, 1995{\natexlab{a}}.

\bibitem[Dioh et~al.(1995{\natexlab{b}})Dioh, Ivankovi\'{c}, Leevers, and
  Williams]{Dioh1995b}
N.~N. Dioh, A.~Ivankovi\'{c}, P.~S. Leevers, and J.~G. Williams.
\newblock Stress wave propagation effects in split hopkinson pressure bar
  tests.
\newblock In \emph{Proceedings of The Royal Society; Proceedings: Mathematical
  and Physical Sciences}, pages 187--204. 449, 1995{\natexlab{b}}.

\bibitem[Oosterkamp et~al.(2000)Oosterkamp, Ivankovic, and
  Venizelos]{Oosterkamp2000}
L.~Djapic Oosterkamp, A.~Ivankovic, and G.~Venizelos.
\newblock High strain rate properties of selected aluminium alloys.
\newblock \emph{Materials Science and Engineering: A}, 278\penalty0
  (1):\penalty0 225 -- 235, 2000.
\newblock ISSN 0921-5093.
\newblock \doi{https://doi.org/10.1016/S0921-5093(99)00570-5}.
\newblock URL
  \url{http://www.sciencedirect.com/science/article/pii/S0921509399005705}.

\bibitem[Weller et~al.(1998)Weller, Tabor, Jasak, and Fureby]{Weller1998}
H.~G. Weller, G.~Tabor, H.~Jasak, and C.~Fureby.
\newblock A tensorial approach to computational continuum mechanics using
  object orientated techniques.
\newblock \emph{Computers in Physics}, 12:\penalty0 620--631, 1998.

\bibitem[Jasak and Weller(2000{\natexlab{b}})]{Jasak2000a}
H.~Jasak and H.~G. Weller.
\newblock Application of the finite volume method and unstructured meshes to
  linear elasticity.
\newblock \emph{International Journal for Numerical Methods in Engineering},
  48:\penalty0 267--287, 2000{\natexlab{b}}.

\bibitem[Suvanjumrat and Chaichanasiri(2011)]{Suvanjumrat2011}
C.~Suvanjumrat and E.~Chaichanasiri.
\newblock Implementation and validation of finite volume {C++} codes for plane
  stress analysis.
\newblock In \emph{{CST02}, The Second {TSME}}, Krabi, 2011. International
  Conference on Mechanical Engineering.

\bibitem[Haider et~al.(2018)Haider, Lee, Gil, Huerta, and Bonet]{Haider2018}
J.~Haider, C.~H. Lee, A.~J. Gil, A.~Huerta, and J.~Bonet.
\newblock An upwind cell centred total {Lagrangian} finite volume algorithm for
  nearly incompressible explicit fast solid dynamic applications.
\newblock \emph{Computer Methods in Applied Mechanics and Engineering},
  340:\penalty0 684?727, 2018.

\bibitem[Demird{\v{z}}i\'{c} et~al.(1997{\natexlab{a}})Demird{\v{z}}i\'{c},
  Muzaferija, and Peri\'{c}]{Demirdzic1997c}
I.~Demird{\v{z}}i\'{c}, S.~Muzaferija, and M.~Peri\'{c}.
\newblock Benchmark solutions of some structural analysis problems using
  finite-volume method and multigrid acceleration.
\newblock \emph{International Journal for Numerical Methods in Engineering},
  40:\penalty0 1893--1908, 1997{\natexlab{a}}.

\bibitem[Cardiff et~al.(2014{\natexlab{d}})Cardiff, Tukovi\'{c}, Jasak, and
  A.~Ivankovi\'{c}]{Cardiff2014e}
P.~Cardiff, {\v{Z}}.~Tukovi\'{c}, H.~Jasak, and A~A.~Ivankovi\'{c}.
\newblock Block-coupled finite volume methodology for linear elasticity.
\newblock In \emph{9$^{th}$ {OpenFOAM} Workshop}, volume~9, University of
  Zagreb, Croatia, 2014{\natexlab{d}}.

\bibitem[Gonz\'{a}lez et~al.(15 June 2018)Gonz\'{a}lez, Naseri, Chiva, Rigola,
  and P\'{e}rez-Segarra]{Gonzalez2018}
I.~Gonz\'{a}lez, A.~Naseri, J.~Chiva, J.~Rigola, and C.~D. P\'{e}rez-Segarra.
\newblock An enhanced finite volume based solver for thermoelastic materials in
  fluid-structure coupled problems.
\newblock In \emph{{6$^{th}$ European Conference on Computational Mechanics
  (ECCM 6), 7$^{th}$ European Conference on Computational Fluid Dynamics (ECFD
  7)}}, Glasgow, UK, 15 June 2018.

\bibitem[Fowler and Yee(2003)]{Fowler2003}
B.~L. Fowler and R.~K. Yee.
\newblock Application of finite volume method for solid mechanics.
\newblock In \emph{{Proceedings of the ASME IMECE Conference International
  Mechanical Engineering Congress and RD\&D Expo}}, pages 15--21, November
  2003.
\newblock IMECE2003-44297.

\bibitem[Bijelonja(2011{\natexlab{a}})]{Bijelonja2011a}
I.~Bijelonja.
\newblock A numerical method for almost incompressible body problem.
\newblock In Katalini\'{c}, editor, \emph{Proceedings of the 22$^{nd}$
  International {DAAAM} Symposium}, pages 321--322, Vienna, Austria,
  2011{\natexlab{a}}.

\bibitem[Oliveira and Rente(1999)]{Oliveira1999}
P.~J. Oliveira and C.~J. Rente.
\newblock Development and application of a finite volume method for static and
  transient stress analysis.
\newblock In \emph{Proceedings of {NAFEMS} World Congress'99 on Effective
  Engineering Analysis}, pages 297--309, 1999.

\bibitem[Demird{\v{z}}i\'{c}(2016)]{Demirdzic2016}
I.~Demird{\v{z}}i\'{c}.
\newblock A fourth-order finite volume method for structural analysis.
\newblock \emph{Applied Mathematical Modelling}, 40:\penalty0 3104--3114, 2016.

\bibitem[Fallah(2008)]{Fallah2008a}
N.~Fallah.
\newblock A method for calculation of face gradients in two-dimensional cell
  centred finite volume formulation for stress analysis in solid problems.
\newblock \emph{Scientia Iranica}, 15:\penalty0 286--294, 2008.

\bibitem[Nordbotten(2014)]{Nordbotten2014}
J.~M. Nordbotten.
\newblock Cell-centered finite volume discretizations for deformable porous
  media.
\newblock \emph{International Journal for Numerical Methods in Engineering},
  100:\penalty0 399--418, 2014.

\bibitem[Nordbotten(2015)]{Nordbotten2015}
J.~M. Nordbotten.
\newblock Convergence of a cell-centered finite volume discretization for
  linear elasticity.
\newblock \emph{{SIAM} Journal on Numerical Analysis}, 53:\penalty0 2605--2625,
  2015.

\bibitem[Keilegavlen and Nordbotten(2017)]{Keilegavlen2017}
E.~Keilegavlen and J.~M. Nordbotten.
\newblock Finite volume methods for elasticity with weak symmetry.
\newblock \emph{International Journal for Numerical Methods in Engineering},
  2017.
\newblock \doi{10.1002/nme.5538}.

\bibitem[Tukovi\'{c} et~al.(2018{\natexlab{a}})Tukovi\'{c}, Kara\v{c}, Cardiff,
  Jasak, and Ivankovi\'{c}]{Tukovic2018a}
{\v{Z}}.~Tukovi\'{c}, A.~Kara\v{c}, P.~Cardiff, H.~Jasak, and A.~Ivankovi\'{c}.
\newblock {OpenFOAM} finite volume solver for fluid-solid interaction.
\newblock \emph{Transactions of {FAMENA}}, 42\penalty0 (3):\penalty0 1--31,
  2018{\natexlab{a}}.
\newblock \doi{10.21278/TOF.42301}.

\bibitem[Fallah(Winter 2008)]{Fallah2008b}
N.~Fallah.
\newblock {Finite Volume Based Formulations for the Analysis of Bernouli and
  {Timoshenko} Beams}.
\newblock \emph{Journal of Numerical Simulation in Engineering}, 1\penalty0
  (3):\penalty0 259--268, Winter 2008.

\bibitem[Golubovi\'{c}(2017)]{Golubovic2017b}
A.~Golubovi\'{c}.
\newblock \emph{Finite volume analysis of laminated composite plates}.
\newblock PhD thesis, University of Sarajevo, 2017.

\bibitem[Godunov(1959)]{Godunov1959}
S.~K. Godunov.
\newblock A difference method for numerical calculation of discontinuous
  solutions of the equations of hydrodynamics.
\newblock \emph{Matematicheskii Sbornik}, 47:\penalty0 271--306, 1959.

\bibitem[Godunov(1962)]{Godunov1962}
S.~K. Godunov.
\newblock The problem of a generalized solution in the theory of quasi-linear
  equations and in gas dynamics.
\newblock \emph{Russian Mathematical Surveys}, 17:\penalty0 145--156, 1962.

\bibitem[Trangenstein and Pember(1992)]{Trangenstein1992}
J.~A. Trangenstein and R.~B. Pember.
\newblock Numerical algorithms for strong discontinuities in elastic-plastic
  solids.
\newblock \emph{Journal of Computational Physics}, 103\penalty0 (1):\penalty0
  63--89, 1992.

\bibitem[Trangenstein(1994)]{Trangenstein1994}
J.~A. Trangenstein.
\newblock Second-order {Godunov} algorithm for two-dimensional solid mechanics.
\newblock \emph{Computational Mechanics}, 13:\penalty0 343--359, 1994.

\bibitem[Miller and Puckett(1996)]{Miller1996}
G.~H. Miller and E.~G. Puckett.
\newblock A high-order {Godunov} method for multiple condensed phases.
\newblock \emph{Journal of Computational Physics}, 128:\penalty0 134--164,
  1996.

\bibitem[Tang and Sotiropoulos(1999)]{Tang1999}
H.~Tang and F.~Sotiropoulos.
\newblock A second-order {Godunov} method for wave problems in coupled solid
  water gas systems.
\newblock \emph{Journal of Computational Physics}, 151:\penalty0 790--815,
  1999.

\bibitem[Berezovski and Maugin(2001)]{Berezovski2001}
A.~Berezovski and G.~A. Maugin.
\newblock Simulation of thermoelastic wave propagation by means of a composite
  wave-propagation algorithm.
\newblock \emph{Journal of Computational Physics}, 168:\penalty0 249--264,
  2001.

\bibitem[Howell and Ball(2002)]{Howell2002}
B.~P. Howell and G.~J. Ball.
\newblock A free-{Lagrange} augmented {Godunov} method for the simulation of
  elastic-plastic solids.
\newblock \emph{Journal of Computational Physics}, 175:\penalty0 128--167,
  2002.

\bibitem[Berezovski and Maugin(2003)]{Berezovski2003}
A.~Berezovski and G.~A. Maugin.
\newblock Simulation of wave and front propagation in thermoelastic materials
  with phase transformation.
\newblock \emph{Computational Materials Science}, 28:\penalty0 478--485, 2003.

\bibitem[Kluth and Despr\'{e}s(2008)]{Kluth2008}
G.~Kluth and B.~Despr\'{e}s.
\newblock {F. V.} schemes for hyperelastic-plastic models in finite
  deformations.
\newblock In Reymard and J.-M. H\'{e}rard, editors, \emph{Finite Volumes for
  Complex Applications {V}-Problems \& Perspectives}. John Wiley \& Sons, 2008.

\bibitem[Carr\'{e} et~al.(2009)Carr\'{e}, Pino, Despr\'{e}s, and
  Labourasse]{Carre2009}
G.~Carr\'{e}, S.~Del Pino, B.~Despr\'{e}s, and E.~Labourasse.
\newblock A cell-centered {Lagrangian} hydrodynamics scheme on general
  unstructured meshes in arbitrary dimension.
\newblock \emph{Journal of Computational Physics}, 228:\penalty0 5160--5183,
  2009.

\bibitem[Maire et~al.(2013)Maire, Abgrall, Breil, Loub\`{e}re, and
  Rebourcet]{Maire2013}
P.-H. Maire, R.~Abgrall, J.~Breil, R.~Loub\`{e}re, and B.~Rebourcet.
\newblock A nominally second-order cell-centered {Lagrangian} scheme for
  simulating elastic- plastic flows on two-dimensional unstructured grids.
\newblock \emph{Journal of Computational Physics}, 235:\penalty0 626--665,
  2013.

\bibitem[Sambasivan et~al.(2013)Sambasivan, Shashkov, and
  Burton]{Sambasivan2013}
S.~K. Sambasivan, M.-J. Shashkov, and D.~E. Burton.
\newblock A finite volume cell-centered lagrangian hydrodynamics approach for
  solids in general unstructured grids.
\newblock \emph{International Journal for Numerical Methods in Fluids},
  72:\penalty0 770--810, 2013.

\bibitem[Sijoy and Chaturvedi(2015)]{Sijoy2015}
C.~D. Sijoy and Shashank Chaturvedi.
\newblock An {Eulerian} multi-material scheme for elastic-plastic impact and
  penetration problems involving large material deformations.
\newblock \emph{European Journal of Mechanics - B/Fluids}, 53:\penalty0 85 --
  100, 2015.
\newblock ISSN 0997-7546.
\newblock \doi{10.1016/j.euromechflu.2015.04.004}.
\newblock URL
  \url{http://www.sciencedirect.com/science/article/pii/S0997754615000527}.

\bibitem[Despr\'{e}s and Labourasse(2015)]{Despres2015}
B.~Despr\'{e}s and E.~Labourasse.
\newblock Angular momentum preserving cell-centered {Lagrangian} and {Eulerian}
  schemes on arbitrary grids.
\newblock \emph{Journal of Computational Physics}, 290:\penalty0 28 -- 54,
  2015.
\newblock ISSN 0021-9991.
\newblock \doi{10.1016/j.jcp.2015.02.032}.
\newblock URL
  \url{http://www.sciencedirect.com/science/article/pii/S0021999115001023}.

\bibitem[Ndanou et~al.(2015)Ndanou, Favrie, and Gavrilyuk]{Ndanou2015}
S.~Ndanou, N.~Favrie, and S.~Gavrilyuk.
\newblock Multi-solid and multi-fluid diffuse interface model: Applications to
  dynamic fracture and fragmentation.
\newblock \emph{Journal of Computational Physics}, 295:\penalty0 523 -- 555,
  2015.
\newblock ISSN 0021-9991.
\newblock \doi{10.1016/j.jcp.2015.04.024}.
\newblock URL
  \url{http://www.sciencedirect.com/science/article/pii/S0021999115002806}.

\bibitem[Cheng et~al.(2015)Cheng, Toro, Jiang, Yu, and Tang]{Cheng2015}
Jun-Bo Cheng, Eleuterio~F. Toro, Song Jiang, Ming Yu, and Weijun Tang.
\newblock A high-order cell-centered {Lagrangian} scheme for one-dimensional
  elastic-plastic problems.
\newblock \emph{Computers \& Fluids}, 122:\penalty0 136 -- 152, 2015.
\newblock ISSN 0045-7930.
\newblock \doi{10.1016/j.compfluid.2015.08.029}.
\newblock URL
  \url{http://www.sciencedirect.com/science/article/pii/S0045793015003059}.

\bibitem[Loubere et~al.(2016)Loubere, Maire, and Rebourcet]{Loubere2016}
R.~Loubere, P.-H. Maire, and B.~Rebourcet.
\newblock Chapter 13 - staggered and colocated finite volume schemes for
  {Lagrangian} hydrodynamics.
\newblock In R.~Abgrall and C.-W. Shu, editors, \emph{Handbook of Numerical
  Methods for Hyperbolic Problems}, volume~17 of \emph{Handbook of Numerical
  Analysis}, pages 319 -- 352. Elsevier, 2016.
\newblock \doi{10.1016/bs.hna.2016.07.003}.
\newblock URL
  \url{http://www.sciencedirect.com/science/article/pii/S1570865916300059}.

\bibitem[Boscheri et~al.(2016)Boscheri, Dumbser, and Loub\'{e}re]{Boscheri2016}
W.~Boscheri, M~Dumbser, and R.~Loub\'{e}re.
\newblock {Cell centered direct Arbitrary-Lagrangian-Eulerian ADER-WENO finite
  volume schemes for nonlinear hyperelasticity}.
\newblock \emph{Computers \& Fluids}, 134-135:\penalty0 111 -- 129, 2016.
\newblock ISSN 0045-7930.
\newblock \doi{10.1016/j.compfluid.2016.05.004}.
\newblock URL
  \url{http://www.sciencedirect.com/science/article/pii/S004579301630144X}.

\bibitem[Vilar et~al.(2016)Vilar, Shu, and Maire]{Vilar2016}
F.~Vilar, C.-W. Shu, and P.-H. Maire.
\newblock Positivity-preserving cell-centered {Lagrangian} schemes for
  multi-material compressible flows: From first-order to high-orders. part {I}:
  The one-dimensional case.
\newblock \emph{Journal of Computational Physics}, 312:\penalty0 385 -- 415,
  2016.
\newblock ISSN 0021-9991.
\newblock \doi{10.1016/j.jcp.2016.02.027}.
\newblock URL
  \url{http://www.sciencedirect.com/science/article/pii/S0021999116000917}.

\bibitem[Georges et~al.(2017)Georges, Breil, and Maire]{Georges2017}
G.~Georges, J.~Breil, and P.-H. Maire.
\newblock A {3D} finite volume scheme for solving the updated {Lagrangian} form
  of hyperelasticity.
\newblock \emph{International Journal for Numerical Methods in Fluids},
  84\penalty0 (1):\penalty0 41--54, 2017.
\newblock \doi{10.1002/fld.4336}.
\newblock URL \url{https://onlinelibrary.wiley.com/doi/abs/10.1002/fld.4336}.

\bibitem[Heuz{\'e}(2017)]{Hueze2017}
T.~Heuz{\'e}.
\newblock {Lax-Wendroff and TVD finite volume methods for unidimensional
  thermomechanical numerical simulations of impacts on elastic-plastic solids}.
\newblock \emph{Journal of Computational Physics}, 346:\penalty0 369 -- 388,
  2017.
\newblock ISSN 0021-9991.
\newblock \doi{10.1016/j.jcp.2017.06.027}.
\newblock URL
  \url{http://www.sciencedirect.com/science/article/pii/S0021999117304722}.

\bibitem[Cheng et~al.(2017{\natexlab{a}})Cheng, Jia, Jiang, Toro, and
  Yu]{Cheng2017a}
Jun-Bo Cheng, Yueling Jia, Song Jiang, Eleuterio~F. Toro, and Ming Yu.
\newblock A second-order cell-centered {Lagrangian} method for two-dimensional
  elastic-plastic flows.
\newblock \emph{Communications in Computational Physics}, 22\penalty0
  (5):\penalty0 1224--1257, 2017{\natexlab{a}}.
\newblock \doi{10.4208/cicp.OA-2016-0173}.

\bibitem[Cheng et~al.(2017{\natexlab{b}})Cheng, Huang, Jiang, and
  Tian]{Cheng2017b}
Jun-Bo Cheng, Weizhang Huang, Song Jiang, and Baolin Tian.
\newblock A third-order moving mesh cell-centered scheme for one-dimensional
  elastic-plastic flows.
\newblock \emph{Journal of Computational Physics}, 349:\penalty0 137 -- 153,
  2017{\natexlab{b}}.
\newblock ISSN 0021-9991.
\newblock \doi{10.1016/j.jcp.2017.08.018}.
\newblock URL
  \url{http://www.sciencedirect.com/science/article/pii/S0021999117305892}.

\bibitem[Fridrich et~al.(2017)Fridrich, Liska, and Wendroff]{Fridrich2017}
David Fridrich, Richard Liska, and Burton Wendroff.
\newblock Cell-centered {Lagrangian} {Lax-Wendroff HLL} hybrid method for
  elasto-plastic flows.
\newblock \emph{Computers \& Fluids}, 157:\penalty0 164 -- 174, 2017.
\newblock ISSN 0045-7930.
\newblock \doi{10.1016/j.compfluid.2017.08.030}.
\newblock URL
  \url{http://www.sciencedirect.com/science/article/pii/S0045793017303080}.

\bibitem[Heuz{\'e}(2018)]{Heuze2018}
T.~Heuz{\'e}.
\newblock Simulation of impacts on elastic--viscoplastic solids with the
  flux-difference splitting finite volume method applied to non-uniform
  quadrilateral meshes.
\newblock \emph{Advanced Modeling and Simulation in Engineering Sciences},
  5\penalty0 (1):\penalty0 9, 2018.
\newblock ISSN 2213-7467.
\newblock \doi{10.1186/s40323-018-0101-z}.
\newblock URL \url{10.1186/s40323-018-0101-z}.

\bibitem[Aguirre et~al.(2015)Aguirre, Gil, Bonet, and Lee]{Aguirre2015}
M.~Aguirre, A.~J. Gil, J.~Bonet, and C.~H. Lee.
\newblock An upwind vertex centred finite volume solver for {Lagrangian} solid
  dynamics.
\newblock \emph{Journal of Computational Physics}, 300:\penalty0 387--422,
  2015.

\bibitem[Selim et~al.(June 2016)Selim, Koomullil, and Mcdaniel]{Selim2016}
M.~M. Selim, R.~P. Koomullil, and D.~R. Mcdaniel.
\newblock Linear elasticity finite volume based structural dynamics solver.
\newblock In \emph{{AIAA} Modeling and Simulation Technologies Conference},
  Washington, D.C., USA, June 2016.
\newblock \doi{10.2514/6.2016-4418}.

\bibitem[Selim et~al.(2017)Selim, Koomullil, and McDaniel]{Selim2017}
M.~M. Selim, R.~Koomullil, and D.~R. McDaniel.
\newblock Finite volume based fluid-structure interaction solver.
\newblock In \emph{58$^{th}$ AIAA/ASCE/AHS/ASC Structures, Structural Dynamics,
  and Materials Conference, Grapevine, Texas}, 2017.

\bibitem[Baliga and Patankar(1980)]{Baliga1980}
B.~R. Baliga and S.~V. Patankar.
\newblock A new finite-element formulation for convection-diffusion problems.
\newblock \emph{Numerical Heat Transfer, Part B: Fundamentals}, 3:4:\penalty0
  393--409, 1980.
\newblock \doi{10.1080/01495728008961767}.

\bibitem[Taylor et~al.(2003)Taylor, Bailey, and Cross]{Taylor2003}
G.~A. Taylor, C.~Bailey, and M.~Cross.
\newblock A vertex-based finite volume method applied to non-linear material
  problems in computational solid mechanics.
\newblock \emph{International Journal for Numerical Methods in Engineering},
  56:\penalty0 507--529, 2003.

\bibitem[Wheel(1996{\natexlab{a}})]{Wheel1996a}
M.~A. Wheel.
\newblock A geometrically versatile finite volume formulation for plane
  elastostatic stress analysis.
\newblock \emph{The Journal of Strain Analysis for Engineering Design},
  31:\penalty0 111--116, 1996{\natexlab{a}}.

\bibitem[Wheel(1996{\natexlab{b}})]{Wheel1996b}
M.~A. Wheel.
\newblock A finite-volume approach to the stress analysis of pressurised
  axisymmetric structures.
\newblock \emph{International Journal of Pressure Vessels and Piping},
  68:\penalty0 311--317, 1996{\natexlab{b}}.

\bibitem[Wheel(1998)]{Wheel1998}
M.~A. Wheel.
\newblock Applying the finite volume approach to structural analysis.
\newblock In P.~E.~O'Donoghue S.~N.~Atluri, editor, \emph{Modelling and
  Simulation in Engineering}, pages 229--234. 1, 1998.

\bibitem[Wheel(1999)]{Wheel1999a}
M.~A. Wheel.
\newblock A mixed finite volume formulation for determining the small strain
  deformation of incompressible materials.
\newblock \emph{International Journal for Numerical Methods in Engineering},
  44:\penalty0 1843--1861, 1999.

\bibitem[Costa and Sousa(2000)]{Costa2000}
V.~A.~F. Costa and A.~C.~M. Sousa.
\newblock A control volume-based fem for the solution of the three-dimensional,
  elastic stress-strain equations.
\newblock In \emph{WSEAS Proceedings 2$^{nd}$ World MCME (Mathematics \&
  Computers in Mechanical Engineering, MCME 2000)}, Vouliagmeni, Greece, 2000.

\bibitem[Fallah(2005)]{Fallah2005}
N.~Fallah.
\newblock \emph{Using shape function in cell centred finite volume formulation
  for two dimensional stress analysis, Lecture series on computer and
  computational sciences ({ICCMSE} 2005)}, volume~4.
\newblock Brill Academic Publishers, The Netherlands, 2005.

\bibitem[Xia et~al.(2007)Xia, Zhao, Yeo, and Lv]{Xia2007}
G.~H. Xia, Y.~Zhao, J.~H. Yeo, and X.~Lv.
\newblock A {3D} implicit unstructured-grid finite volume method for structural
  dynamics.
\newblock \emph{Computational Mechanics}, 40:\penalty0 299--312, 2007.

\bibitem[Tsui et~al.(2013)Tsui, Huang, Huang, and Lin]{Tsui2013}
Y.-Y. Tsui, Y.-C. Huang, C.-L. Huang, and S.-W. Lin.
\newblock A finite-volume-based approach for dynamic fluid-structure
  interaction.
\newblock \emph{Numerical Heat Transfer, Part B: Fundamentals}, 64\penalty0
  (4):\penalty0 326--349, 2013.
\newblock \doi{10.1080/10407790.2013.806691}.

\bibitem[Wu et~al.(2013)Wu, Xie, and Chen]{Wu2013}
Y.~Wu, X.~Xie, and L.~Chen.
\newblock Hybrid stress finite volume method for linear elasticity problems.
\newblock \emph{International Journal of Numerical Analysis and Modeling},
  10:\penalty0 634--656, 2013.

\bibitem[Taylor et~al.(1995)Taylor, Bailey, and Cross]{Taylor1995a}
G.~A. Taylor, C.~Bailey, and M.~Cross.
\newblock Solution of the elastic/visco-plastic constitutive equations: A
  finite volume approach.
\newblock \emph{Applied Mathematical Modelling}, 19:\penalty0 746--760, 1995.

\bibitem[Ferguson(1998)]{Ferguson1998}
W.~J. Ferguson.
\newblock The control volume finite element numerical solution technique
  applied to creep in softwoods.
\newblock \emph{International Journal of Solids and Structures}, 35:\penalty0
  1325--1338, 1998.

\bibitem[Hambleton et~al.(2011)Hambleton, Sloan, Pyatigorets, and
  Voller]{Hambleton2011}
J.P. Hambleton, S.W. Sloan, A.V. Pyatigorets, and V.R. Voller.
\newblock Lower bound limit analysis using the control volume finite element
  method.
\newblock In \emph{13th International Conference of the {IACMAG}}, pages
  88--93, Melbourne, Australia, 2011.

\bibitem[Fallah et~al.(1999)Fallah, Bailey, and Cross]{Fallah1999}
N.~Fallah, C.~Bailey, and M.~Cross.
\newblock Finite volume method for stress analysis.
\newblock In \emph{Proceedings {ASME} The 7$^{th}$ annual conference of the
  association for computational mechanics in engineering}, volume~99, pages
  135--138, Durham, UK, 1999.

\bibitem[Fallah et~al.(2000{\natexlab{a}})Fallah, Bailey, Cross, and
  Taylor]{Fallah2000a}
N.~Fallah, C.~Bailey, M.~Cross, and G.~A. Taylor.
\newblock Comparison of finite element and finite volume methods application in
  geometrically nonlinear stress analysis.
\newblock \emph{Applied Mathematical Modelling}, 24:\penalty0 439--455,
  2000{\natexlab{a}}.

\bibitem[Fallah et~al.(2000{\natexlab{b}})Fallah, Bailey, and
  Cross]{Fallah2000b}
N.~Fallah, C.~Bailey, and M.~Cross.
\newblock {CFD} approach for solid mechanics analysis.
\newblock In \emph{{European} Congress on Computational Methods in Applied
  Sciences and Engineering, {ECCOMAS} 2000}, Barcelona, Spain, 11-14 September
  2000{\natexlab{b}}.

\bibitem[Fallah(2000)]{Fallah2000c}
N.~Fallah.
\newblock \emph{Computational stress analysis using finite volume methods}.
\newblock PhD thesis, University of Greenwich, 2000.

\bibitem[Slone et~al.(2002{\natexlab{a}})Slone, Fallah, Bailey, and
  Cross]{Slone2002b}
A.~K. Slone, N.~Fallah, C.~Bailey, and M.~Cross.
\newblock A finite volume approach to geometrically nonlinear stress analysis.
\newblock In \emph{Third International Symposium on Finite Volumes for Complex
  Applications - Problems and Perspectives}, pages 663--670,
  2002{\natexlab{a}}.
\newblock ISBN 1-9039-9634-1.

\bibitem[Teran et~al.(2003)Teran, Blemker, Hing, and Fedkiw]{Teran2003}
J.~Teran, S.~Blemker, V.~N.~T. Hing, and R.~Fedkiw.
\newblock Finite volume methods for the simulation of skeletal muscle.
\newblock In \emph{Eurographics/{SIGGRAPH} Symposium on Computer Animation},
  2003.

\bibitem[Limache and Idelsohn(2007)]{Limache2007}
A.~C. Limache and S.~R. Idelsohn.
\newblock On the development of finite volume methods for computational solid
  mechanics.
\newblock \emph{Meci\'{a}nica Computacional}, 26:\penalty0 827--843, 2007.

\bibitem[Wheel(1997)]{Wheel1997}
M.~A. Wheel.
\newblock A finite volume method for analyzing the bending deformation of thick
  and thin plates.
\newblock \emph{Computer Methods in Applied Mechanics and Engineering},
  147:\penalty0 199--208, 1997.

\bibitem[Beveridge and Wheel(2009)]{Beveridge2009}
A.~J. Beveridge and M.~Wheel.
\newblock A control volume based formulation of the discrete {Kirchoff}
  triangular thin plate bending element.
\newblock In \emph{The 17$^{th}$ UK National Conference on Computational
  Mechanics in Engineering}, pages 287--290, Nottingham, UK, 2009.

\bibitem[Taylor(1996)]{Taylor1996}
G.~A. Taylor.
\newblock \emph{A vertex based discretisation scheme applied to material
  non-linearity within a multi-physics framework}.
\newblock PhD thesis, University of Greenwich, 1996.

\bibitem[Bailey et~al.(1996)Bailey, Chow, Cross, Fryer, and
  Pericleous]{Bailey1996}
C.~Bailey, P.~Chow, M.~Cross, Y.~Fryer, and K.~Pericleous.
\newblock Multiphysics modelling of metals casting process.
\newblock \emph{Proceedings of the Royal Society of London A: Mathematical,
  Physical and Engineering Sciences}, pages 459--486, 1996.

\bibitem[Bailey et~al.(1997)Bailey, Taylor, Bounds, Moran, and
  Cross]{Bailey1997}
C.~Bailey, G.~A. Taylor, S.~M. Bounds, G.~Moran, and M.~Cross.
\newblock {PHYSICA}: a multiphysics computational framework and its application
  to casting.
\newblock \emph{Mineral \& Metal Processing and Power Generation}, pages
  419--425, 1997.

\bibitem[Bailey et~al.(1999{\natexlab{a}})Bailey, Bounds, Cross, Moran,
  Pericleous, and Taylor]{Bailey1999a}
C.~Bailey, S.~Bounds, M.~Cross, G.~Moran, K.~Pericleous, and G.~A. Taylor.
\newblock Multiphysics modeling and its application to the casting process.
\newblock \emph{Computer Modelling and Simulation in Engineering}, 4:\penalty0
  206--212, 1999{\natexlab{a}}.

\bibitem[Bailey et~al.(1999{\natexlab{b}})Bailey, Taylor, Cross, and
  Chow]{Bailey1999b}
C.~Bailey, G.~A. Taylor, M.~Cross, and P.~Chow.
\newblock Discretisation procedures for for multi-physics phenomena.
\newblock \emph{Journal of Computational and Applied Mathematics},
  103:\penalty0 3--17, 1999{\natexlab{b}}.

\bibitem[Oldroyd et~al.(1999)Oldroyd, Wheel, and Scanlon]{Oldroyd1999}
A.~B. Oldroyd, M.~A. Wheel, and T.~J. Scanlon.
\newblock An integrated volume based approach for analysing flow induced
  vibrations.
\newblock In \emph{Proceedings European Conference on Computational Mechanics
  {(ECCM)}}, Munich, Germany, September 1999.

\bibitem[Wheel et~al.(1999)Wheel, Oldroyd, Scanlon, and Wenke]{Wheel1999b}
M.~A. Wheel, A.~Oldroyd, T.~J. Scanlon, and P.~Wenke.
\newblock Integrating finite volume based structural analysis procedures with
  {CFD} software to analyse fluid structure interaction.
\newblock In \emph{{Proceedings of 2$^{nd}$ Int. Conf. on Finite Volumes for
  Complex Applications}}, Duisburg, Germany, 1999.

\bibitem[Slone(2000)]{Slone2000}
A.~K. Slone.
\newblock \emph{A finite volume unstructured mesh approach to dynamic
  fluid-structure interactions between fluids and linear elastic solids}.
\newblock PhD thesis, University of Greenwich, 2000.

\bibitem[Slone et~al.(2001{\natexlab{a}})Slone, Pericleous, Bailey, and
  Cross]{Slone2001}
A.~K. Slone, K.~Pericleous, C.~Bailey, and M.~Cross.
\newblock Details of an integrated approach to three- dimensional dynamic fluid
  structure interaction.
\newblock In \emph{Fluid Structure Interaction}, pages 57--66. {WIT} Press,
  2001{\natexlab{a}}.
\newblock ISBN 1-85312-881-3.
\newblock \url{https://www.witpress.com/books/978-1-85312-881-3}.

\bibitem[Slone et~al.(2002{\natexlab{b}})Slone, Pericleous, Bailey, and
  Cross]{Slone2002}
A.~K. Slone, K.~Pericleous, C.~Bailey, and M.~Cross.
\newblock Dynamic fluid-structure interaction using finite volume unstructured
  mesh procedures.
\newblock \emph{Computers \& Structures}, 80:\penalty0 371--390,
  2002{\natexlab{b}}.

\bibitem[Slone et~al.(2003)Slone, Bailey, and Cross]{Slone2003}
A.~K. Slone, C.~Bailey, and M.~Cross.
\newblock Dynamic solid mechanics using finite volume methods.
\newblock \emph{Applied Mathamatical Modelling}, 27:\penalty0 69--87, 2003.

\bibitem[Slone et~al.(2004)Slone, Pericleous, Bailey, Cross, and
  Bennett]{Slone2004}
A.~K. Slone, K.~Pericleous, C.~Bailey, M.~Cross, and C.~Bennett.
\newblock A finite volume unstructured mesh approach to dynamic fluid-structure
  interaction: An assessment of the challenge of predicting the onset of
  flutter.
\newblock \emph{Applied Mathematical Modelling}, 28:\penalty0 211--239, 2004.

\bibitem[Slone et~al.(2007)Slone, Croft, Williams, and Cross]{Slone2007}
A.~K. Slone, T.~N. Croft, A.~J. Williams, and M.~Cross.
\newblock An alternative mixed eulerian lagrangian approach to high speed
  collision between solid structures on parallel clusters.
\newblock \emph{Advances in Engineering Software}, 38\penalty0 (4):\penalty0
  244 -- 255, 2007.
\newblock ISSN 0965-9978.
\newblock \doi{https://doi.org/10.1016/j.advengsoft.2006.09.015}.
\newblock URL
  \url{http://www.sciencedirect.com/science/article/pii/S0965997806001542}.

\bibitem[Cross et~al.(2007)Cross, Croft, Slone, Williams, Christakis, Patel,
  Bailey, and Pericleous]{Cross2007}
M.~Cross, T.~N. Croft, A.~K. Slone, A.~J. Williams, N.~Christakis, M.~K. Patel,
  C.~Bailey, and K.~Pericleous.
\newblock Computational modelling of multi-physics and multi-scale processes in
  parallel.
\newblock \emph{International Journal for Computational Methods in Engineering
  Science and Mechanics}, 8\penalty0 (2):\penalty0 63--74, 2007.
\newblock \doi{10.1080/15502280601149510}.
\newblock URL \url{https://doi.org/10.1080/15502280601149510}.

\bibitem[Lv et~al.(2007)Lv, Zhao, Huang, Xia, and Su]{Lv2007}
X.~Lv, Y.~Zhao, X.~Y. Huang, G.~H. Xia, and X.~H. Su.
\newblock A matrix-free implicit unstructured multigrid finite volume method
  for simulating structural dynamics and fluid-structure interaction.
\newblock \emph{Journal of Computational Physics}, 225:\penalty0 120--144,
  2007.

\bibitem[Croft et~al.(2008)Croft, Williams, Slone, and Cross]{Croft2008}
T.~N. Croft, A.~J. Williams, A.~K. Slone, and M.~Cross.
\newblock A two-dimensional prototype multi-physics model of the right
  ventricle of the heart.
\newblock \emph{International Journal for Numerical Methods in Fluids},
  57\penalty0 (5):\penalty0 583--600, 2008.
\newblock \doi{10.1002/fld.1764}.
\newblock URL \url{https://onlinelibrary.wiley.com/doi/abs/10.1002/fld.1764}.

\bibitem[Xia and l.~Lin(2008)]{Xia2008}
G.~H. Xia and C.~l.~Lin.
\newblock An unstructured finite volume approach for structural dynamics in
  response to fluid motions.
\newblock \emph{Computers \& Structures}, 86:\penalty0 684--701, 2008.

\bibitem[Hejranfar and Azampour(2016)]{Hejranfar2016}
K.~Hejranfar and M.~H. Azampour.
\newblock Simulation of {2D} fluid--structure interaction in inviscid
  compressible flows using a cell-vertex central difference finite volume
  method.
\newblock \emph{Journal of Fluids and Structures}, 67:\penalty0 190--218, 2016.

\bibitem[Taylor(1995)]{Taylor1995b}
G.~A. Taylor.
\newblock Material non-linearity within a finite volume framework for the
  simulation of a metal casting process.
\newblock \emph{{Computational Plasticity: Fundamentals and Applications II}},
  pages 1459--1470, 1995.

\bibitem[Taylor et~al.(1998)Taylor, Bailey, and Cross]{Taylor1998}
G.~A. Taylor, C.~Bailey, and M.~Cross.
\newblock A three dimensional finite volume approach to the thermomechanical
  modelling of the shape casting of metals.
\newblock In \emph{Proceedings of 8$^{th}$ Int. Conf. on Modelling of Casting,
  Welding and Advanced Solidification Processes}, 1998.

\bibitem[Cross(1996)]{Cross1996}
M.~Cross.
\newblock Computational issues in the modelling of materials-based
  manufacturing processes.
\newblock \emph{Journal of Computer-Aided Materials Design}, 3:\penalty0
  100--116, 1996.

\bibitem[Bounds et~al.(2000)Bounds, Moran, Pericleous, Cross, and
  Croft]{Bounds2000}
S.~Bounds, G.~Moran, K.~Pericleous, M.~Cross, and T.~N. Croft.
\newblock A computational model for defect prediction in shape castings based
  on the interaction of free surface flow, heat transfer, and solidification
  phenomena.
\newblock \emph{Metallurgical and Materials Transactions B}, 31:\penalty0
  515--527, 2000.

\bibitem[Williams et~al.(2001)Williams, Croft, and Cross]{Williams2001}
A.~J. Williams, T.~N. Croft, and M.~Cross.
\newblock Computational modelling of metals extrusion and forging processes.
\newblock In M.~Cross, J.~W. Ewans, and C.~Bailey, editors,
  \emph{{Computational Modelling of Materials, Minerals and Metals Processing,
  TMS}}, pages 481--490. Materials Society, 2001.

\bibitem[Williams et~al.(2002)Williams, Croft, and Cross]{Williams2002}
A.~J. Williams, T.~N. Croft, and M.~Cross.
\newblock Computational modelling of metal extrusion and forging processes.
\newblock \emph{Journal of Materials Processing Technology}, 125:\penalty0
  573--582, 2002.
\newblock \doi{10.1016/S0924-0136(02)00401-6}.
\newblock URL
  \url{http://www.sciencedirect.com/science/article/pii/S0924013602004016}.

\bibitem[Williams et~al.(2010)Williams, Slone, Croft, and Cross]{Williams2010}
A.~J. Williams, A.~K. Slone, T.~N. Croft, and M.~Cross.
\newblock {A mixed Eulerian-Lagrangian method for modelling metal extrusion
  processes}.
\newblock \emph{Computer Methods in Applied Mechanics and Engineering},
  199:\penalty0 2123--2134, 2010.

\bibitem[Taylor et~al.(2000{\natexlab{a}})Taylor, Hughes, Strusevich, and
  Pericleous]{Taylor2000a}
G.~A. Taylor, M.~Hughes, N.~Strusevich, and K.~Pericleous.
\newblock The application of three dimensional finite volume methods to the
  modelling of welding phenomena.
\newblock In Sahm, P.~N. Hansen, and J.~G. Conley, editors, \emph{Modeling of
  Casting, Welding and Advanced Solidification Processes {IX}},
  2000{\natexlab{a}}.

\bibitem[Taylor et~al.(2002{\natexlab{a}})Taylor, Hughes, Strusevich, and
  Pericleous]{Taylor2002a}
G.~A. Taylor, M.~Hughes, N.~Strusevich, and K.~Pericleous.
\newblock The application of three-dimensinal finite volume methods to the
  modelling of welding phenomena.
\newblock \emph{Applied Mathematical Modelling}, 26:\penalty0 309--320,
  2002{\natexlab{a}}.

\bibitem[Taylor et~al.(2002{\natexlab{b}})Taylor, Hughes, Strusevich, and
  Pericleous]{Taylor2002b}
G.~A. Taylor, M.~Hughes, N.~Strusevich, and K.~Pericleous.
\newblock Finite volume methods applied to the computational modelling of
  welding phenomena.
\newblock \emph{Applied Mathematical Modelling}, 26:\penalty0 309--320,
  2002{\natexlab{b}}.

\bibitem[Taylor et~al.(2000{\natexlab{b}})Taylor, Breiguine, Bailey, and
  Cross]{Taylor2000b}
G.~A. Taylor, V.~Breiguine, C.~Bailey, and M.~Cross.
\newblock An augmented {Lagrangian} contact algorithm employing a vertex-based
  finite volume method.
\newblock In \emph{{ACME}}, 2000{\natexlab{b}}.
\newblock Available at:
  \url{http://www.brunel.ac.uk/~eesrgat/research/ps_pubs}.

\bibitem[Gong et~al.(2013)Gong, Xuan, Ming, and Zhang]{Gong2013}
J.-F. Gong, L.-K. Xuan, P.-J. Ming, and W.-P. Zhang.
\newblock Thermoelastic analysis of functionally graded solids using a
  staggered finite volume method.
\newblock \emph{Composite Structures}, 104:\penalty0 134--143, 2013.

\bibitem[Gong et~al.(2014)Gong, Ming, Xuan, and Zhang]{Gong2014}
J.-F. Gong, P.-J. Ming, L.-K. Xuan, and W.-P. Zhang.
\newblock Thermoelastic analysis of three-dimensional functionally graded
  rotating disks based on finite volume method.
\newblock \emph{Journal of Mechanical Engineering Science}, 2014.

\bibitem[Perr\'{e} and Passard(1995)]{Perre1995}
P.~Perr\'{e} and J.~Passard.
\newblock A control-volume procedure compared with the finite-element method
  for calculating stress and strain during wood drying.
\newblock \emph{Drying Technology}, 13:\penalty0 635--660, 1995.

\bibitem[Salinas et~al.(2011)Salinas, Ch\'{a}vez, Gatica, and
  Ananias]{Salinas2011}
C.~Salinas, C.~Ch\'{a}vez, Y.~Gatica, and R.~Ananias.
\newblock Two-dimensional wood drying stress simulation using control
  control-volume mixed finite element methods ({CVFEM}).
\newblock \emph{Ingener\'{i}a e Investigacion}, 31:\penalty0 171--183, 2011.

\bibitem[Wheel(2008)]{Wheel2008}
M.~A. Wheel.
\newblock A control volume-based finite element method for plane micropolar
  elasticity.
\newblock \emph{International Journal for Numerical Methods in Engineering},
  75:\penalty0 992--1006, 2008.

\bibitem[Beveridge et~al.(2013)Beveridge, Wheel, and Nash]{Beveridge2013}
A.~J Beveridge, M.~Wheel, and D.~Nash.
\newblock A higher order control volume based finite element method to predict
  the deformation of heterogeneous materials.
\newblock \emph{Computers \& Structures}, 129:\penalty0 56--62, 2013.

\bibitem[Zhu et~al.(2012)Zhu, Ming, Xuan, and Zhang]{Zhu2012}
M.~Zhu, P.-J. Ming, L.~Xuan, and W.-P. Zhang.
\newblock An unstructured finite volume time domain method for structural
  dynamics.
\newblock \emph{Applied Mathematical Modelling}, 36:\penalty0 183--192, 2012.

\bibitem[Xuan et~al.(2014{\natexlab{a}})Xuan, Ming, Gong, Zheng, and
  Zhang]{Xuan2014a}
L.~Xuan, P.-J. Ming, J.~Gong, D.~Zheng, and W.-P. Zhang.
\newblock A finite volume time domain method for in-plane vibration on mixed
  grids.
\newblock \emph{Journal of Vibration and Acoustics}, 136, 2014{\natexlab{a}}.

\bibitem[Xuan et~al.(2014{\natexlab{b}})Xuan, Ming, Zhang, Jin, and
  Gong]{Xuan2014b}
L.~Xuan, P.-J. Ming, W.-P. Zhang, G.~Jin, and J.~Gong.
\newblock Time domain finite volume method for the transient response and
  natural characteristics of structural-acoustic coupling in an enclosed
  cavity.
\newblock \emph{Shengxue Xuebao/Acta Acustica}, 39:\penalty0 215--225,
  2014{\natexlab{b}}.

\bibitem[Wenke and Wheel(2003)]{Wenke2003}
P.~Wenke and M.~A. Wheel.
\newblock A finite volume method for solid mechanics incorporating rotational
  degrees of freedom.
\newblock \emph{Computers \& Structures}, 81:\penalty0 321--329, 2003.

\bibitem[Pan et~al.(2010)Pan, Wheel, and Qin]{Pan2010}
W.~Pan, M.~Wheel, and Y.~Qin.
\newblock Six-node triangle finite volume method for solids with a rotational
  degree of freedom for incompressible material.
\newblock \emph{Computers \& Structures}, 88:\penalty0 1506--1511, 2010.

\bibitem[Pan and Wheel(2011)]{Pan2011}
W.~Pan and M.~Wheel.
\newblock A finite-volume method for solids with a rotational degrees of
  freedom based on the 6-node triangle.
\newblock \emph{International Journal of Numerical Methods in Biomedical
  Engineering}, 27:\penalty0 1411--1426, 2011.

\bibitem[McManus et~al.(2000)McManus, Cross, Walshaw, Johnson, Leggett, and
  Scalable]{McManus2000}
K.~McManus, M.~Cross, C.~Walshaw, S.~Johnson, P.~Leggett, and A.~Scalable.
\newblock Strategy for the parallelization of multiphysics unstructured
  mesh-iterative codes on distributed-memory systems.
\newblock \emph{International Journal of High Performance Computing
  Applications}, 14:\penalty0 137--174, 2000.

\bibitem[McManus et~al.(2002)McManus, Cross, Walshaw, Croft, and
  Williams]{McManus2002}
K.~McManus, M.~Cross, C.~Walshaw, T.~N. Croft, and A.~J. Williams.
\newblock Parallel performance in multi-physics simulation.
\newblock In \emph{{ICCS '02 Proceedings of the International Conference on
  Computational Science-Part II}}, pages 806--815, London, 2002.
  Springer-Verlag.

\bibitem[Maitre et~al.(2002)Maitre, Rezgui, Souhail, and Zine]{Maitre2002}
J.~F. Maitre, A.~Rezgui, H.~Souhail, and A.~M. Zine.
\newblock High order finite volume schemes: Application to non-linear
  elasticity problems.
\newblock In \emph{Finite volumes for complex applications, {III}
  (Porquerolles), Hermes Sci. Publ., Paris}, pages 391--398, 2002.

\bibitem[Souhail(2004)]{Souhail2004}
H.~Souhail.
\newblock \emph{Sch\'{e}ma volumes finis: Estimation derreur a posteriori
  hi\'{e}rarchique par \'{e}l\'{e}ments finis mixtes. R\'{e}solution de
  probl\'{e}mes d'\'{e}lasticit\'{e} non-lin\'{e}aire}.
\newblock PhD thesis, Ecole Centrale de Lyon, 2004.

\bibitem[Zhang and Liu(1999)]{Zhang1999}
J.~Zhang and T.~Liu.
\newblock {P-SV} wave propagation in heterogeneous media: Grid method.
\newblock \emph{Geophysical Journal International}, 136:\penalty0 431--438,
  1999.

\bibitem[Zhang and Liu(2002)]{Zhang2002}
J.~Zhang and T.~Liu.
\newblock {Elastic wave modelling in 3-D heterogeneous media: 3-D grid method}.
\newblock \emph{Geophysical Journal International}, 150:\penalty0 780--799,
  2002.

\bibitem[Zhang(2004)]{Zhang2004}
J.~Zhang.
\newblock Wave propagation across fluid-solid interfaces: A grid method
  approach.
\newblock \emph{Geophysical Journal International}, 159:\penalty0 240--252,
  2004.

\bibitem[Liu et~al.(2004)Liu, Liu, and Zhang]{Liu2004}
T.~Liu, K.~Liu, and J.~Zhang.
\newblock Unstructured grid method for stress wave propagation in elastic
  media.
\newblock \emph{Computer Methods in Applied Mechanics and Engineering},
  193:\penalty0 2427--2452, 2004.

\bibitem[Liu et~al.(2005)Liu, Liu, and Zhang]{Liu2005}
T.~Liu, K.~Liu, and J.~Zhang.
\newblock Triangular grid method for stress-wave propagation in {2-D}
  orthotropic materials.
\newblock \emph{Archive of Applied Mechanics}, 74:\penalty0 477--488, 2005.

\bibitem[Gao and Zhang(2006)]{Gao2006}
H.~Gao and J.~Zhang.
\newblock Parallel {3-D} simulation of seismic wave propagation in
  heterogeneous anisotropic media: A grid method approach.
\newblock \emph{Geophysical Journal International}, 165:\penalty0 875--888,
  2006.

\bibitem[Dormy and Tarantola(1995)]{Dormy1995}
E.~Dormy and A.~Tarantola.
\newblock Numerical simulation of elastic wave propagation using a finite
  volume method.
\newblock \emph{Journal of geophysical research}, 100:\penalty0 2123--2133,
  1995.
\newblock \doi{10.1029/94JB02648}.

\bibitem[Harlow and Welch(1965)]{Harlow1965}
F.~H. Harlow and J.~E. Welch.
\newblock Numerical calculation of the time-dependent viscous incompressible
  flow of fluid with free surface.
\newblock \emph{Physics of Fluids}, 8, 1965.

\bibitem[Spalding(1993)]{Spalding1993}
D.~B. Spalding.
\newblock Simulation of fluid flow, heat transfer and solid deformation
  simultaneously.
\newblock In \emph{{NAFEMS} Conference No 4}, Brighton, UK, 1993.

\bibitem[Spalding(1997)]{Spalding1997}
D.~B. Spalding.
\newblock Simultaneous fluid-flow, heat-transfer and solid-stress computation
  in a single computer code.
\newblock In \emph{Keynote lecture 4$^{th}$ International Colloquium on Process
  Simulation}, Helsinki University of Technology, Espoo, Finland, 1997.

\bibitem[Spalding(1998{\natexlab{a}})]{Spalding1998a}
D.~B. Spalding.
\newblock Fluid-structure interaction in the presence of heat transfer and
  chemical reaction.
\newblock In \emph{ASME/JSME oint Pressure Vessels and Piping Conference}, San
  Diego, CA, USA, 1998{\natexlab{a}}.

\bibitem[Spalding(1998{\natexlab{b}})]{Spalding1998b}
D.~B. Spalding.
\newblock Fluid-structure interaction in the presence of heat transfer and
  chemical reaction.
\newblock In V.~Kudriavtsev and W.~Cheng, editors, \emph{V}. Computational
  Technologies for Fluid/Thermal/Structural/Chemical Systems With Industrial
  Applications, 1998{\natexlab{b}}.

\bibitem[Spalding(2002)]{Spalding2002}
D.~B. Spalding.
\newblock Simultaneous prediction of solid stress, heat transfer and fluid flow
  by a single algorithm.
\newblock In \emph{{ASME} Pressure Vessels and Piping Conference}, Vancouver,
  British Columbia, Canada, March 2002.

\bibitem[Spalding(2006)]{Spalding2004}
D.~B. Spalding.
\newblock Extending the boundaries of heat transfer.
\newblock In \emph{The 13th International Heat Transfer Conference, Sydney,
  Australia, August 16}, 2006.

\bibitem[Spalding(2008)]{Spalding2008}
D.~B. Spalding.
\newblock Enlarging the frontiers of computational fluid dynamics.
\newblock In \emph{International Symposium Heat \& Mass Transfer \&
  Hydrodynamics in Swirling Flow}, Moscow, Russia, 2008.

\bibitem[Patanker and Spalding(1972)]{Patankar1972}
S.~V. Patanker and D.~B. Spalding.
\newblock A calculation procedure for heat, mass, and momentum transfer in
  three-dimensional parabolic flows.
\newblock \emph{International Journal of Heat and Mass Transfer}, 15:\penalty0
  1787, 1972.

\bibitem[Patankar(1980)]{Patankar1980}
S.~V. Patankar.
\newblock \emph{Numerical Hear Transfer and Fluid Flow}.
\newblock Hemisphere Publishing Corporation, McGraw-Hill Book Company,
  Washington, 1980.

\bibitem[Hattel et~al.(1993{\natexlab{b}})Hattel, Hansen, and
  Hansen]{Hattel1993b}
J.~H. Hattel, P.~Hansen, and L.~F. Hansen.
\newblock Analysis of thermally induced stresses in die casting using a novel
  control volume technique.
\newblock In T.~Piwonka, editor, \emph{Modelling of Casting and Welding and
  Advanced Solidification Processes Advanced Solidification Processes Minerals
  Advanced Solidification Processes Minerals Metals and Materials Society}.
  Modelling of Casting, Welding and Advanced Solidification Processes,
  Minerals, Metals \& Materials Society, {TMS}, 1993{\natexlab{b}}.

\bibitem[Hattel et~al.(1993{\natexlab{c}})Hattel, Hansen, and
  Andersen]{Hattel1993c}
J.~H. Hattel, P.~N. Hansen, and S.~Andersen.
\newblock Modeling of thermal induced stresses in high pressure die casting
  dies.
\newblock In \emph{{NADCA} 19$^{th}$ International die casting congress,
  {NADC}, Transactions}, Rosemant, IL, USA, 1993{\natexlab{c}}.

\bibitem[Hattel(1993{\natexlab{a}})]{Hattel1993d}
J.~H. Hattel.
\newblock Stress calculations using a control volume based finite difference
  method.
\newblock \emph{Revnedannelse og Brudmekanik, {DMS} Vintermderbog},
  1993{\natexlab{a}}.

\bibitem[Hattel(1993{\natexlab{b}})]{Hattel1993e}
J.~H. Hattel.
\newblock \emph{Control volume based finite difference method - Numerical
  modeling of thermal and mechanical conditions in casting and heat treatment
  processes}.
\newblock PhD thesis, Institute of Manufacturing Engineering, Technical
  University of Denmark, 1993{\natexlab{b}}.

\bibitem[Hattel and Hansen(1994)]{Hattel1994}
J.~H. Hattel and P.~N. Hansen.
\newblock {1-D} analytical model for the thermally induced stresses in the mold
  surface during die casting.
\newblock \emph{Applied Mathematical Modeling}, 18:\penalty0 550--559, 1994.

\bibitem[Hattel and Hansen(1995)]{Hattel1995}
J.~H. Hattel and P.~N. Hansen.
\newblock A control volume-based finite difference method for solving the
  equilibrium equations in terms of displacements.
\newblock \emph{Applied Mathematical Modelling}, 19:\penalty0 210--243, 1995.

\bibitem[Hattel(1997)]{Hattel1997}
J.~H. Hattel.
\newblock Numerical modelling of stresses and deformations in casting
  processes.
\newblock In \emph{Proceedings of CASTING 1997}, Helsinki University of
  Technology, 1997. {Int. ADI and Simulation Conf. Helsinki}.

\bibitem[Pryds and Hattel(1997)]{Pryds1997}
N.~Pryds and J.~H. Hattel.
\newblock Numerical modelling of rapid solidification.
\newblock \emph{Modelling and Simulation in Materials Science and Engineering},
  5\penalty0 (5):\penalty0 451--472, 1997.

\bibitem[Hattel et~al.(1998)Hattel, Thorborg, and Andersen]{Hattel1998}
J.~H. Hattel, J.~Thorborg, and S.~Andersen.
\newblock Stress/strain modelling of casting processes in the framework of the
  control-volume method.
\newblock In \emph{{Modeling of Casting and Advanced Solidification Processes
  VIII, Warrendale, USA : TMS, The Minerals, Metals and Materials Society}},
  pages 763--770, 1998.

\bibitem[Hattel and Pryds(2001)]{Hattel2001}
J.~H. Hattel and N.~Pryds.
\newblock Modelling rapid solidification with the control volume method.
\newblock In \emph{A.R. Dinesen, M. Eldrup, D. Juul Jensen, S. Linderoth, T.B.
  Pedersen, N.H. Pryds, A. Schrder Pedersen, J.A. Wert, Proceedings of Science
  of Metastable and Nanocrystalline Alloys - Structure, Properties and
  Modelling, Roskilde, Ris National Laboratory}, pages 241--247, 2001.

\bibitem[Thorborg(2001)]{Thorborg2001}
J.~Thorborg.
\newblock \emph{Nonlinear constitutive modelling in thermomechanical processes
  with the control volume method}.
\newblock PhD thesis, Department of Manufacturing Engineering, Technical
  University of Denmark, 2001.

\bibitem[Hattel and Thorborg(2003)]{Hattel2003}
J.~H. Hattel and J.~Thorborg.
\newblock A numerical model for predicting the thermomechanical conditions
  during hydration of early-age concrete.
\newblock \emph{Applied Mathematical Modelling}, 27:\penalty0 1--26, 2003.

\bibitem[Thorborg and Hattel(2003)]{Thorborg2003}
J.~Thorborg and J.~H. Hattel.
\newblock Thermo-elasto-plasticity in solidification processes using the
  control volume method on staggered grid.
\newblock In Stefanescu et. al., editor, \emph{{Modelling of Casting, Welding
  and Advanced Solidification Processes, Warrendale: TMS - The Minerals, Metals
  \& Materials Society}}, 2003.

\bibitem[Wang and Melnik(2007)]{Wang2007}
L.~Wang and R.~Melnik.
\newblock Finite volume analysis of nonlinear thermo-mechanical dynamics of
  shape memory alloys.
\newblock \emph{Heat and Mass Transfer}, 43, 2007.

\bibitem[Rajagopal et~al.(2014)Rajagopal, Srinivasa, and
  Ponnalagu]{Rajagopal2014}
K.~R. Rajagopal, A.~R. Srinivasa, and A.~Ponnalagu.
\newblock Thermo-inelastic response of polymeric solids.
\newblock Technical report, Final Report, Texas Engineering Experiment Station,
  Harvey Mitchell Parkway South, Suite 300, College Station, TX, 2014.

\bibitem[Aboudi et~al.(1999)Aboudi, Pindera, and Arnold]{Aboudi1999}
J.~Aboudi, M.-J. Pindera, and S.M. Arnold.
\newblock Higher-order theory for functionally graded materials.
\newblock \emph{Composites Part B: Engineering}, 30:\penalty0 777--832, 1999.

\bibitem[Aboudi(2001)]{Aboudi2001a}
J.~Aboudi.
\newblock Micromechanical analysis of fully coupled
  electro-magneto-thermo-elastic multiphase composites.
\newblock \emph{Smart Materials and Structures}, 10:\penalty0 867--877, 2001.

\bibitem[Aboudi et~al.(2001)Aboudi, Pindera, and Arnold]{Aboudi2001b}
J.~Aboudi, M.-J. Pindera, and S.M. Arnold.
\newblock Linear thermoelastic higher-order theory for periodic multiphase
  materials.
\newblock \emph{{ASME} Journal of Applied Mechanics}, 68:\penalty0 697--707,
  2001.

\bibitem[Haj-Ali and Aboudi(2009)]{Haj-Ali2009}
R.~Haj-Ali and J.~Aboudi.
\newblock Nonlinear micromechanical formulation of the high fidelity
  generalized method of cells.
\newblock \emph{International Journal of Solids and Structures}, 46:\penalty0
  2577--2592, 2009.

\bibitem[Haj-Ali and Aboudi(2012)]{Haj-Ali2012}
R.~Haj-Ali and J.~Aboudi.
\newblock Discussion paper: Has renaming the high fidelity generalized method
  of cells been justified?
\newblock \emph{International Journal of Solids and Structures}, 49:\penalty0
  2051--2058, 2012.

\bibitem[Bansal and Pindera(2005)]{Bansal2005}
Y.~Bansal and M.-J. Pindera.
\newblock A second look at the higher-order theory for periodic multiphase
  materials.
\newblock \emph{{ASME} Journal of Applied Mechanics}, 72:\penalty0 177--195,
  2005.

\bibitem[Bansal and Pindera(2006)]{Bansal2006}
Y.~Bansal and M.-J. Pindera.
\newblock Finite-volume direct averaging micromechanics of heterogeneous
  materials with elastic-plastic phases.
\newblock \emph{International Journal of Plasticity}, 22:\penalty0 775--825,
  2006.

\bibitem[Aboudi(1982)]{Aboudi1982}
J.~Aboudi.
\newblock A continuum theory for fiber-reinforced elastic viscoplastic
  composites.
\newblock \emph{International Journal of Engineering Science}, 20:\penalty0
  605--621, 1982.

\bibitem[Aboudi(1991)]{Aboudi1991}
J.~Aboudi.
\newblock \emph{Mechanics of Composite Materials: A Unified Micromechanical
  Approach}.
\newblock Elsevier, Amsterdam, Netherlands, 1991.

\bibitem[Paley and Aboudi(1992)]{Paley1992}
M.~Paley and J.~Aboudi.
\newblock Micromechanical analysis of composites by the generalized cells
  model.
\newblock \emph{Mechanics of Materials}, 14:\penalty0 127--139, 1992.

\bibitem[Cavalcante and Pindera(2016)]{Cavalcante2016}
M.~A.~A. Cavalcante and M.-J. Pindera.
\newblock Generalized {FVDAM} theory for elastic-plastic periodic materials.
\newblock \emph{International Journal of Plasticity}, 77:\penalty0 90--117,
  2016.

\bibitem[Aboudi(2004)]{Aboudi2004}
J.~Aboudi.
\newblock The generalized method of cells and high-fidelity generalized method
  of cells micromechanical models: a review.
\newblock \emph{Mechanics of Advanced Materials and Structures}, 11:\penalty0
  329--366, 2004.

\bibitem[Pindera et~al.(2009)Pindera, Khatam, Drago, and Bansal]{Pindera2009}
M.-J. Pindera, H.~Khatam, A.~S. Drago, and Y.~Bansal.
\newblock Micromechanics of spatially uniform heterogeneous media: A critical
  review and emerging approaches.
\newblock \emph{Composites Part B: Engineering}, 40:\penalty0 349--378, 2009.

\bibitem[Charalambakis and Murat(2006)]{Charalambakis2006}
N.~Charalambakis and F.~Murat.
\newblock Homogenization of stratified thermoviscoplastic materials.
\newblock \emph{Quarterly of Applied Mathematics}, 64(2):\penalty0 359--99,
  2006.

\bibitem[Aboudi et~al.(2007)Aboudi, Arnold, and Bednarcyk]{Aboudi2007}
J.~Aboudi, S.~M. Arnold, and B.~A. Bednarcyk.
\newblock \emph{Micromechanical analyses of smart composite materials}.
\newblock Nova Science Publishers, New York, USA, 2007.

\bibitem[Aboudi(2008)]{Aboudi2008}
J.~Aboudi.
\newblock Finite strain micromechanical modeling of multiphase composites.
\newblock \emph{International Journal for Multiscale Computational
  Engineering}, 6:\penalty0 411--434, 2008.

\bibitem[Atluri and Shen(2002)]{Atluri2002}
S.~N. Atluri and S.~Shen.
\newblock The meshless local {Petrov-Galerkin} {(MLPG)} method: a simple \&
  less-costly alternative to the finite element and boundary element methods.
\newblock \emph{Computer Modeling in Engineering and Sciences}, 3:\penalty0
  11--51, 2002.

\bibitem[Atluri and Zhu(1998)]{Atluri1998}
S.~N. Atluri and T.~Zhu.
\newblock A new meshless local {Petrov-Galerkin} {(MLPG)} approach in
  computational mechanics.
\newblock \emph{Computational Mechanics}, 22:\penalty0 117--127, 1998.
\newblock \doi{10.1007/s004660050346}.

\bibitem[Ching and Batra(2001)]{Ching2001}
H.-K. Ching and R.~C. Batra.
\newblock {Determination of Crack Tip Fields in Linear Elastostatics by the
  Meshless Local Petrov-Galerkin (MLPG) Method}.
\newblock \emph{Computer Modeling in Engineering and Sciences}, 2:\penalty0
  273--289, 2001.

\bibitem[Warlock et~al.(2002)Warlock, Ching, Kapila, and Batra]{Warlock2002}
A.~Warlock, H.-K. Ching, A.~K. Kapila, and R.~C. Batra.
\newblock Plane strain deformations of an elastic material compressed in a
  rough rectangular cavity.
\newblock \emph{International Journal of Engineering Science}, 40:\penalty0
  991--1010, 2002.

\bibitem[Qian et~al.(2003)Qian, Batra, and Chen]{Qian2003}
L.~Qian, R.~Batra, and L.~Chen.
\newblock Elastostatic deformations of a thick plate by using a high- er-order
  shear and normal deformable plate theory and two meshless local
  {Petrov-Galerkin} {(MLPG)} methods.
\newblock \emph{Computer Modeling in Engineering and Sciences}, 4\penalty0
  (1):\penalty0 161--76, 2003.

\bibitem[Raju and Phillips(2003)]{Raju2003}
I.~Raju and D.~Phillips.
\newblock Further developments in the {MLPG} method for beam problems.
\newblock \emph{Computer Modeling in Engineering and Sciences}, 4\penalty0
  (1):\penalty0 141--60, 2003.

\bibitem[Atluri et~al.(2004)Atluri, Han, and Rajendran]{Atluri2004}
S.~N. Atluri, Z.~Han, and A.~Rajendran.
\newblock A new implementation of the meshless finite volume method, through
  the {MLPG} ``mixed'' approach.
\newblock \emph{Computer Modeling in Engineering and Sciences}, 6\penalty0
  (6):\penalty0 491--514, 2004.

\bibitem[Han and Atluri(2004)]{Han2004}
Z.~Han and S.~N. Atluri.
\newblock Meshless local {Petrov-Galerkin} {(MLPG)} approaches for solving {3D}
  problems in elasto-statics.
\newblock \emph{Computer Modeling in Engineering and Sciences}, 6:\penalty0
  169--88, 2004.

\bibitem[Batra et~al.(2004)Batra, Porfiri, and Spinello]{Batra2004}
R.~C. Batra, M.~Porfiri, and D.~Spinello.
\newblock Treatment of material discontinuity in two meshless local
  {Petrov-Galerkin} {(MLPG)} formulations of axisymmetric transient heat
  conduction.
\newblock \emph{International Journal for Numerical Methods in Engineering},
  61:\penalty0 2461--2479, 2004.

\bibitem[Han et~al.(2005)Han, Rajendran, and Atluri]{Han2005}
Z.~Han, A.~Rajendran, and S.~N. Atluri.
\newblock Meshless local {Petrov-Galerkin} {(MLPG)} approaches for solving
  nonlinear problems with large deformations and rotations.
\newblock \emph{Computer Modeling in Engineering and Sciences}, 10\penalty0
  (1):\penalty0 1, 2005.

\bibitem[Sladek et~al.(2008)Sladek, Sladek, Solek, and Saez]{Sladek2008}
Jan Sladek, Vladimir Sladek, Peter Solek, and Andres Saez.
\newblock Dynamic {3D} axisymmetric problems in continuously non-homogeneous
  piezoelectric solids.
\newblock \emph{International Journal of Solids and Structures}, 45\penalty0
  (16):\penalty0 4523 -- 4542, 2008.
\newblock ISSN 0020-7683.
\newblock \doi{10.1016/j.ijsolstr.2008.03.027}.
\newblock URL
  \url{http://www.sciencedirect.com/science/article/pii/S0020768308001376}.

\bibitem[Moosavi and Khelil(2008)]{Moosavi2008}
M.~R. Moosavi and A.~Khelil.
\newblock Accuracy and computational efficiency of the finite volume method
  combined with the meshless local {Petrov-Galerkin} in comparison with the
  finite element method in elasto-static problem.
\newblock \emph{{ICCES}}, 5:\penalty0 211--38, 2008.

\bibitem[Moosavi and Khelil(2009)]{Moosavi2009}
M.~R. Moosavi and A.~Khelil.
\newblock Finite volume meshless local {Petrov-Galerkin} method in
  elastodynamic problems.
\newblock \emph{Engineering Analysis with Boundary Elements}, 33:\penalty0
  1016--1021, 2009.

\bibitem[Moosavi et~al.(2011{\natexlab{a}})Moosavi, Delfanian, and
  Khelil]{Moosavi2011a}
M.~R. Moosavi, F.~Delfanian, and A.~Khelil.
\newblock Orthogonal meshless finite volume method in elasticity.
\newblock \emph{Thin-Walled Structures}, 49:\penalty0 708--712,
  2011{\natexlab{a}}.

\bibitem[Moosavi et~al.(2011{\natexlab{b}})Moosavi, Delfanian, and
  Khelil]{Moosavi2011b}
M.~R. Moosavi, F.~Delfanian, and A.~Khelil.
\newblock The orthogonal meshless finite volume method for solving
  {Euler-Bernoulli} beam and thin plate problems.
\newblock \emph{Thin-Walled Structures}, 49:\penalty0 923--932,
  2011{\natexlab{b}}.

\bibitem[Hosseini et~al.(2011)Hosseini, Sladek, and Sladek]{Hosseini2011}
Seyed~Mahmoud Hosseini, Jan Sladek, and Vladimir Sladek.
\newblock Meshless local {Petrov-Galerkin} method for coupled thermoelasticity
  analysis of a functionally graded thick hollow cylinder.
\newblock \emph{Engineering Analysis with Boundary Elements}, 35\penalty0
  (6):\penalty0 827 -- 835, 2011.
\newblock ISSN 0955-7997.
\newblock \doi{10.1016/j.enganabound.2011.02.001}.
\newblock URL
  \url{http://www.sciencedirect.com/science/article/pii/S0955799711000233}.

\bibitem[Soares et~al.(2012)Soares, Sladek, and Sladek]{Soares2012}
D.~Soares, V.~Sladek, and J.~Sladek.
\newblock Modified meshless local {Petrov-Galerkin} formulations for
  elastodynamics.
\newblock \emph{International Journal for Numerical Methods in Engineering},
  90\penalty0 (12):\penalty0 1508--1828, 2012.
\newblock \doi{10.1002/nme.3373}.
\newblock URL \url{https://onlinelibrary.wiley.com/doi/abs/10.1002/nme.3373}.

\bibitem[Moosavi et~al.(2012{\natexlab{a}})Moosavi, Delfanian, and
  Khelil]{Moosavi2012a}
M.~R. Moosavi, F.~Delfanian, and A.~Khelil.
\newblock Orthogonal meshless finite volume method applied to crack problems.
\newblock \emph{Thin-Walled Structures}, 52:\penalty0 61--65,
  2012{\natexlab{a}}.

\bibitem[Moosavi et~al.(2012{\natexlab{b}})Moosavi, Delfanian, and
  Khelil]{Moosavi2012b}
M.~Moosavi, F.~Delfanian, and A.~Khelil.
\newblock Orthogonal meshless finite volume method in shell analysis.
\newblock \emph{Finite Elements in Analysis and Design}, 62:\penalty0 1--7,
  2012{\natexlab{b}}.

\bibitem[Moosavi(2013)]{Moosavi2013}
M.~R. Moosavi.
\newblock Orthogonal meshless finite volume method applied to elastodynamic
  crack problems.
\newblock \emph{International Journal of Fracture}, 179:\penalty0 1--7, 2013.
\newblock \doi{10.1007/s10704-012-9752-9}.

\bibitem[Ebrahimnejad et~al.(2015)Ebrahimnejad, Fallah, and
  Khoei]{Ebrahimnejad2015}
M.~Ebrahimnejad, N.~Fallah, and A.~R. Khoei.
\newblock Adaptive refinement in the meshless finite volume method for
  elasticity problems.
\newblock \emph{Computers \& Mathematics with Applications}, 69, 2015.

\bibitem[Ebrahimnejad et~al.(2017)Ebrahimnejad, Fallah, and
  Khoei]{Ebrahimnejad2017}
M.~Ebrahimnejad, N.~Fallah, and A.R. Khoei.
\newblock Three types of meshless finite volume method for the analysis of
  two-dimensional elasticity problems.
\newblock \emph{Computational and Applied Mathematics}, 36\penalty0
  (2):\penalty0 971--990, 2017.

\bibitem[Fallah(2018)]{Fallah2018c}
N.~Fallah.
\newblock Mesh-free and mesh based finite volume methods for solid mechanics
  analysis.
\newblock In \emph{41$^{st}$ Solid Mechanics Conference ({SOLMECH} 2018)},
  Warsaw, Poland, 2018.

\bibitem[Davoudi-Kia and Fallah(2017)]{Davoudi-Kia2017}
A.~Davoudi-Kia and N.~Fallah.
\newblock Comparison of enriched meshless finite-volume and element-free
  {Galerkin} methods for the analysis of heterogeneous media.
\newblock \emph{Engineering with Computers}, Dec 2017.
\newblock ISSN 1435-5663.
\newblock \doi{10.1007/s00366-017-0573-3}.
\newblock URL \url{10.1007/s00366-017-0573-3}.

\bibitem[Davoudi-Kia and Fallah(April 2018)]{Davoudi-Kia2018}
A.~Davoudi-Kia and N.~Fallah.
\newblock An enriched meshless finite volume method for the modeling of
  material discontinuity problems in {2D} elasticity.
\newblock \emph{Latin American Journal of Solids and Structures}, 15\penalty0
  (2):\penalty0 209--219, April 2018.
\newblock \doi{10.1590/1679-78254121}.

\bibitem[Fallah and Nikraftar(2018)]{Fallah2018d}
N.~Fallah and N.~Nikraftar.
\newblock Meshless finite volume method for the analysis of fracture problems
  in orthotropic media.
\newblock \emph{Engineering Fracture Mechanics}, 2018.
\newblock \doi{10.1016/j.engfracmech.2018.09.029}.
\newblock In press.

\bibitem[Ba\v{s}i\'{c} et~al.(2001)Ba\v{s}i\'{c}, Demird{\v{z}}i\'{c}, and
  Muzaferija]{Basic2001}
H.~Ba\v{s}i\'{c}, I.~Demird{\v{z}}i\'{c}, and S.~Muzaferija.
\newblock Analysis of plastic flow of metals during extrusion processes using
  finite volume method.
\newblock In \emph{Proceedings of 3$^{rd}$ International Conference on
  Industrial Tools}, pages 22--26, Slovenia, 2001. Rogaska Slatina.

\bibitem[Ba\v{s}i\'{c}(2008)]{Basic2008}
H.~Ba\v{s}i\'{c}.
\newblock Friction models comparison in finite volume method simulation of bulk
  metal forming technologies.
\newblock \emph{Journal for Technology of Plasticity}, 33:\penalty0 113--122,
  2008.

\bibitem[Ba\v{s}i\'{c}(2009)]{Basic2009}
H.~Ba\v{s}i\'{c}.
\newblock The constitutive models in numerical simulation of steady-state metal
  forming processes.
\newblock \emph{Journal for Technology of Plasticity}, 34:\penalty0 27--36,
  2009.

\bibitem[Chen et~al.(2007)Chen, Lou, and Ruan]{Chen2007}
Z.Z. Chen, Z.L. Lou, and X.Y. Ruan.
\newblock Finite volume simulation and mould optimization of aluminum profile
  extrusion.
\newblock \emph{Journal of Materials Processing Technology}, 190\penalty0
  (1-3):\penalty0 382 -- 386, 2007.
\newblock \doi{10.1016/j.jmatprotec.2007.01.032}.

\bibitem[Jafari et~al.(2007)Jafari, Zebarjad, and Kolahan]{Jafari2007}
M.~R. Jafari, S.~M. Zebarjad, and F.~Kolahan.
\newblock Simulation of thixoformability of {A356} aluminum alloy using finite
  volume method.
\newblock \emph{Materials Science and Engineering A}, 454:\penalty0 558--563,
  2007.
\newblock \doi{10.1016/j.msea.2006.11.124}.

\bibitem[Lou et~al.(2008)Lou, Zhao, Wang, and Wu]{Lou2008}
S.~Lou, G.~Zhao, R.~Wang, and X.~Wu.
\newblock Modeling of aluminum alloy profile extrusion process using finite
  volume method.
\newblock \emph{Journal of Materials Processing Technology}, 206:\penalty0
  481--490, 2008.

\bibitem[Al-Athel and Gadala(2011)]{Al-Athel2011}
K.S. Al-Athel and M.S. Gadala.
\newblock {Eulerian} volume of solid ({VOS}) approach in solid mechanics and
  metal forming.
\newblock \emph{Computer methods in applied mechanics and engineering},
  200:\penalty0 2145--2159, 2011.

\bibitem[Wang and Li(2011)]{Wang2011}
R.~Wang and H.~Z. Li.
\newblock Modeling of aluminum extrusion process using non-orthogonal block
  structured grids based {FVM}.
\newblock \emph{Advanced Materials Research}, 189:\penalty0 1749--1752, 2011.

\bibitem[Wang(2012)]{Wang2012}
R.~Wang.
\newblock Body fitted grids based {FVM} simulation of aluminum extrusion
  process.
\newblock \emph{Advanced Materials Research}, 418:\penalty0 2102--2105, 2012.

\bibitem[Bressan et~al.(2013)Bressan, Martins, and Button]{Bressan2013}
J.~D. Bressan, M.~M. Martins, and S.~T. Button.
\newblock Aluminium extrusion analysis by the finite volume method.
\newblock In O\~{n}ate, D.~R.~J. Owen, D.~Peric, and B.~Su{\'a}rez, editors,
  \emph{XII International Conference on Computational Plasticity. Fundamentals
  and Applications COMPLAS XII}, 2013.

\bibitem[Zhang et~al.(2016)Zhang, Chen, Zhao, Zhang, and Lou]{Zhang2016}
C.~Zhang, H.~Chen, G.~Zhao, L.~Zhang, and S.~Lou.
\newblock Optimization of porthole extrusion dies with the developed algorithm
  based on finite volume method.
\newblock \emph{The International Journal of Advanced Manufacturing
  Technology}, pages 1--13, 2016.

\bibitem[de~Brauer et~al.(2016)de~Brauer, Iollo, and Milcent]{DeBrauer2016}
A.~de~Brauer, A.~Iollo, and T.~Milcent.
\newblock A {Cartesian} scheme for compressible multimaterial models in {3D}.
\newblock \emph{Journal of Computational Physics}, 313:\penalty0 121 -- 143,
  2016.
\newblock ISSN 0021-9991.
\newblock \doi{10.1016/j.jcp.2016.02.032}.
\newblock URL
  \url{http://www.sciencedirect.com/science/article/pii/S0021999116000966}.

\bibitem[de~Brauer et~al.(2017)de~Brauer, Iollo, and Milcent]{DeBrauer2017}
A.~de~Brauer, A.~Iollo, and T.~Milcent.
\newblock A {Cartesian} scheme for compressible multimaterial hyperelastic
  models with plasticity.
\newblock \emph{Communications in Computational Physics}, 22\penalty0
  (5):\penalty0 1362?1384, 2017.
\newblock \doi{10.4208/cicp.OA-2017-0018}.

\bibitem[Teng et~al.(1999)Teng, Chen, and Hu]{Teng1999}
J.~G. Teng, S.~F. Chen, and J.~L. Hu.
\newblock A finite volume method for deformation analysis of woven fabrics.
\newblock \emph{International Journal for Numerical Methods in Engineering},
  46:\penalty0 2061--2098, 1999.

\bibitem[Chen et~al.(2001)Chen, Hu, and Teng]{Chen2001}
S.~F Chen, J.~L Hu, and J.~G Teng.
\newblock A finite-volume method for contact drape simulation of woven fabrics
  and garments.
\newblock \emph{Finite Elements in Analysis and Design}, 37:\penalty0 513--531,
  2001.

\bibitem[Martin and Pascal(2011)]{Martin2011}
B.~Martin and F.~Pascal.
\newblock Discrete duality finite volume method applied to linear elasticity.
\newblock \emph{{Finite Volumes for Complex Applications VI: Problems \&
  Perspectives}}, 4:\penalty0 663--671, 2011.

\bibitem[Martin(2012)]{Martin2012}
B.~Martin.
\newblock \emph{Elaboration de solveurs volumes finis {2D}/{3D} pour
  r\'{e}soudre le probl\,{e}me de l'elasticit\'{e} lin\'{e}aire}.
\newblock PhD thesis, Ecole normale sup\'{e}rieure de Cachan - {ENS} Cachan,
  Francais, 2012.

\bibitem[Pietro et~al.(2011)Pietro, Eymard, Lemaire, and Masson]{Pietro2011}
D.~A.~Di Pietro, R.~Eymard, S.~Lemaire, and R.~Masson.
\newblock Hybrid finite volume discretization of linear elasticity models on
  general meshes.
\newblock In J.~F{\"{u}}rst, J.~Halama, R.~Herbin, and F.~Hubert, editors,
  \emph{{Finite Volumes for Complex Applications: VI Problems \& Perspectives,
  Volume 4 of Springer Proceedings in Mathematics, Springer Berlin
  Heidelberg}}, pages 331--339, 2011.

\bibitem[Wilkins(1963)]{Wilkins1963}
M.~L. Wilkins.
\newblock Calculation of elastic-plastic flow, {T3 - UCRL; 7322}.
\newblock Technical report, Lawrence Radiation Laboratory, Lawrence Livermore
  Laboratory, University of California, Berkeley, 1963.
\newblock \url{https://catalog.hathitrust.org/Record/007293160}.

\bibitem[Wilkins(1964)]{Wilkins1964}
M.~L. Wilkins.
\newblock Calculations of elastic-plastic flow.
\newblock \emph{Methods of Computational Physics}, 3, 1964.
\newblock Edited by B. Adler, S. Fernback, and M. Rotenberg.

\bibitem[Wilkins(1999)]{Wilkins1999}
M.~L. Wilkins.
\newblock \emph{Computer simulation of dynamic phenomena}.
\newblock Springer-Verlag, Berlin Heidelberg New York, 1999.

\bibitem[Bessonov et~al.(2009)Bessonov, Golovashchenko, and
  Volpert]{Bessonov2009}
N.~M. Bessonov, S.~F. Golovashchenko, and V.~A. Volpert.
\newblock Numerical modelling of contact elastic-plastic flows.
\newblock \emph{Mathematical Modelling of Natural Phenomena}, 4\penalty0
  (1):\penalty0 44--87, 2009.
\newblock URL \url{http://eudml.org/doc/222347}.

\bibitem[{Rhie} and {Chow}(1983)]{Rhie1983}
C.~M. {Rhie} and W.~L. {Chow}.
\newblock Numerical study of the turbulent flow past an airfoil with trailing
  edge separation.
\newblock \emph{{AIAA} Journal}, 21:\penalty0 1525--1532, November 1983.
\newblock \doi{10.2514/3.8284}.

\bibitem[Jameson et~al.(1981)Jameson, Schmidt, and Turkel]{Jameson1981}
A.~Jameson, W.~Schmidt, and E.~Turkel.
\newblock Numerical solution of the {Euler} equations by finite volume methods
  using {Runge-Kutta} time-stepping schemes.
\newblock In \emph{{AIAA} 5$^{th}$ Computational Fluid Dynamics Conference},
  volume~81, page 1259, 1981.

\bibitem[Jacobs(1980)]{Jacobs1980}
D.~A.~H. Jacobs.
\newblock Preconditioned conjugate gradient methods for solving systems of
  algebraic equations.
\newblock Technical report, Central Electricity Research Laboratories Report
  {(RD/L/N193/80)}, 1980.

\bibitem[Hassan(2019)]{Hassan2019}
O.~I.~I. Hassan.
\newblock \emph{A vertex centred Finite Volume algorithm for fast solid
  dynamics: Total and Updated Lagrangian descriptions}.
\newblock PhD thesis, University of Swansea, 2019.

\bibitem[Belytschko et~al.(2014)Belytschko, Liu, Moran, and
  Elkhodary]{Belytschko2014}
T.~Belytschko, W.~K. Liu, B.~Moran, and K.~I. Elkhodary.
\newblock \emph{Nonlinear Finite Elements for Continua and Structures}.
\newblock Wiley, Chichester, West Sussex, 2$^{nd}$ edition, 2014.

\bibitem[Bathe(1996)]{Bathe1996}
K.~J. Bathe.
\newblock \emph{Finite element procedures}.
\newblock Prentice Hall, New Jersey, 1996.

\bibitem[Sch\"{a}fer(2006)]{Schafer2006}
M.~Sch\"{a}fer.
\newblock \emph{Computational Engineering - Introduction to Numerical Methods}.
\newblock Springer-Verlag, Berlin, Heidelberg, 2006.
\newblock ISBN 3540306854.

\bibitem[Zienkiewicz and Taylor(2000)]{Zienkiewicz2000}
O.~C. Zienkiewicz and R.~L. Taylor.
\newblock \emph{The finite element method, volume 2, solid mechanics}.
\newblock Butterworth Heinemann, Oxford, 5$^{th}$ edition, 2000.

\bibitem[Laursen(2002)]{Laursen2002}
T.~A. Laursen.
\newblock \emph{Computational Contact and Impact Mechanics}.
\newblock Springer-Verlag Berlin Heidelberg {GmBH}, 2002.

\bibitem[Corp.(2018)]{Abaqus}
Dassault Syst\'{e}mes~Simulia Corp.
\newblock Abaqus 6.14.
\newblock \url{http://www.simulia.com/products/abaqus_fea.html}, 2018.

\bibitem[O{\~{n}}ate and Cervera(1993)]{Onate1993}
E.~O{\~{n}}ate and M.~Cervera.
\newblock Derivation of thin plate bending elements with one degree of freedom
  per node.
\newblock \emph{Engineering Computations}, 10:\penalty0 543--561, 1993.

\bibitem[O{\~{n}}ate(1998)]{Onate1998}
E.~O{\~{n}}ate.
\newblock Elementos finitos y volumenes finitos puntos de encuentro y
  posibilidad de nuevas aplicaciones.
\newblock Technical report, {CIMNE, Barcelona}, 1998.

\bibitem[O{\~{n}}ate and Zarate(2000)]{Onate2000}
E.~O{\~{n}}ate and F.~Zarate.
\newblock Rotation-free triangular plate and shell elements.
\newblock \emph{International Journal for Numerical Methods in Engineering},
  47:\penalty0 557--603, 2000.

\bibitem[Zienkiewicz(1995)]{Zienkiewicz1995}
O.~C. Zienkiewicz.
\newblock Origins, milestones and directions of the finite element method: A
  personal view.
\newblock \emph{Archives of Computational Methods in Engineering State of the
  art reviews}, 2\penalty0 (1):\penalty0 1--48, 1995.

\bibitem[Lahrmann(1992)]{Lahrmann1992}
A.~Lahrmann.
\newblock An element formulation for the classical finite difference and finite
  volume method applied to arbitrarily shaped domains.
\newblock \emph{International Journal for Numerical Methods in Engineering},
  35:\penalty0 893--913, 1992.

\bibitem[Harrild and Henriquez(1997)]{Harrild1997}
D.~M. Harrild and C.~S. Henriquez.
\newblock A finite volume model of cardiac propagation.
\newblock \emph{Annals of Biomedical Engineering}, 25:\penalty0 315--334, 1997.

\bibitem[Fang et~al.(2002)Fang, Tsuchiya, and Yamamoto]{Fang2002}
Q.~Fang, T.~Tsuchiya, and T.~Yamamoto.
\newblock Finite difference, finite element and finite volume methods applied
  to two-point boundary value problems.
\newblock \emph{Journal of Computational and Applied Mathematics},
  139:\penalty0 9--19, 2002.

\bibitem[Yamamoto et~al.(2002)Yamamoto, Fang, Tsuchiya, and
  Difference]{Yamamoto2002}
T.~Yamamoto, Q.~Fang, T.~Tsuchiya, and Finite Difference.
\newblock Finite element and finite volume methods applied to two-point
  boundary value problems.
\newblock \emph{Journal of Computational and Applied Mathematics},
  139:\penalty0 9--19, 2002.

\bibitem[Jacquemet and Henriquez(2005)]{Jacquemet2005a}
V.~Jacquemet and C.~S. Henriquez.
\newblock {Finite Volume Stiffness Matrix for Solving Anisotropic Cardiac
  Propagation in 2-D and 3-D Unstructured Meshes}.
\newblock \emph{{IEEE} Transactions on Biomedical Engineering}, 52:\penalty0
  1490--1492, 2005.

\bibitem[Jacquemet(2005)]{Jacquemet2005b}
V.~Jacquemet.
\newblock {Link between the FEM and FVM formulations of anisotropic cardiac
  propagation in unstructured meshes}.
\newblock Technical report, {ITS Technical Report, TR-ITS 2005.021}, 2005.

\bibitem[Filippini et~al.(2014)Filippini, Maliska, and {Vaz
  Jr.}]{Filippini2014}
G.~Filippini, C.~R. Maliska, and M.~{Vaz Jr.}
\newblock A physical perspective of the element-based finite volume method and
  {FEM-Galerkin} methods within the framework of the space of finite elements.
\newblock \emph{International Journal for Numerical Methods in Engineering},
  98:\penalty0 24--43, 2014.

\bibitem[Demird\v{z}i\'{c}(2020)]{Demirdzic2020}
I.~Demird\v{z}i\'{c}.
\newblock Finite volumes vs finite elements. there is a choice.
\newblock \emph{Coupled Systems Mechanics}, 9\penalty0 (1):\penalty0 5--28,
  2020.
\newblock \doi{10.12989/csm.2020.9.1.005}.

\bibitem[Hassan et~al.(2019)Hassan, Ghavamian, Lee, Gil, Bonet, and
  Auricchio]{Hassan2019paper}
Osama~I. Hassan, Ataollah Ghavamian, Chun~Hean Lee, Antonio~J. Gil, Javier
  Bonet, and Ferdinando Auricchio.
\newblock An upwind vertex centred finite volume algorithm for nearly and truly
  incompressible explicit fast solid dynamic applications: Total and updated
  lagrangian formulations.
\newblock \emph{Journal of Computational Physics: X}, 3:\penalty0 100025, 2019.
\newblock ISSN 2590-0552.
\newblock \doi{https://doi.org/10.1016/j.jcpx.2019.100025}.
\newblock URL
  \url{http://www.sciencedirect.com/science/article/pii/S2590055219300411}.

\bibitem[Smith et~al.(2013)Smith, Griffiths, and Margetts]{ParaFEM}
I.~M. Smith, D.~V. Griffiths, and L.~Margetts.
\newblock \emph{Programming the Finite Element Method}.
\newblock Wiley, 5$^{th}$ edition, 2013.

\bibitem[Reed and Hill(1973)]{Reed1973}
W.~H. Reed and T.~R. Hill.
\newblock Triangular mesh methods for the neutron transport equation.
\newblock Technical report, Los Alamos Scientific Lab., N.Mex. (USA), 1973.
\newblock Technical Report 836, {LA-UR-73-479; CONF-730414-2}.

\bibitem[Cockburn(2003)]{Cockburn2003}
B.~Cockburn.
\newblock Discontinuous galerkin methods.
\newblock Technical report, School of Mathematics, Univeristy of Minnesota,
  2003.

\bibitem[Cockburn et~al.(2009)Cockburn, Gopalakrishnan, and
  Lazarov]{Cockburn2009}
B.~Cockburn, J.~Gopalakrishnan, and R.~Lazarov.
\newblock Unified hybridization of discontinuous galerkin, mixed, and
  continuous galerkin methods for second order elliptic problems.
\newblock \emph{{SIAM} Journal on Numerical Analysis}, 47:\penalty0 1319--1365,
  2009.
\newblock \doi{10.1137/070706616}.

\bibitem[Fu et~al.(2015)Fu, Cockburn, and Stolarski]{Fu2015}
G.~Fu, B.~Cockburn, and H.~Stolarski.
\newblock Analysis of an hdg method for linear elasticity.
\newblock \emph{International Journal for Numerical Methods in Engineering},
  102\penalty0 (3-4):\penalty0 551--575, 2015.
\newblock \doi{10.1002/nme.4781}.
\newblock URL \url{https://onlinelibrary.wiley.com/doi/abs/10.1002/nme.4781}.

\bibitem[Qiu et~al.(2018)Qiu, Shen, and Shi]{Qiu2017}
W.~Qiu, J.~Shen, and K.~Shi.
\newblock An hdg method for linear elasticity with strong symmetric stresses.
\newblock \emph{Math. Comp.}, 87:\penalty0 69--93, 2018.
\newblock \doi{10.1090/mcom/3249}.

\bibitem[Sevilla et~al.(2018{\natexlab{c}})Sevilla, Giacomini, Karkoulias, and
  Huerta]{Sevilla2018HDG}
R.~Sevilla, M.~Giacomini, A.~Karkoulias, and A.~Huerta.
\newblock A superconvergent hybridisable discontinuous galerkin method for
  linear elasticity.
\newblock \emph{International Journal for Numerical Methods in Engineering},
  116:\penalty0 91--116, 2018{\natexlab{c}}.
\newblock \doi{10.1002/nme.5916}.

\bibitem[Hesthaven and Warburton(2007)]{Hesthaven2007}
J.~S. Hesthaven and T.~Warburton.
\newblock \emph{Nodal Discontinuous {Galerkin} Methods: Algorithms, Analysis,
  and Applications (Texts in Applied Mathematics)}.
\newblock Springer, New York, USA, 2007.

\bibitem[Slone et~al.(June 1997)Slone, Pericleous, Bailey, and
  Cross]{Slone1997}
A.~K. Slone, K.~Pericleous, C.~Bailey, and M.~Cross.
\newblock Dynamic fluid structure interactions using finite volume unstructured
  mesh procedures.
\newblock In \emph{International Forum on Aero-elasticity and Structural
  Dynamics}, pages 417--424, June 1997.

\bibitem[Demird{\v{z}}i\'{c}(1998)]{Demirdzic1998}
I.~Demird{\v{z}}i\'{c}.
\newblock Finite volume approach to multi-physics problems.
\newblock In S.N. Alturi and P.~E. O'Donoghue, editors, \emph{Modelling and
  Simulation Based Engineering}, pages 1757--1762. Tech Science Press,
  Palmdale, 1998.

\bibitem[Sch{\"{a}}fer et~al.(2000)Sch{\"{a}}fer, Meynen, Sieber, and
  Teschauer]{Schafer2000}
M.~Sch{\"{a}}fer, S.~Meynen, R.~Sieber, and I.~Teschauer.
\newblock Multigrid methods for coupled fluid-solid problems.
\newblock In \emph{European Congress on Computational Methods in Applied
  Sciences and Engineering, {ECCOMAS} 2000}, Barcelona, Spain, 2000.

\bibitem[Slone et~al.(2000{\natexlab{a}})Slone, Cross, Pericleous, and
  Bailey]{Slone2000b}
A.~K. Slone, M.~Cross, K.~Pericleous, and C.~Bailey.
\newblock A finite volume approach to dynamic fluid structure interaction.
\newblock In \emph{8th Annual Conference of the Association for Computational
  Mechanics in Engineering ({ACME} 2000)}, pages 218--221, Greenwich
  University, London, UK, 2000{\natexlab{a}}.

\bibitem[Slone et~al.(2000{\natexlab{b}})Slone, Pericleous, Bailey, and
  Cross]{Slone2000c}
A.~K. Slone, K.~Pericleous, C.~Bailey, and M.~Cross.
\newblock Dynamic fluid structure interactions using finite volume unstructured
  mesh procedures.
\newblock In \emph{8th Symposium on Multidisciplinary Analysis and
  Optimization, Multidisciplinary Analysis Optimization Conferences}. American
  Institute of Aeronautics and Astronautics, 2000{\natexlab{b}}.
\newblock \doi{10.2514/6.2000-4788}.

\bibitem[Slone et~al.(2001{\natexlab{b}})Slone, Pericleous, Bailey, and
  Cross]{Slone2001b}
A.~K. Slone, K.~Pericleous, C.~Bailey, and M.~Cross.
\newblock Dynamic fluid structure interactions using finite volume unstructured
  mesh procedures.
\newblock In \emph{{ECCOMAS} Computational Fluid Dynamics Conference},
  2001{\natexlab{b}}.
\newblock ISBN 0 905 091 12 4.

\bibitem[Slone et~al.(2001{\natexlab{c}})Slone, Cross, Pericleous, Bailey, and
  Cross]{Slone2001c}
A.~K. Slone, M.~Cross, K.~Pericleous, C.~Bailey, and M.~Cross.
\newblock Using finite volume unstructured mesh approach to dynamic fluid
  structure interaction: An assessment of the challenge of flutter analysis.
\newblock In W.~A. Wall, K.~U. Bletzinger, and K.~Schweizerhof, editors,
  \emph{Trends in Computational Structural Mechanics}, pages 741--750. {CIMNE}:
  International Centre for Numerical Methods in Engineering,
  2001{\natexlab{c}}.
\newblock ISBN 84-89925-77-1.
\newblock
  \url{www.tib.eu/de/suchen/id/BLCP%3ACN050065390/A-finite-volume-unstructured-mesh-approach-to-dynamic/}.

\bibitem[D{\v{z}}aferovi\'{c}(2002)]{Dzaferovic2002}
E.~D{\v{z}}aferovi\'{c}.
\newblock \emph{Interaction of viscoplastic fluid and viscoelastic solid -
  Numerical modelling}.
\newblock PhD thesis, University of Sarajevo, 2002.
\newblock In Bosnian.

\bibitem[Kara\v{c} and Ivankovi\'{c}(2002)]{Karac2002}
A.~Kara\v{c} and A.~Ivankovi\'{c}.
\newblock Drop impact of fluid-filled plastic containers: finite volume method
  for coupled fluid-structure-fracture problems.
\newblock In H.A. Mang, F.~G. Rammerstorfer, and J.~Eberhardsteiner, editors,
  \emph{{WCCM V, Fifth World Congres on Computational Mechanics}}, Vienna,
  Austria, 2002.

\bibitem[Slone et~al.(2002{\natexlab{c}})Slone, Croft, Williams, and
  Cross]{Slone2002c}
A.~K. Slone, T.cN. Croft, A.~J. Williams, and M.~Cross.
\newblock A two fluid approach to high impact interaction amongst solid
  structures.
\newblock In H.~A. Mang, F.~G. Rammerstorfer, and J.~Eberhardsteiner, editors,
  \emph{Fifth World Congress on Computational Mechanics Proceedings}, Vienna,
  Austria, 2002{\natexlab{c}}. Vienna University of Technology, Austria.
\newblock ISBN 3-9501554-0-6.

\bibitem[Slone et~al.(2002{\natexlab{d}})Slone, Grossman, Williams, Pericleous,
  Bailey, and Cross]{Slone2002d}
A.~K. Slone, D.~Grossman, A.J. Williams, K.~Pericleous, C.~Bailey, and
  M.~Cross.
\newblock A time and space accurate numerical approach to closely coupled fluid
  structure interaction problems.
\newblock In \emph{11th International Colloquium on Numerical Analysis and
  Computer Science With Applications}, Plovdiv, Bulgaria, 2002{\natexlab{d}}.

\bibitem[Kara\v{c} and Ivankovi\'{c}(2003{\natexlab{a}})]{Karac2003a}
A.~Kara\v{c} and A.~Ivankovi\'{c}.
\newblock Fully predictive model of the drop impact and fracture of
  fluid-filled plastic containers.
\newblock In \emph{Proceedings of 11th {ACME} conference on computational
  mechanics in engineering}, pages 113--116, Glasgow, 2003{\natexlab{a}}.
  University of Strathclyde Publishing.

\bibitem[Kara\v{c} and Ivankovi\'{c}(2003{\natexlab{b}})]{Karac2003b}
A.~Kara\v{c} and A.~Ivankovi\'{c}.
\newblock Modelling the drop impact behaviour of fluid-filled polyethylene
  containers.
\newblock In B.~R.~K. Blackman, A.~Pavan, and J.~G. Williams, editors,
  \emph{Fracture of Polymers}, pages 253--264. {ESIS} publication 32,
  2003{\natexlab{b}}.
\newblock Composites and Adhesives.

\bibitem[Kara\v{c}(2003)]{Karac2003c}
A.~Kara\v{c}.
\newblock \emph{Drop impact of fluid-filled polyethylene containers}.
\newblock PhD thesis, Imperial College London, 2003.

\bibitem[Cross et~al.(2004)Cross, Slone, Croft, and Williams]{Cross2004}
M.~Cross, A.~K. Slone, T.~N. Croft, and A.~J. Williams.
\newblock Computational modelling of thermal fluid structure interaction
  processes.
\newblock In B.~H.~V. Topping, editor, \emph{Progress in Engineering
  Computational Technology}, page 111 – 126. Saxe-Coberg Publications, 2004.

\bibitem[Kara\v{c} and Ivankovi\'{c}(2004)]{Karac2004}
A.~Kara\v{c} and A.~Ivankovi\'{c}.
\newblock Modelling drop impact and fracture of fluid-filled plastic
  containers.
\newblock In \emph{Proceedings of The 15$^{th}$ European Conference on Fracture
  - Advanced Fracture Mechanics for Life and Safety Assessments}, Stockholm,
  Sweden, 2004.

\bibitem[Giannopapa(2004)]{Giannopapa2004a}
C.~G. Giannopapa.
\newblock \emph{Fluid-Structure interaction in flexible vessels}.
\newblock PhD thesis, University of London, 2004.

\bibitem[Giannopapa and Papadakis(2004)]{Giannopapa2004b}
C.~G. Giannopapa and G.~Papadakis.
\newblock A new formulation for solids suitable for a unified solution method
  for fluid-structure interaction problems.
\newblock In \emph{{ASME PVP}}, volume 491-1, pages 111--117, San Diego
  California, 2004.

\bibitem[Slone and Cross(2006)]{Slone2004b}
A.~K. Slone and M.~Cross.
\newblock A comparison of finite element and finite volume methods for
  computational structural mechanics and their application in multi-physics
  problems.
\newblock In \emph{5th International Conference on Engineering Computational
  Technology}. Civil-Comp Press, 2006.
\newblock ISBN 0-948749-96-2.
\newblock paper number 71.

\bibitem[Slone et~al.(2006)Slone, Croft, Williams, and Cross]{Slone2004c}
A.~K. Slone, T.~N. Croft, A.~J. Williams, and M.~Cross.
\newblock A mixed eulerian-lagrangian approach to high speed collision between
  solid structures on parallel clusters.
\newblock In \emph{5th International Conference on Engineering Computational
  Technology}. Civil-Comp Press, 2006.
\newblock ISBN 1-905088-10-8.
\newblock paper number 118.

\bibitem[Cross et~al.(2006)Cross, Croft, McBride, Slone, and
  Williams]{Cross2006}
M.~Cross, T.~N. Croft, D.~McBride, A.~K. Slone, and A.~J. Williams.
\newblock Using mixed discretisation schemes in multi-physics simulation.
\newblock In \emph{Innovation in Engineering Computational Technology}, pages
  309--324. Saxe-Coburg Publications, 2006.
\newblock \url{www.ctresources.info/csets/download/csets.315.pdf}.

\bibitem[Papadakis and Giannopapa(2006)]{Papadakis2006}
G.~Papadakis and C.~G. Giannopapa.
\newblock Towards a unified solution method for fluid-structure interaction
  problems: Progress and challenges.
\newblock In \emph{{Proceedings of PVP 2006-ICPVT11 10$^{th}$ International
  Symposium on Emerging Technology in Fluids}}, Vancouver, Canada, 2006.

\bibitem[Giannopapa and Papadakis(2007)]{Giannopapa2007}
C.~G. Giannopapa and G.~Papadakis.
\newblock Indicative results and progress on the development of the unified
  single solution method for fluid-structure interaction problems.
  {(CASA-report, No. 0711)}.
\newblock Technical report, Technische Universiteit Eindhoven, Eindhoven, 2007.

\bibitem[Kara\v{c} and Ivankovi\'{c}(2009{\natexlab{a}})]{Karac2009a}
A.~Kara\v{c} and A.~Ivankovi\'{c}.
\newblock Investigating the behaviour of fluid-filled polyethylene containers
  under base drop impact: A combined experimental/numerical approach.
\newblock \emph{International Journal of Impact Engineering}, 36:\penalty0
  621--631, 2009{\natexlab{a}}.

\bibitem[Kara\v{c} and Ivankovi\'{c}(2009{\natexlab{b}})]{Karac2009b}
A.~Kara\v{c} and A.~Ivankovi\'{c}.
\newblock Behaviour of fluid filled {PE} containers under impact: Theoretical
  and numerical investigation.
\newblock \emph{International Journal of Impact Engineering}, 36:\penalty0
  621--631, 2009{\natexlab{b}}.

\bibitem[Safari et~al.(2009)Safari, Ivankovi\'{c}, Tukovi\'{c}, and
  M.~Walter]{Safari2009}
A.~Safari, A.~Ivankovi\'{c}, {\v{Z}}.~Tukovi\'{c}, and E.~Casey M.~Walter.
\newblock A fluid-structure interaction study of biofilm detachment.
\newblock In \emph{1$^{st}$ International Conference on Mathematical and
  Computational Biomedical Engineering - {CMBE2009}, June 29 -July 1 Swansea,
  UK}. P. Nithiarasu , R. L, 2009.

\bibitem[Slone et~al.(2009)Slone, Williams, Croft, and Cross]{Slone2009}
A.~K. Slone, A.~J. Williams, T.~N. Croft, and M.~Cross.
\newblock Dynamic fluid structure interaction in parallel: a challenge for
  scalability.
\newblock In B.~H.~V. Topping and P.~Ivanyi, editors, \emph{Parallel,
  Distributed and Grid Computing for Engineering}, pages 329--350. Saxe-Coburg
  Publications, 2009.
\newblock based upon an invited keynote presentation.

\bibitem[Das et~al.(2010)Das, Mathur, and Murthy]{Das2010}
S.~Das, S.~R. Mathur, and J.~Y. Murthy.
\newblock {An unstructured finite-volume method for structure-electrostatic
  interactions in MEMS}.
\newblock In \emph{{Proceedings of IMECE2010 2010 ASME International Mechanical
  Engineering Congress and Exposition}}. Vancouver, Canada, 2010.

\bibitem[Kelly and O'Rourke(2010)]{Kelly2010}
A.~Kelly and M.-J. O'Rourke.
\newblock Two system, single analysis, fluid-structure interaction modelling of
  the abdominal aortic aneurysms.
\newblock \emph{Proceedings of the Institution of Mechanical Engineers Part H -
  Journal of Engineering in Medicine}, 224\penalty0 (H8):\penalty0 955--970,
  2010.

\bibitem[Kelly and O'Rourke(2012)]{Kelly2012}
A.~Kelly and M.-J. O'Rourke.
\newblock Fluid, solid and fluid-structure interaction simulations on
  patient-based abdominal aortic aneurysm models.
\newblock \emph{Proceedings of the Institution of Mechanical Engineers Part H -
  Journal of Engineering in Medicine}, 226\penalty0 (4):\penalty0 288--304,
  2012.

\bibitem[Das(2013)]{Das2013}
S.~Das.
\newblock \emph{Fluid-Structure Interactions in Microstructures}.
\newblock PhD thesis, University of Texas at Austin, 2013.

\bibitem[de~Oliveira et~al.(2017)de~Oliveira, Gasche, Militzer, and
  Baccin]{Oliveira2017a}
I.~L. de~Oliveira, J.~L. Gasche, J.~Militzer, and C.~E. Baccin.
\newblock {Using {FOAM-Extend} to Assess the Influence of Fluid-Solid
  Interaction on the Flow in Intracranial Aneurismus}.
\newblock In \emph{COBEM-2017-0851, 24$^{th}$ ABCM International Congress of
  Mechanical Engineering}, Curitiba, PR, Brazil, 2017.

\bibitem[de~Oliveira(2017)]{Oliveira2017b}
I.~L. de~Oliveira.
\newblock \emph{Using {FOAM-Extend} to assess the influence of fluid-structure
  interaction on the rupture of intracranial aneurysms}.
\newblock PhD thesis, Sao Paulo State University, J\'{u}lio De Mesquita Filho,
  2017.

\bibitem[Cardiff et~al.(2018{\natexlab{a}})Cardiff, Kara\v{c}, Jaeger, Jasak,
  Nagy, Ivankovi\'{c}, and Tukovi\'{c}]{Cardiff2018a}
P.~Cardiff, A.~Kara\v{c}, P.~De Jaeger, H.~Jasak, J.~Nagy, A.~Ivankovi\'{c},
  and {\v{Z}}.~Tukovi\'{c}.
\newblock Towards the development of an extendable solid mechanics and
  fluid-solid interactions toolbox for {OpenFOAM}.
\newblock \emph{preprint}, 2018{\natexlab{a}}.
\newblock arXiv:1808.10736 [math.NA], available at
  \url{https://arxiv.org/abs/1808.10736}.

\bibitem[Tukovi\'{c} et~al.(2018{\natexlab{b}})Tukovi\'{c}, Buka\v{c}, Cardiff,
  Jasak, and Ivankovi\'{c}]{Tukovic2018b}
{\v{Z}}.~Tukovi\'{c}, M.~Buka\v{c}, P.~Cardiff, H.~Jasak, and A.~Ivankovi\'{c}.
\newblock Added mass partitioned fluid-structure interaction solver based on a
  robin boundary condition for pressure.
\newblock In \emph{{OpenFOAM} Selected papers of the 11th Workshop}, pages
  1--23. Springer International, 2018{\natexlab{b}}.

\bibitem[Leevers et~al.(1993)Leevers, Venizelos, and
  Ivankovi\'{c}]{Leevers1993}
P.~S. Leevers, G.~Venizelos, and A.~Ivankovi\'{c}.
\newblock Rapid crack propagation along pressurized pipe: small-scale testing
  and numerical modelling.
\newblock \emph{Construction and Building Materials}, 1993.

\bibitem[Demird{\v{z}}i\'{c} et~al.(1996)Demird{\v{z}}i\'{c}, Ivankovi\'{c},
  MacGillivray, and Maneeratana]{Demirdzic1996b}
I.~Demird{\v{z}}i\'{c}, A.~Ivankovi\'{c}, H.~J. MacGillivray, and
  K.~Maneeratana.
\newblock Numerical modelling of high-rate tensile tests using finite volume
  formulation.
\newblock In \emph{Proceedings of {IUTAM} Symposium on Innovative Computational
  Methods for Fracture and Damage}, Dublin, Ireland, 1996.

\bibitem[Murphy and Ivankovi\'{c}(1999)]{Murphy1999}
N.~Murphy and A.~Ivankovi\'{c}.
\newblock Dynamic fracture simulation of brittle material characterised by
  microcrack-dominated failure mechanisms.
\newblock In \emph{{Proceedings of 7$^{th}$ ACME Conf Comp Mech in Eng.}},
  pages 99--102, Durham, UK, 1999.

\bibitem[Stylianou(1999)]{Stylianou1999}
V.~Stylianou.
\newblock \emph{Finite volume modelling of rapid crack propagation ({RCP}) in
  brittle polymers}.
\newblock PhD thesis, Imperial College London, 1999.

\bibitem[Pandya et~al.(2000{\natexlab{a}})Pandya, Ivankovi\'{c}, and
  Williams]{Pandya2000a}
K.~C. Pandya, A.~Ivankovi\'{c}, and J.~G. Williams.
\newblock Cohesive zone modelling of crack growth in polymers - part 2 -
  numerical simulation of crack growth.
\newblock \emph{Plastics, Rubber and Composites}, 29:\penalty0 447--452,
  2000{\natexlab{a}}.

\bibitem[Pandya et~al.(2000{\natexlab{b}})Pandya, Ivankovi\'{c}, and
  Williams]{Pandya2000b}
K.~C. Pandya, A.~Ivankovi\'{c}, and J.~G. Williams.
\newblock Predicting crack growth in tough polyethylene from measured cohesive
  zone traction-separation curves.
\newblock In \emph{Proceedings of 11th Int. Conf. on Deformation, Yield and
  Fracture of Polymers}, 2000{\natexlab{b}}.

\bibitem[Pandya et~al.(2000{\natexlab{c}})Pandya, Ivankovi\'{c}, and
  Williams]{Pandya2000c}
K.~C. Pandya, A.~Ivankovi\'{c}, and J.~G. Williams.
\newblock Predictive fracture modelling in tough polyethylenes using
  experimentally measured cohesive zone traction curves.
\newblock In \emph{{Proceedings of 13th European Conference on Fracture -
  ECF13}}, San Sebastian, Spain, 2000{\natexlab{c}}.

\bibitem[Ivankovi\'{c} et~al.(2002{\natexlab{b}})Ivankovi\'{c}, Jasak,
  Kara\v{c}, Trop\v{s}a, and Leevers]{Ivankovic2002a}
A.~Ivankovi\'{c}, H.~Jasak, A.~Kara\v{c}, V.~Trop\v{s}a, and P.~Leevers.
\newblock Fully predictive model of {RCP} in plastic pipes.
\newblock In \emph{Proceedings of 14$^{th}$ European Conference on Fracture},
  Kracow, Poland, 2002{\natexlab{b}}.

\bibitem[Ivankovi\'{c} et~al.(2002{\natexlab{c}})Ivankovi\'{c}, Jasak,
  Kara\v{c}, and Trop\v{s}a]{Ivankovic2002c}
A.~Ivankovi\'{c}, H.~Jasak, A.~Kara\v{c}, and V.~Trop\v{s}a.
\newblock Prediction of dynamic fracture in pressurised plastic pipes.
\newblock In \emph{The Annual Conference of the Association of Computational
  Mechanics in Engineering}, volume~10, pages 173--176, 2002{\natexlab{c}}.

\bibitem[Ivankovi\'{c} et~al.(2004)Ivankovi\'{c}, Murphy, and
  Hillmansen]{Ivankovic2004b}
A.~Ivankovi\'{c}, N.~Murphy, and S.~Hillmansen.
\newblock Evolution of dynamic fractures in {PMMA}: experimental and numerical
  investigations.
\newblock In M.~H. Aliabadi and A.~Ivankovi\'{c}, editors, \emph{Advances in
  fracture mechanics, vol. 9.} WIT Press/Computational Mechanics Publications,
  Southampton, UK, 2004.

\bibitem[Murphy(2007)]{Murphy2007}
N.~Murphy.
\newblock \emph{Dynamic Fracture of PMMA: A combined experimental and numerical
  investigation}.
\newblock PhD thesis, University College Dublin, 2007.

\bibitem[McAuliffe et~al.(2011)McAuliffe, Kara\v{c}, Murphy, and
  Ivankovi\'{c}]{McAuliffe2011}
D.~McAuliffe, A.~Kara\v{c}, N.~Murphy, and A.~Ivankovi\'{c}.
\newblock Transferability of adhesive fracture toughness measurements between
  peel and {TDCB} test methods for a nano-toughened epoxy.
\newblock \emph{Adhesion Society}, 2011.
\newblock \url{http://hdl.handle.net/10197/4765}.

\bibitem[McAuliffe et~al.(2012)McAuliffe, Kara\v{c}, Murphy, and
  Ivankovi\'{c}]{McAuliffe2012a}
D.~McAuliffe, A.~Kara\v{c}, N.~Murphy, and A.~Ivankovi\'{c}.
\newblock Determination of the cohesive strength and toughening mechanisms of a
  nano-modified adhesive under a triaxial stress.
\newblock \emph{Adhesion Society}, 2012.
\newblock \url{http://hdl.handle.net/10197/4785}.

\bibitem[McAuliffe(2012)]{McAuliffe2012b}
D.~McAuliffe.
\newblock \emph{Fracture Toughness Characterisation of a Nano-modified
  Structural Adhesives}.
\newblock PhD thesis, University College Dublin, 2012.

\bibitem[Cooper et~al.(2012)Cooper, Ivankovic, Kara\v{c}, McAuliffe, and
  Murphy]{Cooper2012}
V.~Cooper, A.~Ivankovic, A.~Kara\v{c}, D.~McAuliffe, and N.~Murphy.
\newblock Effects of bond gap thickness on the fracture of nano-toughened epoxy
  adhesive joints.
\newblock \emph{Polymer}, 53\penalty0 (24):\penalty0 5540 -- 5553, 2012.
\newblock \doi{10.1016/j.polymer.2012.09.049}.

\bibitem[Georgiou et~al.(2003{\natexlab{a}})Georgiou, Ivankovi\'{c}, Kinloch,
  and Trop\v{s}a]{Georgiou2003a}
I.~Georgiou, A.~Ivankovi\'{c}, A.J. Kinloch, and V.~Trop\v{s}a.
\newblock Rate dependent fracture behaviour of adhesively bonded joints.
\newblock In A.~Pavan B.~R. K.~Blackman and J.~G. Williams, editors,
  \emph{Fracture of Polymers, Composites and Adhesives {II}}, volume~32 of
  \emph{European Structural Integrity Society}, pages 317 -- 328. Elsevier,
  2003{\natexlab{a}}.

\bibitem[Georgiou et~al.(2003{\natexlab{b}})Georgiou, Hadavinia, Ivankovi\'{c},
  Kinloch, Trop\v{s}a, and Williams]{Georgiou2003b}
I.~Georgiou, H.~Hadavinia, A.~Ivankovi\'{c}, A.~J. Kinloch, V.~Trop\v{s}a, and
  J.~G. Williams.
\newblock Cohesive zone models and the plastically-deforming peel test.
\newblock \emph{The Journal of Adhesion}, 79:\penalty0 239--265,
  2003{\natexlab{b}}.

\bibitem[Cooper et~al.(2008)Cooper, Ivankovi\'{c}, and Kara\v{c}]{Cooper2008}
V.~Cooper, A.~Ivankovi\'{c}, and A.~Kara\v{c}.
\newblock A mode {I} fracture behaviour analysis of adhesively bonded joints.
\newblock In \emph{European conference on fracture}, volume~17, Brno, Czech
  Republic, 2008.

\bibitem[Cooper(2010)]{Cooper2010}
V.~J. Cooper.
\newblock \emph{The Fracture Behaviour of Nano-toughened Structural Epoxy
  Adhesives}.
\newblock PhD thesis, University College Dublin, 2010.

\bibitem[Tabakovi\'{c} et~al.(2010)Tabakovi\'{c}, Kara\v{c}, Ivankovi\'{c},
  Gibney, McNally, and Gilchrist]{Tabakovic2010}
A.~Tabakovi\'{c}, A.~Kara\v{c}, A.~Ivankovi\'{c}, A.~Gibney, C.~McNally, and
  M.~D. Gilchrist.
\newblock Modelling the quasi-static behaviour of bituminous material using a
  cohesive zone model.
\newblock \emph{Engineering Fracture Mechanics}, 77:\penalty0 2403--2418, 2010.

\bibitem[Carolan(2011)]{Carolan2011}
D.~Carolan.
\newblock \emph{Mechanical and Fracture Properties of PCBN as a Function of
  Rate and Temperature}.
\newblock PhD thesis, University College Dublin, 2011.

\bibitem[Petrovi\'{c}(2011)]{Petrovic2011}
M.~Petrovi\'{c}.
\newblock \emph{The Behaviour of Polycrystalline Diamonds as a Function of Rate
  and Temperature}.
\newblock PhD thesis, University College Dublin, 2011.

\bibitem[Carolan et~al.(2012{\natexlab{a}})Carolan, Ivankovi\'{c}, and
  Murphy]{Carolan2012b}
D.~Carolan, A.~Ivankovi\'{c}, and N.~Murphy.
\newblock Numerical investigation into dynamic fracture of pcbn.
\newblock \emph{Key Engineering Materials}, 488:\penalty0 553--556,
  2012{\natexlab{a}}.

\bibitem[Carolan et~al.(2013{\natexlab{b}})Carolan, Ivankovi\'{c}, and
  Murphy]{Carolan2013b}
D.~Carolan, A.~Ivankovi\'{c}, and N.~Murphy.
\newblock A combined experimental-numerical investigation of fracture of
  polycrystalline cubic boron nitride.
\newblock \emph{Engineering Fracture Mechanics}, 99:\penalty0 101--117,
  2013{\natexlab{b}}.

\bibitem[Alveen et~al.(2014)Alveen, McNamara, Carolan, Murphy, and
  Ivankovi\'{c}]{Alveen2014}
P.~Alveen, D.~McNamara, D.~Carolan, N.~Murphy, and A.~Ivankovi\'{c}.
\newblock Analysis of two-phase ceramic composites using micromechanical
  models.
\newblock \emph{Computational Materials Science}, 92:\penalty0 318--324, 2014.

\bibitem[McNamara et~al.(2014)McNamara, Alveen, Carolan, Murphy, and
  Ivankovi\'{c}]{McNamara2014}
D.~McNamara, P.~Alveen, D.~Carolan, N.~Murphy, and A.~Ivankovi\'{c}.
\newblock Micromechanical study of strength and toughness of advanced ceramics.
\newblock \emph{Procedia Materials Science}, 3:\penalty0 1810--1815, 2014.

\bibitem[Alveen(2015)]{Alveen2015}
P.~Alveen.
\newblock \emph{An experimental-numerical investigation into the properties of
  polycrystalline cubic boron nitride towards materials optimisation}.
\newblock PhD thesis, University College Dublin, 2015.

\bibitem[general poro-elastic model for pad-scale fracturing of~horizontal
  wells(2015)]{Manchanda2015}
R.~A general poro-elastic model for pad-scale fracturing of~horizontal wells.
\newblock \emph{A general poro-elastic model for pad-scale fracturing of
  horizontal wells}.
\newblock PhD thesis, University of Texas at Austin, 2015.

\bibitem[McNamara et~al.(2015)McNamara, Alveen, Carolan, Murphy, and
  Ivankovi\'{c}]{McNamara2015}
D.~McNamara, P.~Alveen, D.~Carolan, N.~Murphy, and A.~Ivankovi\'{c}.
\newblock Numerical analysis of the strength of polycrystalline diamond as a
  function of microstructure.
\newblock \emph{International Journal of Refractory Metals and Hard Materials},
  52:\penalty0 195--202, 2015.

\bibitem[McNamara(2015)]{McNamara2015b}
D.~McNamara.
\newblock \emph{The mechanical and fracture properties of polycrystalline
  diamond as a function of microstructure}.
\newblock PhD thesis, University College Dublin, 2015.

\bibitem[Lee(2017)]{Lee2017}
D.~Lee.
\newblock \emph{A Model for Hydraulic Fracturing and Proppant Placement in
  Unconsolidated Sands}.
\newblock PhD thesis, University of Texas at Austin, 2017.

\bibitem[Yi(2018)]{ShitingYi2018}
S.~Shiting Yi.
\newblock \emph{Development of Computationally Efficient 2D and Pseudo-3D
  Multi- Fracture Models with Applications to Fracturing and Refracturing}.
\newblock PhD thesis, University of Texas at Austin, 2018.

\bibitem[Sabbagh-Yazdi et~al.(2018)Sabbagh-Yazdi, Farhoud, and
  Gharebaghi]{Sabbagh-Yazdi2018}
S.~R. Sabbagh-Yazdi, A.~Farhoud, and S.~A. Gharebaghi.
\newblock Simulation of {2D} linear crack growth under constant load using
  {GFVM} and two-point displacement extrapolation method.
\newblock \emph{Applied Mathematical Modelling}, 2018.
\newblock \doi{10.1016/j.apm.2018.05.022}.

\bibitem[Pindera(1991)]{Pindera1991}
M.-J. Pindera.
\newblock Local/global stiffness matrix formulation for composite materials and
  structures.
\newblock \emph{Composites Engineering}, 1(2):\penalty0 69--83, 1991.

\bibitem[Aboudi et~al.(1994)Aboudi, Pindera, and Arnold]{Aboudi1994}
J.~Aboudi, M.-J. Pindera, and S.M. Arnold.
\newblock Elastic response of metal matrix composites with tailored
  microstructures to thermal gradients.
\newblock \emph{International Journal of Solids and Structures}, 31:\penalty0
  1393--1428, 1994.

\bibitem[Aboudi(2002)]{Aboudi2002a}
J.~Aboudi.
\newblock Micromechanical analysis of the fully coupled finite thermoelastic
  response of rubberlike matrix composites.
\newblock \emph{International Journal of Solids and Structures}, 39:\penalty0
  2587--2612, 2002.

\bibitem[Aboudi et~al.(2002)Aboudi, Pindera, and Arnold]{Aboudi2002b}
J.~Aboudi, M.-J. Pindera, and S.M. Arnold.
\newblock High-fidelity generalized method of cells for inelastic periodic
  multiphase materials.
\newblock Technical report, NASA TM-2002-211469, 2002.

\bibitem[Zhong et~al.(2002)Zhong, Pindera, and Arnold]{Zhong2002}
Y.~Zhong, M.-J. Pindera, and S.~Arnold.
\newblock Efficient reformulation of {HOTFGM}: Heat conduction with variable
  thermal conductivity.
\newblock \emph{NASA CR 2002-211910}, November 2002.

\bibitem[Aboudi et~al.(2003)Aboudi, Pindera, and Arnold]{Aboudi2003}
J.~Aboudi, M.-J. Pindera, and S.M. Arnold.
\newblock Higher-order theory for periodic multiphase materials with inelastic
  phases.
\newblock \emph{International Journal of Plasticity}, 19:\penalty0 805--847,
  2003.

\bibitem[Bansal and Pindera(2003)]{Bansal2003}
Y.~Bansal and M.-J. Pindera.
\newblock Efficient reformulation of the thermo-elastic higher-order theory for
  {FGMs}.
\newblock \emph{Journal of Thermal Stresses}, 26(11-12):\penalty0 1055--1092,
  2003.

\bibitem[Arnold et~al.(2004)Arnold, Bednarcyk, and Aboudi]{Arnold2004}
S.M. Arnold, B.~Bednarcyk, and J.~Aboudi.
\newblock Comparison of the computational efficiency of the original versus
  reformulated high-fidelity generalized method of cells.
\newblock Technical report, NASA/TM-2004-213438., 2004.

\bibitem[Bansal and M.-J.(2004)]{Bansal2004}
Y.~Bansal and Pindera M.-J.
\newblock Testing the predictive capability of the high-fidelity generalized
  method of cells using an efficient reformulation.
\newblock Technical report, NASA/CR-2004-213043, 2004.

\bibitem[Bednarcyk et~al.(2004)Bednarcyk, Arnold, Aboudi, and
  Pindera]{Bednarcyk2004}
B.A. Bednarcyk, S.M. Arnold, J.~Aboudi, and M.-J. Pindera.
\newblock Local field effects in titanium matrix composites subject to
  fiber-matrix debonding.
\newblock \emph{International Journal of Plasticity}, 20:\penalty0 1707--1737,
  2004.

\bibitem[Pindera et~al.(2004)Pindera, Bansal, and Zhong]{Pindera2004}
M.-J. Pindera, Y.~Bansal, and Y.~Zhong.
\newblock {Finite-Volume Direct Averaging Theory for Functionally Graded
  Materials (FVDAT-FGM)}.
\newblock Technical report, {NASA} Disclosure of Invention and New Technology
  Form 1679, 2004.

\bibitem[Zhong et~al.(2004)Zhong, Bansal, and Pindera]{Zhong2004}
Y.~Zhong, Y.~Bansal, and M.-J. Pindera.
\newblock Efficient reformulation of the thermal higher-order theory for
  {FGM}'s with variable thermal conductivity.
\newblock \emph{International Journal of Computational Engineering Science},
  5\penalty0 (4):\penalty0 795--831, 2004.

\bibitem[Aboudi(2005)]{Aboudi2005a}
J.~Aboudi.
\newblock Micromechanically established constitutive equations for multiphase
  materials with viscoelastic-viscoplastic phases.
\newblock \emph{Mechanics of Time-Dependent Materials}, 9:\penalty0 121--145,
  2005.

\bibitem[Aboudi and Gilat(2005)]{Aboudi2005b}
J.~Aboudi and R.~Gilat.
\newblock Micromechanical analysis of lattice blocks.
\newblock \emph{International Journal of Solids and Structures}, 42:\penalty0
  4372--4392, 2005.

\bibitem[Cavalcante(2006)]{Cavalcante2006}
M.~A.~A. Cavalcante.
\newblock Modelling of the transient thermo-mechanical behavior of composite
  material structures by the finite-volume theory.
\newblock Master's thesis, Federal University of Alagoas, Maceio, Alagoas,
  Brazil, 2006.

\bibitem[Pindera and Bansal(2006)]{Pindera2006}
M.-J. Pindera and Y.~Bansal.
\newblock Finite volume direct averaging micromechanics of heterogeneous
  materials with elastic-plastic phases.
\newblock \emph{International Journal of Plasticity}, 22:\penalty0 775--825,
  2006.

\bibitem[Bruck et~al.(2007)Bruck, Gilat, Aboudi, and Gershon]{Bruck2007}
H.A. Bruck, R.~Gilat, J.~Aboudi, and A.L. Gershon.
\newblock A new approach for optimizing the mechanical behavior of porous
  microstructures for porous materials by design.
\newblock \emph{Modelling and Simulation in Materials Science and Engineering},
  15:\penalty0 653--674, 2007.

\bibitem[Cavalcante et~al.(2007{\natexlab{a}})Cavalcante, Marques, and
  Pindera]{Cavalcante2007a}
M.~A.~A. Cavalcante, S.~P.~C. Marques, and M.-J. Pindera.
\newblock Parametric formulation of the finite-volume theory for functionally
  graded materials - part {I}: Analysis.
\newblock \emph{{ASME} Journal of Applied Mechanics}, 74:\penalty0 935--945,
  2007{\natexlab{a}}.

\bibitem[Cavalcante et~al.(2007{\natexlab{b}})Cavalcante, Marques, and
  Pindera]{Cavalcante2007b}
M.~A.~A. Cavalcante, S.~P.~C. Marques, and M.-J. Pindera.
\newblock Parametric formulation of the finite-volume theory for functionally
  graded materials - part {II}: Numerical results.
\newblock \emph{{ASME} Journal of Applied Mechanics}, 74:\penalty0 946--957,
  2007{\natexlab{b}}.

\bibitem[Drago and Pindera(2007)]{Drago2007}
A.~S. Drago and M.-J. Pindera.
\newblock Micro-macromechanical analysis of heterogeneous materials:
  macroscopically homogeneous vs periodic microstructures.
\newblock \emph{Composites Science and Technology}, 67(6):\penalty0 1243--63,
  2007.

\bibitem[Gattu(2007)]{Gattu2007}
M.~Gattu.
\newblock Parametric finite volume theory for periodic heterogeneous materials.
\newblock Master's thesis, University of Virginia, Charlottesville, Virginia,
  USA, 2007.

\bibitem[Pindera and Bansal(2007)]{Pindera2007}
M.-J. Pindera and Y.~Bansal.
\newblock On the micromechanics-based simulation of metal matrix composite
  response.
\newblock \emph{Journal of Engineering Materials and Technology},
  129(3):\penalty0 468--82, 2007.

\bibitem[Ryvkin and Aboudi(2007)]{Ryvkin2007}
M.~Ryvkin and J.~Aboudi.
\newblock The effect of fiber loss in periodic composites.
\newblock \emph{International Journal of Solids and Structures}, 44:\penalty0
  3497--3513., 2007.

\bibitem[Bednarcyk et~al.(2008)Bednarcyk, Aboudi, Arnold, and
  Sullivan]{Bednarcyk2008}
B.A. Bednarcyk, J.~Aboudi, S.M. Arnold, and R.M. Sullivan.
\newblock Analysis of space shuttle external tank spray-on foam insulation with
  internal pore pressure.
\newblock \emph{Journal of Engineering Materials Technology}, 130:\penalty0
  041005--0410016, 2008.

\bibitem[Cavalcante et~al.(2008)Cavalcante, Marques, and
  Pindera]{Cavalcante2008}
M.~A.~A. Cavalcante, S.~P.~C. Marques, and M.-J. Pindera.
\newblock Computational aspects of the parametric finite-volume theory for
  functionally graded materials.
\newblock \emph{Computational Materials Science}, 44:\penalty0 422--438, 2008.

\bibitem[Gattu et~al.(2008)Gattu, Khatam, Drago, and Pindera]{Gattu2008}
M.~Gattu, H.~Khatam, A.~S. Drago, and M.-J. Pindera.
\newblock Parametric finite-volume micromechanics of uniaxial,
  continuously-reinforced periodic materials with elastic phases.
\newblock \emph{Journal of Engineering Materials and Technology}, 130:\penalty0
  31015--31030, 2008.

\bibitem[Paulino et~al.(2008)Paulino, Pindera, Dodds, Rochinha, Dave, and
  Chen]{Paulino2008}
G.~H. Paulino, M.-J. Pindera, R.~H. Dodds, F.~E. Rochinha, E.~V. Dave, and
  L.~Chen.
\newblock Multiscale and functionally graded materials.
\newblock In \emph{{AIP} conference proceedings}, volume 973, Melville, New
  York, 2008.

\bibitem[Cavalcante et~al.(2009)Cavalcante, Marques, and
  Pindera]{Cavalcante2009}
M.~A.~A. Cavalcante, S.~P.~C. Marques, and M.-J. Pindera.
\newblock Transient thermo-mechanical analysis of a layered cylinder by the
  parametric finite-volume theory.
\newblock \emph{Journal of Thermal Stresses}, 32:\penalty0 112--134, 2009.

\bibitem[Gao et~al.(2009)Gao, Song, and Sun]{Gao2009}
X.~Gao, Y.~Song, and Z.~Sun.
\newblock Quadrilateral subcell based finite volume micromechanics theory for
  multiscale analysis of elastic periodic materials.
\newblock \emph{{ASME} Journal of Applied Mechanics}, 76:\penalty0 011013--1,
  2009.

\bibitem[Khatam and Pindera(2009{\natexlab{a}})]{Khatam2009a}
H.~Khatam and M.-J. Pindera.
\newblock Parametric finite-volume micromechanics of periodic materials with
  elastoplastic phases.
\newblock \emph{International Journal of Plasticity}, 25:\penalty0 1386--1411,
  2009{\natexlab{a}}.

\bibitem[Khatam and Pindera(2009{\natexlab{b}})]{Khatam2009b}
H.~Khatam and M.-J. Pindera.
\newblock Thermo-elastic moduli of lamellar composites with wavy architectures.
\newblock \emph{Composites Part B: Engineering}, 40(1):\penalty0 50--64,
  2009{\natexlab{b}}.

\bibitem[Khatam et~al.(2009)Khatam, Chen, and Pindera]{Khatam2009c}
H.~Khatam, L.~Chen, and M.-J. Pindera.
\newblock Elastic and plastic response of perforated plates with different
  porosity architectures.
\newblock \emph{Journal of Engineering Materials and Technology, Transactions
  of the {ASME}}, 131(3):\penalty0 031015--14, 2009.

\bibitem[Aboudi and Freed(2010)]{Aboudi2010}
J.~Aboudi and Y.~Freed.
\newblock \emph{Shape Memory Alloys: Manufacture Properties and Applications,
  Micromechanical modeling of shape memory alloy composites}.
\newblock Nova Science Publishers, 2010.

\bibitem[Bednarcyk et~al.(2010)Bednarcyk, Aboudi, Arnold, and
  Sullivan]{Bednarcyk2010}
B.A. Bednarcyk, J.~Aboudi, S.M. Arnold, and R.M. Sullivan.
\newblock Micromechanics modeling of composites subjected to multiaxial
  progressive damage in the constituents.
\newblock \emph{AIAA Journal}, 48:\penalty0 1367--1378, 2010.

\bibitem[Haj-Ali and Aboudi(2010)]{Haj-Ali2010}
R.~Haj-Ali and J.~Aboudi.
\newblock Formulation of the high-fidelity generalized method of cells with
  arbitrary cell geometry for refined micromechanics and damage in composites.
\newblock \emph{International Journal of Solids and Structures}, 47:\penalty0
  3447--3461, 2010.

\bibitem[Khatam and Pindera(2010)]{Khatam2010}
H.~Khatam and M.-J. Pindera.
\newblock Plasticity-triggered architectural effects in periodic multilayers
  with wavy microstructures.
\newblock \emph{International Journal of Plasticity}, 26(2):\penalty0 273--287,
  2010.

\bibitem[Aboudi(2011)]{Aboudi2011}
J.~Aboudi.
\newblock The effect of anisotropic damage evolution on the behavior of ductile
  and brittle matrix composites.
\newblock \emph{International Journal of Solids and Structures}, 48:\penalty0
  2102--2119, 2011.

\bibitem[Cavalcante et~al.(2011{\natexlab{a}})Cavalcante, Marques, and
  Pindera]{Cavalcante2011a}
M.~A.~A. Cavalcante, S.~P.~C. Marques, and M.-J. Pindera.
\newblock Transient finite-volume analysis of a graded cylindrical shell under
  thermal shock loading.
\newblock \emph{Mechanics of Advanced Materials and Structures}, 18:\penalty0
  53--67, 2011{\natexlab{a}}.

\bibitem[Cavalcante et~al.(2011{\natexlab{b}})Cavalcante, Khatam, and
  Pindera]{Cavalcante2011b}
M.~A.~A. Cavalcante, H.~Khatam, and M.-J. Pindera.
\newblock {Homogenization of elastic-plastic periodic materials by FVDAM and
  FEM approaches - An assessment}.
\newblock \emph{Composites Part B: Engineering}, 42:\penalty0 1713--1730,
  2011{\natexlab{b}}.

\bibitem[Chareonsuk and Vessakosol(2011)]{Chareonsuk2011}
J.~Chareonsuk and P.~Vessakosol.
\newblock Numerical solution for functionally graded solids under thermal and
  mechanical loads using a high-order control volume finite element method.
\newblock \emph{Applied Thermal Engineering}, 31:\penalty0 213--27, 2011.

\bibitem[Khatam and Pindera(2011)]{Khatam2011}
H.~Khatam and M.-J. Pindera.
\newblock Plastic deformation modes in perforated sheets and their relation to
  yield and limit surfaces.
\newblock \emph{International Journal of Plasticity}, 27(10):\penalty0
  1537--59, 2011.

\bibitem[Carolan et~al.(2012{\natexlab{b}})Carolan, Tukovi\'{c}, McNamara,
  Alveen, Murphy, and Ivankovi\'{c}]{Carolan2012a}
D.~Carolan, {\v{Z}}.~Tukovi\'{c}, D.~McNamara, P.~Alveen, N.~Murphy, and
  A.~Ivankovi\'{c}.
\newblock Effect of microstructure on the fracture toughness of polycrystalline
  cubic boron nitride.
\newblock In \emph{7$^{th}$ {OpenFOAM} Workshop}, Darmstadt, Germany,
  2012{\natexlab{b}}.

\bibitem[Cavalcante and Pindera(2012{\natexlab{a}})]{Cavalcante2012a}
M.~A.~A. Cavalcante and M.-J. Pindera.
\newblock Generalized finite-volume theory for elastic stress analysis in solid
  mechanics - part {I}: Framework.
\newblock \emph{{ASME} Journal of Applied Mechanics}, 79, 2012{\natexlab{a}}.

\bibitem[Cavalcante and Pindera(2012{\natexlab{b}})]{Cavalcante2012b}
M.~A.~A. Cavalcante and M.-J. Pindera.
\newblock Generalized finite-volume theory for elastic stress analysis in solid
  mechanics - part {II}: Results.
\newblock \emph{{ASME} Journal of Applied Mechanics}, 79, 2012{\natexlab{b}}.

\bibitem[Cavalcante(2012)]{Cavalcante2012d}
M.~A.~A. Cavalcante.
\newblock \emph{Generalized finite-volume micromechanics theory for
  heterogeneous materials}.
\newblock PhD thesis, University of Virginia, 2012.

\bibitem[Khatam and Pindera(2012)]{Khatam2012}
H.~Khatam and M.-J. Pindera.
\newblock Microstructural scale effects in the nonlinear elastic response of
  bio-inspired wavy multilayers undergoing finite deformation.
\newblock \emph{Composites Part B: Engineering}, 43(3):\penalty0 869--84, 2012.

\bibitem[Cavalcante and Pindera(2013)]{Cavalcante2013}
M.~A.~A. Cavalcante and M.-J. Pindera.
\newblock Generalized {FVDAM} theory for periodic materials with elastic-
  plastic phases.
\newblock In \emph{{CILAMCE 2013 Proceedings of the XXXIV Iberian
  Latin-American Congress on Computational Methods in Engineering}}, {ABMEC},
  Pirenp\'{o}lis, {GO}, Brazil, 2013.

\bibitem[Cavalcante and Pindera(2014{\natexlab{a}})]{Cavalcante2014a}
M.A.A. Cavalcante and M.-J. Pindera.
\newblock Generalized {FVDAM} theory for periodic materials undergoing finite
  deformations - part {I}: Framework.
\newblock \emph{{ASME} Journal of Applied Mechanics}, 81(2):\penalty0
  021005--021010., 2014{\natexlab{a}}.

\bibitem[Cavalcante and Pindera(2014{\natexlab{b}})]{Cavalcante2014b}
M.A.A. Cavalcante and M.-J. Pindera.
\newblock Generalized {FVDAM} theory for periodic materials undergoing finite
  deformations - part {II}: numerical results.
\newblock \emph{{ASME} Journal of Applied Mechanics}, 81(2):\penalty0
  021006--021012, 2014{\natexlab{b}}.

\bibitem[Cardiff et~al.(2014{\natexlab{e}})Cardiff, Leonard, Murphy, and
  Ivankovi\'{c}]{Cardiff2014g}
P.~Cardiff, M.~Leonard, N.~Murphy, and A.~Ivankovi\'{c}.
\newblock Fracture toughness optimization of nano-toughened structural
  adhesives: A representative volume element approach.
\newblock In \emph{Proceedings of the 37$^{th}$ Annual Meeting of the Adhesion
  Society}, San Diego, CA, USA, 2014{\natexlab{e}}.

\bibitem[Leonard(2014)]{Leonard2014}
M.~Leonard.
\newblock \emph{Micro-Mechanical Modelling of Toughening Mechanisms in
  Nano-Toughened Structural Adhesives}.
\newblock PhD thesis, University College Dublin, 2014.

\bibitem[Tu and Pindera(2014)]{Tu2014}
W.~Tu and M.-J. Pindera.
\newblock Cohesive zone-based damage evolution in periodic materials via finite
  volume homogenization.
\newblock \emph{{ASME} Journal of Applied Mechanics}, 81 (10):\penalty0 1--12,
  01005, 2014.

\bibitem[Carolan et~al.(2015)Carolan, H.~M.~Chong, Kinloch, and
  Taylor]{Carolan2015}
D.~Carolan, A.~Ivankovic H.~M.~Chong, A.~J. Kinloch, and A.~C. Taylor.
\newblock Co-continuous polymer systems: A numerical investigation.
\newblock \emph{Computational Materials Science}, 98:\penalty0 24--33, 2015.

\bibitem[Tu(2016)]{Tu2016}
W.~Tu.
\newblock \emph{{CZM}-based Finite-Volume Homogenization and Optimization of
  Periodic Composites}.
\newblock PhD thesis, University of Virginia, 2016.

\bibitem[Chen et~al.(2017)Chen, Wang, Chen, and Geng]{Chen2017}
Q.~Chen, G.~Wang, X.~Chen, and J.~Geng.
\newblock Finite-volume homogenization of elastic/viscoelastic periodic
  materials.
\newblock \emph{Composite Structures}, 182:\penalty0 457--470, 2017.

\bibitem[Chen et~al.(2018{\natexlab{a}})Chen, Wang, and Pindera]{Chen2018a}
Q.~Chen, G.~Wang, and M.-J. Pindera.
\newblock Finite-volume homogenization and localization of nanoporous materials
  with cylindrical voids. part 1: Theory and validation.
\newblock \emph{European Journal of Mechanics / A Solids}, 2018{\natexlab{a}}.
\newblock \doi{10.1016/j.euromechsol.2018.02.004}.

\bibitem[Chen et~al.(2018{\natexlab{b}})Chen, Tu, Liu, and Chen]{Chen2018b}
Q.~Chen, W.~Tu, R.~Liu, and X.~Chen.
\newblock Parametric multiphysics finite-volume theory for periodic composites
  with thermo-electro-elastic phases.
\newblock \emph{Journal of Intelligent Material Systems and Structures},
  29:\penalty0 530--552, 03 2018{\natexlab{b}}.

\bibitem[Ye et~al.(2018)Ye, Hong, Cai, Wang, Zhai, and Shi]{Ye2018}
J.~Ye, Y.~Hong, H.~Cai, Y.~Wang, Z.~Zhai, and B.~Shi.
\newblock A new three-dimensional parametric {FVDAM} for investigating the
  effective elastic moduli of particle-reinforced composites with interphase.
\newblock \emph{Mechanics of Advanced Materials and Structures}, 0\penalty0
  (0):\penalty0 1--11, 2018.
\newblock \doi{10.1080/15376494.2018.1452321}.
\newblock URL \url{10.1080/15376494.2018.1452321}.

\bibitem[Bijelonja and Demird{\v{z}}i\'{c}(2000)]{Bijelonja2000}
I.~Bijelonja and S.~Muzaferija Demird{\v{z}}i\'{c}.
\newblock Some computational aspects of finite volume analysis of solid body
  deformation.
\newblock In \emph{Proceedings of 3$^{rd}$ Congres of Croatian Society of
  Mechanics}, pages 261--267, Dubrovnik, Croatia, 2000.

\bibitem[Cross et~al.(2003)Cross, Walshaw, Williams, Slone, Croft, and
  McManus]{Cross2003}
M.~Cross, C.~Walshaw, A.~J. Williams, A.~K. Slone, T.~N. Croft, and K.~McManus.
\newblock Parallel processing for nonlinear problems.
\newblock In E.~O{\~{n}}ate and D.~R.~J. Owen, editors, \emph{{VII}
  International Conference on Computational Plasticity}, Barcelona, Spain,
  2003. {CIMNE}.

\bibitem[Demird{\v{z}}i\'{c} et~al.(2003)Demird{\v{z}}i\'{c},
  D{\v{z}}aferovi\'{c}, and Ivankovi\'{c}]{Demirdzic2003}
I.~Demird{\v{z}}i\'{c}, E.~D{\v{z}}aferovi\'{c}, and A.~Ivankovi\'{c}.
\newblock Predicting residual stresses due to solidification in cast plastic
  plates.
\newblock In \emph{4$^{th}$ International Congress of Croatian Society of
  Mechanics}, 2003.

\bibitem[Williams et~al.(April 2003)Williams, Slone, Croft, and
  Cross]{Williams2003}
A.~J. Williams, A.~K. Slone, T.~N. Croft, and M.~Cross.
\newblock A mixed eulerian- lagrangian approach for metal forming.
\newblock In \emph{6th {ESAFORM} conference on Metal Forming}, Salerno, Italy,
  April 2003.

\bibitem[Kalkan(2011)]{Kalkan2011}
H.~Kalkan.
\newblock A combined experimental-numerical investigation on aluminium
  extrusion.
\newblock Master's thesis, Atilim University, Ankara, Turkey, 2011.

\bibitem[Mohan et~al.(2011)Mohan, Kara\v{c}, Murphy, and
  Ivankovi\'{c}]{Mohan2011}
J.~Mohan, A.~Kara\v{c}, N.~Murphy, and A.~Ivankovi\'{c}.
\newblock An experimental and numerical investigation of the mixed-mode
  fracture toughness and lap shear strength of aerospace grade composite
  joints.
\newblock \emph{Key Engineering Materials}, 488:\penalty0 549--552, 2011.

\bibitem[Cardiff et~al.(2014{\natexlab{f}})Cardiff, Jaeger, Tukovi\'{c}, and
  Ivankovi\'{c}]{Cardiff2014f}
P.~Cardiff, P.~De Jaeger, {\v{Z}}.~Tukovi\'{c}, and A.~Ivankovi\'{c}.
\newblock A finite approach to simulation of wire rolling.
\newblock In \emph{Joint Symposium of Irish Mechanics Society \& Irish Society
  for Scientific \& Engineering Computation, Galway}, 2014{\natexlab{f}}.

\bibitem[Cardiff et~al.(2018{\natexlab{b}})Cardiff, \v{Z}. Tukovi\'{c},
  Ivankovic, and Jaeger]{Cardiff2018b}
P.~Cardiff, \v{Z}. Tukovi\'{c}, A.~Ivankovic, and P.~De Jaeger.
\newblock Development of an arbitrary {Lagrangian-Eulerian} finite volume
  method for metal forming simulation in {OpenFOAM}.
\newblock In \emph{The 13th {OpenFOAM} Workshop {(OFW13)}}, Shanghai, China,
  June 24-29 2018{\natexlab{b}}.

\bibitem[Clancy et~al.(2018)Clancy, Cardiff, Jaeger, and Ivankovic]{Clancy2018}
M.~Clancy, P.~Cardiff, P.~De Jaeger, and A.~Ivankovic.
\newblock Implementation of advanced plasticity models in {OpenFOAM}.
\newblock In \emph{The 13th {OpenFOAM} Workshop {(OFW13)}}, Shanghai, China,
  June 24-29 2018.

\bibitem[Grossman et~al.(2002)Grossman, Bailey, Pericleous, and
  Slone]{Grossman2002}
D.~Grossman, C.~Bailey, K.~Pericleous, and A.~K. Slone.
\newblock Computational modelling of blood flow and artery wall interaction.
\newblock In \emph{9th Workshop on the Finite Element Methods in Biomedical
  Engineering Biomechanics and Related Fields}, University of Ulm, Germany,
  2002.
\newblock ISBN 3-9806183-5-8.

\bibitem[Anthony(2003)]{Anthony2003}
C.~M. Anthony.
\newblock Finite volume modelling of the human leg.
\newblock Master's thesis, Imperial College London, 2003.

\bibitem[Alakija et~al.(2005)Alakija, Ivankovi\'{c}, and
  Kara\v{c}]{Alakija2005}
O.~Alakija, A.~Ivankovi\'{c}, and A.~Kara\v{c}.
\newblock Finite volume solution to high rate wave propagation through a lung
  alveoli stack.
\newblock In \emph{IUTAM Symposium on Impact Biomechanics: From Fundamental
  Insights to Applications}, pages 281--288, 2005.

\bibitem[Quinn et~al.(2007)Quinn, Ivankovi\'{c}, and Kara\v{c}]{Quinn2007a}
N.~M. Quinn, A.~Ivankovi\'{c}, and A.~Kara\v{c}.
\newblock An experimental and numerical investigation into the deformation
  profiles of mock arteries.
\newblock In \emph{{ASME Summer Bioengineering Conference}}, 2007.

\bibitem[Quinn et~al.(2-3 April, 2007)Quinn, Ivankovi\'{c}, and
  Kara\v{c}]{Quinn2007b}
N.~M. Quinn, A.~Ivankovi\'{c}, and A.~Kara\v{c}.
\newblock A combined experimental and numerical investigation into early
  atherosclerosis.
\newblock In \emph{The Fifteenth {UK} Conference of the Association of
  Computational Mechanics in Engineering}, Glasgow, United Kingdom, 2-3 April,
  2007. Civil-Comp Press.
\newblock \doi{10.4203/ccp.85.11}.
\newblock Paper 11.

\bibitem[Safari et~al.(2008)Safari, Ivankovi\'{c}, and Tukovi\'{c}]{Safari2008}
A.~Safari, A.~Ivankovi\'{c}, and {\v{Z}}.~Tukovi\'{c}.
\newblock Numerical modelling of viscoelastic response of bacterial biofilm to
  mechanical stress.
\newblock In \emph{Bioengineering In Ireland Conference}, pages 25--26,
  Radisson Hotel, Sligo, January, 2008.

\bibitem[Kanyanta et~al.(2009{\natexlab{b}})Kanyanta, Ivankovi\'{c}, and
  Kara\v{c}]{Kanyanta2009b}
V.~Kanyanta, A.~Ivankovi\'{c}, and A.~Kara\v{c}.
\newblock Accurate prediction of blood flow transients: A fluid-structure
  interaction approach in hemodynamic wall shear stress.
\newblock In P.~Nithiarasu, editor, \emph{1$^{st}$ International Conference on
  Mathematical and Computational Biomedical Engineering - {CMBE2009}}, Swansea,
  UK, 2009{\natexlab{b}}.

\bibitem[Kanyanta(2009)]{Kanyanta2009c}
V.~Kanyanta.
\newblock \emph{Towards early diagnosis of atherosclerosis: Wall shear
  prediction}.
\newblock PhD thesis, University College Dublin, Ireland, 2009.

\bibitem[Kelly(2009)]{Kelly2009}
S.~Kelly.
\newblock \emph{Thrombus growth and its influence on the stress distribution in
  patient-based abdominal aortic aneurysm models}.
\newblock PhD thesis, University College Dublin, 2009.

\bibitem[Cardiff et~al.(2010)Cardiff, Kara\v{c}, Flavin, FitzPatrick, and
  Ivankovi\'{c}]{Cardiff2010}
P.~Cardiff, A.~Kara\v{c}, R.~Flavin, D.~FitzPatrick, and A.~Ivankovi\'{c}.
\newblock The development of a numerical model of the hip joint for complex
  soft tissue reconstructions around the hip joint.
\newblock In \emph{13$^{th}$ Annual Sir Bernard Crossland Symposium, University
  College Dublin}, Dublin, Ireland, 2010.

\bibitem[Cardiff et~al.(2011{\natexlab{c}})Cardiff, Kara\v{c}, Flavin,
  FitzPatrick, and Ivankovi\'{c}]{Cardiff2011a}
P.~Cardiff, A.~Kara\v{c}, R.~Flavin, D.~FitzPatrick, and A.~Ivankovi\'{c}.
\newblock The development of a numerical model of the hip joint.
\newblock In \emph{17$^{th}$ Bioengineering In Ireland}, Galway, Ireland,
  2011{\natexlab{c}}.

\bibitem[Cardiff et~al.(2011{\natexlab{d}})Cardiff, Kara\v{c}, Flavin,
  FitzPatrick, and Ivankovi\'{c}]{Cardiff2011b}
P.~Cardiff, A.~Kara\v{c}, R.~Flavin, D.~FitzPatrick, and A.~Ivankovi\'{c}.
\newblock Numerical analysis of the hip joint bones in contact.
\newblock In \emph{ACME-UK, Heriott-Watt University}, Edinburgh, Scotland,
  2011{\natexlab{d}}.

\bibitem[Quinn(2011)]{Quinn2011}
N.~Quinn.
\newblock \emph{Towards early diagnosis of atherosclerosis: combined
  experimental and numerical investigation into the deformation of mock
  arterial models}.
\newblock PhD thesis, University College Dublin, 2011.

\bibitem[Cardiff et~al.(2012{\natexlab{c}})Cardiff, Kara\v{c}, Flavin,
  FitzPatrick, and Ivankovi\'{c}]{Cardiff2012b}
P.~Cardiff, A.~Kara\v{c}, R.~Flavin, D.~FitzPatrick, and A.~Ivankovi\'{c}.
\newblock Modelling the muscles for hip joint stress analysis using a finite
  volume methodology.
\newblock In \emph{18$^{th}$ Bioengineering In Ireland}, Belfast, Northern
  Ireland, 2012{\natexlab{c}}.

\bibitem[Cardiff et~al.(2014{\natexlab{g}})Cardiff, Kara\v{c}, FitzPatrick,
  Flavin, and Ivankovi\'{c}]{Cardiff2014c}
P.~Cardiff, A.~Kara\v{c}, D.~FitzPatrick, R.~Flavin, and A.~Ivankovi\'{c}.
\newblock Development of mapped stress-field boundary conditions based on a
  {Hill}-type muscle model.
\newblock \emph{International Journal for Numerical Methods in Biomedical
  Engineering}, 2014{\natexlab{g}}.
\newblock \doi{10.1002/cnm}.

\bibitem[Parsa(2014)]{Parsa2014}
H.~Khalili Parsa.
\newblock \emph{Compression tests on fluid-filled gelatine microcapsules: A
  combined experimental/numerical study}.
\newblock PhD thesis, University College Dublin, 2014.

\bibitem[Safari(2015)]{Safari2015}
A.~Safari.
\newblock \emph{A combined experimental and numerical study of biofilm
  detachment}.
\newblock PhD thesis, University College Dublin, 2015.

\bibitem[Fitzgerald et~al.(2016)Fitzgerald, Cardiff, Flavin, and
  Ivankovic]{Fitzgerald2016}
K.~Fitzgerald, P.~Cardiff, R.~Flavin, and A.~Ivankovic.
\newblock Calculation of hip joint contact pressures using a high resolution
  finite volume model with {CT}-based properties.
\newblock In \emph{Bioengineering in Ireland}, 2016.

\bibitem[Fitzgerald et~al.(2017)Fitzgerald, Cardiff, Flavin, and
  Ivankovic]{Fitzgerald2017}
K.~Fitzgerald, P.~Cardiff, R.~Flavin, and A.~Ivankovic.
\newblock Towards in silico analysis of total hip arthroplasty mechanics.
\newblock In \emph{{XXVI Congress of the International Society of
  Biomechanics}}, Brisbane, Australia, 2017.

\bibitem[Muralidharan et~al.(2017)Muralidharan, Cardiff, Flavin, and
  Ivankovic]{Muralidharan2017}
L.~Muralidharan, P.~Cardiff, R.~Flavin, and A.~Ivankovic.
\newblock A numerical model for the calculation of ankle joint stresses.
\newblock In \emph{{XXVI Congress of the International Society of
  Biomechanics}}, Brisbane, Australia, 2017.

\bibitem[Kova\v{c}evi\'{c} et~al.(2002{\natexlab{a}})Kova\v{c}evi\'{c},
  Sto\v{s}i\'{c}, and Smith]{Kovacevic2002a}
A.~Kova\v{c}evi\'{c}, N.~Sto\v{s}i\'{c}, and I.~K. Smith.
\newblock Solid-fluid interaction in screw compressors.
\newblock In \emph{{XVI} International Compressor Engineering Conference at
  Purdue}, 2002{\natexlab{a}}.

\bibitem[Kova\v{c}evi\'{c} et~al.(2002{\natexlab{b}})Kova\v{c}evi\'{c},
  Sto\v{s}i\'{c}, and Smith]{Kovacevic2002b}
A.~Kova\v{c}evi\'{c}, N.~Sto\v{s}i\'{c}, and I.~K. Smith.
\newblock Three-dimensional modelling of solid-fluid interaction as a design
  tool in screw compressors.
\newblock In \emph{International Design Conference - DESIGN 2002}, Dubrovnik,
  Croatia, 2002{\natexlab{b}}.

\bibitem[Kova\v{c}evi\'{c} et~al.(2002{\natexlab{c}})Kova\v{c}evi\'{c},
  Sto\v{s}i\'{c}, and Smith]{Kovacevic2002c}
A.~Kova\v{c}evi\'{c}, N.~Sto\v{s}i\'{c}, and I.~K. Smith.
\newblock The influence of rotor deflection upon the screw compressor process.
\newblock In \emph{Schraubencompressor tagung (Screw compressor meeting)},
  Dortmund, Germany, 2002{\natexlab{c}}.

\bibitem[Kova\v{c}evi\'{c} et~al.(2002{\natexlab{d}})Kova\v{c}evi\'{c},
  Sto\v{s}i\'{c}, and Smith]{Kovacevic2002d}
A.~Kova\v{c}evi\'{c}, N.~Sto\v{s}i\'{c}, and I.~K. Smith.
\newblock Numerical simulation of fluid flow and solid structure in screw
  compressors.
\newblock In \emph{Proceedings of 2002 {ASME} Congress}, New Orleans, USA,
  2002{\natexlab{d}}. Symposium on the Analysis and Applications of Heat Pump
  and Refrigeration Systems.

\bibitem[Kova\v{c}evi\'{c} et~al.(2003)Kova\v{c}evi\'{c}, Sto\v{s}i\'{c}, and
  Smith]{Kovacevic2003}
A.~Kova\v{c}evi\'{c}, N.~Sto\v{s}i\'{c}, and I.~K. Smith.
\newblock Fluid-solid interaction for extension of range in screw machine
  application.
\newblock In \emph{Advances of {CFD} in Fluid Machinery Design, {ImechE}
  Seminar}, London, UK, 2003.

\bibitem[Kova\v{c}evi\'{c} et~al.(2004{\natexlab{b}})Kova\v{c}evi\'{c},
  Sto\v{s}i\'{c}, Smith, and Muji\'{c}]{Kovacevic2004a}
A.~Kova\v{c}evi\'{c}, N.~Sto\v{s}i\'{c}, I.~K. Smith, and E.~Muji\'{c}.
\newblock Fluid-solid interaction in the design of multifunctional screw
  machines.
\newblock In \emph{8$^{th}$ International Design Conference - Design 2004},
  volume~2, pages 1289--1295, Dubrovnik, Croatia, 2004{\natexlab{b}}.

\bibitem[Kova\v{c}evi\'{c} et~al.(2011)Kova\v{c}evi\'{c}, Muji\'{c},
  Sto\v{s}i\'{c}, and Smith]{Kovacevic2011}
A.~Kova\v{c}evi\'{c}, E.~Muji\'{c}, N.~Sto\v{s}i\'{c}, and I.~K. Smith.
\newblock Extending the role of computational fluid dynamics in screw machines.
\newblock \emph{Journal of Process Mechanical Engineering, Part {E},
  Proceedings of the Institution of Mechanical Engineers}, pages 83--97, 2011.

\bibitem[Bijelonja(2011{\natexlab{b}})]{Bijelonja2011b}
I.~Bijelonja.
\newblock A finite volume method for a geomechanics problem.
\newblock In \emph{Proceedings of the 22$^{nd}$ International {DAAAM}
  Symposium}, pages 323--324, 2011{\natexlab{b}}.

\bibitem[Demird{\v{z}}i\'{c} et~al.(1997{\natexlab{b}})Demird{\v{z}}i\'{c},
  Muzaferija, and Peri\'{c}]{Demirdzic1997a}
I.~Demird{\v{z}}i\'{c}, S.~Muzaferija, and M.~Peri\'{c}.
\newblock Advances in computation of heat transfer, fluid flow, and solid body
  deformation using finite volume approaches.
\newblock In E.~M.~Sparrow W.~J.~Minkowycz, editor, \emph{Advances in Numerical
  Heat Transfer, Chapter 2}, pages 59--96. Taylor \& Francis,
  1997{\natexlab{b}}.

\bibitem[Bijelonja(2005)]{Bijelonja2005b}
I.~Bijelonja.
\newblock Finite volume method analysis of large strain elasto-plastic
  deformation.
\newblock In \emph{The 16$^{th}$ {DAAAM} International Symposium}, Opatia,
  Croatia, 2005.

\bibitem[Sabbagh-Yazdi et~al.(2008{\natexlab{a}})Sabbagh-Yazdi, Alkhamis,
  Mastorakis, and Esmaili]{Sabbagh-Yazdi2008a}
S.~R. Sabbagh-Yazdi, M.~T. Alkhamis, N.~E. Mastorakis, and M.~Esmaili.
\newblock Finite volume analysis of two-dimensional strain in a thick pipe with
  internal fluid pressure.
\newblock \emph{International Journal of Mathematical Models and Methods in
  Applied Sciences}, 2:\penalty0 162--167, 2008{\natexlab{a}}.

\bibitem[Sabbagh-Yazdi et~al.(2008{\natexlab{b}})Sabbagh-Yazdi, Mastorakis, and
  Esmaili]{Sabbagh-Yazdi2008b}
S.~R. Sabbagh-Yazdi, N.~E. Mastorakis, and M.~Esmaili.
\newblock Explicit {2D} matrix free {Galerkin} finite volume solution of plane
  strain structural problems on triangular meshes.
\newblock \emph{International Journal of Mathematics and Computers in
  Simulations}, 2:\penalty0 1--8, 2008{\natexlab{b}}.

\bibitem[Alkhamis et~al.(2008)Alkhamis, Sabbagh-Yazdi, Esmaeili, and
  Wegian]{Alkhamis2008}
M.~T. Alkhamis, S.~R. Sabbagh-Yazdi, M.~Esmaeili, and F.~M. Wegian.
\newblock Utilizing {NASIR} {Galerkin} finite volume analyzer for {2D} plane
  strain problems under static and vibrating concentrated loads.
\newblock \emph{Jordan Journal of Civil Engineering}, 2:\penalty0 335--343,
  2008.

\bibitem[Sabbagh-Yazdi et~al.(2009)Sabbagh-Yazdi, Alimohammadi, and
  Mastorakis]{Sabbagh-Yazdi2009}
S.~R. Sabbagh-Yazdi, S.~Alimohammadi, and N.~E. Mastorakis.
\newblock Comparison of finite element and finite volume solvers results for
  plane-stress displacements in plate with oval hole.
\newblock In \emph{{Proceedings of the 4$^{th}$ IASME/WSEAS International
  Conference on Continuum Mechanics CM'09}}, pages 168--173, 2009.

\bibitem[Sabbagh-Yazdi and Amiri-Saadatabadi(2011)]{Sabbagh-Yazdi2011a}
S.~R. Sabbagh-Yazdi and T.~Amiri-Saadatabadi.
\newblock Sequential computations of two-dimensional temperature profiles and
  thermal stresses on an unstructured triangular mesh by {GFVM} method.
\newblock \emph{International Journal of Civil Engineering}, 9:\penalty0
  171--182, 2011.

\bibitem[Sabbagh-Yazdi et~al.(2011)Sabbagh-Yazdi, Esmaili, and
  Alkhamis]{Sabbagh-Yazdi2011b}
S.~R. Sabbagh-Yazdi, M.~Esmaili, and M.~T. Alkhamis.
\newblock Symmetric conditions for strain analysis in a long thick cylinder
  under internal pressure using {NASIR} unstructured {GFVM} solver.
\newblock \emph{Jordan Journal of Civil Engineering}, 5:\penalty0 258--267,
  2011.

\bibitem[Sabbagh-Yazdi and Ali-Mohammadi(2011)]{Sabbagh-Yazdi2011c}
S.~R. Sabbagh-Yazdi and S.~Ali-Mohammadi.
\newblock Performance evaluation of iterative {GFVM} on coarse unstructured
  triangular meshes and comparison with matrix manipulation based solution
  methods.
\newblock \emph{Scientia Iranica}, 18:\penalty0 131--138, 2011.

\bibitem[Sabbagh-Yazdi et~al.(2012)Sabbagh-Yazdi, AliMohammadi, and
  Pipelzadeh]{Sabbagh-Yazdi2012a}
S.~R. Sabbagh-Yazdi, S.~AliMohammadi, and M.~K. Pipelzadeh.
\newblock Unstructured finite volume method for matrix-free explicit solution
  of stress-strain fields in two-dimensional problems with curved boundaries in
  equilibrium condition.
\newblock \emph{Applied Mathematical Modelling}, 36:\penalty0 2224--2236, 2012.

\bibitem[Sabbagh-Yazdi and Bayatlou(2012)]{Sabbagh-Yazdi2012b}
S.~R. Sabbagh-Yazdi and M.~Bayatlou.
\newblock Equilibrium condition nonlinear modeling of a cracked concrete beam
  using a {2D} {Galerkin} finite volume solver.
\newblock \emph{Computational Methods in Civil Engineering}, 3:\penalty0
  63--76, 2012.

\bibitem[Bailey et~al.(1999{\natexlab{c}})Bailey, Chow, Cross, Pericleous,
  Taylor, Croft, Wheeler, and Lu]{Bailey1999c}
C.~Bailey, P.~Chow, M.~Cross, K.~Pericleous, G.~A. Taylor, T.~N. Croft,
  D.~Wheeler, and H.~Lu.
\newblock Finite volume methods for multiphysics problems.
\newblock In D.~Haenel, R.~Vilsmeirer, and F.~Benkhaldoun, editors,
  \emph{Finite Volumes for Complex Applications, II - Problems and
  Perspectives}, Duisburg, Germany, 1999{\natexlab{c}}. Hermes Science.

\bibitem[Hitchings et~al.(1987)Hitchings, Davies, and
  Kamoulakos]{Hitchings1987}
D.~Hitchings, G.~A.~O. Davies, and A.~Kamoulakos.
\newblock \emph{Linear static benchmarks}.
\newblock {International Association for the Engineering Analysis Community \&
  National Agency for Finite Element Methods \& Standards (NAFEMS)}, Glasgow,
  UK, 1987.

\bibitem[Voller(2009)]{Voller2009}
V.~R. Voller.
\newblock \emph{Basic Control Volume Finite Element Methods for Fluids and
  Solids}.
\newblock World Scientific, Singapore, 2009.

\end{thebibliography}

\end{document}